\documentclass{article}
\usepackage{mathtools}
\usepackage[export]{adjustbox}
\usepackage{amsfonts}
\usepackage{amsmath}
\usepackage{amsthm}
\usepackage{amsbsy}
\usepackage{amssymb}
\usepackage{authblk}
\usepackage{algcompatible}
\usepackage{algpseudocodex}
\usepackage{algorithm}
\newfloat{algorithm}{t}{lop}
\usepackage{array}
\usepackage{bm}
\usepackage{bbm}
\usepackage{booktabs}
\usepackage{color,soul}
\usepackage{comment}
\usepackage{caption}
\usepackage{dsfont}
\usepackage{diagbox}
\usepackage{empheq}
\usepackage[shortlabels]{enumitem}
\usepackage{etoolbox}
\usepackage{float}
\usepackage{geometry}
\usepackage{graphicx} % Required for inserting images
\usepackage{hhline}
\usepackage{hyperref}
\usepackage{makecell}
\usepackage{mathrsfs}
\usepackage{mathtools}
\usepackage{multirow}
\usepackage{natbib}
\usepackage{romannum}
\usepackage{caption}
\usepackage{subcaption}
\usepackage{titling}
\usepackage{tikz}
\tikzset{algpxIndentLine/.style={draw=gray,very thin}}
\usepackage[toc,page]{appendix}
\newcommand{\subtitle}[1]{%
  \posttitle{%
    \par\end{center}
    \begin{center}\large#1\end{center}
    \vskip0.5em}%
}
\DeclareMathOperator{\C}{\mathcal{C}}

\DeclareMathOperator{\K}{\mathcal{K}}
\DeclareMathOperator{\M}{\mathcal{M}}

\DeclareMathOperator{\R}{\mathbb{R}}
\DeclareMathOperator{\s}{\mathcal{S}}
\DeclareMathOperator{\T}{\mathcal{T}}

\DeclareMathOperator{\W}{\mathcal{W}}

\DeclareMathOperator{\sse}{\subseteq}
\DeclareMathOperator{\del}{\partial}

\DeclareMathOperator{\on}{\text{ on }}

\DeclareMathOperator{\sign}{sign}
\DeclareMathOperator{\prob}{\mathbb{P}}

\newtheorem{theorem}{Theorem}[section]
\newtheorem*{theorem*}{Theorem}

\newtheorem*{lemma*}{Lemma}

\newtheorem*{proposition*}{Proposition}

\newtheorem{innercustomgeneric}{\customgenericname}
\providecommand{\customgenericname}{}
\newcommand{\newcustomtheorem}[2]{%
  \newenvironment{#1}[1]
  {%
   \renewcommand\customgenericname{#2}%
   \renewcommand\theinnercustomgeneric{##1}%
   \innercustomgeneric
  }
  {\endinnercustomgeneric}
}
\newcustomtheorem{customthm}{Theorem}
\newcustomtheorem{customlemma}{Lemma}
\newcustomtheorem{customproposition}{Proposition}

\newcommand{\bdr}[1]{\partial{#1}}
\newcommand{\adj}[1]{{#1}^*}

\newcommand{\lint}[1]{d \ell_{\bm{#1}}}
\newcommand{\n}[1]{\bm{n}(\bm{#1})}
\newcommand{\bdrint}[1]{\int_{\bdr{#1}}}
\newcommand{\slp}[1]{\mathcal{S}_{\bdr{\Omega_{#1}}}[\adj{\gamma}_{#1}]}
\title{Walk-on-Interfaces: A Monte Carlo Estimator for an Elliptic Interface Problem with Nonhomogeneous Flux Jump Conditions and a Neumann Boundary Condition}
\author{Xinwen Ding, Adam R Stinchcombe}
\date{\today}

\begin{document}
\pagenumbering{arabic}
\fontsize{12pt}{15pt}\selectfont
\maketitle
%Journals of interest: JCP, SISC, JSC
\tableofcontents
\newpage
\section*{Abstract}
Elliptic interface problems arise in numerous scientific and engineering applications, modeling heterogeneous materials in which physical properties change discontinuously across interfaces. In this paper, we present \textit{Walk-on-Interfaces} (WoI), a grid-free Monte Carlo estimator for a class of Neumann elliptic interface problems with nonhomogeneous flux jump conditions. Our Monte Carlo estimators maintain consistent accuracy throughout the domain and, thus, do not suffer from the well-known close-to-source evaluation issue near the interfaces. We also presented a simple modification with reduced variance. Estimation of the gradient of the solution can be performed, with almost no additional cost, by simply computing the gradient of the Green's function in WoI. Taking a scientific machine learning approach, we use our estimators to provide training data for a deep neural network that outputs a continuous representation of the solution. This regularizes our solution estimates by removing the high-frequency Monte Carlo error. All of our estimators are highly parallelizable, have a $\mathcal{O}(1 / \sqrt{\mathcal{W}})$ convergence rate in the number of samples, and generalize naturally to higher dimensions. We solve problems with many interfaces that have irregular geometry and in up to dimension six. Numerical experiments demonstrate the effectiveness of the approach and to highlight its potential in solving problems motivated by real-world applications.
\section{Introduction}\label{intro}
Elliptic interface problems arise in a wide range of scientific and engineering applications, from electrostatics of biomolecules~\cite{ho2012fast} and tissues with packed cells~\cite{Ying_Beale_2013}, to porous media flow~\cite{philip1970flow, Stastna2023}, geophysics~\cite{sagar2018handbook}, and composite materials~\cite{WANG2019117}. Elliptic interface problems are often the forward problem in important inverse problems. In particular, electrical impedance tomography~\cite{TyniStinchcombeAlexakis2024} and 
electroencephalography~\cite{FEMforEEG} can be formulated as such an inverse problem.

Mathematically, interface problems are characterized by partial differential equations (PDEs) with spatially varying coefficients that are discontinuous across interfaces, reflecting sudden changes in physical properties. The canonical elliptic PDE is 
\begin{equation}\label{interface_problem_PDE}
    \nabla \cdot (\sigma(\bm{x}) \nabla u(\bm{x})) = 0, \quad \bm{x} \in \Omega.
\end{equation}
We assume that $\sigma(\bm{x}) \geq 0$ is piecewise constant. The boundaries between the regions of constant $\sigma$ give rise to interfaces. Since $\sigma(\bm{x})$ is not differentiable across these interfaces, it is common to derive interface conditions consistent with Eq.~\eqref{interface_problem_PDE}.

%cite Shravan Veerapaneni https://dept.math.lsa.umich.edu/~shravan/Publications.html

\begin{figure}[!ht]
\centering
\includegraphics[width=0.8\columnwidth]{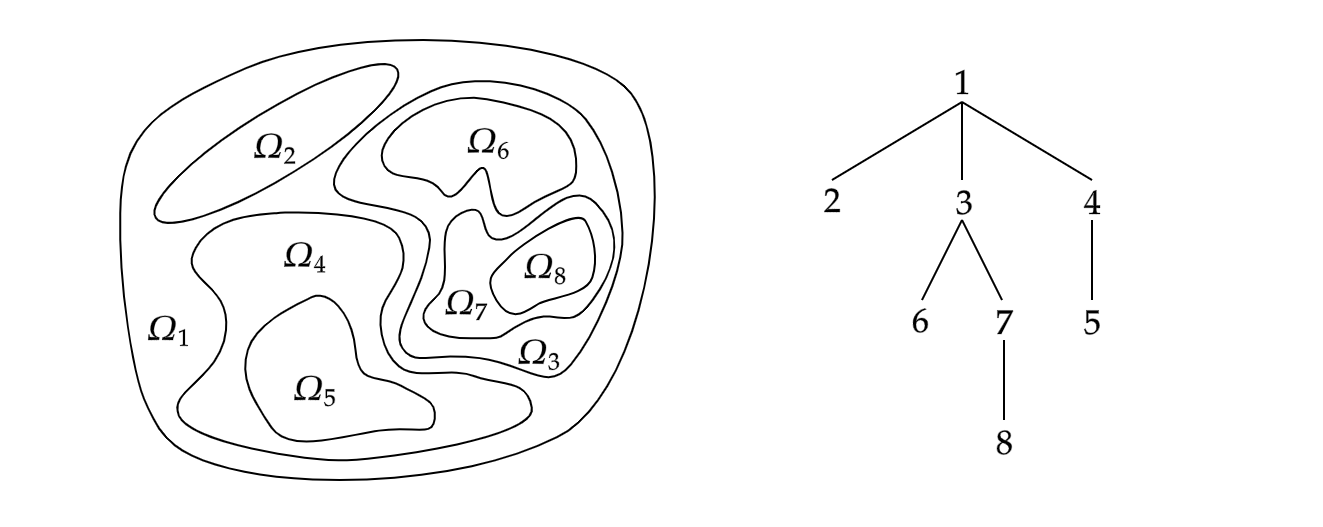}
\caption{An example of tree-structured domain defined by \cite{Bower_thesis}. The outer
boundary of the root region of the tree, $\Omega_1$, is also the boundary of the domain. Each interface is shared between a region and its unique parent.}\label{ex_tree_structured_domain}
\end{figure}

To organize the regions of constant $\sigma$, we employ a tree definition of domains with subregions that have different physical properties separated by interfaces. Consider a bounded, simply connected domain $\Omega \sse \R^m$ with a smooth boundary, $\partial \Omega$. Let $N < \infty$ be the total number of non-overlapping regions inside $\Omega$. For each $i \in \{1,\ldots,N\}$, we denote the $i^{\text{th}}$ region and all its descendants as $\Omega_i \subseteq \Omega$, and we require the interface $\bdr{\Omega_i} \in \mathcal{C}^2$. The tree definition from \cite{Bower_thesis}, where $\bdr{\Omega_1} = \bdr{\Omega}$, is adopted here to describe the hierarchical relationships among the regions. An example of a domain layout and its associated tree definition is provided in Figure~\ref{ex_tree_structured_domain}. 

Let $\n{x}$ denote the unit outward normal on the boundary. Correspondingly, we denote $\del_{\bm{n}}$ or $\frac{\del}{\del \bm{n}}$ as the normal derivative. We further define the jump of a function $f: \bdr{\Omega_i} \to \R$ across some interface $\bdr{\Omega_i}$ at some point $\bm{x}_0$ by
\[
[f(\bm{x}_0)] \coloneqq f_{+}(\bm{x}_0) - f_{-}(\bm{x}_0), 
\]
where $f_{+}(\bm{x}_0)$ is the limit point of $f(\bm{x})$ as $\bm{x} \to \bm{x}_0$ from the exterior of $\Omega_i$ (i.e. from $\Omega \backslash \Omega_i$), while $f_{-}(\bm{x}_0)$ is the limit point of $f(\bm{x})$ as $\bm{x} \to \bm{x}_0$ from the interior of $\Omega_i$ (i.e. from $\Omega_i$).

Under the assumption that $\sigma(\bm{x})$ is piecewise constant, the PDE Eq.~\eqref{interface_problem_PDE} can be cast into the interface problem with Neumann boundary condition,
\begin{subequations}\label{interface_problem}
\begin{empheq}[left=\empheqlbrace]{align}
 \Delta u(\bm{x}) &= 0 \quad \text{in $\Omega \ \backslash \ \cup_i \bdr{\Omega_i}$}  \label{Laplace_eq }\\
 % ========================
 \sigma_1 \partial_{\bm{n}}u(\bm{x}) &= b_1(\bm{x}) \quad \text{on $\bdr{\Omega_1}$} \label{boundary_condition} \\
 % ========================
 [u](\bm{x}) &= 0 \quad \text{on $\bdr{\Omega_i}$, $\forall i \in \{2,3, \dots, N\}$} \label{no_jump} \\
 % ========================
 [\sigma \partial_{\bm{n}} u](\bm{x}) &= b_i(\bm{x}) \quad \text{on $\bdr{\Omega_i}$, $\forall i \in \{2,3, \dots, N\}$} \label{interface_condition},
\end{empheq}
\end{subequations}
where $b_i \in \mathcal{C}^{0,1}$ and $\sigma(\bm{x})|_{\Omega_i} \equiv \sigma_i$ for each $i = 1, 2, \dots, N$. For the $i^{th}$ region, we let $p_i$ denote the index of its unique parent in the tree definition of the domain. Then, $\sigma(\bm{x})$ jumps across $\bdr{\Omega_i}$ from $\sigma_i$ to $\sigma_{p_i}$. We additionally require $\int_{\bdr{\Omega_i}} b_i(\bm{x}) \ dA = 0, \forall i \in \{1,\ldots,N\}$ to ensure the well-posedness of the interface problem, up to an additive constant \cite{evansPDE}. 

Classical numerical methods for solving elliptic interface problems have been extensively studied over the past several decades. These numerical methods typically require a high-quality underlying mesh to define the boundary and the interface. Methods that explicitly compute discretized derivatives, such as the finite element method~\cite{TAVARES20147535,GUO2020109478}, the finite difference method~\cite{leveque1994immersed,LiEllipticInterfaceProblems, FENG2024112635}, and the finite volume method \cite{dong2005comparison}, often lead to a large and ill-conditioned linear system. Integral equation methods, such as the boundary element method \cite{deMunck} and the boundary integral equation method \cite{Ying_Beale_2013}, result in well-conditioned dense linear systems after discretization. The linear systems can be solved rapidly when coupled with the fast multiple method \cite{greengard1987fast, martinsson2019fast}. One well-known challenge for integral equation methods is accurately evaluating a nearly singular kernel at points close to the interface, which is usually tackled by carefully designed quadratures \cite{Helsing2008, Klockner2013}. Another common challenge is handling domains that contain a large number (often close-to-touching) of interfaces \cite{Barnett2018}. %\xinwen{\href{https://www.sciencedirect.com/science/article/abs/pii/S0021999124004662}{this paper} listed lots of methods in intro for Dirichlet problem, not sure whether we need to look at them.}

Monte Carlo methods are a powerful class of grid-free computational tools that estimate solutions to PDEs through random sampling. There are two major types of algorithms in the Monte Carlo family, \textit{Walk-on-Sphere} (WoS) and \textit{Walk-on-Boundary} (WoB). WoS~\cite{Muller1956} was first proposed to solve Laplace equations by simulating the paths of Brownian motion. Later, it was generalized to other problems \cite{bossy2010, kyprianou2017, sawhneyseyb22gridfree}, and more complex domain geometries \cite{Sawhney2023_WoSt, DeLambilly2023, Miller2024_WoSt}. Algorithms in the class of WoS rely on an $\epsilon$-shell near the boundary to stop the Brownian motion, which results in inaccurate estimation near the boundary. The WoB methods, rooted in potential theory, do not suffer from this limitation. It was first introduced to solve the Lam\'e equation \cite{Sabelfeld1982}. Later, WoB was extended to solve Laplace and Poisson problems with Dirichlet, Neumann, and Robin boundary conditions \cite{Sabelfeld_Simonov_2016, sabelfeld1991, sabelfeld1994}. One recent work \cite{Sugimoto_2023} summarized existing WoB estimators, proposed a new WoB estimator for the Neumann problem, and emphasized particular applications of WoB in computer graphics. 

In this paper, we present \textit{Walk-on-Interfaces} (WoI), an extension of the family of WoB methods designed to solve the Neumann elliptic interface problem with nonhomogeneous flux jump conditions in Eq.~\eqref{interface_problem} with $N \geq 1$. To the best of our knowledge, this is the first Monte Carlo estimator for an interface problem.

The proposed Monte Carlo estimators offer several important advantages. They are grid-free, highly parallelizable, and can be easily generalized to solve high-dimensional problems. Our estimators maintain consistent accuracy near and even on the boundary and the interfaces, as well as in the interior of the domain. By computing the gradient of our Monte Carlo estimators, we obtain an estimation to the gradient of the solution to the interface problem almost without extra cost. The estimators and their associated gradient estimators have $\mathcal{O}(1/\sqrt{\mathcal{W}})$ convergence rate in the number of samples $\mathcal{W}$.

As an application, these estimators can be combined with a neural network to form a WoI framework, in which the estimators play the role of data generators. The output of this framework provides a continuous approximation of the solution to the interface problem, free from high-frequency Monte Carlo error. This data-driven approach is especially beneficial in domains with multiple interfaces.

The rest of this paper is organized as follows: In Section \ref{wob}, we review and summarize the WoB estimator for the interior Neumann problem. In section \ref{bie}, we present the boundary integral formulation for the elliptic interface problem, which is foundational for the WoI estimator. Building on this, we present the WoI estimator and its variance-reduced counterpart in Section \ref{woi}. Then, in Section \ref{numerical_results}, we investigate the performance of our estimators and the WoI framework via both synthetic examples with known ground truth and examples motivated by real-world applications. Finally, we conclude this paper and discuss possible directions for future work in Section \ref{conclusion}.
\section{Summary of the Walk-on-Boundary Method}\label{wob}
In this section, we provide a brief summary of the WoB method \cite{Sabelfeld_Simonov_2016} for the interior Neumann problem for the Laplace equation,
\begin{equation}\label{interior_Neumann}
\begin{cases}
    \Delta v(\bm{x}) = 0, \quad &\bm{x} \in \Omega \\
    \frac{\del v}{\del \n{x}} =  g(\bm{x}), \quad &\bm{x} \on \bdr{\Omega}
\end{cases}.
\end{equation}
The Neumann boundary data $g$ must satisfy the solvability condition,
\[
\bdrint{\Omega} g(\bm{x}) \ dA = 0.
\]
%\adam{this should be in a generic dimension, so replace $dl$ with $dA$}
It is well-known that the solution to Eq.~\eqref{interior_Neumann} is unique up to a constant \cite{evansPDE}. The WoB estimator estimates one particular solution, while the full family of solutions is obtained by shifting the particular solution according to a selected reference point.
The estimator relies on the fact that solutions, up to an additive constant, to the interior Neumann problem can be represented with a single-layer potential,
\[
v(\bm{x}) = \bdrint{\Omega} \Phi(\bm{x}, \bm{y}) \adj{\gamma}(\bm{y}) \ d A_{\bm{y}}, \quad \forall \bm{x} \in \Omega,
\]
in which \[
\Phi(\bm{x}, \bm{y}) = 
\begin{cases}
-\frac{1}{2\pi} \ln|\bm{x} - \bm{y}| & d = 2 \\
\frac{1}{d(d-2)\alpha(d)}\frac{1}{|\bm{x} - \bm{y}|^{d-2}} & d \geq 3
\end{cases}\quad
\]
is the Green's function, $\alpha(d) = \frac{\pi^{d/2}}{\Gamma\left(\frac{d}{2} + 1\right)}$ is the volume of the unit ball in $\R^d$, and $\adj{\gamma}(\bm{y})$ is a continuous density on the boundary $\bdr{\Omega}$. 
Let us denote 
\[
\adj{K}(\bm{x}, \bm{y}) 
\coloneqq  -\frac{\del \Phi (\bm{x}, \bm{y})}{\del \n{x}} 
= \frac{1}{d\alpha(d)} \frac{(\bm{x} - \bm{y}) \cdot \n{x}}{|\bm{x} - \bm{y}|^d}.
\]

\begin{figure}[!ht]
    \centering
    \begin{subfigure}[b]{0.45\columnwidth}
    \includegraphics[width=\textwidth]{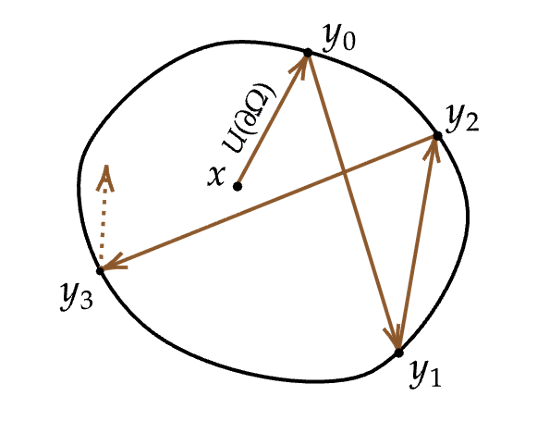}
    \caption{}
    \label{wob_convex_trajectory}
    \end{subfigure}
    \begin{subfigure}[b]{0.45\columnwidth}
    \includegraphics[width=\textwidth]{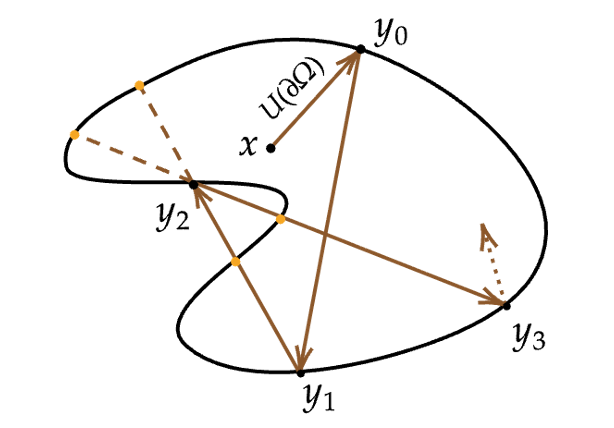}
    \caption{}
    \label{wob_general_trajectory}
    \end{subfigure}
    \caption{A example diagram of the WoB method in a $d=2$ domain. (a) One possible WoB trajectory with a convex domain with $M = 3$. Notice that when $\Omega$ is convex, $q(\bm{y}_{i-1}, \bm{y}_i) = 1$ always holds when $\ell(t)$ is not in the tangent plane. (b) One possible WoB trajectory with a general domain with $M = 3$. In both situations, $\bm{y}_0$ is sampled uniformly on $\bdr{\Omega}$, annotated by $U(\bdr{\Omega})$ alongside the arrow from $\bm{x}$ to $\bm{y}_0$. For all other arrows without annotations, $\bm{y}_i$'s are determined by finding an intersection of $\bdr{\Omega}$ and a ray with random direction with origin $\bm{y}_{i-1}$.}
    \label{wob_trajectory}
\end{figure}

Let $\{Y_0, Y_1, \dots, Y_i \dots \}$ be a Markov chain with transitional probability density functions (PDFs),
\begin{equation}\label{wob_neumann_p}
    p(\bm{x}, \bm{y}_0) = \frac{1}{|\bdr{\Omega}|} \quad \text{and} \quad
    p(\bm{y}_{i-1}, \bm{y}_i) = \frac{2}{\textit{d}\alpha(\textit{d})} \frac{|(\bm{y}_i - \bm{y}_{i-1}) \cdot \bm{n}(\bm{y}_i)|}{q(\bm{y}_{i-1}, \bm{y}_i)|\bm{y}_i -  \bm{y}_{i-1}|^m}, \notag
\end{equation}
such that $\bm{y}_i$ is a realization of the state $Y_i$, and that $q(\bm{y}_{i-1}, \bm{y}_i)$ is the number of intersections, distinct from $\bm{y}_{i-1}$, between $\Omega$ and a line,
\[
\ell(t) = t \bm{d} + \bm{y}_{i-1},
\]
whose direction $\bm{d}$ is uniformly randomly selected in $\R^d$.
The WoB estimator is \begin{align}\label{WoB_neumann}
v(\bm{x}) &\approx \mathbb{E}_{Y_i} \bigg[ 2 \sum_{i=0}^{M-1} 2^i T_i Q_i \Phi(\bm{x}, Y_i) + 2^M T_M Q_M \Phi(\bm{x}, Y_M) \bigg],
\end{align}
in which 
\begin{equation} \label{wob_Q}
Q_i = 
\begin{cases}
    \frac{g(Y_0)}{p(\bm{x}, Y_0)} \quad & \text{if $i = 0$}\\
     Q_{i-1} \frac{\adj{K}(Y_i, Y_{i-1})}{p(Y_{i-1}, Y_i)} = \frac{1}{2} Q_{i-1} \quad & \text{if $i \geq 1$}\\
\end{cases},
\end{equation}
and
\begin{equation}\label{WoB_T}
    T_i = 
    \begin{cases}
        1 \quad & \text{if $i = 0$} \\
        \prod_{m=1}^{i} q(Y_{m-1}, Y_m) \sign[(Y_m - Y_{m-1}) \cdot \bm{n}(Y_m)] \quad & \text{if $i \geq 1$}
    \end{cases}.
\end{equation}
 We provide concrete examples in Figure \ref{wob_trajectory} to visualize WoB trajectories. The Markov chain moves around the boundary by finding the intersections between the boundary and lines in uniformly random directions. For nonconvex domains, where multiple intersections may occur, one is selected uniformly at random.
%\section{WoB Archive}
%\input{content/2-wob_archive}
\section{Boundary Integral Formulation}\label{bie}
As the first step towards our \textit{Walk-on-Interfaces} estimator, we present the boundary integral formulation for the elliptic interface problem. This formulation is based on the work in a recent thesis by Kyle Waller Bower~\cite{Bower_thesis}.

Given a query point $\bm{x} \in \Omega$, the solution $u(\bm{x})$ of the elliptic interface problem can be represented as the sum of $N$ single-layer potentials,
\begin{equation}\label{interface_problem_soln_single_layer_potential}
u(\bm{x}) = \sum_{j=1}^{N} \slp{j}(\bm{x}) = \sum_{j=1}^{N} \bdrint{\Omega_{j}} \Phi(\bm{x}, \bm{y}) \adj{\gamma}_j(\bm{y}) \ d A_{\bm{y}}, \notag
\end{equation}
in which $\Phi (\bm{x}, \bm{y})$ is the Green's function, and $\adj{\gamma}_j \in \mathcal{C}(\bdr{\Omega_j})$ is the charge density function on $\bdr{\Omega_j}$. Furthermore, since we assumed the interfaces do not overlap with each other, the solution is then
\begin{equation}\label{interface_problem_soln_sum_of_density}
    u(\bm{x})
    = \int_{\Gamma} \Phi(\bm{x}, \bm{y}) \bigg( 
    \sum_{j=1}^{N} \adj{\gamma}_j(\bm{y}) \mathds{1}_{\bdr{\Omega_j}}
    \bigg) \ d A_{\bm{y}},
\end{equation}
in which $\Gamma = \bigsqcup_{j=1}^N \bdr{\Omega}_j$ is the \textit{disjoint} union of the boundary and all of the interfaces.

We need to solve for the charge densities, $\adj{\gamma}_j(\bm{y})$, to solve for $u(\bm{x})$. A classic theorem in potential theory~\cite{Kress_2014} characterizes the jump of $\frac{\del \slp{j}(\bm{x})}{\del \n{x}}$ across the interface,
\begin{equation}\label{interface_jump}
    \bigg(\frac{\del \slp{j}(\bm{x})}{\del \n{x}}\bigg)_{\pm}
    =
    \mp \frac{1}{2}\adj{\gamma}_j (\bm{x}) \delta_{ij} + \int_{\bdr{\Omega_j}} \frac{\partial \Phi(\bm{x}, \bm{y})}{\partial \n{x}} \adj{\gamma}_j (\bm{y}) \ d A_{\bm{y}}, \quad \forall \bm{x} \in \bdr{\Omega_j},
\end{equation}
in which $\delta_{ij}$ is the  Kronecker-delta function, and the subscript $+$ (or $-$) on the left-hand-side denotes the outward (or inward) normal direction.

We introduce the following notation that encodes the location of $\bm{x}$ and $\bm{y}$ using subscript indices $i$ and $j$. In particular,
\[
\adj{K}_{i,j}(\bm{x}, \bm{y}) \coloneqq -\frac{\del \Phi(\bm{x}, \bm{y})}{\del \n{x}} = \frac{1}{d\alpha(d)} \frac{(\bm{x} - \bm{y}) \cdot \n{x}}{|\bm{x} - \bm{y}|^d}, \quad \text{for } \bm{x} \in \bdr{\Omega_i}, \bm{y} \in \bdr{\Omega_j},
\]
in which $d$ is the dimension of the problem, and $\alpha(d) = \frac{\pi^{d/2}}{\Gamma\left(\frac{d}{2} + 1\right)}$ is the volume of the unit ball in $\R^d$.
Correspondingly, we define an integral operator $\adj{\K}_{i,j}$ such that
\[
\adj{\K}_{i,j} \adj{\gamma}_j = \bdrint{\Omega_j} \adj{K}_{i,j} (\bm{x}, \bm{y})  \adj{\gamma}_j(\bm{y}) \ d A_{\bm{y}}.
\]

Imposing the boundary condition \eqref{boundary_condition} and the interface condition \eqref{interface_condition} using Eq.~\eqref{interface_jump} yields a system of $N$ integral equations, one for each charge density function $\adj{\gamma_i}, \forall i \in \{1, \dots, N\}$,
\begin{subequations}\label{init_charge_density_sys}
\begin{empheq}[left=\empheqlbrace]{align}
\frac{1}{2} \sigma_1 \adj{\gamma}_1(\bm{x}) - \sigma_1 \sum_{j=1}^{N} \adj{\K}_{1,j} \adj{\gamma}_j &= b_1(\bm{x}), \quad \bm{x} \in \bdr{\Omega_1} \\
-\frac{1}{2} (\sigma_{p_i} + \sigma_i) \adj{\gamma}_i(\bm{x}) - (\sigma_{p_i} - \sigma_i) \sum_{j=1}^{N} \adj{\K}_{i,j} \adj{\gamma}_j &= b_i(\bm{x}), \quad 
\begin{aligned}[t]
&\text{$\bm{x} \in \bdr{\Omega_i}$}, \\
&\text{for $i = 2, 3, \dots, N$},
\end{aligned}                   
\end{empheq}
\end{subequations}
in which the region indexed by $p_i \in \{1, 2, \dots, N\}$ is the parent of region $i$, under the tree-structure definition of the domain. Rearranging the system so that the coefficients in front of $\adj{\gamma}_i$ are equal to 1, we obtain
\begin{equation}\label{unit_coeff_charge_density_sys}
    \adj{\gamma}_i(\bm{x}) - \alpha_i \sum_{j=1}^{N} \adj{\K}_{i,j} \adj{\gamma}_j = \beta_i(\bm{x}), \quad \text{$\bm{x} \in \bdr{\Omega_i}$, for $i = 1,\dots, N$}, 
\end{equation}
in which
$\alpha_i =
\begin{cases}
2 & \text{if $i = 1$} \\
\frac{2(\sigma_i - \sigma_{p_i})}{\sigma_{p_i} + \sigma_i} & \text{if $i = 2, 3, \dots, N$}
\end{cases}$
are constants, and $\beta_i: \bdr{\Omega_i} \to \R$ are functions defined as
$\beta_i(\bm{x}) = 
\begin{cases}
\frac{2}{\sigma_1} b_1(\bm{x}) & \text{if $i = 1$} \\
-\frac{2}{\sigma_{p_i} + \sigma_i} b_i(\bm{x}) & \text{if $i = 2, 3, \dots, N$}.
\end{cases}$

We can write Eq.~\eqref{unit_coeff_charge_density_sys} in matrix form,
\begin{equation}\label{matrix_form_sys}
(I_{N \times N} - \M) \adj{\bm{\gamma}} = \bm{\beta},
\end{equation}
where
\[
\M = \begin{bmatrix}
\alpha_1 \adj{\K}_{1,1} & \alpha_1 \adj{\K}_{1,2} & \hdots & \alpha_1 \adj{\K}_{1,N} \\
\alpha_2 \adj{\K}_{2,1} & \alpha_2 \adj{\K}_{2,2} & \hdots & \alpha_2 \adj{\K}_{2,N} \\
\vdots & \vdots & \ddots & \vdots \\
\alpha_N \adj{\K}_{N,1} & \alpha_N \adj{\K}_{N,2} & \hdots & \alpha_N \adj{\K}_{N,N}
\end{bmatrix}, \quad
 \adj{\bm{\gamma}} = \begin{bmatrix}
\adj{\gamma}_1 \\ \adj{\gamma}_2 \\ \vdots \\ \adj{\gamma}_N
\end{bmatrix},  \text{ and }
\bm{\beta} = \begin{bmatrix}
\beta_1 \\
\beta_2 \\
\vdots \\
\beta_N
\end{bmatrix}.
\]
We can express the solution to Eq.~\eqref{matrix_form_sys} as a Neumann series,
\begin{align}
    \adj{\bm{\gamma}} &= \sum_{i=0}^{\infty} \M^i \bm{\beta} \\
    &=
    \begin{bmatrix}
        \beta_1 + \sum_{h_0=1}^{N} (\alpha_1 \adj{\K}_{1, h_0}) \beta_{h_0} + \sum_{h_1=1}^{N}\sum_{h_0=1}^{N} (\alpha_1 \adj{\K}_{1, h_1}) (\alpha_{h_1} \adj{\K}_{h_1, h_0}) \beta_{h_0} + \dots \\
        \beta_2 + \sum_{h_0=1}^{N} (\alpha_2 \adj{\K}_{2, h_0}) \beta_{h_0} + \sum_{h_1=1}^{N}\sum_{h_0=1}^{N} (\alpha_2 \adj{\K}_{2, h_1}) (\alpha_{h_1} \adj{\K}_{h_1, h_0}) \beta_{h_0} + \dots\\
        \vdots \\
        \beta_N + \sum_{h_0=1}^{N} (\alpha_N \adj{\K}_{N, h_0}) \beta_{h_0} + \sum_{h_1=1}^{N}\sum_{h_0=1}^{N} (\alpha_N \adj{\K}_{N, h_1}) (\alpha_{h_1} \adj{\K}_{h_1, h_0}) \beta_{h_0} + \dots
    \end{bmatrix} \notag \\
    &= 
    \begin{bmatrix}
    \beta_1 \\
    \beta_2 \\
    \vdots \\
    \beta_N
    \end{bmatrix} + 
    \sum_{i=1}^{\infty}
    \begin{bmatrix}
    \sum_{h_{i-1}=1}^{N} \dots \sum_{h_1=1}^{N} \sum_{h_0=1}^{N} (\alpha_1\adj{\K}_{1, h_{i-1}}) \cdots (\alpha_{h_2}\adj{\K}_{h_2, h_1} ) (\alpha_{h_1}\adj{\K}_{h_1, h_0}) \beta_{h_0} \\
    \sum_{h_{i-1}=1}^{N} \dots \sum_{h_1=1}^{N} \sum_{h_0=1}^{N} (\alpha_2\adj{\K}_{2, h_{i-1}}) \cdots (\alpha_{h_2}\adj{\K}_{h_2, h_1} ) (\alpha_{h_1}\adj{\K}_{h_1, h_0}) \beta_{h_0} \\
    \\ 
    \vdots
    \\
    \sum_{h_{i-1}=1}^{N} \dots \sum_{h_1=1}^{N} \sum_{h_0=1}^{N} (\alpha_N\adj{\K}_{N, h_{i-1}}) \cdots (\alpha_{h_2}\adj{\K}_{h_2, h_1} ) (\alpha_{h_1}\adj{\K}_{h_1, h_0}) \beta_{h_0} \\
    \end{bmatrix}.
\end{align}
As for the WoB estimator, we obtain a finite-term approximation of $\adj{\bm{\gamma}}$ using a truncation method~\cite{Sabelfeld_Simonov_2016} as in the WoB estimator. For some fixed $0 \leq M < \infty$,
\begin{align}\label{approx_vector_form_gamma}
\adj{\bm{\gamma}} \approx 
\begin{bmatrix}
\beta_1 \\
\beta_2 \\
\vdots \\
\beta_N
\end{bmatrix} &+ 
\sum_{i=1}^{M-1}
\begin{bmatrix}
\sum_{h_{i-1}=1}^{N} \dots \sum_{h_1=1}^{N} \sum_{h_0=1}^{N} (\alpha_1\adj{\K}_{1, h_{i-1}}) \cdots (\alpha_{h_2}\adj{\K}_{h_2, h_1} ) (\alpha_{h_1}\adj{\K}_{h_1, h_0}) \beta_{h_0} \\
\sum_{h_{i-1}=1}^{N} \dots \sum_{h_1=1}^{N} \sum_{h_0=1}^{N} (\alpha_2\adj{\K}_{2, h_{i-1}}) \cdots (\alpha_{h_2}\adj{\K}_{h_2, h_1} ) (\alpha_{h_1}\adj{\K}_{h_1, h_0}) \beta_{h_0} \\
\\ 
\vdots
\\
\sum_{h_{i-1}=1}^{N} \dots \sum_{h_1=1}^{N} \sum_{h_0=1}^{N} (\alpha_N\adj{\K}_{N, h_{i-1}}) \cdots (\alpha_{h_2}\adj{\K}_{h_2, h_1} ) (\alpha_{h_1}\adj{\K}_{h_1, h_0}) \beta_{h_0} \\
\end{bmatrix} \notag \\
&+ \frac{1}{2}
\begin{bmatrix}
\sum_{h_{M-1}=1}^{N} \dots \sum_{h_1=1}^{N} \sum_{h_0=1}^{N} (\alpha_1\adj{\K}_{1, h_{M-1}}) \cdots (\alpha_{h_2}\adj{\K}_{h_2, h_1} ) (\alpha_{h_1}\adj{\K}_{h_1, h_0}) \beta_{h_0} \\
\sum_{h_{M-1}=1}^{N} \dots \sum_{h_1=1}^{N} \sum_{h_0=1}^{N} (\alpha_2\adj{\K}_{2, h_{M-1}}) \cdots (\alpha_{h_2}\adj{\K}_{h_2, h_1} ) (\alpha_{h_1}\adj{\K}_{h_1, h_0}) \beta_{h_0} \\
\\ 
\vdots
\\
\sum_{h_{M-1}=1}^{N} \dots \sum_{h_1=1}^{N} \sum_{h_0=1}^{N} (\alpha_N\adj{\K}_{N, h_{M-1}}) \cdots (\alpha_{h_2}\adj{\K}_{h_2, h_1} ) (\alpha_{h_1}\adj{\K}_{h_1, h_0}) \beta_{h_0} \\
\end{bmatrix}.
\end{align}
Substitute Eq.~\eqref{approx_vector_form_gamma} into Eq.~\eqref{interface_problem_soln_sum_of_density}, the solution to the interface problem at some query point $\bm{x}$ can be estimated by
\begin{align}
    u(\bm{x})
    &\approx 
    \begin{aligned}[t]
        & \int_{\Gamma} \Phi(\bm{x}, \bm{y}) 
        \left( 
        \sum_{i=0}^{M} w_i 
        \sum_{h_{i}, \dots, h_1, h_0=1}^{N}
        (\alpha_{h_i} \adj{\K}_{h_i, h_{i-1}}) \cdots (\alpha_{h_2}\adj{\K}_{h_2, h_1} ) (\alpha_{h_1}\adj{\K}_{h_1, h_0}) \beta_{h_0}
        \mathds{1}_{\bdr{\Omega}_{h_i}} 
        \right) \ d A_{\bm{y}}.
    \end{aligned} \notag
\end{align}
The weights are $w_i = 1 -\frac{1}{2} \delta_{i,M}$. Exchanging sums and integrals gives
\begin{align}\label{finite_sum_approx_ux}
    u(\bm{x})
    &\approx
    \begin{aligned}[t]
        & \sum_{i=0}^{M} w_i
        \sum_{h_{i}, \dots, h_1, h_0=1}^{N}
        \int_{\bdr{\Omega}_{h_i}} \Phi(\bm{x}, \bm{y}) 
        \bigg(
        (\alpha_{h_i} \adj{\K}_{h_i, h_{i-1}}) \cdots (\alpha_{h_2}\adj{\K}_{h_2, h_1} ) (\alpha_{h_1}\adj{\K}_{h_1, h_0}) \beta_{h_0}
        \bigg) \ d A_{\bm{y}}.
    \end{aligned}
\end{align}
Note that $\bigg((\alpha_{h_i} \adj{\K}_{h_i, h_{i-1}}) \cdots (\alpha_{h_2}\adj{\K}_{h_2, h_1} ) (\alpha_{h_1}\adj{\K}_{h_1, h_0}) \beta_{h_0} \bigg)$ is a function of $\bm{y}$.

\section{The Walk-on-Interfaces Method}\label{woi}
In this section, we derive and present our \textit{Walk-on-Interfaces} (WoI) estimator for the elliptic interface problem Eq.~\eqref{interface_problem}. We begin with a theorem serving as the backbone of our WoI estimators. Next, we propose a naive yet computationally intensive estimator. For better time performance, we exploit the linearity inherited in the naive estimator to develop a time-efficient alternative. Finally, we close this section by discussing the variance-reduced estimator and the gradient estimator.

\subsection{A Theorem Towards Walk-on-Interfaces (WoI) Estimator} \label{theorem_towards_woi}
Observe that Eq.~\eqref{finite_sum_approx_ux} writes $u(\bm{x})$ in a finite sum of boundary integrals.
To establish our first WoI estimator for $u(\bm{x})$, it is sufficient to develop a Monte Carlo estimator for each inner product. We summarize an estimator for a boundary integral of the form
\[
\int_{\bdr{\Omega}_{h_i}} \Phi(\bm{x}, \bm{y}) 
(
\adj{\K}_{h_i, h_{i-1}} \cdots \adj{\K}_{h_2, h_1}\adj{\K}_{h_1, h_0} \beta_{h_0}
) \ d A_{\bm{y}}
\]
in Theorem \ref{naive_woi_theorem}.

We'll first introduce a notation for random variables before stating the theorem. Let $(h_0, h_1, \dots, h_i, \dots)$ be a sequence of interface or boundary indices. Let $Y_{(h_i, \dots, h_1, h_0)} \in \bdr{\Omega_{h_i}}$ be a random variable in a Markov tree whose previous random state, $Y_{(h_{i-1}, \dots, h_1, h_0)}$, is located on $\bdr{\Omega_{h_{i-1}}}$. Let $\bm{y}_{(h_i, \dots, h_1, h_0)} \in \bdr{\Omega_{h_i}}$ be a realization of $Y_{(h_i, \dots, h_1, h_0)}$. As a shorthand notation, we denote $Y_{(h_i, \dots, h_1, h_0)}$ as $Y_i$, and $\bm{y}_{(h_i, \dots, h_1, h_0)}$ as $\bm{y}_i$, respectively. In the rest of this work, we will use the shorthand notation when the sequence of the interface indices can be inferred from equations or expressions.

With this notation, we introduce the following theorem.

\begin{theorem}\label{naive_woi_theorem}
Given a query point $\bm{x}$ and a fixed set of $h_0, h_1, \dots, h_i \in \{1, \dots, N\}$. Let $\{Y_0, Y_1, \dots, Y_i, \dots\}$ be a Markov chain starting from $\bm{x}$ such that $Y_i \in \bdr{\Omega_{h_i}}$. Let $p_0(\bm{x}, \bm{y})$ and $p(\bm{x}, \bm{y})$ be any probability density functions (PDFs) denoting the transitional distribution of $\bm{y}$ given $\bm{x}$. Define a family of functions
\begin{equation}\label{naive_woi_Q*}
    \adj{Q}_{(h_i, \dots, h_1, h_0)} = 
    \begin{cases}
        \frac{\beta_{h_0}(Y_0)}{p_0(\bm{x}, Y_0)} \quad & \text{if $i = 0$}\\
         \adj{Q}_{(h_{i-1}, \dots, h_1, h_0)} \frac{\adj{K}_{h_i, h_{i-1}}(Y_i, Y_{i-1})}{p(Y_{i-1}, Y_i)}\quad & \text{if $i \geq 1$}\\
    \end{cases}.
\end{equation}
Then, $\adj{Q}_{(h_i, \dots, h_1, h_0)} \Phi(\bm{x}, Y_i)$ is an unbiased estimator for 
\[
\int_{\bdr{\Omega}_{h_i}} \Phi(\bm{x}, \bm{y}) 
(
\adj{\K}_{h_i, h_{i-1}} \cdots \adj{\K}_{h_2, h_1}\adj{\K}_{h_1, h_0} \beta_{h_0}
) \ d A_{\bm{y}}.
\]
i.e. 
\[
\int_{\bdr{\Omega}_{h_i}} \Phi(\bm{x}, \bm{y}) 
(\adj{\K}_{h_i, h_{i-1}} \cdots \adj{\K}_{h_2, h_1}\adj{\K}_{h_1, h_0} \beta_{h_0}
) \ d A_{\bm{y}} 
= \mathbb{E} \bigg[\adj{Q}_{(h_i, \dots, h_1, h_0)} \Phi(\bm{x}, Y_i) \bigg].
\]
\end{theorem}
This theorem is proved in Appendix~\ref{Append-proof_of_naive_woi_theorem}.
\subsection{Naive Walk-on-Interfaces} \label{naive_woi}
To translate Theorem~\ref{naive_woi_theorem} into practical algorithms, we begin by specifying key choices and essential components that contribute to our basic WoI estimators. We begin by presenting our choice of transitional PDF in Section \ref{transitional_probability_function}. Then, we state the resulting naive WoI estimator with the selected transitional PDF in Section \ref{the_naive_woi_estimator}. In Section \ref{sampling_strategy} , we explain our approach to sample the transitional PDF. Finally, we conclude our discussion of the naive WoI estimator by addressing the reuse of samples and Markov chains in Section \ref{reuse_markov_chains}. 

\subsubsection{Transitional Probability Function}\label{transitional_probability_function}
The selection of appropriate probability density functions $p_0(\bm{x}, \bm{y}_0)$ and $p(\bm{y}_{i-1}, \bm{y}_i)$ is crucial for our algorithm. Ideal choices of transitional probability should achieve both efficiency, allowing easy sampling with fast algorithms, and accuracy, ensuring that the PDF closely approximates the integration kernel. While there can be many better choices, we use
\begin{align}
    p_0(\bm{x}, \bm{y}_0) &= \frac{1}{|\bdr{\Omega_{h_0}}|} \quad \text{and} \label{naive_woi_p0} \\
    p(\bm{y}_{i-1}, \bm{y}_i) &= 
    \begin{cases}
        \frac{2}{m\alpha(m)} \frac{|(\bm{y}_i - \bm{y}_{i-1}) \cdot \bm{n}(\bm{y}_i)|}{q(\bm{y}_{i-1}, \bm{y}_i)|\bm{y}_i -  \bm{y}_{i-1}|^m} & \quad \text{if $h_{i-1} = h_{i}$} \\
        \frac{1}{m\alpha(m)} \frac{|(\bm{y}_i - \bm{y}_{i-1}) \cdot \bm{n}(\bm{y}_i)|}{q(\bm{y}_{i-1}, \bm{y}_i)|\bm{y}_i - \bm{y}_{i-1}|^m} & \quad \text{if $h_{i-1}$ is a descendant of $h_{i}$} \\
        \frac{1}{|\bdr{\Omega_{h_i}}|} & \quad \text{otherwise}
    \end{cases} \label{naive_woi_p}.
\end{align}
in all our implementations, where $q(\bm{y}_{i-1}, \bm{y}_i)$ is again the number of intersections between $\bdr{\Omega}_i$ and a line (or a ray) passing through $\bm{y}_{i-1}$ and $\bm{y}_{i}$. We will further explain this idea in Section \ref{sampling_strategy}.

\subsubsection{The Naive WoI Estimator}\label{the_naive_woi_estimator}
Applying Theorem \ref{naive_woi_theorem} to Eq.~\eqref{finite_sum_approx_ux}, we can write $u(\bm{x})$ by sums of expected values
\begin{align}
    u(\bm{x}) 
    &\approx
    \begin{aligned}[t]
        &\sum_{i=0}^{M} w_i 
        \sum_{h_{i}, \dots, h_1, h_0=1}^{N}
        \left(\prod_{m=1}^i \alpha_{h_m} \right)
        \mathbb{E} \left[
        \adj{Q}_{(h_i, \dots, h_1, h_0)}\Phi(x, Y_i)
        \right],
    \end{aligned}
\end{align}
where $\adj{Q}_{(h_i, \dots, h_1, h_0)}$ is defined in Eq.~\eqref{naive_woi_Q*}, with its corresponding $p_0(\bm{x}, \bm{y})$ and $p(\bm{x}, \bm{y})$ defined in Eq.~\eqref{naive_woi_p0} and Eq.~\eqref{naive_woi_p}, respectively. 

Notice that we can write part of the transitional probability $p(\bm{y}_{i-1}, \bm{y}_i)$ as multiples of $\adj{K}_{h_i, h_{i-1}}(\bm{y}_{i-1}, \bm{y}_i) = \frac{1}{m\alpha(m)} \frac{(\bm{y}_{i-1} - \bm{y}_i) \cdot \bm{n}(\bm{y}_{i-1})}{|\bm{y}_{i-1} - \bm{y}_i|^m}$
\begin{align}
    p(\bm{y}_{i-1}, \bm{y}_i) 
    &= 
    \begin{cases}
        \frac{2 \sign[(\bm{y}_i - \bm{y}_{i-1}) \cdot \bm{n}(\bm{y}_i)]}{q(\bm{y}_{i-1}, \bm{y}_i)} \adj{K}_{h_i, h_{i-1}}(\bm{y}_{i-1}, \bm{y}_i) & \quad \text{if $h_{i-1} = h_{i}$} \\
        \frac{\sign[(\bm{y}_i - \bm{y}_{i-1}) \cdot \bm{n}(\bm{y}_i)]}{q(\bm{y}_{i-1}, \bm{y}_i)} \adj{K}_{h_i, h_{i-1}}(\bm{y}_{i-1}, \bm{y}_i) & \quad \text{if $h_{i-1}$ is a descendant of $h_{i}$} \\
        \frac{1}{|\bdr{\Omega_{h_i}}|} & \quad \text{otherwise}
    \end{cases}.
\end{align}
This observation allows us to decompose $\adj{Q}_{(h_i, \dots, h_1, h_0)}$ into a product of a straightforward component, 
\begin{equation}\label{naive_woi_Q}
    Q_{(h_i, \dots, h_1, h_0)} = 
    \begin{cases}
        \frac{\beta_{h_0}(Y_0)}{p_0(\bm{x}, Y_0)}  \quad & \text{if $i = 0$} \\
        \frac{1}{2} Q_{(h_{i-1}, \dots, h_1, h_0)}
        \quad & \text{if $i \geq 1$ and $h_{i-1} = h_i$} \\
        Q_{(h_{i-1}, \dots, h_1, h_0)}  
        \quad & \text{if $i \geq 1$ and $h_{i-1}$ is a descendant of $h_i$} \\
        \frac{\adj{K}_{h_i, h_{i-1}}(Y_i, Y_{i-1})}{p(Y_{i-1}, Y_i)} Q_{(h_{i-1}, \dots, h_1, h_0)}
        \quad & \text{if $i \geq 1$ and for any other $h_{i-1}$ and $h_i$}
    \end{cases},
\end{equation}
and another more complicated component,
\begin{equation}\label{naive_woi_T}
    T_{(h_i, \dots, h_1, h_0)} = 
    \begin{cases}
        1 \quad & \text{if $i = 0$} \\
       \prod_{m=1}^{i} q(Y_{m-1}, Y_m) \sign[(Y_m - Y_{m-1}) \cdot \bm{n}(Y_m)]  & \text{if $i \geq 1$}
    \end{cases},
\end{equation}
so that we don't need to evaluate $\adj{K}_{h_i, h_{i-1}}(Y_i, Y_{i-1})$ in every update.

Hence, substituting the construction $\adj{Q}_{(h_i, \dots, h_1, h_0)} = Q_{(h_i, \dots, h_1, h_0)} T_{(h_i, \dots, h_1, h_0)}$ gives
\begin{align}\label{naive_woi_estimator_sum_of_expectation}
    u(\bm{x}) 
    &\approx
    \begin{aligned}[t]
        &\sum_{i=0}^{M} w_i
        \sum_{h_{i}, \dots, h_1, h_0=1}^{N}
        \left(\prod_{m=1}^i \alpha_{h_m} \right)
        \mathbb{E} \left[
        Q_{(h_i, \dots, h_1, h_0)} T_{(h_i, \dots, h_1, h_0)} \Phi(x, Y_i)
        \right],
    \end{aligned}
\end{align}
where
$Q_{(h_i, \dots, h_1, h_0)}$ and $T_{(h_i, \dots, h_1, h_0)}$ are defined in Eq.~\eqref{naive_woi_Q} and Eq.~\eqref{naive_woi_T}, respectively. By linearity of expected value, the naive WoI estimator can be summarized in one single expected value,
\begin{align}\label{naive_woi_estimator}
    u(\bm{x}) 
    &\approx
    \mathbb{E}
        \left[
        \sum_{i=0}^{M} w_i
        \sum_{h_{i}, \dots, h_1, h_0=1}^{N}
        \left(\prod_{m=1}^i \alpha_{h_m} \right)
        Q_{(h_i, \dots, h_1, h_0)} T_{(h_i, \dots, h_1, h_0)} \Phi(x, Y_i) \right].
\end{align}

\subsubsection{Sampling Strategy} \label{sampling_strategy}
We now discuss our motivation for designing the transitional probabilities outlined in Eq.~\eqref{naive_woi_p0} and Eq.~\eqref{naive_woi_p}, and describe our method to sample $\bm{y}_{i+1}$ given $\bm{y}_0, \bm{y}_1, \dots, \bm{y}_i$.

The transitional probability $p_0(\bm{x}, \bm{y}_0)$ is proposed to estimate the integral
\[
\int_{\bdr{\Omega_{h_0}}} \Phi(\bm{x}, \bm{y}) \beta_{h_0}(\bm{y}) \ d A_{\bm{y}},
\]
whose kernel is $\Phi(\bm{x}, \bm{y})$. Since $\Phi(\bm{x}, \bm{y})$ is not proportional to any easy-to-sample distributions defined on  $\bdr{\Omega_{h_0}}$, we simply use a uniform distribution over $\bdr{\Omega_{h_0}}$ for efficiency. Therefore, we sample a random point uniformly on $\bdr{\Omega_{h_0}}$ to obtain $\bm{y}_0$.

We define $p(\bm{x}, \bm{y})$ as the transitional probability to estimate integrals of the form
\[
\int_{\bdr{\Omega}_{h_i}} \Phi(\bm{x}, \bm{y}) 
(
\adj{\K}_{h_i, h_{i-1}} \cdots \adj{\K}_{h_2, h_1}\adj{\K}_{h_1, h_0} \beta_{h_0}
) \ d A_{\bm{y}},
\]
given $\bm{y}_0, \bm{y}_1, \dots, \bm{y}_{i-1}$, and when $i \geq 1$.
Let $f = \adj{\K}_{h_{i-1}, h_{i-2}}\adj{\K}_{h_1, h_0} \beta_{h_0}$, then the target integral can be written as
\begin{align}\label{woi_p_rearrange_target}
    \int_{\bdr{\Omega}_{h_i}} \Phi(\bm{x}, \bm{y}) 
    (\adj{\K}_{h_i, h_{i-1}} f) \ d A_{\bm{y}}
    &=
    \int_{\bdr{\Omega}_{h_i}}  \Phi(\bm{x}, \bm{y}_i) 
    \int_{\bdr{\Omega}_{h_{i-1}}}
    \adj{K}_{h_i, h_{i-1}}(\bm{y}_i, \bm{y}_{i-1}) f(\bm{y}_{i-1}) \ d A_{\bm{y}_{i-1}} \ d A_{\bm{y}_i} \notag \\
    &= \int_{\bdr{\Omega}_{h_{i-1}}} f(\bm{y}_{i-1}) \int_{\bdr{\Omega}_{h_i}} \Phi(\bm{x}, \bm{y}_i) \adj{K}_{h_i, h_{i-1}}(\bm{y}_i, \bm{y}_{i-1})
    \ d A_{\bm{y}_i} \ d A_{\bm{y}_{i-1}}.
\end{align}
It suffices to estimate $f$ from the given points $\bm{y}_0, \bm{y}_1, \dots, \bm{y}_{i-1}$, which leaves the inner integral,
\[
\int_{\bdr{\Omega}_{h_i}} \Phi(\bm{x}, \bm{y}_i) \adj{K}_{h_i, h_{i-1}}(\bm{y}_i, \bm{y}_{i-1})
\ d A_{\bm{y}_i},
\]
to be the only known part in Eq.~\eqref{woi_p_rearrange_target}.

We select the isotropic distribution, which is proportional to $\adj{K}_{h_i, h_{i-1}}(\bm{y}_{i}, \bm{y}_{i-1})$. Also, $h_{i-1}$ is a descendant of $h_{i}$ (i.e. $\Omega_{h_{i-1}} \sse \Omega_{h_i}$), which promotes accuracy. In these cases, we obtain $\bm{y}_i$ by first sampling a random direction $\bm{d}$, which defines a ray, $r(t) = t_{\geq 0} \bm{d} + \bm{y}_{i-1}$, originating at $\bm{y}_{i-1}$. Then we pick any intersection uniformly randomly between $r(t)$ and $\bdr{\Omega_i}$  to be $\bm{y}_i$. 

For a similar reason, we also use the isotropic distribution for the cases when $h_{i-1} = h_{i}$ (i.e. $\Omega_{h_{i-1}} = \Omega_{h_i}$). Notice that sampling a random ray this time does not guarantee the existence of an intersection between $r(t)$ and $\bdr{\Omega_i}$. So, we must experience a time-wasting rejection sampling of ray direction before obtaining $\bm{y}_i$.
As an alternative~\cite{Sabelfeld_Simonov_2016}, we sample a line, $\ell(t) = t\bm{d} + \bm{y}_{i-1}$, in a random direction $\bm{d}$ originating at $\bm{y}_{i-1}$, and set $\bm{y}_i$ to be an arbitrary intersection between $r(t)$ and $\bdr{\Omega}_{h_i}$.

However, for any other relation between $h_{i-1}$ and $h_{i}$, it is not easy to sample $Y_i \in \bdr{\Omega_{h_i}}$ under an isotropic distribution, since sampling a line or a ray following the same procedure as previous does not always guarantee a hit on $\bdr{\Omega_i}$. Rejection sampling is time-consuming, especially when $\Omega_{h_i}$ only occupies a tiny portion of the domain $\Omega$. Therefore, even though the isotropic distribution represents $\adj{K}_{h_i, h_{i-1}}(\bm{y}_i, \bm{y}_{i-1})$ precisely, we still use a uniform distribution over $\bdr{\Omega_{h_i}}$ for the purpose of efficiency.

\subsubsection{Reuse of Markov Chains}\label{reuse_markov_chains}
For each $i \in \{1, \dots, M\}$ and for each tuple of $(h_0, h_1, \dots, h_i)$, the naive WoI estimator Eq.~\eqref{naive_woi_estimator} requires sampling the Markov chain, $\{Y_0, Y_1, \dots, Y_i\}$, to compute the expected value. While it is possible to initialize new Markov chains for each expected value, this approach is suboptimal in both memory storage and computational efficiency. The WoB estimator can achieve better time performance by reusing Markov chains generated for terms with $Q_{(h_i, \dots, h_1, h_0)}$ as part of the Markov chain for terms with $Q_{(h_j \dots h_i, \dots, h_1, h_0)}$, wherever $j > i$. 

Let us fix a tuple $(h_0, h_1, \dots, h_i) \in \{1, 2, \dots, N\}^i$ and assume $\{\bm{y}_0, \bm{y}_1, \dots, \bm{y}_i\}$ to be a realization of the Markov chain for the term 
\[\mathbb{E} \bigg[Q_{(h_i, \dots, h_1, h_0)} T_{(h_i, \dots, h_1, h_0)}\Phi(x, Y_i)\bigg].
\]
Instead of sampling the Markov chain $\{Y_0, Y_1, \dots, Y_{i+1}\}$ with brand new samples, the next point $\bm{y}_{i+1} \in \bdr{\Omega_{h_{i+1}}}$ can be sampled under the transitional probability $p(\bm{y}_{i}, \bm{y}_{i+1})$ so that $\{\bm{y}_0, \dots, \bm{y}_i, \bm{y}_{i+1}\}$ is a set of samples for the term \[
\mathbb{E} \bigg[Q_{(h_{i+1}, h_i, \dots, h_1, h_0)}T_{(h_{i+1}, h_i, \dots, h_1, h_0)}\Phi(x, Y_{i+1})\bigg],
\]
for any $h_{i+1} \in \{1, 2, \dots, N\}$. 

\begin{figure}[!ht]
    \centering
    \begin{subfigure}[b]{0.65\columnwidth}
    \includegraphics[width=\textwidth]{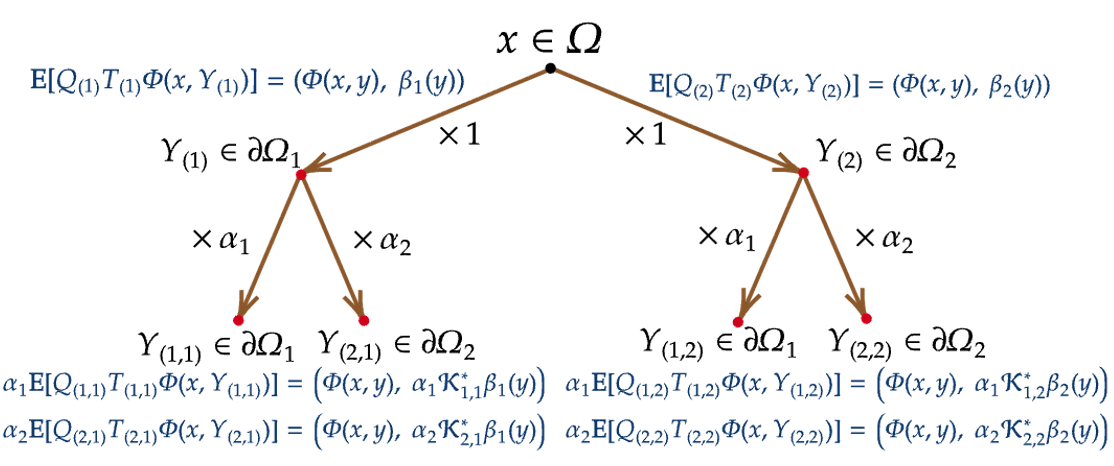}
    \caption{}
    \label{naive_woi_tree}
    \end{subfigure}
    \begin{subfigure}[b]{0.25\columnwidth}
    \includegraphics[width=\textwidth]{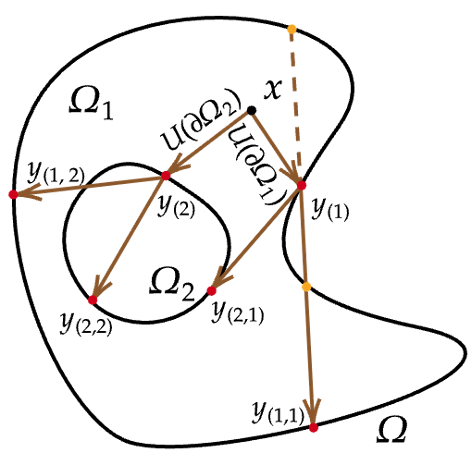}\caption{}\label{a_realization_of_naive_woi_tree}
    \end{subfigure}
    \caption{
    (a) The Markov tree, in which each branch is a Markov chain. We explicitly state which expected value approximates which boundary integral term in $u(\bm{x})$. Note that the boundary integrals are now written in the notation of $L_2$ inner products. The random variables are marked with full subscripts. (b) A realization of the Markov tree, where we reuse existing Markov chains. The red dots are selected to make statistical contribution to the expected value. The dashed line is the line $\bm{\ell}(t)$ we sampled, and the yellow dots are the unselected intersections. $\bm{y}_{(1)}$ and $\bm{y}_{(2)}$ are uniformly sampled on $\bdr{\Omega_1}$ and $\bdr{\Omega_2}$, respectively, since $p_0(\bm{x}, \bm{y}_{h_0})$ is uniform. Given $\bm{y}_{(1)}$, $\bm{y}_{(1,1)}$ is sampled under isotropic distribution, since $h_1 = h_0 = 1$, while $\bm{y}_{(2,1)}$ is sampled uniformly on $\bdr{\Omega_2}$, as $h_0 = 1$ is not a descendant of $h_1 = 2$. Given $\bm{y}_{(2)}$, $\bm{y}_{(1,2)}$ and $\bm{y}_{(2,2)}$ are sampled on $\bdr{\Omega_1}$ and $\bdr{\Omega_2}$, respectively, under isotropic distribution.} 
    \label{naive_woi_tree_with_realization}
\end{figure}

Therefore, the points we need to sample for each observation of the expectation form a $N$-ary tree (with root $\bm{x}$) of height $M + 2$, where each branch of the tree is a Markov chain. We refer to such a tree as \textit{Markov tree}. Figure \ref{naive_woi_tree_with_realization} provides one toy example with $M = 1$, $N = 2$, for which
\begin{align} \label{toy_example_u}
    u(\bm{x})
    &\approx
    \begin{aligned}[t]
        &\int_{\bdr{\Omega_1}} \Phi(\bm{x}, \bm{y}) \beta_1(\bm{y})  \ d A_{\bm{y}} + \int_{\bdr{\Omega_2}} \Phi(\bm{x}, \bm{y}) \beta_2(\bm{y}) \ d A_{\bm{y}}  + \\
        &\frac{1}{2} \int_{\bdr{\Omega_1}} \Phi(\bm{x}, \bm{y}) (\alpha_1 \adj{\K}_{1, 1}\beta_1)  \ d A_{\bm{y}} + 
        \frac{1}{2} \int_{\bdr{\Omega_2}} \Phi(\bm{x}, \bm{y}) (\alpha_2 \adj{\K}_{2, 1}\beta_1)  \ d A_{\bm{y}} + \\
        &\frac{1}{2} \int_{\bdr{\Omega_1}} \Phi(\bm{x}, \bm{y}) (\alpha_1 \adj{\K}_{1, 2}\beta_2)  \ d A_{\bm{y}} + 
        \frac{1}{2} \int_{\bdr{\Omega_2}} \Phi(\bm{x}, \bm{y}) (\alpha_2 \adj{\K}_{2, 2}\beta_2)  \ d A_{\bm{y}},
    \end{aligned} \notag 
\end{align}
as a special case of Eq.~\eqref{finite_sum_approx_ux}. The binary Markov tree of height three and an example of its realization over a nonconvex domain are visualized.
\subsection{Acceleration} \label{acceleration}
Generating samples for the naive WoI requires traversing the $N$-ary Markov tree, which yields exponential time complexity. In each traversal, the number of random variables we need to sample is $\sum_{i=1}^M N^i$. If we sample $\T$ Markov trees for the expected value in Eq.~\eqref{naive_woi_estimator}, the naive WoI estimator relies on $\T \sum_{i=0}^{M} N^i$ points on interfaces to complete the estimation. This quantity grows rapidly with increasing $N$ and escalates exponentially with respect to $M$. Therefore, it is desirable to avoid the tree traversal to maintain a reasonable estimation time while keeping the quality of the estimator. 

\subsubsection{The WoI Estimator}\label{the_woi_estimator}
We avoid tree traversal by leveraging the coefficient, $\alpha_{h_i}$, in the naive WoI estimator, which is formalized in Theorem \ref{woi_theorem} and proven in Appendix \ref{Append-proof_of_woi_theorem}.

\begin{theorem}\label{woi_theorem}
    Let $H_0, H_1, \dots, H_i$ be independent discrete random variables with PDF
    \begin{equation}\label{scheduled_interface_pdf}
        \prob(H_i = n) = 
            \begin{cases}
                \frac{1}{N} & \quad \text{$\forall n \in \{1, \ldots, N\}$, when $i = 0$} \\
                \frac{|\alpha_n|}{||\bm{\alpha}||_1} & \quad \text{$\forall n \in \{1, \ldots, N\}$, when $i \geq 1$}
            \end{cases},
    \end{equation}
    where $||\cdot||_1$ is the $\ell_1$ norm. It can be shown that 
    \begin{align}
        &\sum_{h_{i}, \dots, h_1, h_0=1}^{N}
        \int_{\bdr{\Omega}_{h_i}} \Phi(\bm{x}, \bm{y}) 
        \bigg(
        (\alpha_{h_i} \adj{\K}_{h_i, h_{i-1}}) \cdots (\alpha_{h_2}\adj{\K}_{h_2, h_1} ) (\alpha_{h_1}\adj{\K}_{h_1, h_0}) \beta_{h_0}
        \bigg) \ d A_{\bm{y}} \notag \\
        =
        & N ||\bm{\alpha}||_1^i
        \mathbb{E}_{H_i, \dots, H_0} \bigg[
        \prod_{m=1}^i \sign(\alpha_{H_m})
         \int_{\bdr{\Omega}_{H_i}}
         \Phi(\bm{x}, \bm{y})
         \bigg(\adj{\K}_{H_i, H_{i-1}} \dots \adj{\K}_{H_1, H_0} \beta_{H_0} d A_{\bm{y}}
        \bigg) \bigg]. \notag
    \end{align}
\end{theorem}

Applying the result in Theorem \ref{woi_theorem} to Eq.~\eqref{finite_sum_approx_ux}, we have
\begin{equation}
    u(\bm{x}) 
    \approx
         \sum_{i=0}^{M} w_i N ||\bm{\alpha}||_1^i
        \mathbb{E}_{H_i, \dots, H_0} \left[
        \prod_{m=1}^i \sign(\alpha_{H_m})
         \int_{\bdr{\Omega}_{H_i}}
         \Phi(\bm{x}, \bm{y})
         \left(\adj{\K}_{H_i, H_{i-1}} \dots \adj{\K}_{H_1, H_0} \beta_{H_0} d A_{\bm{y}}
        \right) \right]. \notag
\end{equation}
Then by Theorem \ref{naive_woi_theorem} again, the WoI estimator is
\begin{equation}\label{woi_estimator}
    u(\bm{x}) 
    \approx
    \sum_{i=0}^{M} w_i
        N ||\bm{\alpha}||_1^i \mathbb{E}_{H_i, \dots, H_0} \left[
        \prod_{m=1}^i \sign(\alpha_{H_m})
         \mathbb{E}_{Y_i} \left[Q_{(H_i, \dots H_1, H_0)} T_{(H_i, \dots H_1, H_0)} \Phi(\bm{x}, Y_i) \right] \right],
\end{equation}
for which we substituted $Q_{(H_i, \dots, H_1, H_0)} T_{(H_i, \dots, H_1, H_0)}$ for $\adj{Q}_{(H_i, \dots, H_1, H_0)}$.

\subsubsection{A Summary of the WoI Procedure}
Based on our design of the WoI estimator, we outline the process for estimating $u(\bm{x})$ as presented in Eq.~\eqref{woi_estimator}. After fixing $M$, we first sample $\s$ sets of discrete random variables $H_0, H_1, \dots, H_M$ under the distribution defined in Eq.~\eqref{scheduled_interface_pdf}
for the outer expected value, $\mathbb{E}_{H_M, \dots H_1, H_0}$.
By doing so, we \textit{schedule} the interface or the boundary where the random variables $Y_0, Y_1, \dots Y_M$ are located. Then, for each set of scheduled locations, $\{h_0, h_1, \dots, h_M\}$, we sample $\W$ \textit{walkers} that walk $M + 1$ steps each. The \textit{trajectory} of each walker is a Markov chain $\{Y_0, Y_1, \dots Y_M\}$ such that $Y_i \in \bdr{\Omega_{h_i}}$. The transitional probabilities for the Markov chain $\{Y_0, Y_1, \dots Y_M\}$ with is defined in Eq.~\eqref{naive_woi_p0} and Eq.~\eqref{naive_woi_p}, and the sampling strategy for $Y_i$'s is detailed in Section \ref{sampling_strategy}. 

Hence, we compute the inner expected value, $\mathbb{E}_{Y_i}[\cdot]$, to complete the estimation of $u(\bm{x})$. We present the pseudocode of the WoI estimator in Algorithm \ref{alg_WoI_wrapper} that computes the expected value, and Algorithm \ref{alg_WoI} that samples the states, $Y_i$'s.
\begin{algorithm}[!ht]
\caption{Wrapper Function for WoI Estimator}
\label{alg_WoI_wrapper}
\vspace{-0.2cm}
\begin{flushleft}
\textbf{Input:} Query point $\bm{x}$, number of interfaces $N$, problem dimension $d$,
\\ \hspace*{1.15cm}
integration kernel $\adj{K}_{i,j}(\bm{x}, \bm{y}) = \frac{1}{d\alpha(d)} \frac{(\bm{x} - \bm{y}) \cdot \n{x}}{|\bm{x} - \bm{y}|^m}$, where $m$ is the dimension, 
\\ \hspace*{1.15cm}
constant array $\bm{\alpha} = \begin{bmatrix} \alpha_1 & \alpha_2 & \cdots & \alpha_N \end{bmatrix}$,
\\ \hspace*{1.15cm}
boundary and interface conditions $\bm{\beta} = \begin{bmatrix} \beta_1 & \beta_2 & \cdots & \beta_N \end{bmatrix}$,
\\ \hspace*{1.15cm}
fundamental solution $\Phi(\bm{x}, \bm{y})$,
\\ \hspace*{1.15cm}
$\s$ sets of interface schedules, $\W$ walkers, $M + 1$ steps for each walker. \\
\textbf{Output:} Estimated solution $u(\bm{x})$.
\end{flushleft}
\vspace{-0.4cm}
\begin{algorithmic}[1]
\Function{WoIEstimator}{$\bm{x}$, $\adj{K}_{i,j}(\bm{x}, \bm{y})$, $\bm{\beta}$, $\Phi(\bm{x}, \bm{y})$, $\bm{\alpha}$, $S$, $M$, $\C$}
\State  $u(\bm{x}) = 0$
\For{$s = 1$ to $\s$}
\State $S \gets [h_0, h_1, \dots, h_M]$
\Comment{sample $H_0, H_1, \dots H_M$ to obtain the $s^\text{th}$ interface schedule}
\State $u_s = 0$
\For{$w = 1$ to $\W$}\Comment{walkers start walking}
\State $u_w$ = \textsc{WoI}($\bm{x}$, $\adj{K}_{i,j}(\bm{x}, \bm{y})$,  $\Phi(\bm{x}, \bm{y})$, $\bm{\alpha}$, $\bm{\beta}$, $M$, $S$)
\Comment{Algorithm \ref{alg_WoI}}
\State $u_s = u_s + (u_w - u_s) / w$ \Comment{Welfords algorithm for $\mathbb{E}_{Y_i}$}
\EndFor
\State $u(\bm{x}) = u(\bm{x}) + (u_s - u(\bm{x})) / s$ \Comment{Welfords algorithm for expectation for $\mathbb{E}_{H_i, \dots H_0}$}
\EndFor
\State \textbf{return} $N * u(\bm{x})$
\EndFunction
\end{algorithmic}
\end{algorithm}
\begin{algorithm}[!ht]
\caption{The WoI Estimator}
\label{alg_WoI}
\begin{algorithmic}[1]
\Function{WoI}{$\bm{x}$, $\adj{K}_{i,j}(\bm{x}, \bm{y})$,  $\Phi(\bm{x}, \bm{y})$, $\bm{\alpha}$, $\bm{\beta}$, $M$, $S$}
\State $u_w \gets 0$, $T \gets 1$, $\texttt{weight} \gets 1$.
\For{$i = 0$ to $M$}
\If{$i == M$}
\State $\texttt{weight} \gets 0.5$ \Comment{truncation}
\EndIf
\State $h_i \gets S(i)$
\If{$i == 0$} \Comment{transitional PDF is $p_0(\bm{x}, \bm{y})$}
\State $\bm{y}_i \gets \text{Unif}(\bdr{\Omega_{h_i}})$
\State $Q = \beta_{h_i}(\bm{y}_i) / \frac{1}{|\bdr{\Omega_{h_i}}|}$
\State $u_w = u_w + \texttt{weight} * ||\alpha||_1^{i} * Q * T * \Phi(\bm{x}, \bm{y}_i)$
\Else\Comment{transitional PDF is $p(\bm{x}, \bm{y})$}
\State $h_{i-1} \gets S(i-1)$ \Comment{$i \geq 1$}
\If{$h_{i-1} = h_i$}
\State $\bm{d} \gets \text{an isotropic direction}$
\State $q(\bm{y}_{i-1}, \bm{y}_i) \gets$ \# of intersections between $\bdr{\Omega_{h_i}}$ and a line $\ell(t) = t\bm{d} + \bm{y}_{i-1}$
\State $\bm{y}_i \gets$ any of the $q(\bm{y}_{i-1}, \bm{y}_i)$ intersections
\State $Q \gets Q * 0.5$
\ElsIf{$h_{i-1}$ a descendant of $h_i$}
\State $\bm{d} \gets \text{an isotropic direction}$
\State $q(\bm{y}_{i-1}, \bm{y}_i) \gets$ \# of intersections between $\bdr{\Omega_{h_i}}$ and a ray $r(t) = t_{\geq 0} \bm{d} + \bm{y}_{i-1}$
\State $\bm{y}_i \gets$ any of the $q(\bm{y}_{i-1}, \bm{y}_i)$ intersections
\State $Q \gets Q * 1$
\Else
\State $q(\bm{y_{i-1}}, \bm{y}_i) \gets 1$
\State $\bm{y}_i \gets \text{Unif}(\bdr{\Omega_{h_i}})$
\State $Q = Q * \adj{K}_{h_i, h_{i-1}}(\bm{y}_i, \bm{y}_{i-1})/ \frac{1}{|\bdr{\Omega_{h_i}}|}$ 
\EndIf
\State $T = T * q(\bm{y_{i-1}}, \bm{y}_i) * \sign[(\bm{y}_i - \bm{y}_{i-1})\cdot \bm{n}(\bm{y}_i)]$
\State $u_w = u_w + \texttt{weight} * ||\alpha||_1^{i} * [\prod_{m=0}^{i} \sign(\alpha_{h_m})] * Q * T * \Phi(\bm{x}, \bm{y}_i)$
\EndIf
\EndFor
\State return $u_w$
\EndFunction
\end{algorithmic}
\end{algorithm}

\begin{figure}[!ht]
    \centering
    \begin{subfigure}[b]{0.5\columnwidth}
    \includegraphics[width=\textwidth]{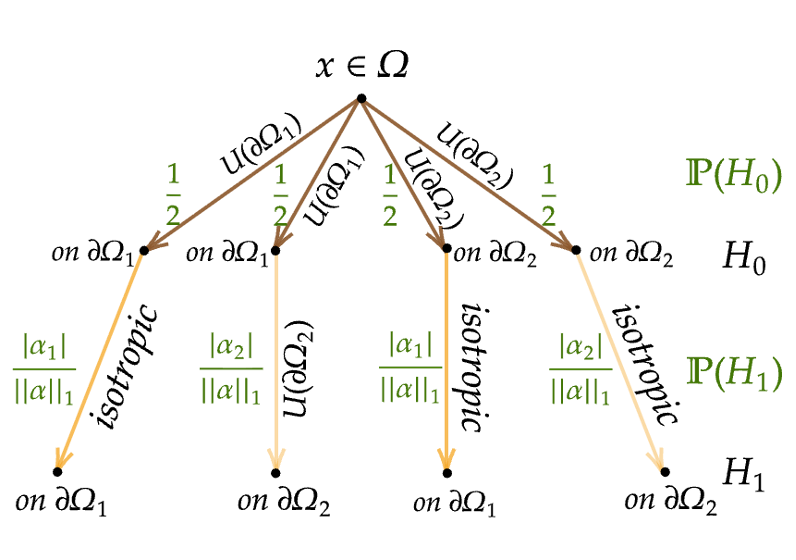}
    \caption{}
    \label{woi_tree}
    \end{subfigure}
    \quad \quad
    \begin{subfigure}[b]{0.3\columnwidth}
    \includegraphics[width=\textwidth]{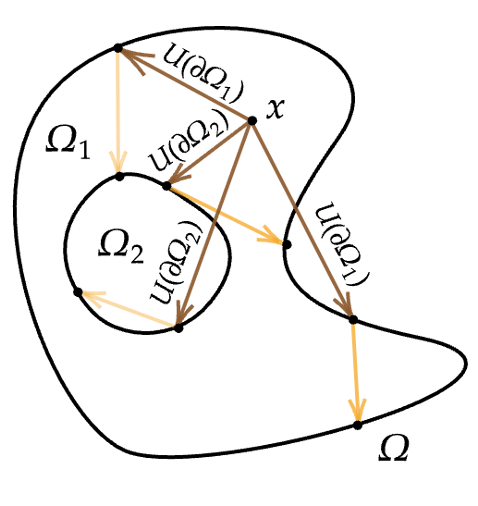}
    \caption{}
    \label{woi_trajectory}
    \end{subfigure}
    \caption{Consider an domain with one interface ($N=2$) and each walker walks two steps $(M = 1)$. Assume $|\alpha_1| > |\alpha_2|$. Then, $\bdr{\Omega_1}$ is more likely to be scheduled than $\bdr{\Omega_2}$ starting from $H_1$. By construction, we have $\prob(H_0 = n) = \frac{1}{2}$, and $\prob(H_1 = n) = \frac{|\alpha_n|}{||\bm{\alpha}||_1}$, $\forall n \in \{1,2\}$. 
    (a) A diagram for all possible WoI schedules. The black dots indicate $Y_i$'s. For points at the same height, a darker-colored arrow pointing toward them indicates a higher likelihood of the interface being scheduled.
    (b) One example trajectory for each possible WoI schedule.} 
    \label{woi_possibe_schedules}
\end{figure}

In Figure \ref{woi_possibe_schedules}, we show all possible schedules of the WoI estimator, each with one example trajectory, when $N = 2$ and $M = 1$. 
\subsection{Variance Reduction} \label{variance_reduction}
As for the WoB, we want the WoI to be an effective estimator with small variance. We follow a similar strategy~\cite{Sabelfeld_Simonov_2016}  to enhance the performance of WoI. Meanwhile, we must point out that there are many other importance sampling methods that can improve the quality of the WoI estimator.

For some fixed $M$, the variance-reduced WoI estimator depends on two simultaneous walkers for each query point. We divide the first interface (or boundary), $\bdr{\Omega_{H_0}}$, into two parts, $\bdr{\Omega_{H_0}^{+}}$ and $\bdr{\Omega_{H_0}^{-}}$ such that 
\begin{align}
\bdr{\Omega_{H_0}^{+}} = \{\bm{y} \in \bdr{\Omega}: \beta_{H_0}(\bm{y}) \geq 0\}, \qquad
\bdr{\Omega_{H_0}^{-}} = \{\bm{y} \in \bdr{\Omega}: \beta_{H_0} (\bm{y}) < 0\} \notag.
\end{align}
Instead of sampling $Y_0$ uniformly over $\bdr{\Omega_{H_0}}$ in the first step of WoI, we introduce two random variables, $Y_0^{(1)}$ and $Y_0^{(2)}$, under a uniform distribution over $\bdr{\Omega_{H_0}^{+}}$ and $\bdr{\Omega_{H_0}^{-}}$, respectively. First, we sample two points $\bm{y}_0^{(1)} \in \bdr{\Omega_{H_0}^{+}}$ and $\bm{y}_0^{(2)} \in \bdr{\Omega_{H_0}^{-}}$ as realizations of $Y_0^{(1)}$ and $Y_0^{(2)}$, respectively. Then, starting from $\bm{y}_0^{(1)}$ and $\bm{y}_0^{(2)}$ respectively, we construct two Markov chains, namely $\{Y_0^{(1)}, Y_1^{(1)}, \dots, Y_M^{(1)}\}$ and $\{Y_0^{(2)}, Y_1^{(2)}, \dots, Y_M^{(2)}\}$ such that $Y_i^{(k)} \in \bdr{\Omega_{H_i}}$, $\forall k \in \{1, 2\}$. We repeat the standard WoI sampling strategy as described in Section \ref{sampling_strategy} for both Markov chains. Finally, we sum up the estimation from both Markov chains. 

\begin{figure}[!ht]
    \centering
    \begin{subfigure}[b]{0.8\columnwidth}
    \includegraphics[width=\textwidth]{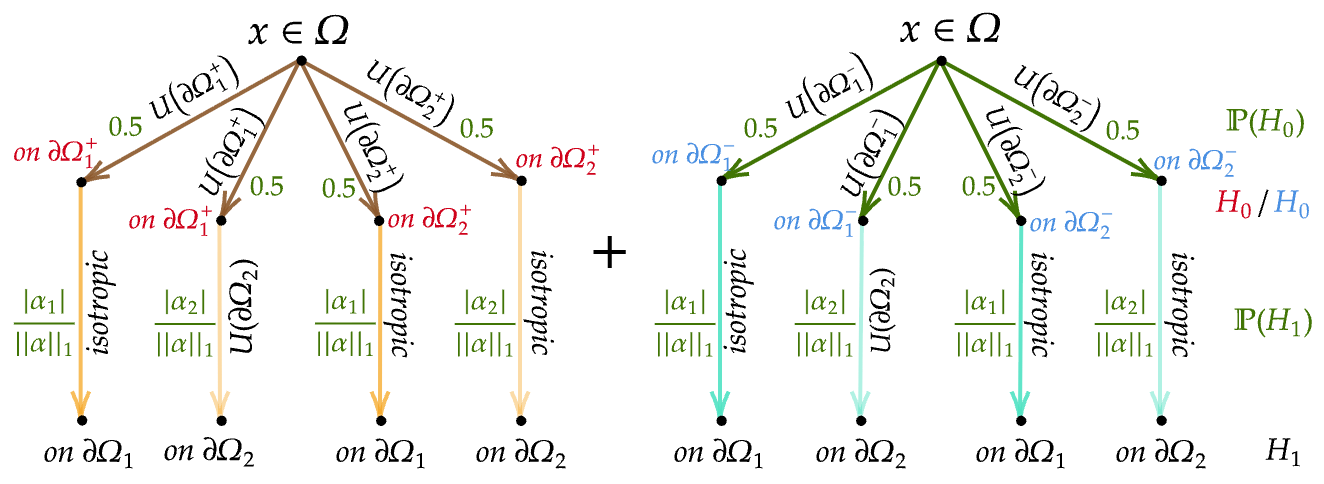}
    \caption{}
    \label{var_reduced_woi_tree}
    \end{subfigure}\\
    \begin{subfigure}[b]{0.8\columnwidth}
    \includegraphics[width=\textwidth]{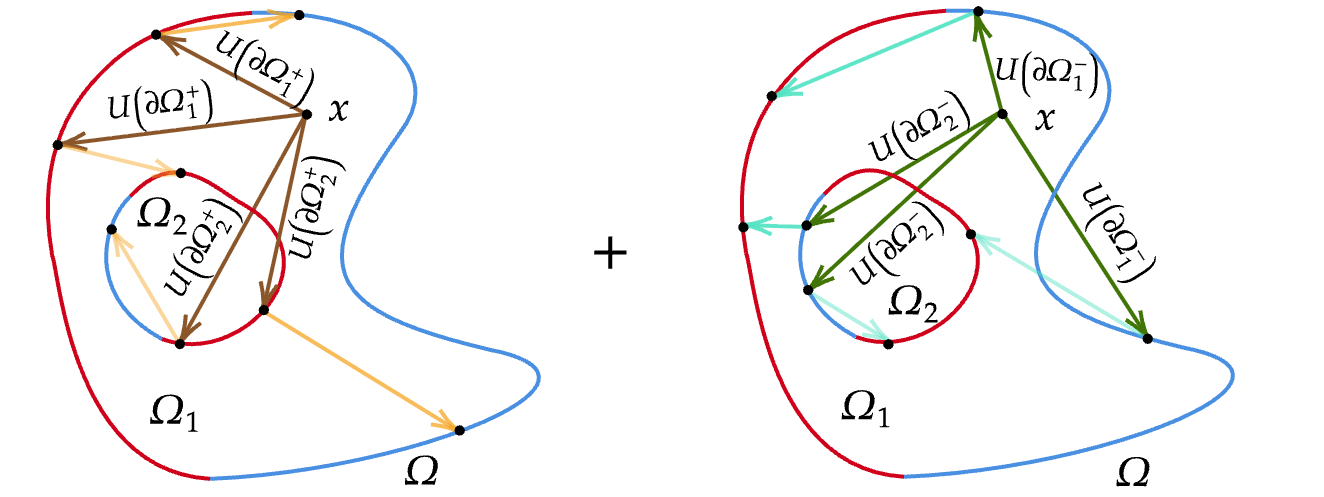}
    \caption{}
    \label{var_reduced_woi_trajectory}
    \end{subfigure} 
    \caption{Consider the domain under the same setup as the one in Figure ~\ref{woi_possibe_schedules}. We still let $|\alpha_1| > |\alpha_2|$, such that $\prob(H_0 = n) = \frac{1}{2}$, and $\prob(H_1 = n) = \frac{|\alpha_n|}{||\bm{\alpha}||_1}$, $\forall n \in \{1,2\}$. (a) A diagram for all possible variance-reduced WoI schedules. The black dots indicate $Y_i^{(k)}$, where $k \in \{1,2\}$. For points at the same height, a darker-colored arrow pointing toward them indicates a higher likelihood of the interface being scheduled.
    (b) One example trajectory for each possible variance-reduced WoI schedules.} 
    \label{var_reduced_woi_possibe_schedules}
\end{figure}

In summary, the variance-reduced WoI estimator is
{\normalsize
\begin{equation}\label{variance_reduced_woi_estimator}
    u(\bm{x}) 
    \approx
    \sum_{i=0}^{M} w_i 
        N ||\bm{\alpha}||_1^i 
        \mathbb{E}_{H_i, \dots, H_0} \left[
        \prod_{m=1}^i \sign(\alpha_{H_m})
         \mathbb{E}_{Y_i} \left[\sum_{k=1}^2
         Q_{(H_i, \dots H_1, H_0)}^{(k)} T_{(H_i, \dots H_1, H_0)}^{(k)}\Phi(\bm{x}, Y_i^{(k)}) \right] \right] ,
\end{equation}
}in which $p_0^{(1)}(\bm{x}, \bm{y}) = \frac{1}{|\bdr{\Omega_{H_0}^{+}}|}$, $p_0^{(2)}(\bm{x}, \bm{y}) = \frac{1}{|\bdr{\Omega_{H_0}^{-}}|}$, $p(\bm{x}, \bm{y})$ is defined in Eq.~\eqref{naive_woi_p},
\begin{equation}\label{variance_reduced_woi_Q}
    Q_{(h_i, \dots, h_1, h_0)}^{(k)} = 
    \begin{cases}
        \frac{\beta_{h_0}(Y_0^{(k)})}{p_0^{(k)}(\bm{x}, Y_0^{(k)})}  \quad & \text{if $i = 0$} \\
        \frac{1}{2} Q_{(h_{i-1}, \dots, h_1, h_0)}^{(k)}
        \quad & \text{if $i \geq 1$ and $h_{i-1} = h_i$} \\
        Q_{(h_{i-1}, \dots, h_1, h_0)}^{(k)}  
        \quad & \text{if $i \geq 1$ and $h_{i-1}$ is a descendant of $h_i$} \\
        \frac{\adj{K}_{h_i, h_{i-1}}(Y_i^{(k)}, Y_{i-1}^{(k)})}{p(Y_{i-1}^{(k)}, Y_i^{(k)})} Q_{(h_{i-1}, \dots, h_1, h_0)}^{(k)}
        \quad & \text{if $i \geq 1$ and for any other $h_{i-1}$ and $h_i$}
    \end{cases},
\end{equation}
and
\begin{equation}\label{variance_reduced_woi_T}
    T_{(h_i, \dots, h_1, h_0)}^{(k)} = 
    \begin{cases}
        1 \quad & \text{if $i = 0$} \\
       \prod_{m=1}^{i} q(Y_{m-1}^{(k)}, Y_m^{(k)}) \sign[(Y_m^{(k)} - Y_{m-1}^{(k)}) \cdot \bm{n}(Y_m^{(k)})]  & \text{if $i \geq 1$}
    \end{cases}.
\end{equation}
Algorithm \ref{alg_var_WoI} provides a detailed explanation of the variance-reduced WoI method, which can replace line 7 of Algorithm \ref{alg_WoI_wrapper} whenever necessary. 
\begin{algorithm}[!ht]
\caption{The Variance Reduced WoI Estimator}
\label{alg_var_WoI}
\begin{algorithmic}[1]
\Function{VarReducedWoI}{$\bm{x}$, $\adj{K}_{i,j}(\bm{x}, \bm{y})$,  $\Phi(\bm{x}, \bm{y})$, $\bm{\alpha}$, $\bm{\beta}$, $M$, $S$}
\State $u_w \gets 0$, $T^{(1)} \gets 1$, $T^{(2)} \gets 1$, $\texttt{weight} \gets 1$.
\For{$i = 0$ to $M$}
\If{$i == M$}
\State $\texttt{weight} \gets 0.5$ \Comment{truncation}
\EndIf
\State $h_i \gets S(i)$
\If{$i == 0$} \Comment{transitional PDF is $p_0(\bm{x}, \bm{y})$}
\State $\bm{y}_i^{(1)} \gets \text{Unif}(\bdr{\Omega_{h_i}^{+}})$, $\bm{y}_i^{(2)} \gets \text{Unif}(\bdr{\Omega_{h_i}^{-}})$
\State $Q^{(1)} = \beta_{h_i}(\bm{y}_i^{(1)}) / \frac{1}{|\bdr{\Omega_{h_i}^{+}}|}$, $Q^{(2)} = \beta_{h_i}(\bm{y}_i^{(2)}) / \frac{1}{|\bdr{\Omega_{h_i}^{-}}|}$
\State $u_w = u_w + \texttt{weight} * ||\alpha||_1^{i} * (\sum_{k=1}^2 Q^{(k)} * T^{(k)} * \Phi(\bm{x}, \bm{y}_i^{(k)}))$
\Else\Comment{transitional PDF is $p(\bm{x}, \bm{y})$}
\State $h_{i-1} \gets S(i-1)$ \Comment{$i \geq 1$}
\For{$k = 1$ to $2$}
\If{$h_{i-1} = h_i$}
\State $\bm{d} \gets \text{an isotropic direction}$
\State $q(\bm{y}_{i-1}^{(k)}, \bm{y}_i^{(k)}) \gets$ \# of intersections between $\bdr{\Omega_{h_i}}$ and line $\ell(t) = t \bm{d} + \bm{y}_{i-1}^{(k)}$
\State $\bm{y}_i^{(k)} \gets$ any of the $q(\bm{y}_{i-1}^{(k)}, \bm{y}_i^{(k)})$ intersections
\State $Q^{(k)} \gets Q^{(k)} * 0.5$
\ElsIf{$h_{i-1}$ a descendant of $h_i$}
\State $\bm{d} \gets \text{an isotropic direction}$
\State $q(\bm{y}_{i-1}^{(k)}, \bm{y}_i^{(k)}) \gets$ \# of intersections between $\bdr{\Omega_{h_i}}$ and ray $r(t) = t_{\geq 0}\bm{d}+ \bm{y}_{i-1}^{(k)}$
\State $\bm{y}_i^{(k)} \gets$ any of the $q(\bm{y}_{i-1}^{(k)}, \bm{y}_i^{(k)})$ intersections
\State $Q^{(k)} \gets Q^{(k)} * 1$
\Else
\State $q(\bm{y}_{i-1}^{(k)}, \bm{y}_i^{(k)}) \gets 1$
\State $\bm{y}_i^{(k)} \gets \text{Unif}(\bdr{\Omega_{h_i}})$
\State $Q^{(k)} = Q^{(k)} * \adj{K}_{h_i, h_{i-1}}(\bm{y}_i^{(k)}, \bm{y}_{i-1}^{(k)})/ \frac{1}{|\bdr{\Omega_{h_i}}|}$ 
\EndIf
\State $T^{(k)} = T^{(k)} * q(\bm{y}_{i-1}^{(k)}, \bm{y}_i^{(k)}) * \sign[(\bm{y}_i^{(k)} - \bm{y}_{i-1}^{(k)})\cdot \bm{n}(\bm{y}_i^{(k)})]$
\EndFor
\State $u_w = u_w + \texttt{weight} * ||\alpha||_1^{i} * [\prod_{m=0}^{i} \sign(\alpha_{h_m})] * (\sum_{k=1}^2 Q^{(k)} * T^{(k)} * \Phi(\bm{x}, \bm{y}_i^{(k)}))$
\EndIf
\EndFor
\State return $u_w$
\EndFunction
\end{algorithmic}
\end{algorithm}
Figure \ref{var_reduced_woi_possibe_schedules} demonstrates all possible schedules that can occur in a variance-reduced WoI estimator. For each type of schedule, a sampled trajectory is also provided within the same figure. 
\subsection{Gradient Estimation} \label{gradient_estimation}
The WoI estimator inherits the property from the WoB that we obtain a gradient estimation with almost no additional cost. This feature enables a fully Monte Carlo approach for computing the gradient of the solution to the elliptic interface problem. In contrast, prior method~\cite{HAN2020109672} relied on applying finite-differences to Monte Carlo estimates of $u(\bm{x})$.

The estimator for $\nabla u(\bm{x})$ is obtained by simply taking the gradient of $\Phi (\bm{x}, \bm{y})$ in the estimator of $u(\bm{x})$. Hence, the WoI gradient estimator is
\begin{equation}\label{naive_woi_grad_estimator}
    \nabla u(\bm{x}) 
    \approx
    \mathbb{E}
        \left[
        \sum_{i=0}^{M} w_i
        \sum_{h_{i}, \dots, h_1, h_0=1}^{N}
        \bigg(\prod_{m=1}^i \alpha_{h_m} \bigg)
        Q_{(h_i, \dots, h_1, h_0)} T_{(h_i, \dots, h_1, h_0)} \nabla_{\bm{x}} \Phi(x, Y_i) \right],
\end{equation}
while a variance-reduced WoI gradient estimator is
{\normalsize
\begin{equation}\label{variance_reduced_woi_grad_estimator}
    \nabla u(\bm{x}) \approx
         \sum_{i=0}^{M} w_i N ||\bm{\alpha}||_1^i 
        \mathbb{E}_{H_i, \dots, H_0} \left[
        \prod_{m=1}^i \sign(\alpha_{H_m})
         \mathbb{E}_{Y_i} \left[\sum_{k=1}^2
         Q_{(H_i, \dots H_1, H_0)}^{(k)} T_{(H_i, \dots H_1, H_0)}^{(k)} 
         \nabla_{\bm{x}} \Phi(\bm{x}, Y_i^{(k)}) \right] \right].
\end{equation}
}
The WoI gradient estimator allows us to use all walkers and trajectories we generated for the WoI estimator to estimate $\nabla u(\bm{x})$, which gives almost no extra cost.
\subsection{Implementation Considerations} \label{mplementation}
We outline some key implementation details of our WoI estimator in this section. The discussion includes proper representation of computational domains, parallelization of walkers, batch estimation of query points, and the management of floating-point precision.

\textit{Domain Definitions.} 
There are many common ways to represent geometries, such as signed distance functions (SDFs), level sets, surface meshes, and implicit functions. Each representation offers trade-offs in terms of accuracy, computational cost, and ease of sampling under certain distributions. According to our design of WoI and variance-reduced WoI estimators, we prioritize representations that facilitate uniform sampling over the domain boundary and efficient computation of ray–mesh and line–mesh intersections. To this end, we adopt surface meshes equipped with axis-aligned bounding box (AABB) tree structures, which are well-suited for both uniform boundary sampling and geometric intersection queries. We additionally accept implicit function representations in cases where established algorithms are available for uniform sampling on the boundary and for performing intersection tests.

\textit{Parallelization of Walkers.} It is well known that Monte Carlo methods are, in general, highly parallelizable. This principle extends naturally to the WoI estimator. Suppose each of the $\mathcal{W}$ walkers performs $M$ steps, and assume that we sampled $H_1, H_2, \dots H_M$ with realizations $h_1, h_2, \dots h_M$ as one schedule of interfaces. When the $i^{th}$ states of several walkers are located on the same boundary, they can be generated in parallel, provided that the hierarchical relationship between $h_{i-1}$ and $h_i$ falls in one of the three cases of $p(\bm{y}_{i-1}, \bm{y}_i)$.

\textit{Batch Estimation with Shared Trajectories.} The WoI estimator, its gradient estimator, and their variance-reduced counterparts all require constructing trajectories for walkers by sampling the interfaces (or the boundary). However, building these stochastic trajectories is computationally expensive. While it is possible to generate independent trajectories for each query point, we instead reuse the same set of trajectories across multiple query points. To simplify the implementation, we randomly shuffle the query points and assign the same set of trajectories to those points with nearby indices in the shuffled list. Hence, we obtain a vectorized evaluation of the WoI estimator, allowing us to compute the solution to the interface problem at a batch of query points simultaneously.

Although this approach introduces correlations in the estimation errors at different spatial locations, it significantly reduces the overall computational cost. With an appropriate choice of batch size, the estimator retains sufficient stochastic properties, and the resulting high-frequency error can still be well-captured. We recommend that the batch size is no larger than 0.5\% of the total number of query points.

\textit{Floating-point Precision.} All the estimators developed so far require the evaluation of either $\Phi(\bm{x}, Y_i)$ or $\nabla \Phi(\bm{x}, Y_i)$. To avoid repeatedly evaluating the Green's function and its gradient, it is desirable to allow vectorized input for both functions. Achieving this level of parallelism requires storing a large number of $(\bm{x}, Y_i)$ pairs in memory. To improve storage efficiency and maximize batch size, we store the query point vector, $\bm{x}$, and the state of the walkers, $Y_i$, in single-precision floating-point format whenever possible. This reduction in memory footprint allows for larger-scale vectorized operations while maintaining sufficient numerical accuracy for estimation.
\section{Results}\label{numerical_results}
In this section, we first evaluate the accuracy and convergence behaviour of the WoI estimator and its variance-reduced counterpart. We provide examples to show how our Monte Carlo estimators solve problems with different domain geometries and how they naturally extend to higher-dimensional domains. Next, we integrate a neural network into our method to obtain a continuous representation of the solution and to smooth out the high-frequency errors of the Monte Carlo estimator. This approach leverages the neural network’s inductive bias toward smooth functions, providing both improved accuracy and a convenient functional form for downstream evaluation. Finally, we conclude this section by demonstrating the practical effectiveness of the WoI framework in solving interface problems under different application scenarios.

\textit{Correctness.} 
We begin by evaluating the accuracy of the WoI estimator. Let $\hat{u}(\bm{x})$ denote the estimated solution and $u(\bm{x})$ represent the ground truth. Given that the solution to the interface problem in Eq.~\eqref{interface_problem} is unique up to a constant, we anticipate that the WoI estimation deviates from the ground truth by a constant offset. To more effectively quantify the estimation error, $e(\bm{x})$, we align the estimated solution by shifting it to match the ground truth at a fixed reference point, $\bm{x}_{\textit{ref}}$,
\[
e(\bm{x}) = |u(\bm{x}) - [\hat{u}(\bm{x}) - \hat{u}(\bm{x}_{\textit{ref}}) + u(\bm{x}_{\textit{ref}})]|.
\]
In our examples, we will consider the $L_2$ norm of the error, $||e(\bm{x})||_2$ and the relative $L_2$ norm of the error, $\frac{||e(\bm{x})||_2}{||u(\bm{x})||_2}$.

\textbf{Example 1. A smooth harmonic function in 2D.} As the first example, we examine the performance of the WoI estimator and the variance-reduced WoI estimator over a two-dimensional domain, for which the ground truth is a $C^{\infty}$ function. Since any harmonic function satisfies the elliptic interface problem, we let the first ground truth be $u(x_1, x_2) = x_1^3 - 3 x_1 x_2^2$. To break the symmetry of the solution, we shift the domain away from the origin and solve the interface problem on a circle of radius 2 centered at $(0.5, -0.5)$. 

\begin{figure}[htbp]
\centering
\begin{tabular}{@{}c@{}c c c c c}
% ============= WOI ==============
\multirow{3}{*}{
  \rotatebox[origin=c]{90}{%
    \parbox[c]{\dimexpr 0.5\columnwidth\relax}{\centering\vfill\textbf{WoI}\vfill}%
  }
} & 
\textbf{Ground Truth} & $\mathcal{W} = 10^6$ & $\mathcal{W} = 5 \times 10^6$ & $\mathcal{W} = 10^7$ & {} \\

% First Row: estimation
{} & \includegraphics[width=0.2\columnwidth]{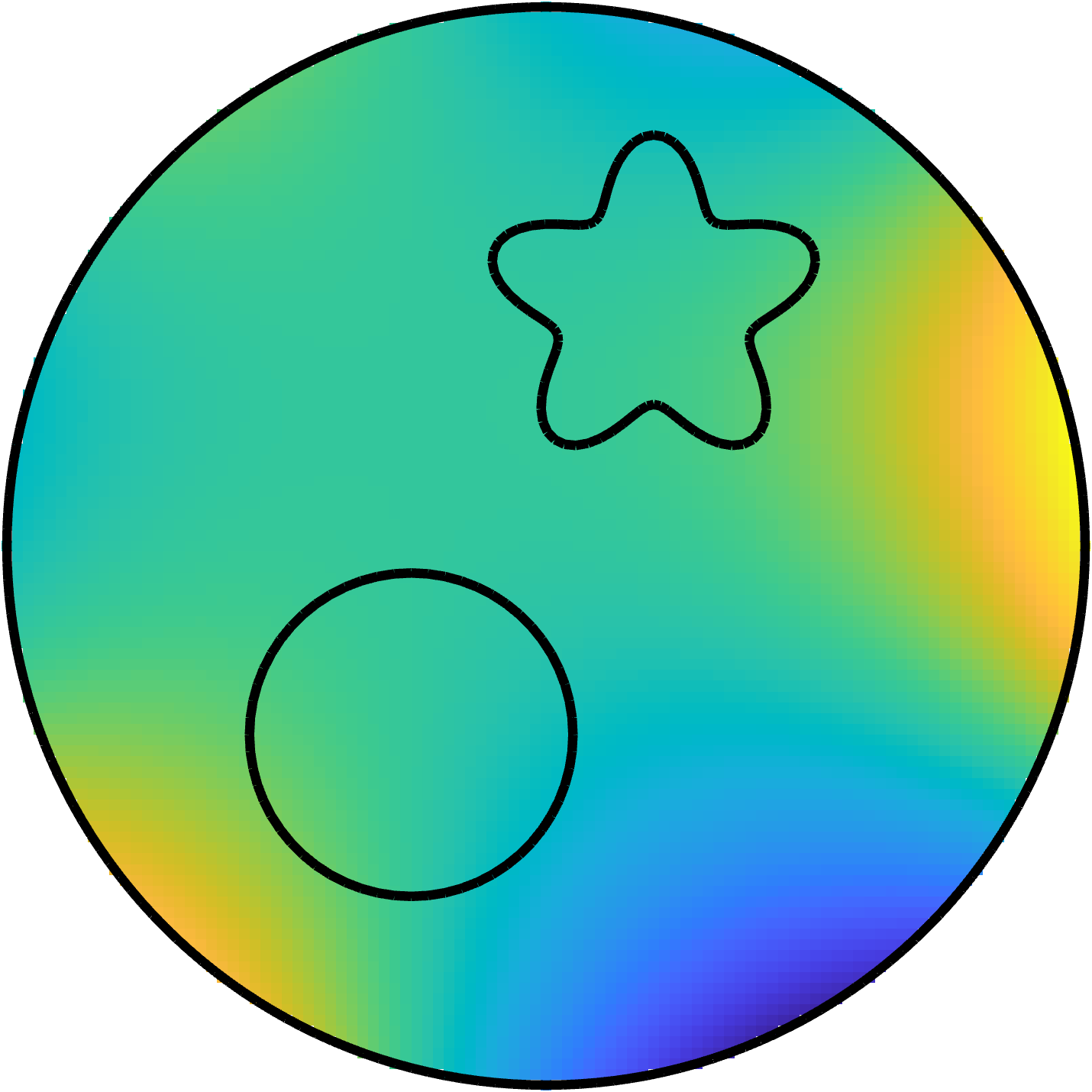} &
\includegraphics[width=0.2\columnwidth]{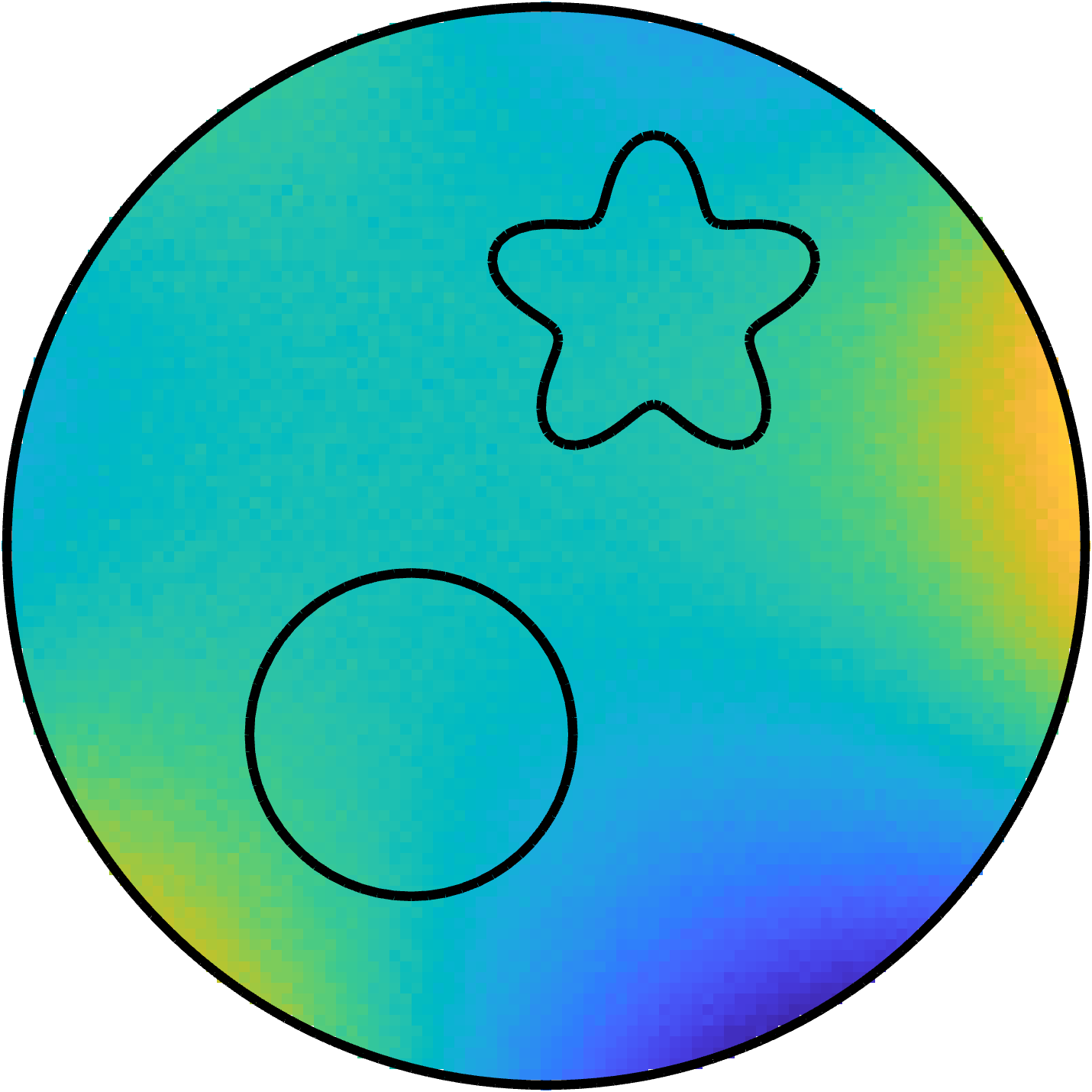}  &
\includegraphics[width=0.2\columnwidth]{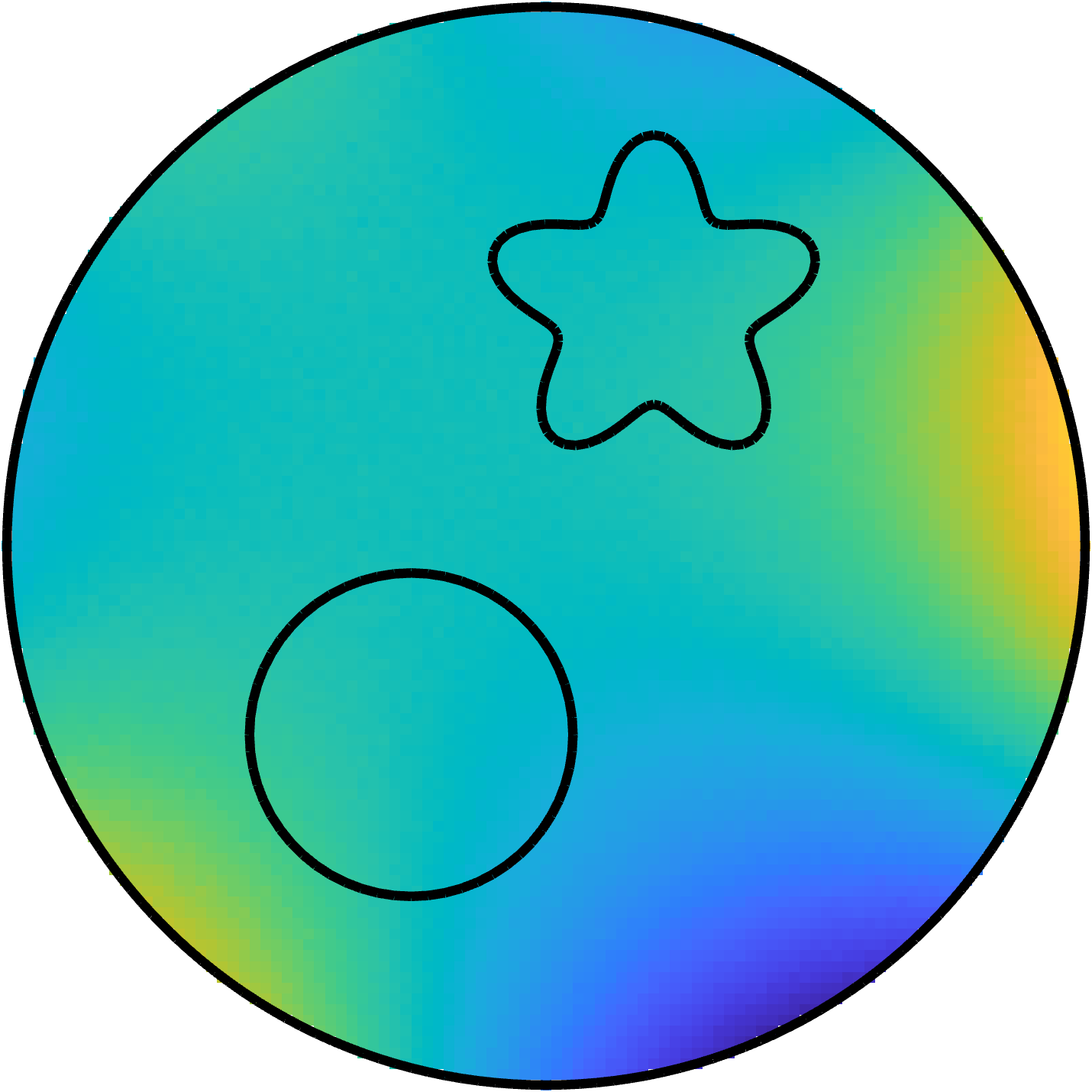} &
\includegraphics[width=0.2\columnwidth]{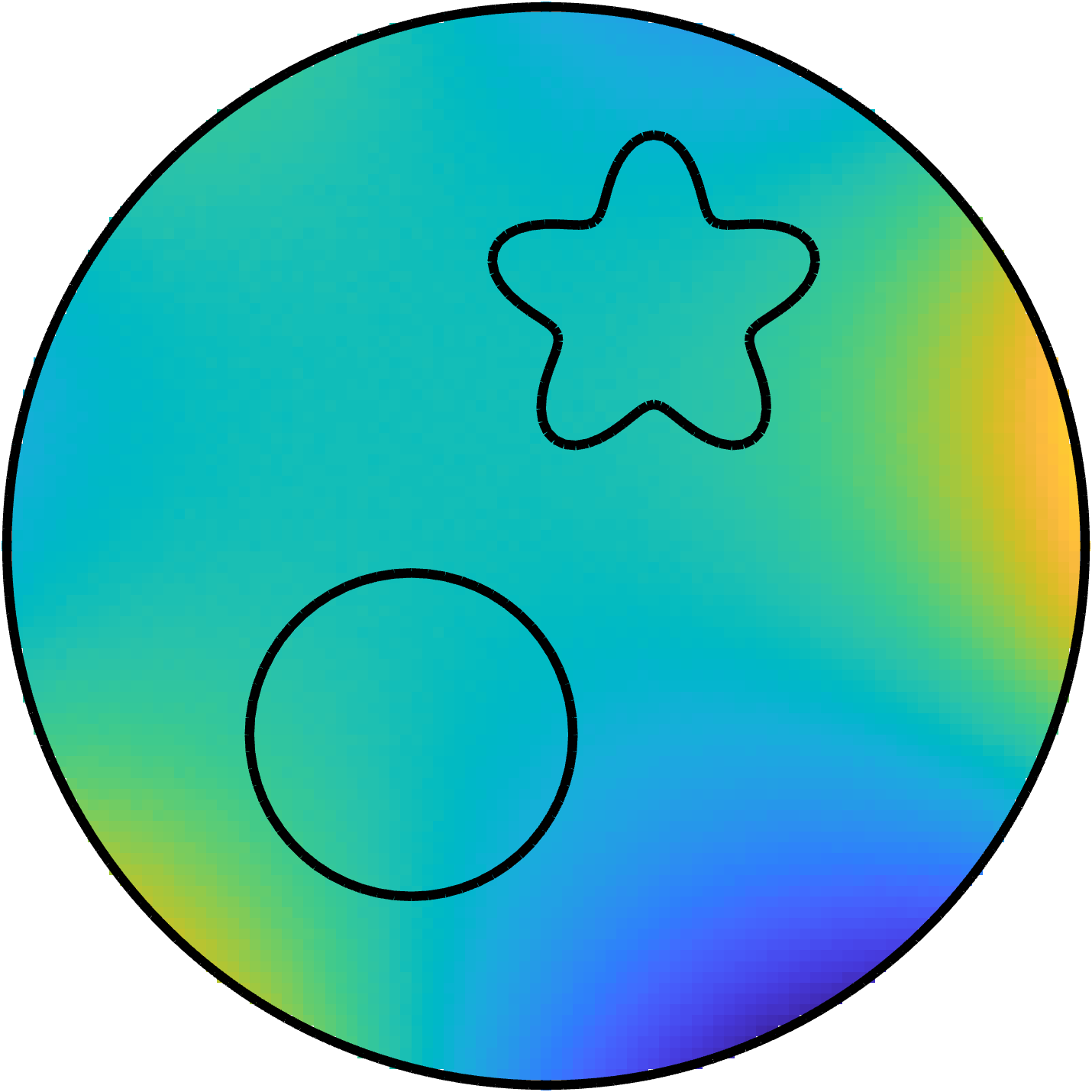} &
\includegraphics[width=0.03\columnwidth]{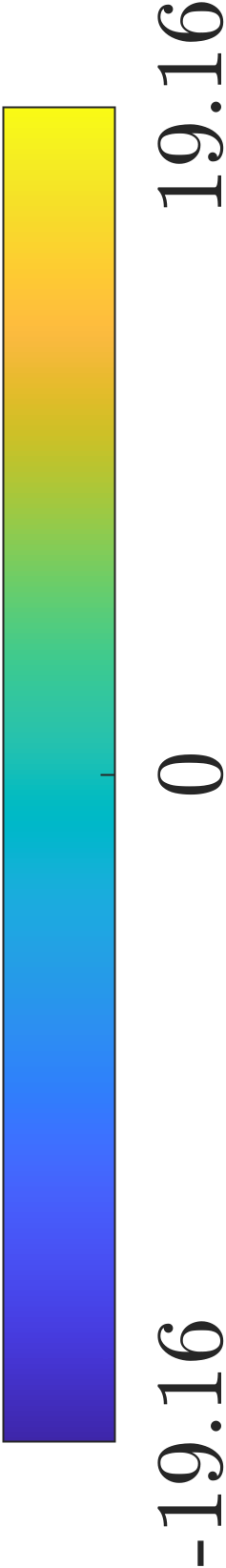} \\

% Second Row: pointwise estimation
{} & {} &
\includegraphics[width=0.2\columnwidth]{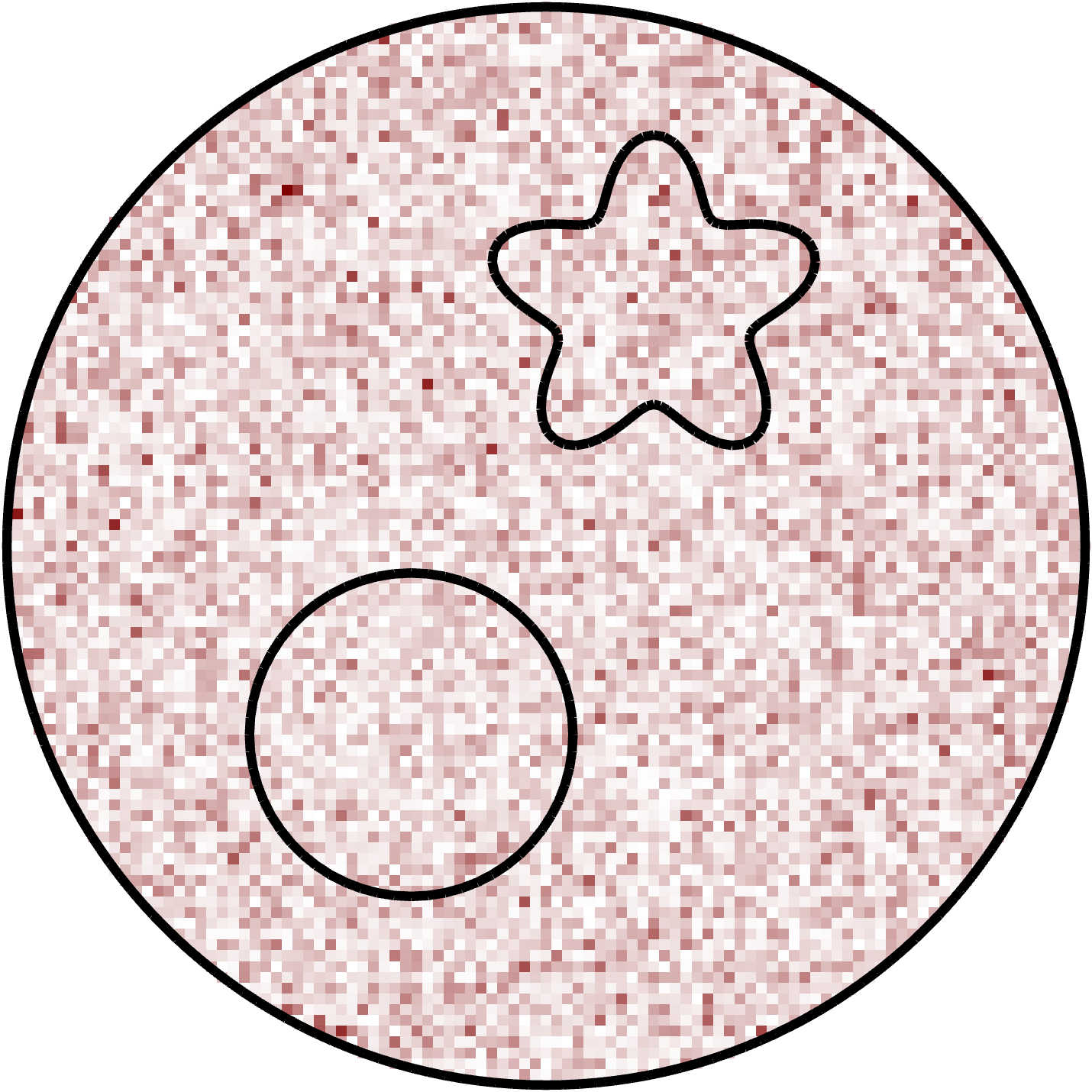} &
\includegraphics[width=0.2\columnwidth]{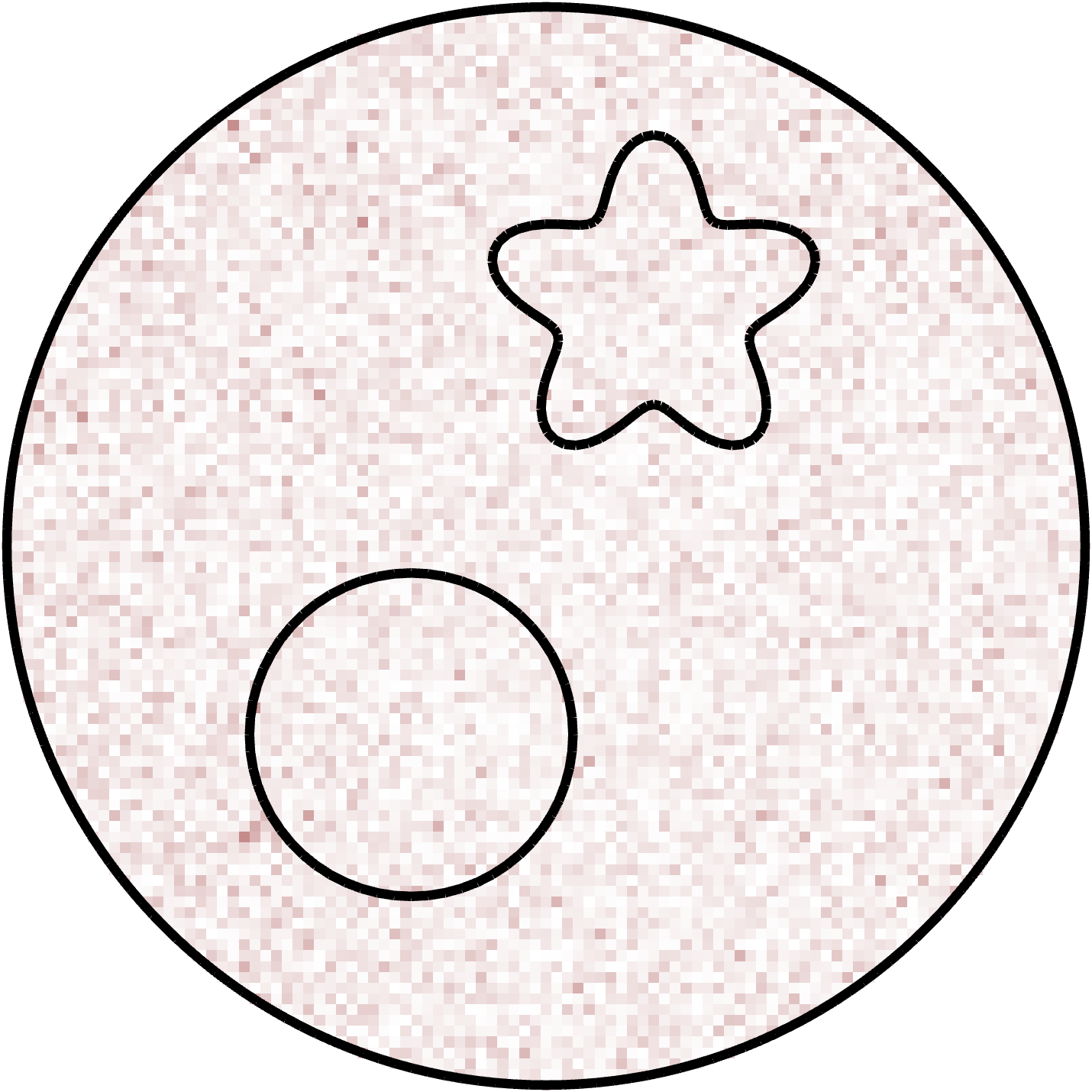} &
\includegraphics[width=0.2\columnwidth]{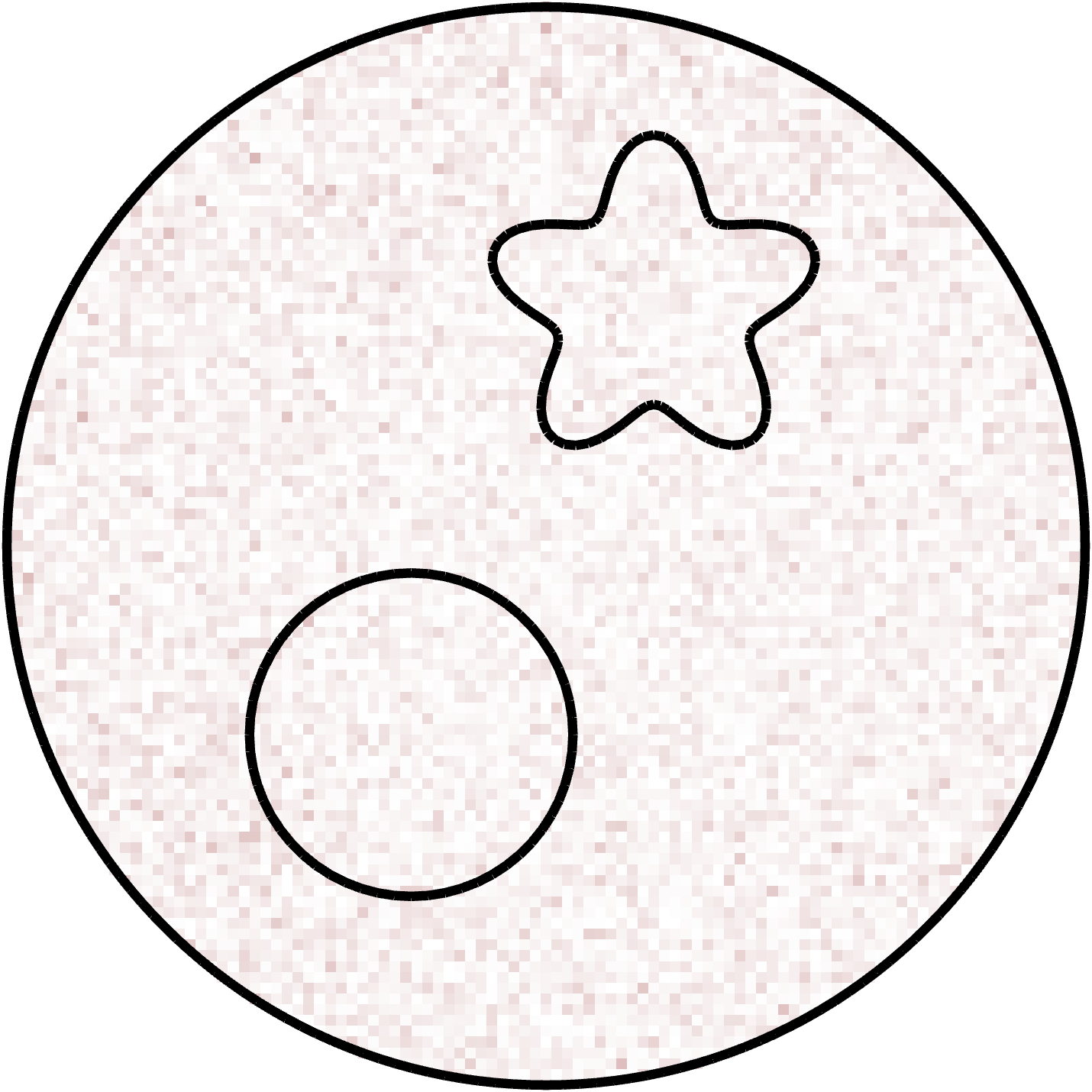} &
\includegraphics[width=0.03\columnwidth]{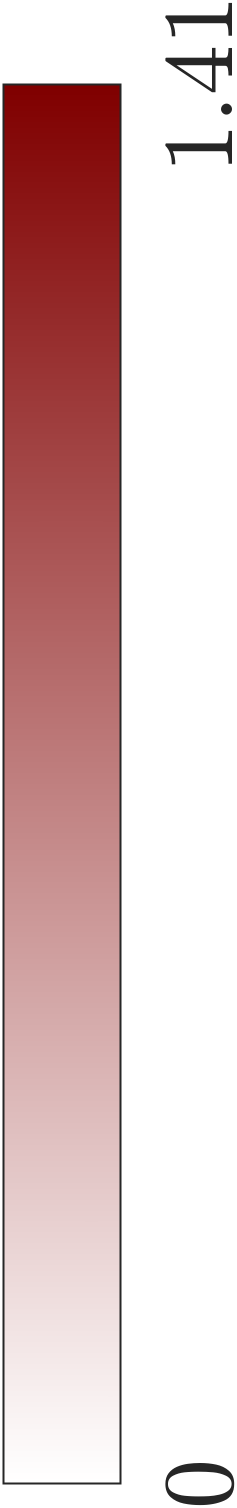}\\
{} & \textbf{$\bm{L_2}$ Error} & $0.319$ & $0.139$ & $0.097$ & {} \\
% relative L2 error
{} & \textbf{Relative $\bm{L_2}$ Error} & $6.77\%$ & $2.94\% $ & $2.07\%$ & {}\\
%===============================================
\hline \hline \\
%===============================================
% =========== variance-reduced WoI =============
\multirow{3}{*}{
  \rotatebox[origin=c]{90}{%
    \parbox[c]{\dimexpr 0.5\columnwidth\relax}{\centering\vfill\textbf{Var-reduced WoI}\vfill}%
  }
} & \textbf{Ground Truth} & $\mathcal{W} = 2 \times 10^6$ & $\mathcal{W} = 10^7$ & $\mathcal{W} = 2 \times 10^7$ & \\
% First Row: estimation
{} & \includegraphics[width=0.2\columnwidth]{figure/N=3_circle_curved_star/N=3_circle_curved_star_ground_truth.png} &
\includegraphics[width=0.2\columnwidth]{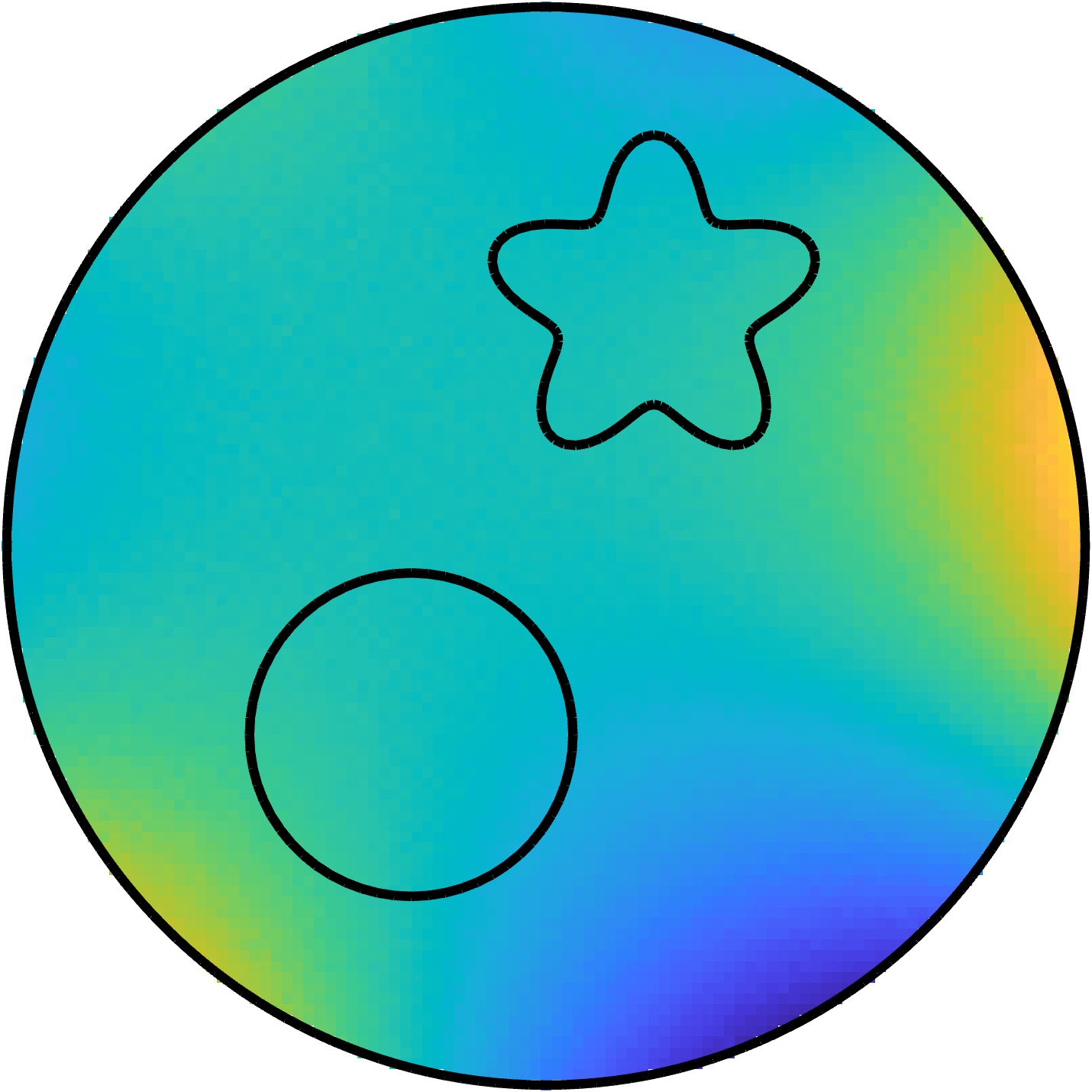}  &
\includegraphics[width=0.2\columnwidth]{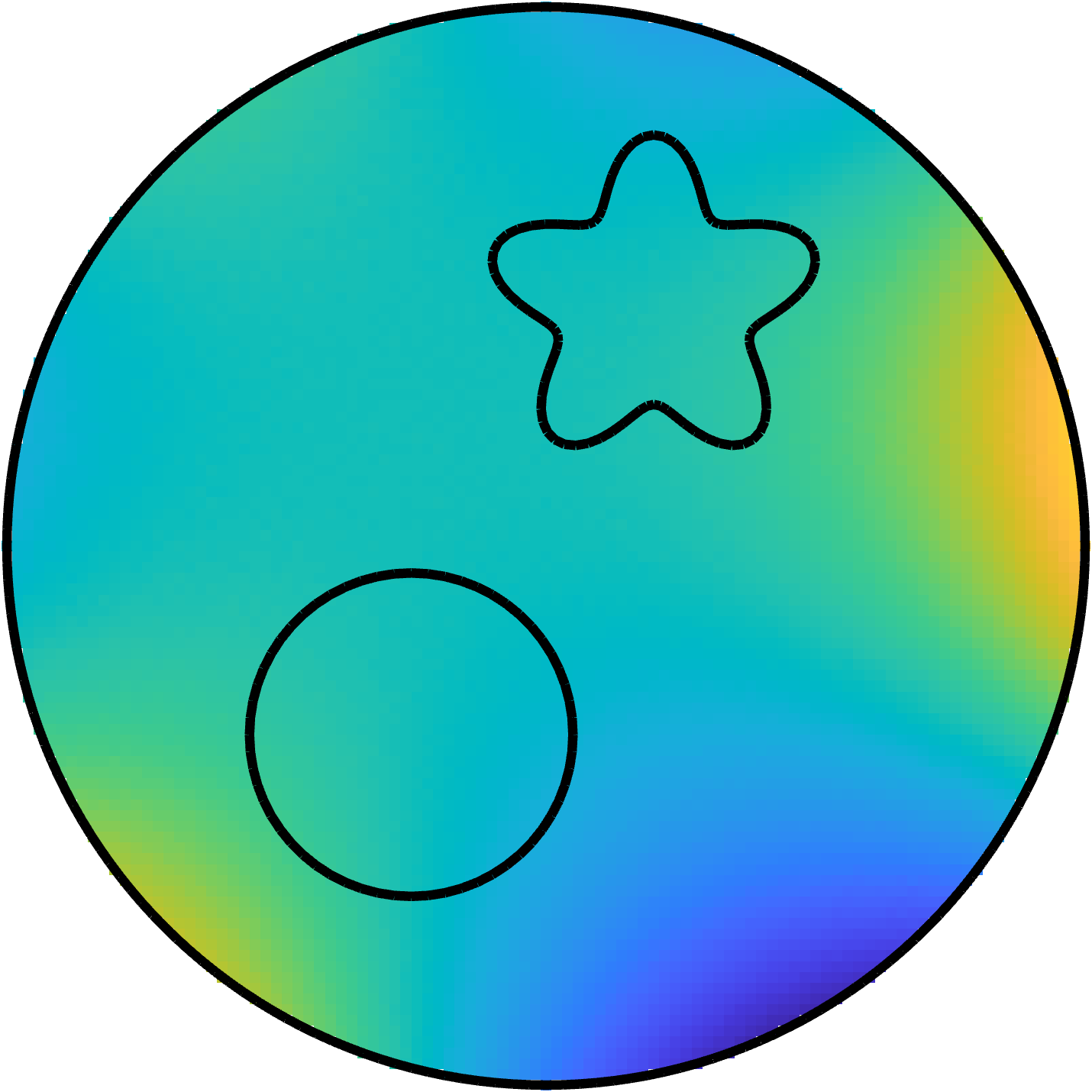} &
\includegraphics[width=0.2\columnwidth]{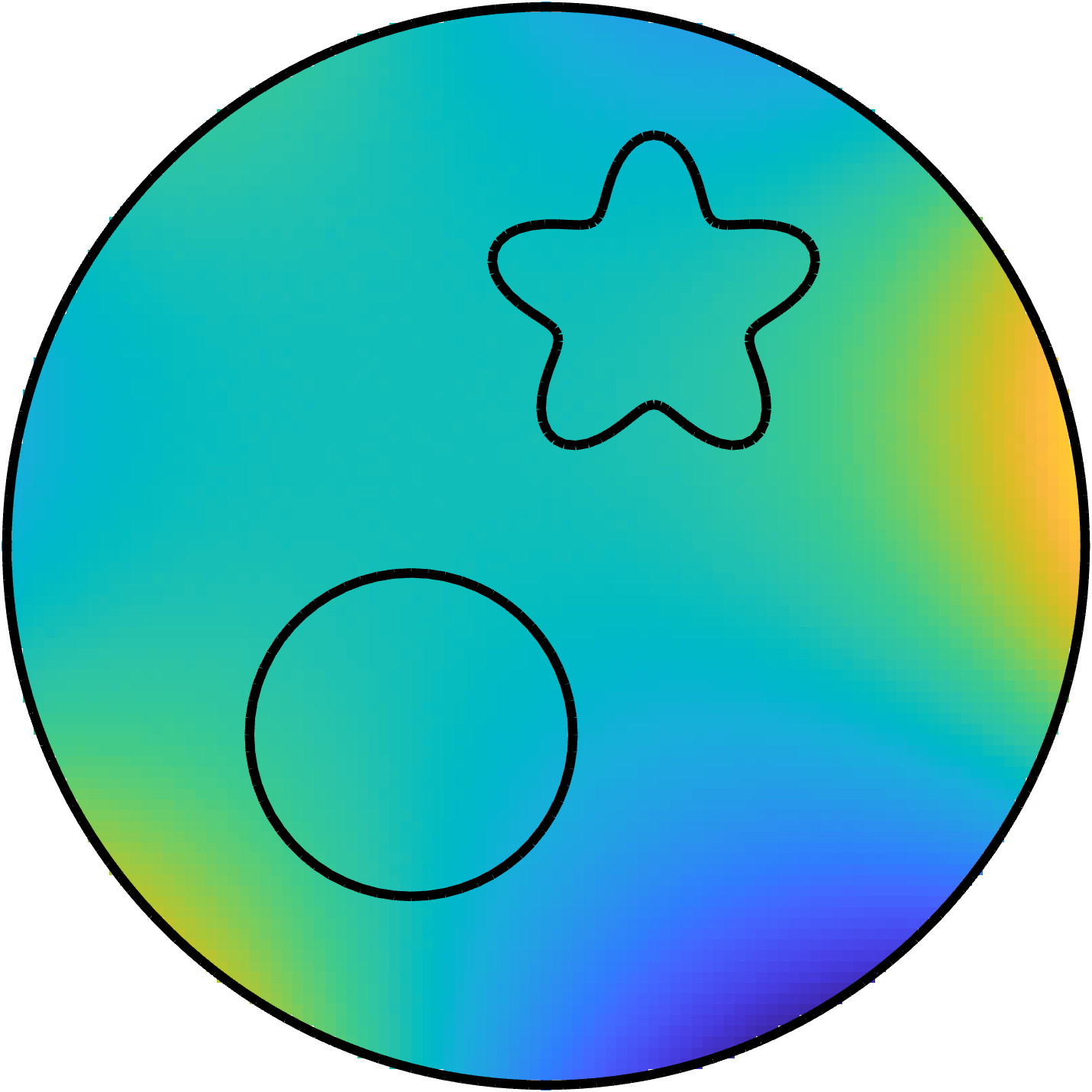} &
\includegraphics[width=0.03\columnwidth]{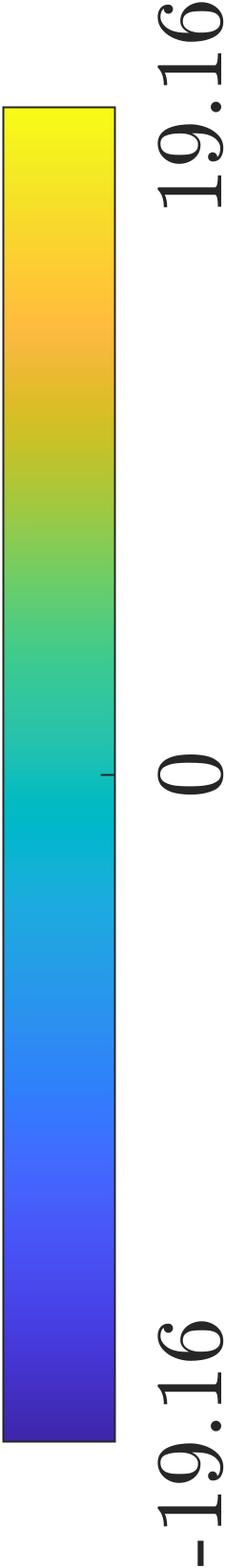}\\
% Second Row: pointwise estimation
{} & {} &
\includegraphics[width=0.2\columnwidth]{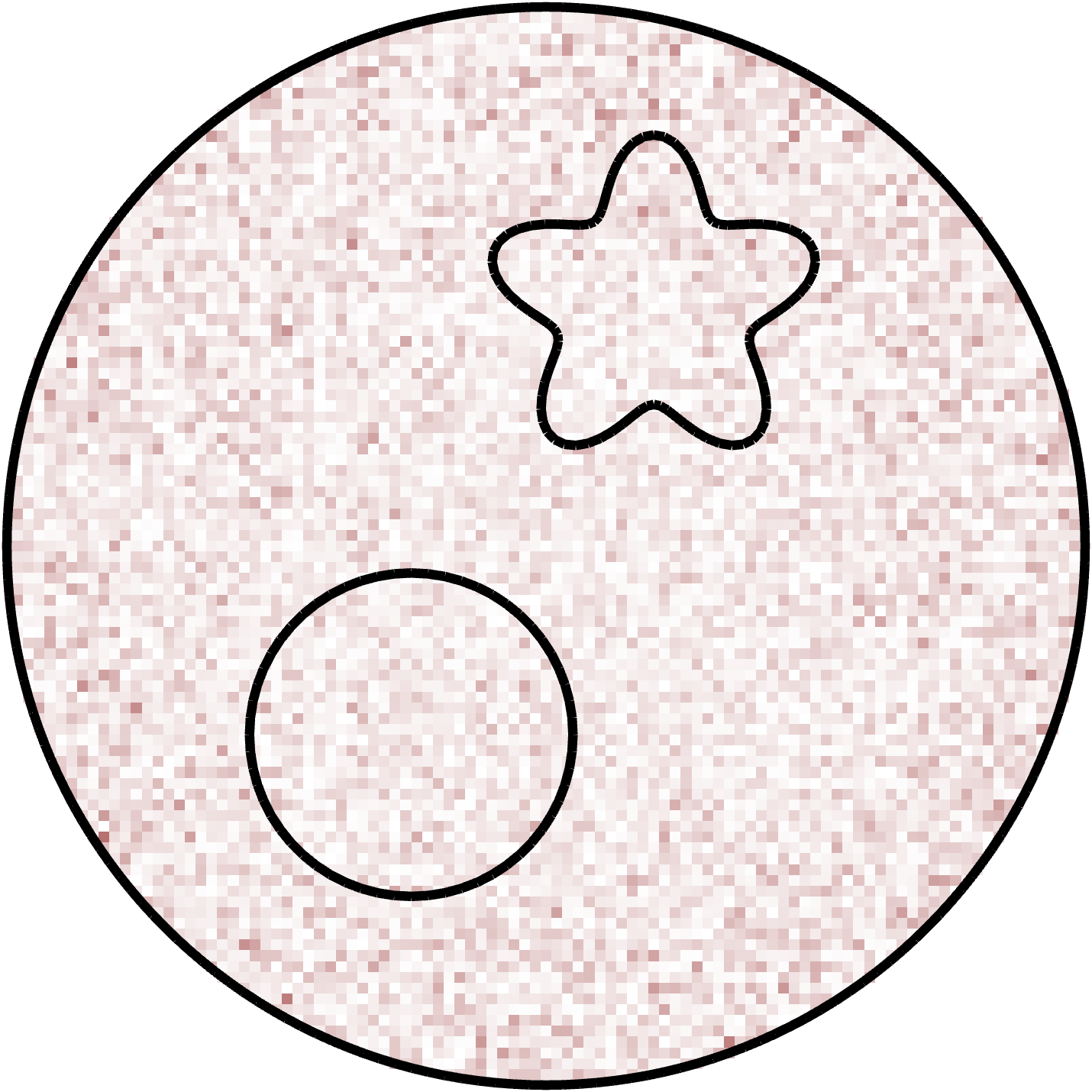} &
\includegraphics[width=0.2\columnwidth]{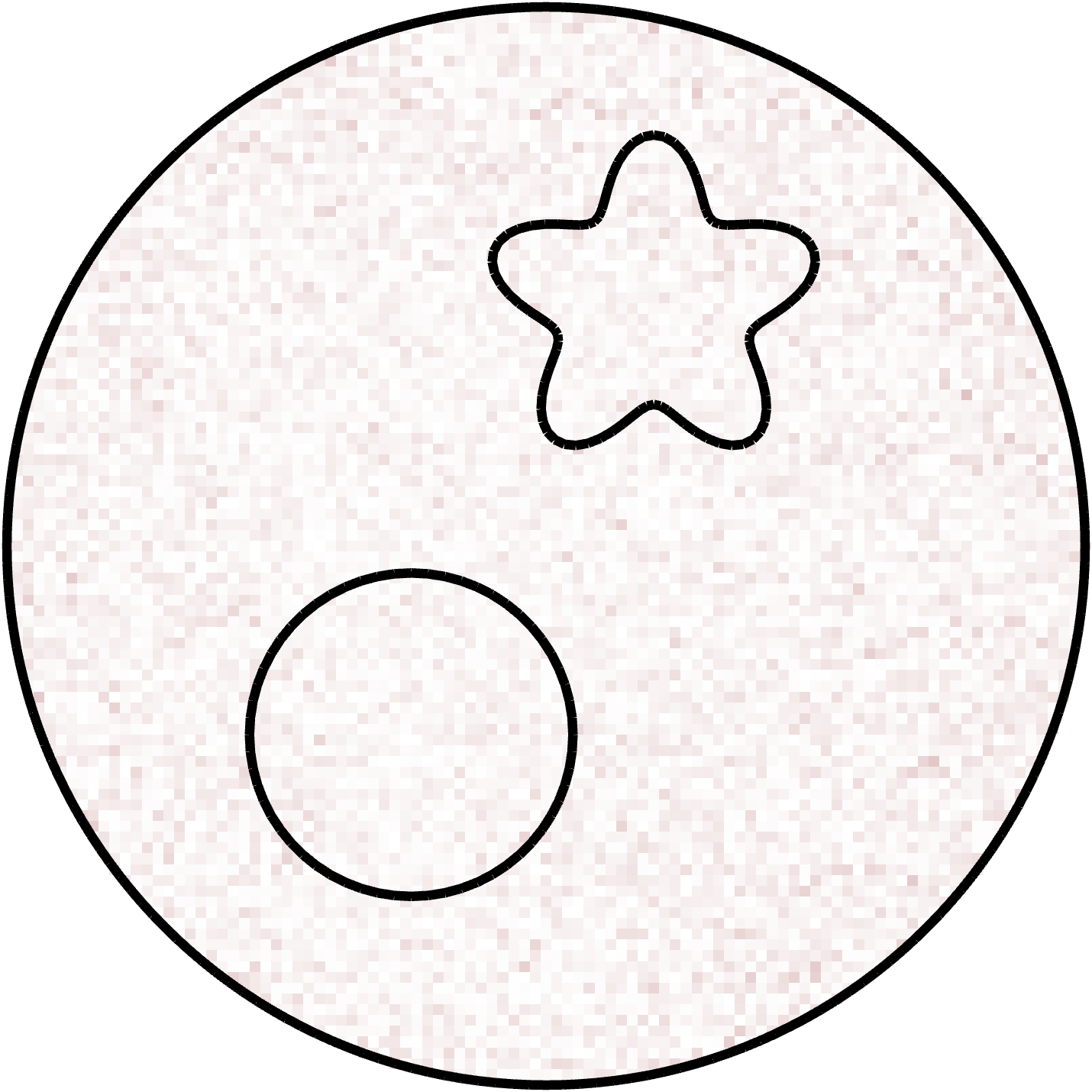} &
\includegraphics[width=0.2\columnwidth]{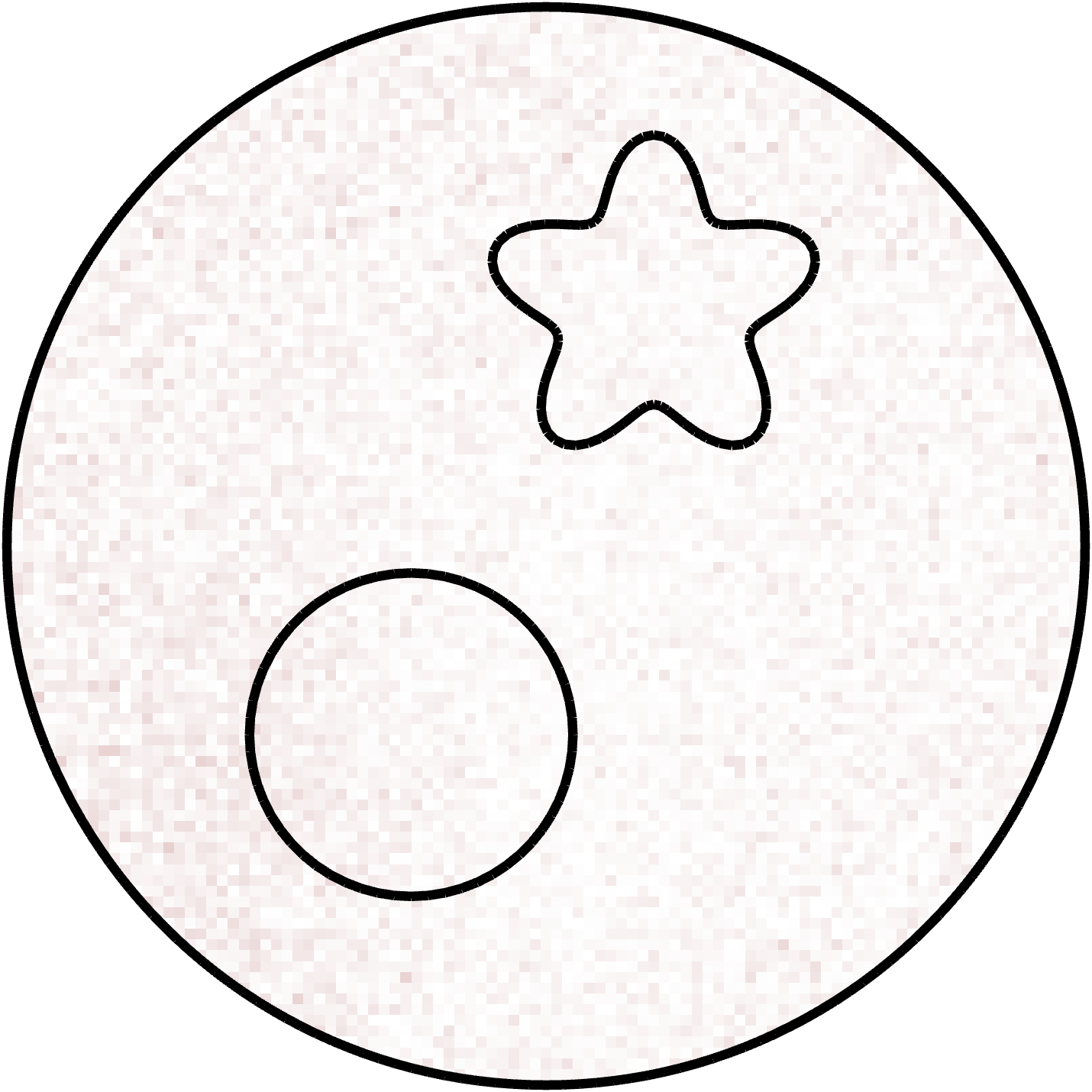} &
\includegraphics[width=0.03\columnwidth]{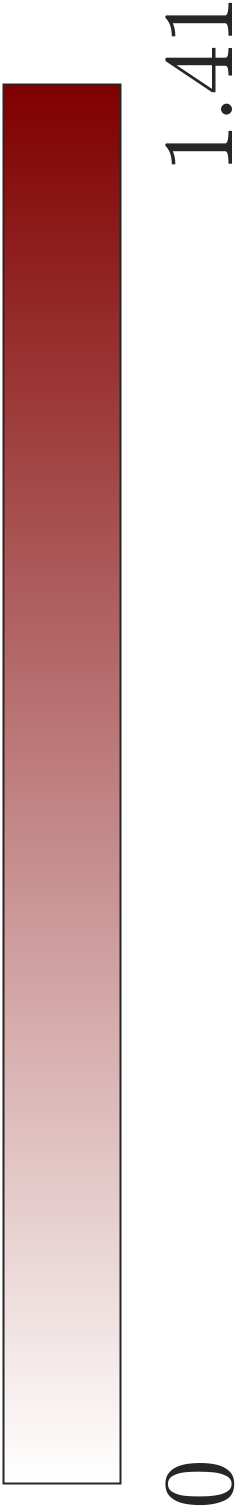}\\

% relative L2 error
{} & \textbf{$\bm{L_2}$ Error} & $0.181$ & $0.081$ & $0.67$ &\\
% relative L2 error
{} & \textbf{Relative $\bm{L_2}$ Error} & $3.83\%$ & $1.72\% $ & $1.41\%$ & \\
\end{tabular}
\caption{The outer-most circle $\Omega_1$ has $\sigma_1 = 1.5$, $\Omega_2$ is the circle inside with $\sigma_2 = 0.5$, and $\Omega_3$ is the smooth star with $\sigma_3 = 1.1$. In both halves of the figure, the first row shows the solution of WoI or variance-reduced WoI with different numbers of walkers. The second row plots the corresponding errors. The ground truth function is $u(x_1, x_2) = x_1^3 - 3x_1 x_2^2$. Both estimators walk $M = 4$ steps. }\label{ex_woi_harmonic2D}
\end{figure}

We allow both estimators to walk $M = 4$ steps. Furthermore, we compute the expectations with respect to $\mathbb{E}_{Y_i}[\cdot]$ using the same total number of samples. In the variance-reduced WoI estimator, each sample of $\mathbb{E}_{Y_i}[\cdot]$ involves contributions from both $Y_i^{(1)}$ and $Y_i^{(2)}$. Therefore, we generate the same number of samples for each of $Y_i^{(k)}$, where $k \in \{1, 2\}$, as the number of $Y_i$ samples used in the standard WoI estimator. Since the variance-reduced estimator relies on two independent walkers per sample, it inherently requires twice as many walkers as the standard WoI estimator. To ensure a fair comparison in terms of computational cost and estimation quality, we match the total number of walkers used in both methods. When $\mathcal{W} = 10^7$, we observe that the variance-reduced WoI estimator outperforms the WoI estimator.
The corresponding results for both methods are shown in Figure~\ref{ex_woi_harmonic2D}, which visualizes the estimated solutions and the pointwise errors relative to the ground truth.

\textbf{Example 2. Potential of a point charge at the origin.} In the second example, we test our estimators with higher-dimensional ellipsoidal domains. The ellipsoids are defined via implicit surfaces of the form $(\bm{x} - \bm{c})^T A (\bm{x} - \bm{c}) = 1$, where $\bm{c}$ is the center, and $A$ is a positive definite matrix. The ground truth function, $u: \R^n \to \R$, is defined as $u(\bm{x}) = ||x||^{2-n}$, where $n \geq 3$. This function arises naturally in physics: in electrostatics, it represents the potential field generated by a unit point charge at the origin in $\R^n$; in gravitational theory, it corresponds to the gravitational potential of a point mass.

Due to the singularity of $u(\bm{x})$ at the origin, we shift the domain to exclude the origin so that the singularity is not approximated. The domains consist of three nested ellipsoids. We name the outer-most ellipsoid $\Omega_1$, the middle ellipsoid $\Omega_2$, and the inner-most ellipsoid $\Omega_3$. The detailed parameters of the domain are recorded in Table \ref{ellipsoid_domain_param}.

\begin{table}
\centering
\captionsetup{justification=centering}
\caption{Parameters for the Nested Ellipsoids in Example 2.}
\label{ellipsoid_domain_param}
\renewcommand{\arraystretch}{1.3}
\begin{tabular}{c | c c c} \toprule
\textbf{Ellipsoid Name} & \textbf{center ($\bm{c}$)} & $\bm{A}$ & $\bm{\sigma}$\\ \hline
$\Omega_1$ & $(1.3, 0, \dots, 0)$ & 
$\text{diag}(A) = [1, 0.9, 1.2, \dots, 1.2]$
& $1.1$ \\ \hline
$\Omega_2$ & $(1.4, 0, \dots, 0)$ & 
\makecell{$\text{diag}(A) = [1.5, 3.2, 2, \dots, 2]$  \\ with $A(1, 2) = A(2, 1) = -1$}
& $1.3$ \\ \hline
$\Omega_3$ & $(1.6, 0, \dots, 0)$ & 
\makecell{$\text{diag}(A) = [8, 6, 4, \dots, 4]$  \\ with $A(1, 2) = A(2, 1) = -2$}
& $0.9$ \\ \bottomrule
\end{tabular}
\end{table}

\begin{figure}[htbp]
\centering
\begin{tabular}{@{}c@{}c c c c c}
% title
\multirow{3}{*}{
  \rotatebox[origin=c]{90}{%
    \parbox[c]{\dimexpr 0.5\columnwidth\relax}{\centering\vfill\textbf{WoI}\vfill}%
  }
} 
& \textbf{Ground Truth} & $\mathcal{W} = 10^6$ & $\mathcal{W} = 5 \times 10^6$ & $\mathcal{W} = 10^7$ & \\
% First Row: estimation
{} & \includegraphics[width=0.20\columnwidth]{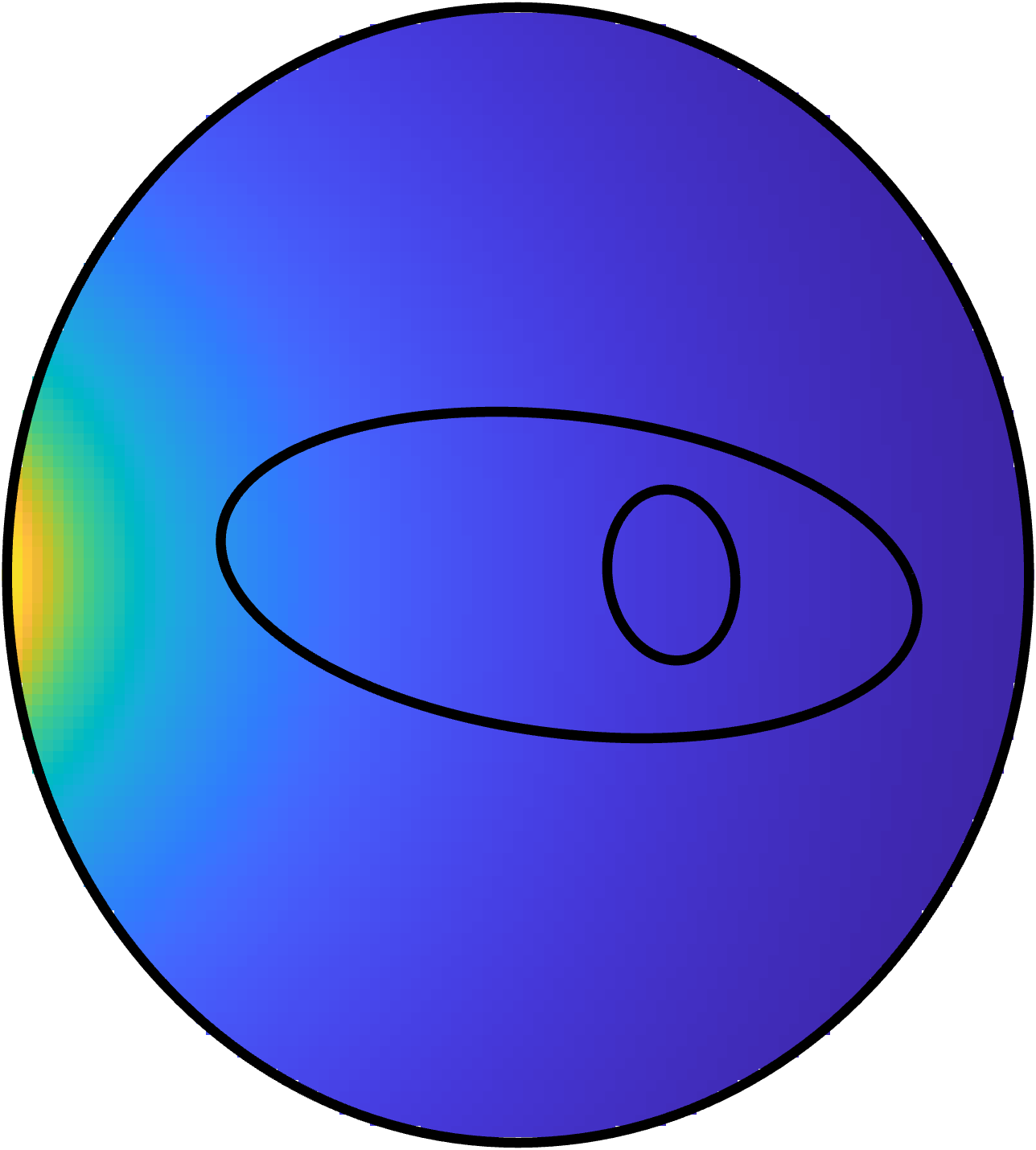} &
\includegraphics[width=0.20\columnwidth]{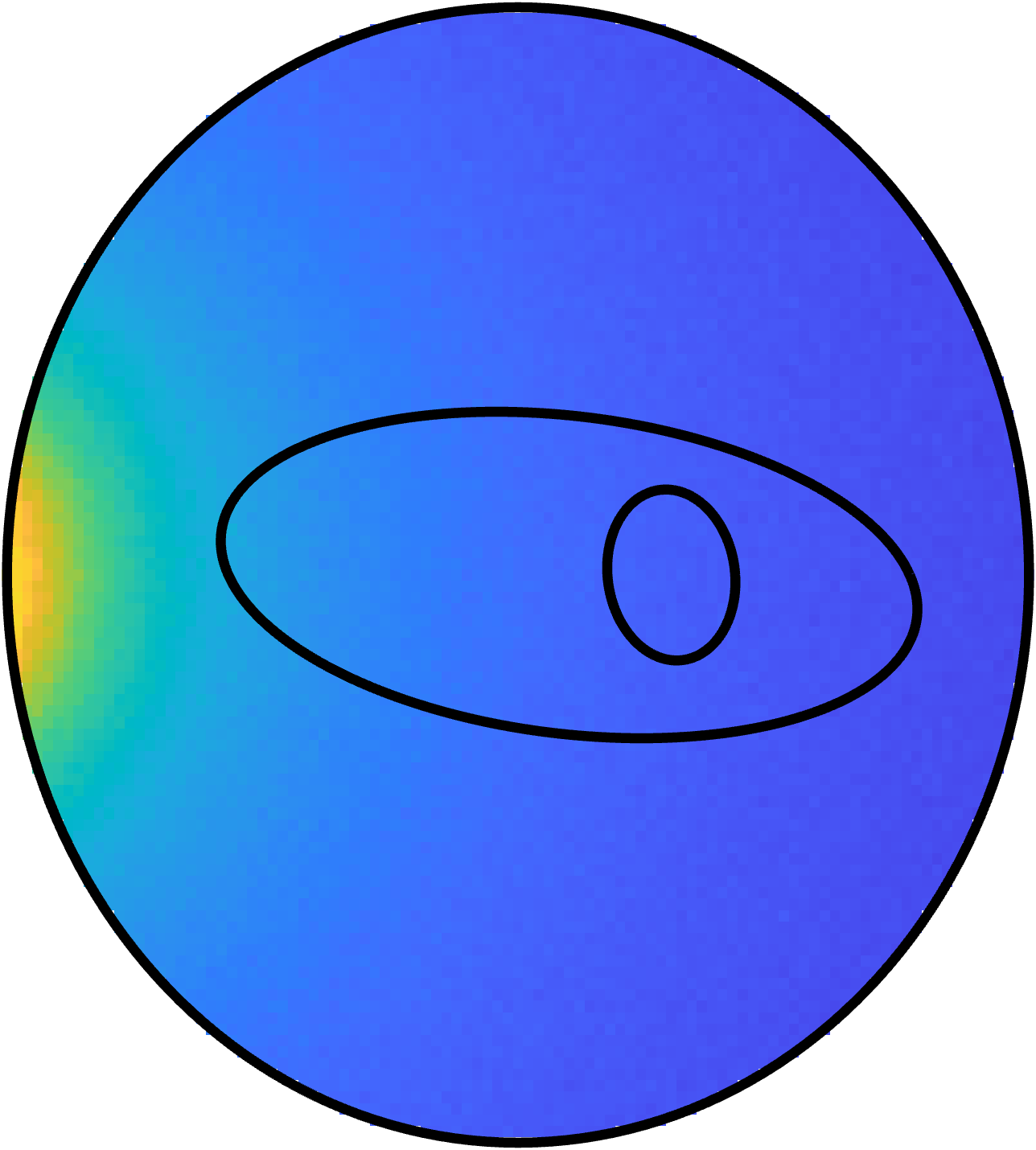}  &
\includegraphics[width=0.20\columnwidth]{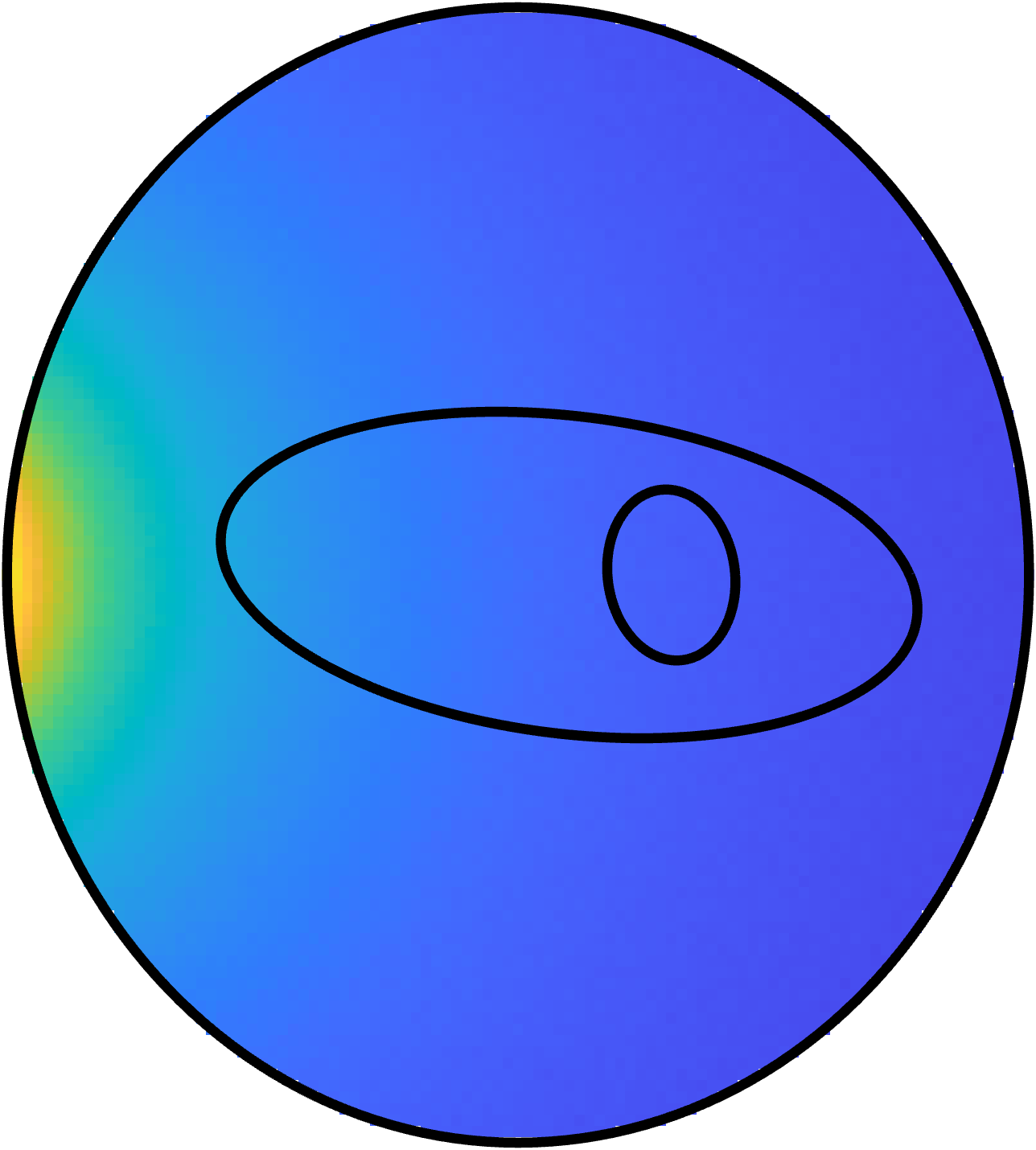} &
\includegraphics[width=0.20\columnwidth]{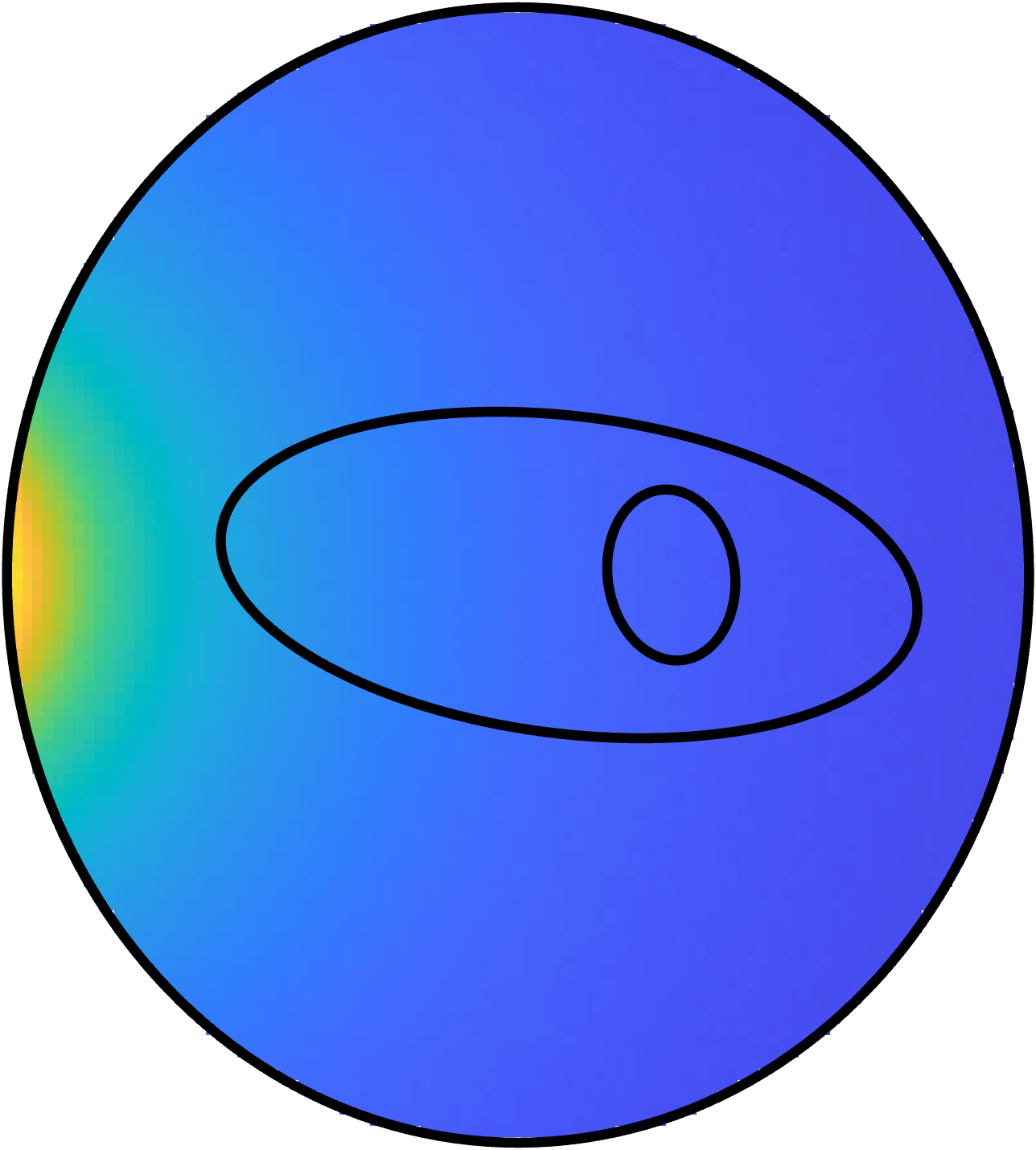} &
\includegraphics[width=0.034\columnwidth]{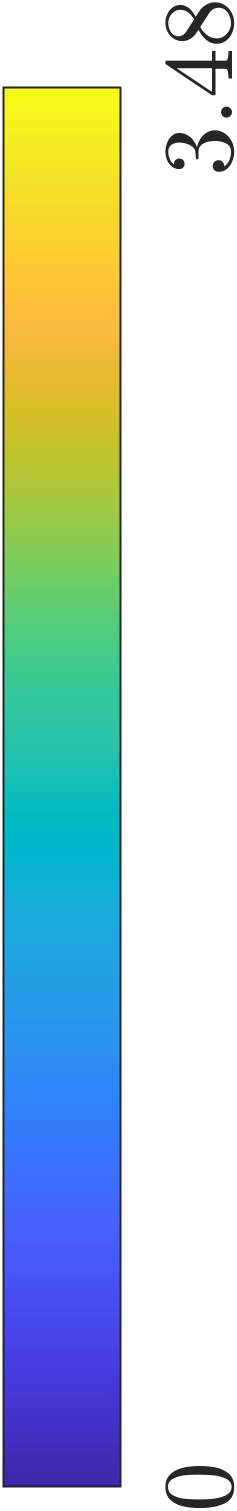}\\

% Second Row: pointwise estimation
{} & \includegraphics[width=0.25\columnwidth]{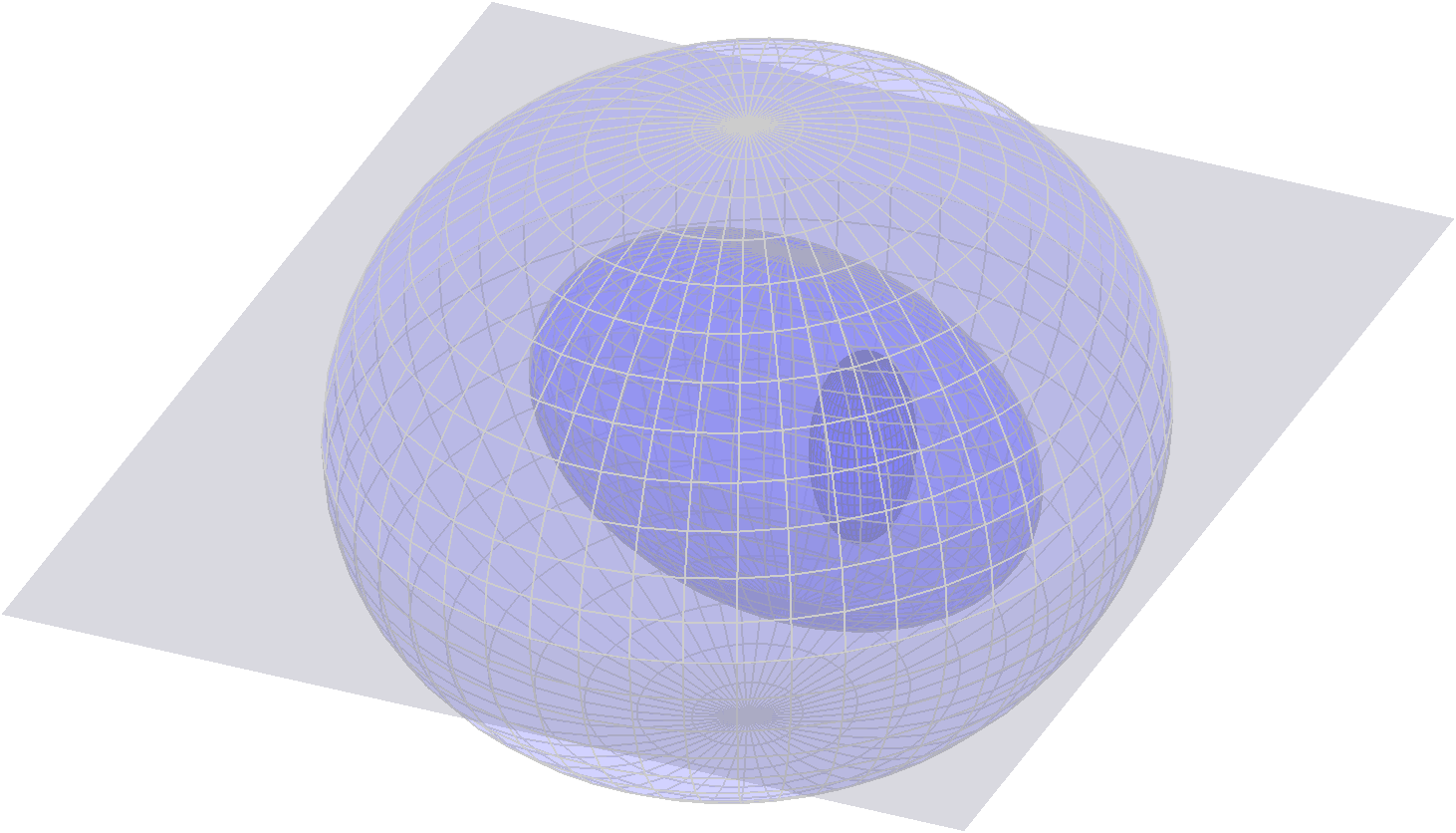} &
\includegraphics[width=0.2\columnwidth]{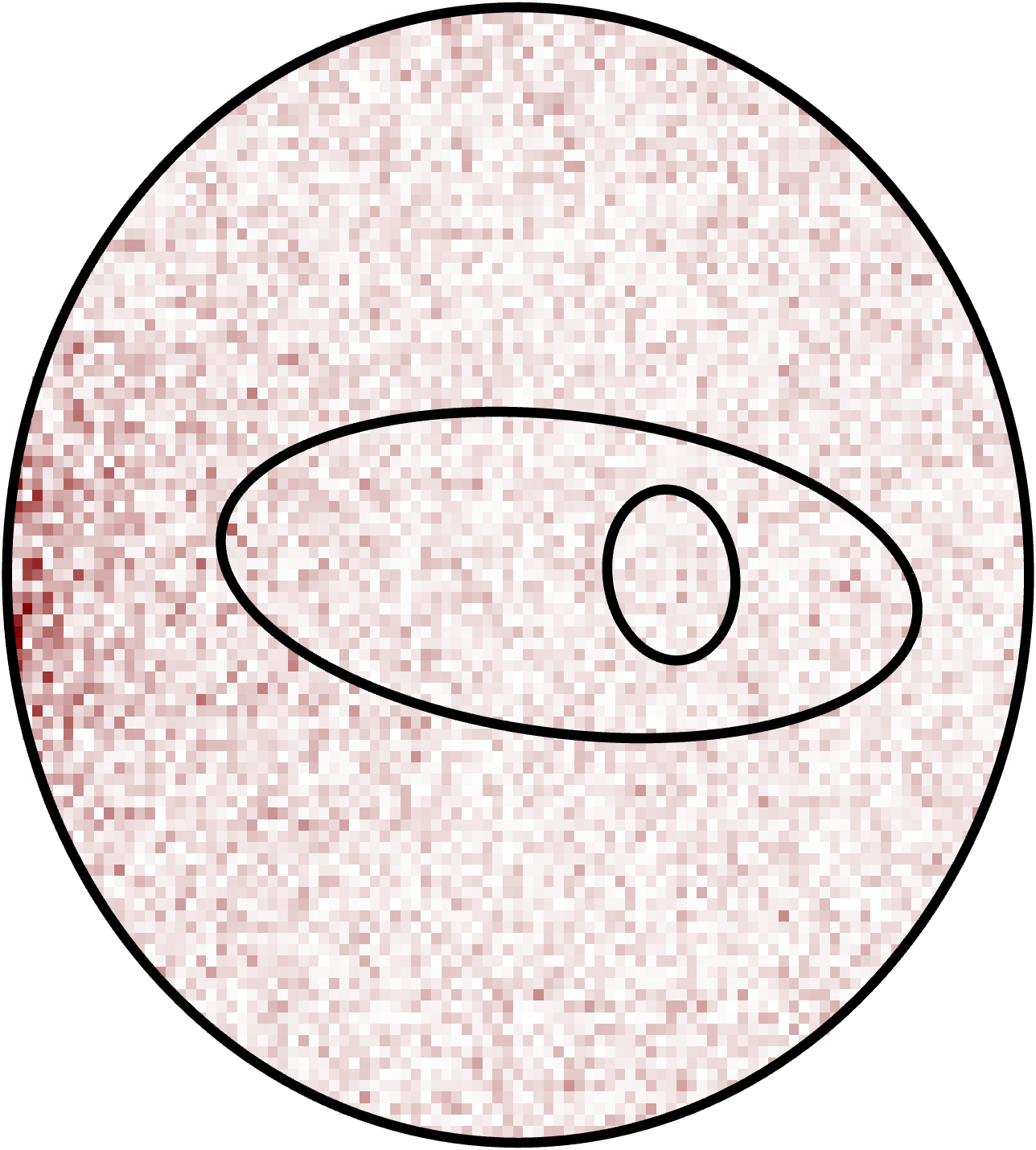} &
\includegraphics[width=0.2\columnwidth]{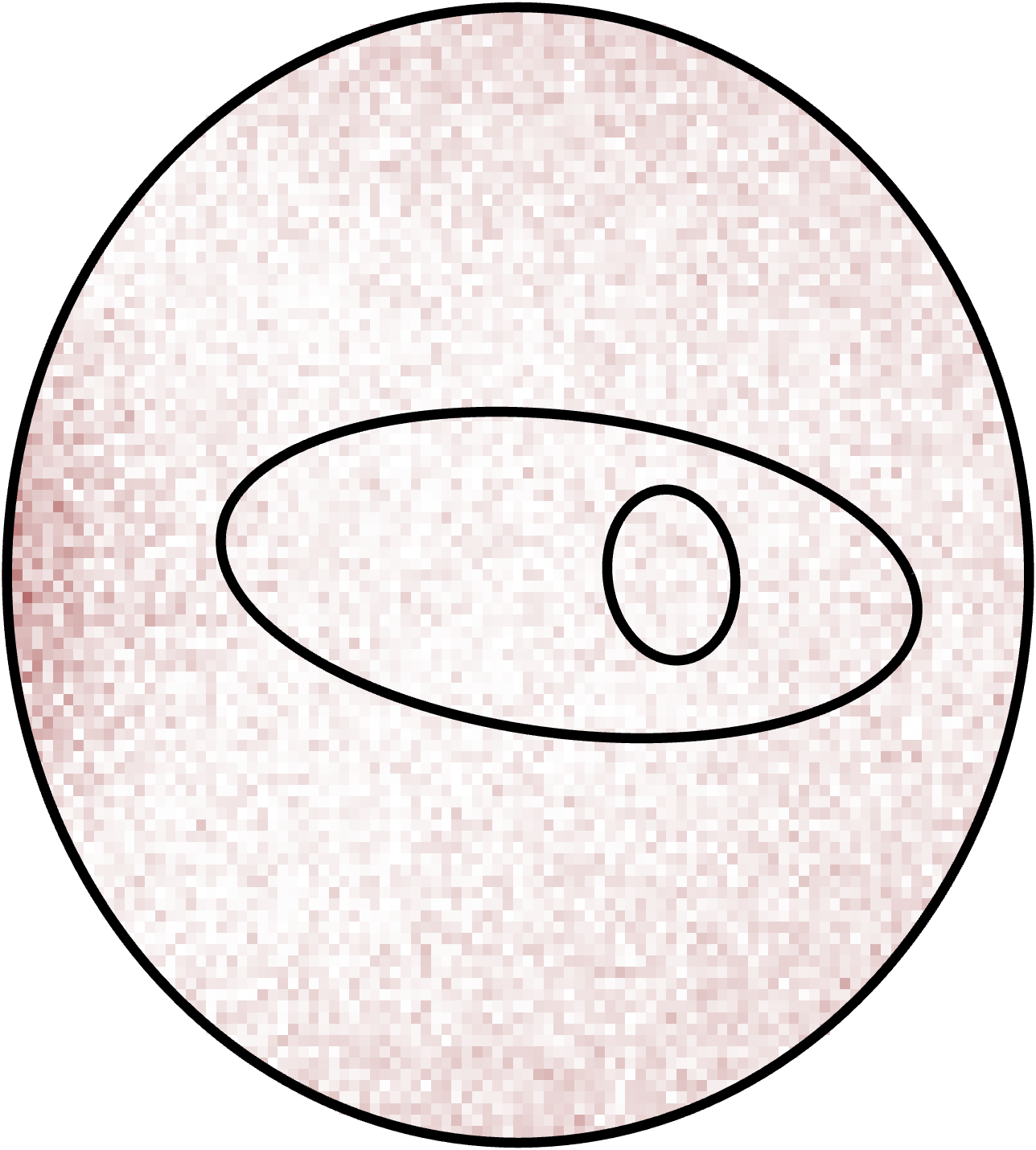} &
\includegraphics[width=0.2\columnwidth]{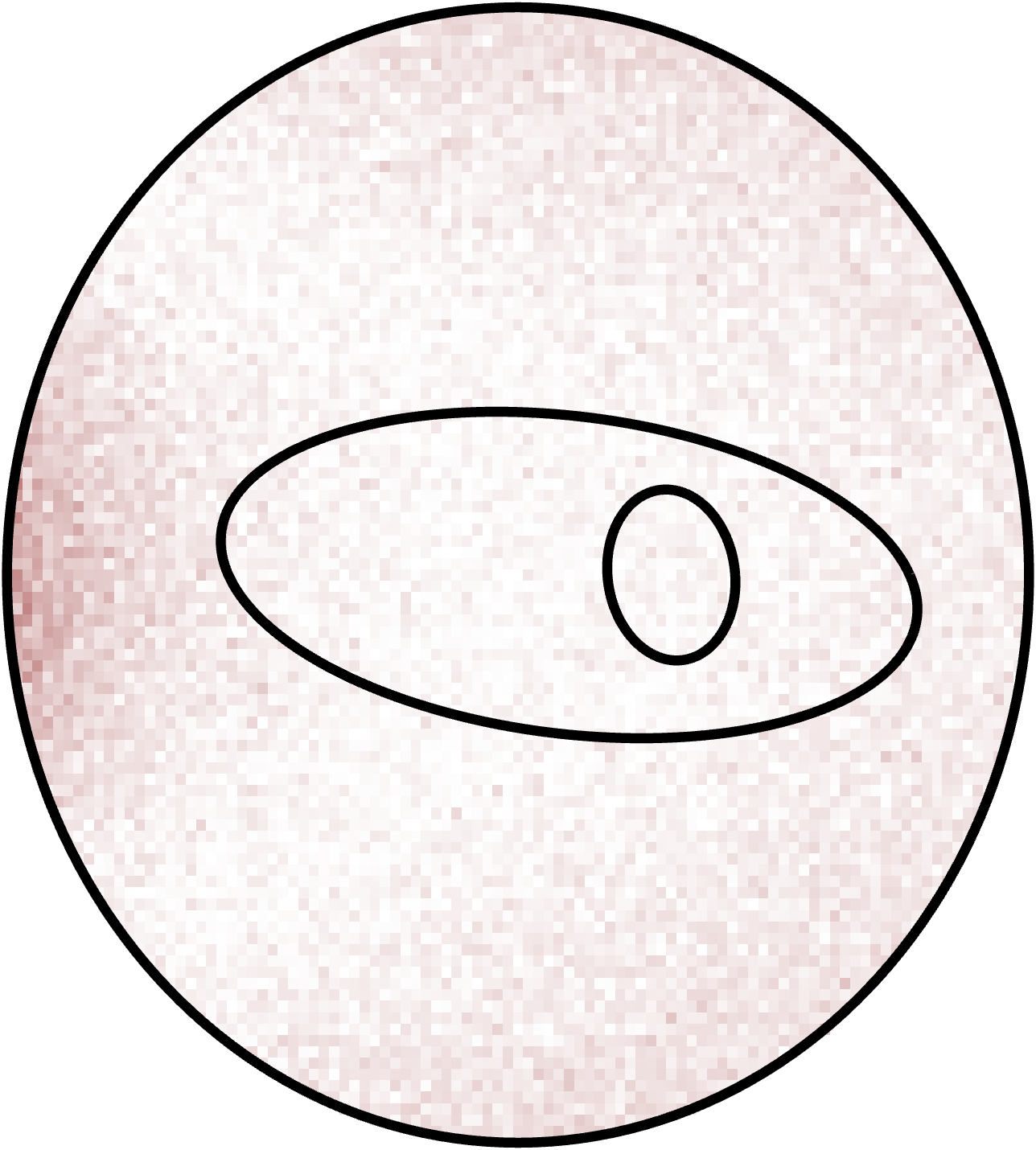} &
\includegraphics[width=0.034\columnwidth]{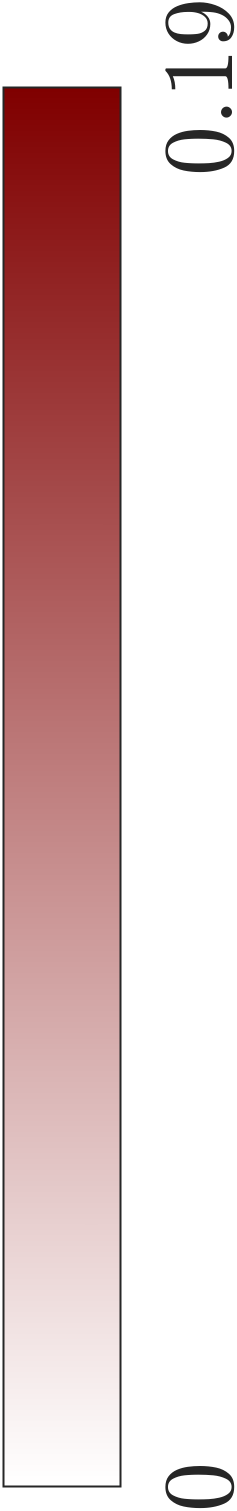}\\

% relative L2 error
{} & \textbf{$\bm{L_2}$ Error} & $0.029$ & $0.019$ & $0.017$ & \\
% relative L2 error
{} & \textbf{Relative $\bm{L_2}$ Error} & $3.10\%$ & $2.09\% $ & $1.79\% $ & \\
%===============================================
\hline \hline \\
%===============================================
% title
\multirow{3}{*}{
  \rotatebox[origin=c]{90}{%
    \parbox[c]{\dimexpr 0.5\columnwidth\relax}{\centering\vfill\textbf{Var-reduced WoI}\vfill}%
  }
} 
& \textbf{Ground Truth} & $\mathcal{W} = 2 \times 10^6$ & $\mathcal{W} = 10^7$ & $\mathcal{W} = 2 \times 10^7$ & \\
% First Row: estimation
{} & \includegraphics[width=0.20\columnwidth]{figure/point_charge_potential/dim=3_grid_ground_truth.png} &
\includegraphics[width=0.20\columnwidth]{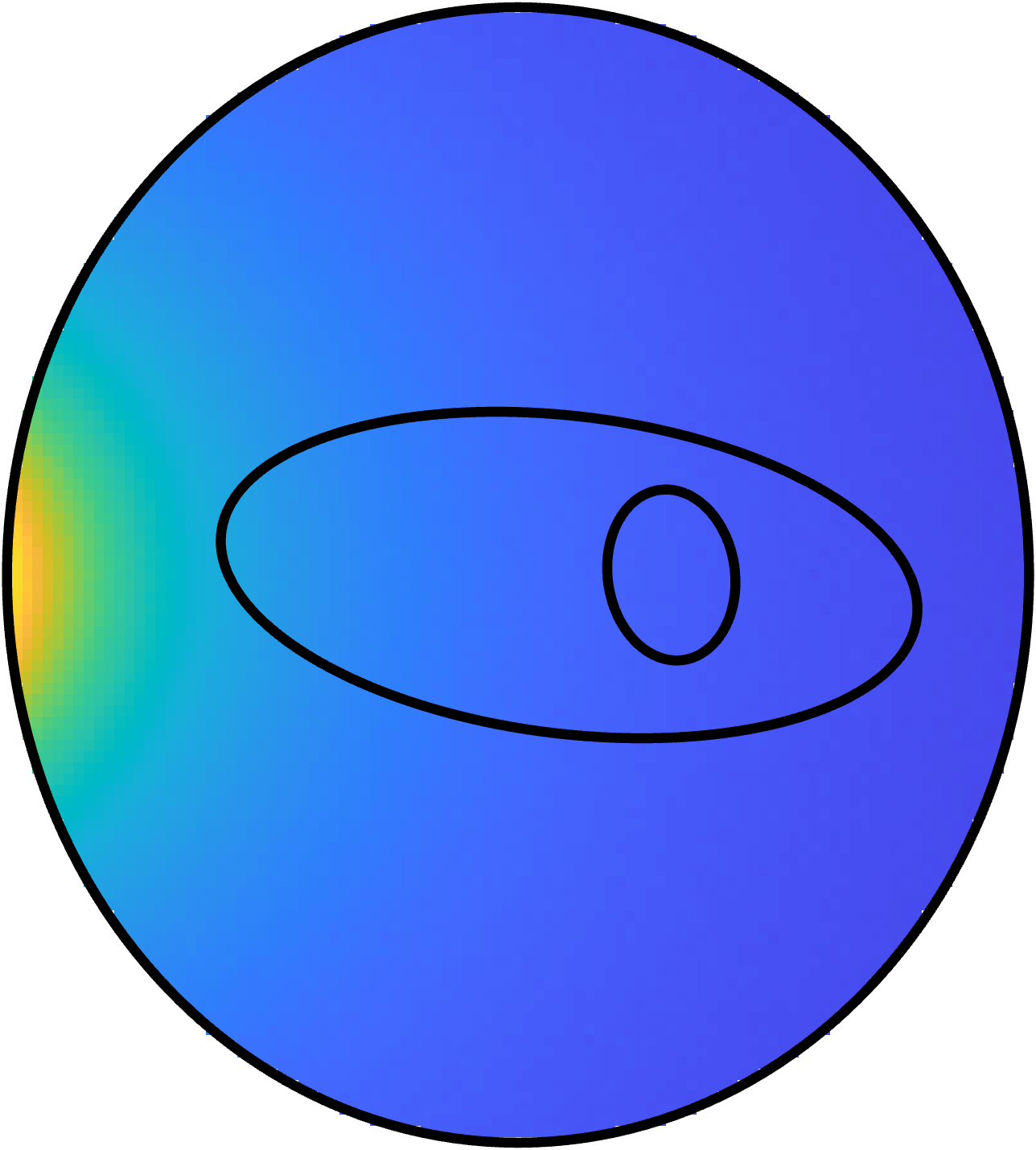}  &
\includegraphics[width=0.20\columnwidth]{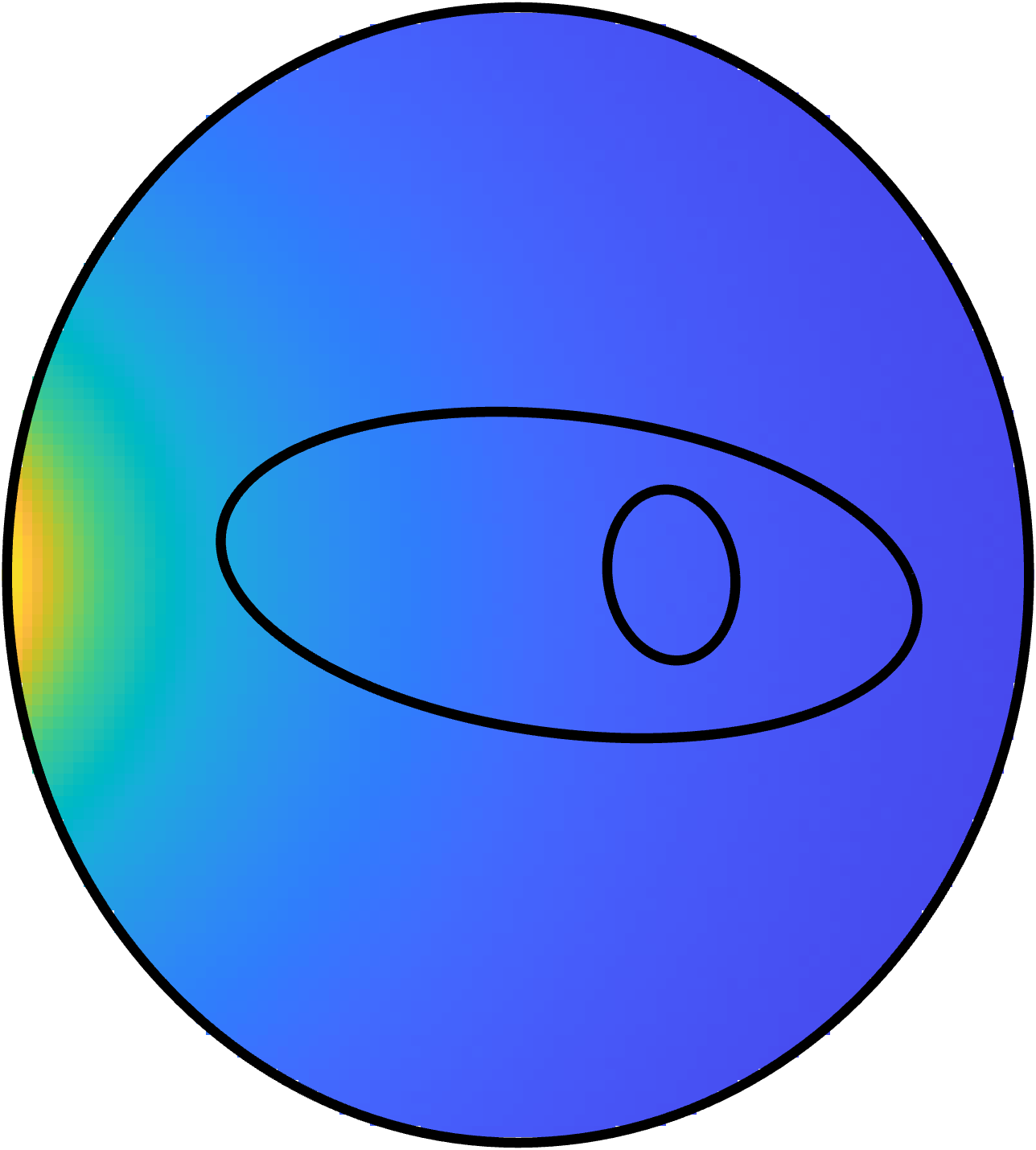} &
\includegraphics[width=0.20\columnwidth]{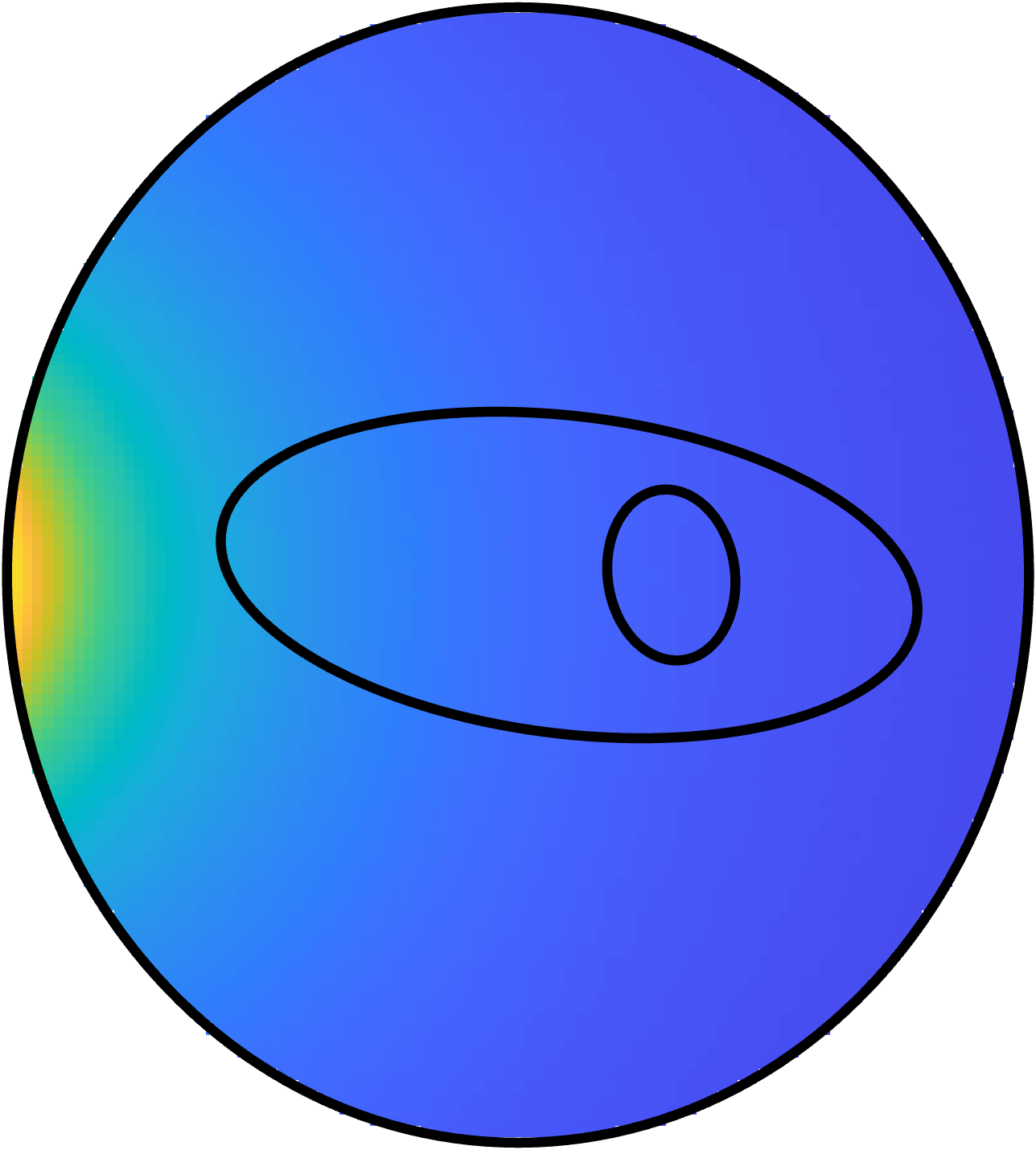} &
\includegraphics[width=0.034\columnwidth]{figure/point_charge_potential/uniform_estimation_cb.png}\\

% Second Row: pointwise estimation
{} & \includegraphics[width=0.25\columnwidth]{figure/point_charge_potential/domain.png} &
\includegraphics[width=0.2\columnwidth]{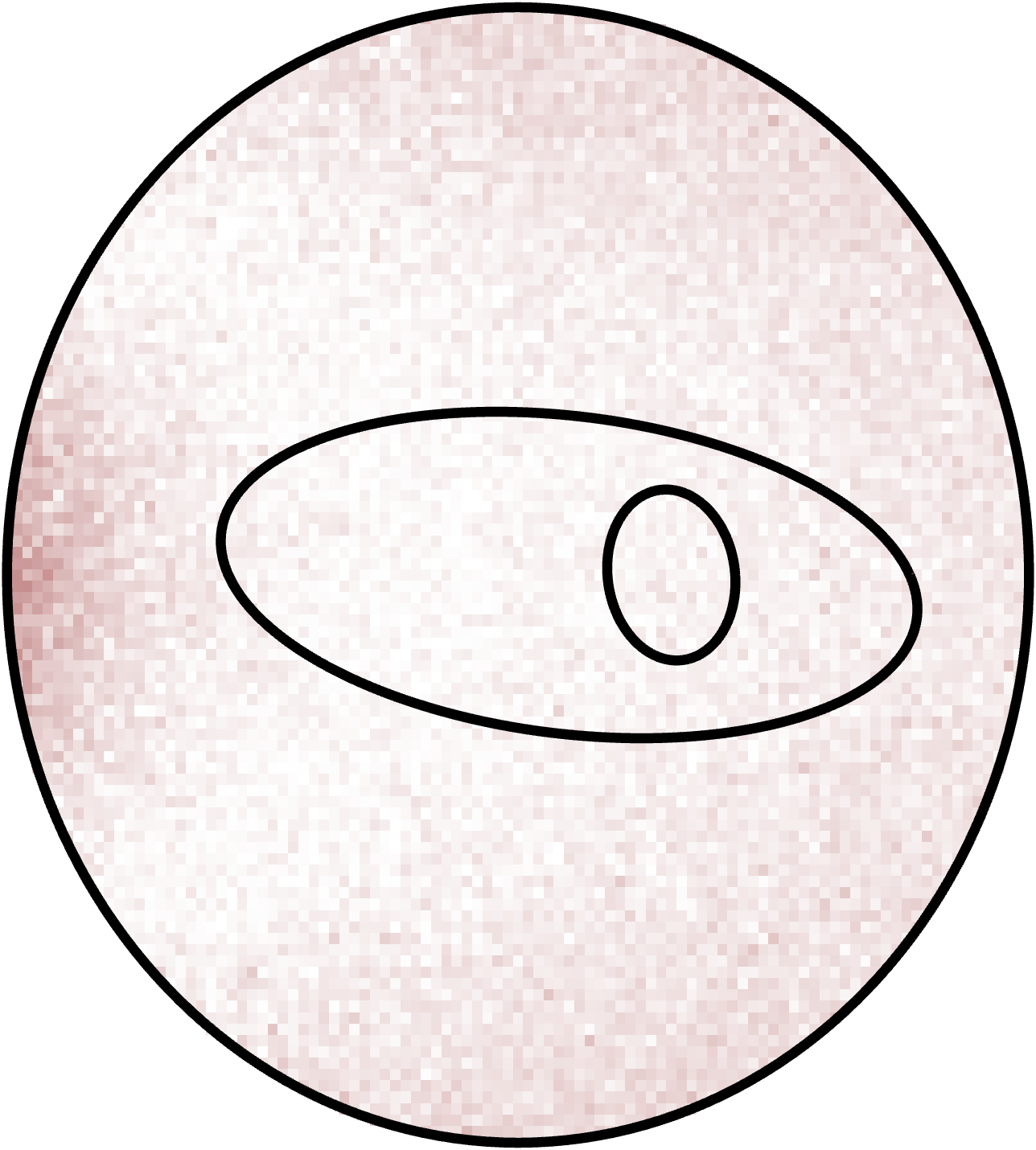} &
\includegraphics[width=0.2\columnwidth]{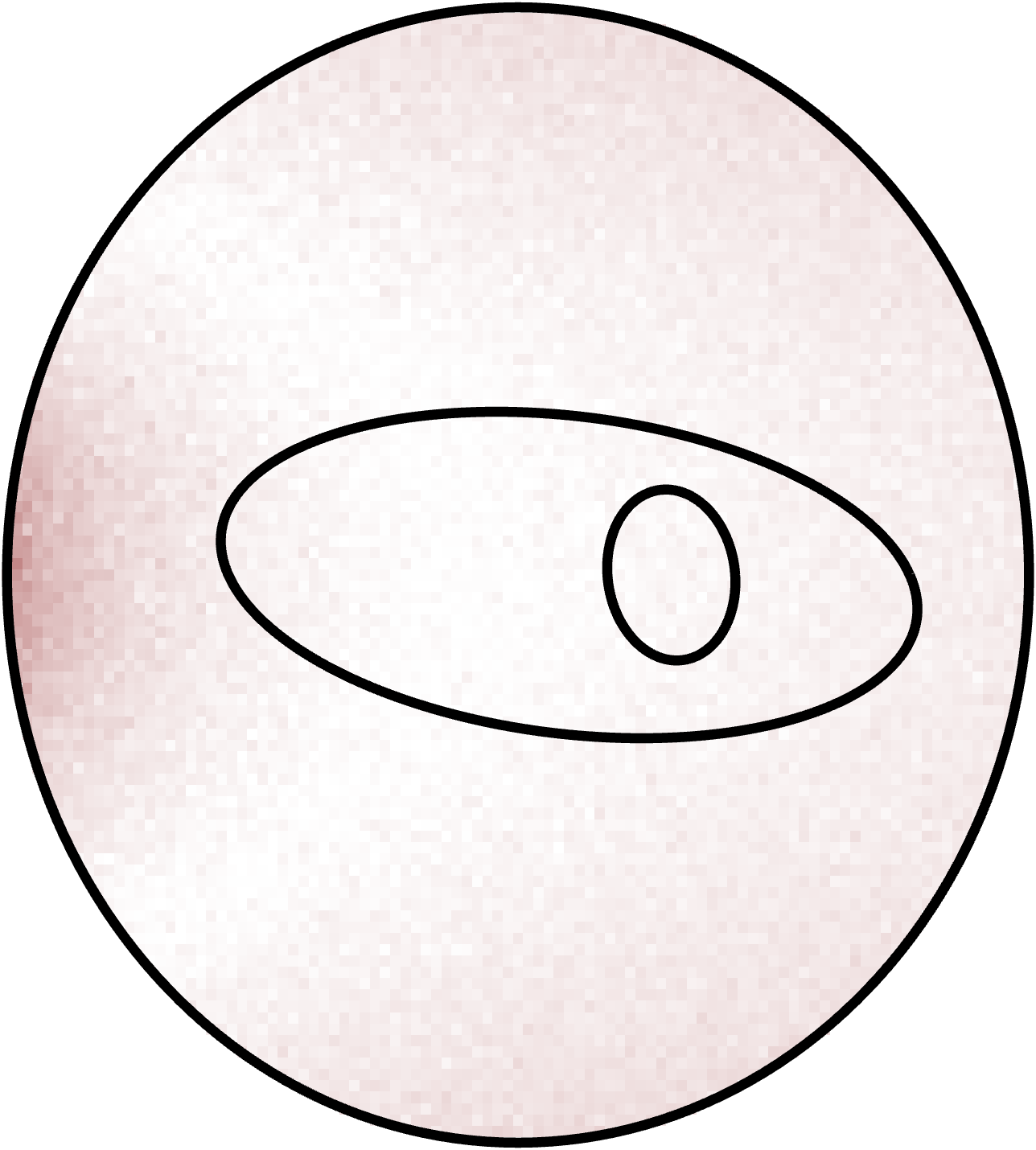} &
\includegraphics[width=0.2\columnwidth]{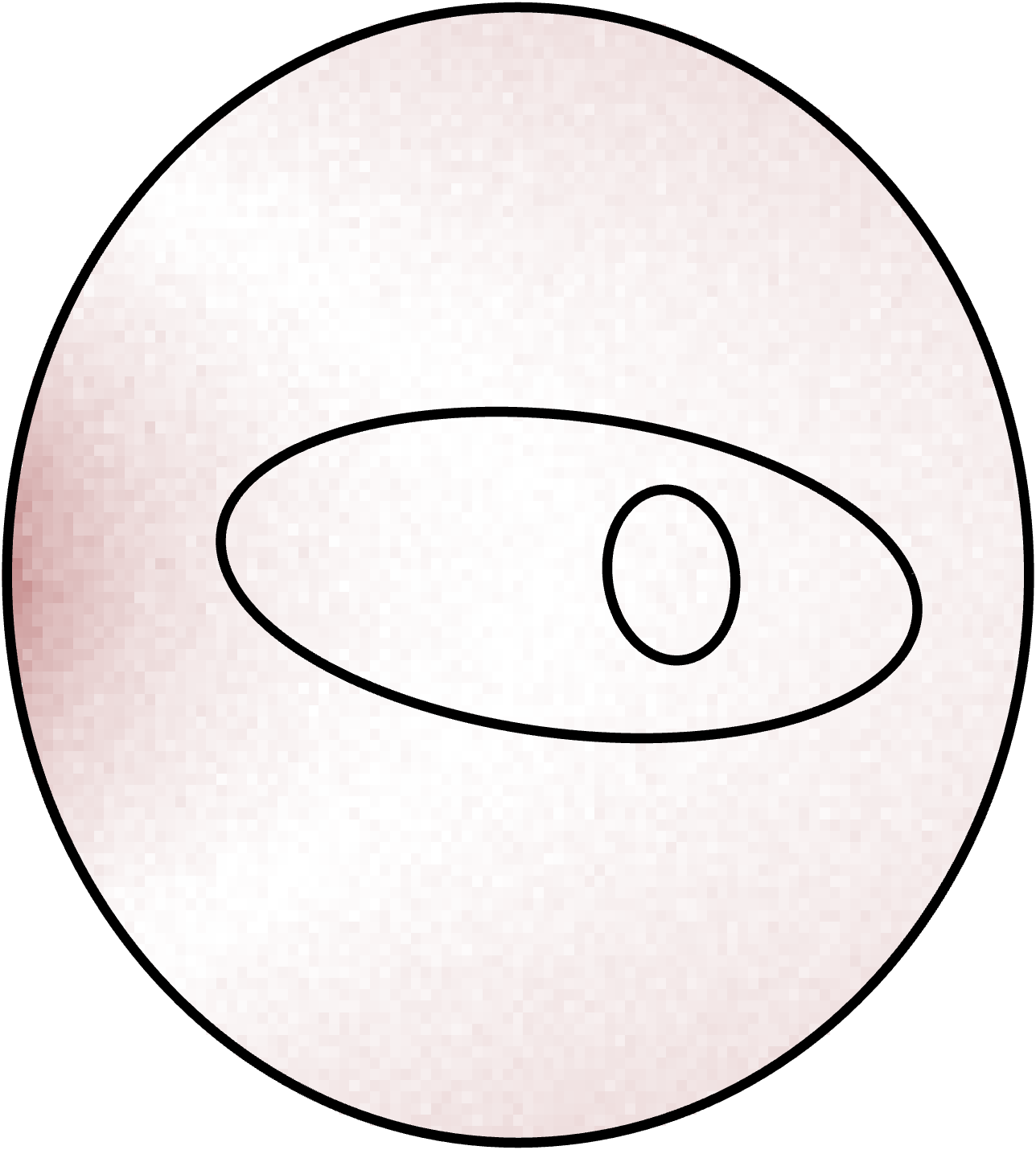} &
\includegraphics[width=0.034\columnwidth]{figure/point_charge_potential/uniform_error_cb.png}\\

% relative L2 error
{} & \textbf{$\bm{L_2}$ Error} & $0.0182$ & $0.0152$ & $0.0148$ & \\
% relative L2 error
{} & \textbf{Relative $\bm{L_2}$ Error} & $1.96\%$ & $1.64\% $ & $1.60\% $ & \\
\end{tabular}
\caption{The outer-most ellipsoid $\Omega_1$ with $\sigma_1 = 1.1$, $\Omega_2$ is the middle ellipsoid with $\sigma_2 = 1.3$, and $\Omega_3$ is the inner-most ellispoid with $\sigma_3 = 0.9$. The first row shows the solution of WoI or variance-reduced WoI with different numbers of walkers. The first figure in the second row shows the domain setup and the equatorial plane. The rest of plots in the second row display corresponding errors. The WoI estimator walks $M = 6$ steps.}\label{unit_point_charge3D}
\end{figure}

We visualize the 3D domain and its associated estimated solution on the equatorial plane in Figure \ref{unit_point_charge3D}. For domains in dimensions 4 and higher, we randomly sample $10^4$ points uniformly within the domain, with the $L_2$ error and the relative $L_2$ error reported in Table \ref{high_dim_ellipsoid_table}.

\begin{table}
\centering
\captionsetup{justification=centering}
\caption{$L_2$ error and relative $L_2$ error of estimation of point charge potential placed at the origin of $d$-dimensional ellipsoidal domains ($d = 3, 4, 5, 6$). Estimator: WoI.}\label{high_dim_ellipsoid_table}
\renewcommand{\arraystretch}{1.8}
\begin{tabular}{c c | c c c c c c} \toprule
\multicolumn{2}{c|}{\diagbox{\textbf{Dim} ($d$)} {\shortstack{\textbf{Walkers} $(\mathcal{W})$}}} & $10^6$ & $5 \times 10^6$ & $10^7$ & $2 \times 10^7$ & $3 \times 10^7$ & $4 \times 10^7$ \\ \hline
\multirow{2}{*}{\shortstack{\textbf{3} \\ $\bm{(M=5)}$}} & \textbf{$\bm{L_2}$ Error} & 0.022 & 0.020 & 0.012 & 0.01062 & 0.01063 & 0.01003 \\
{} & \textbf{Rel. $\bm{L_2}$ Error} & 2.63\% & 2.40\% & 1.44\% & 1.259\% & 1.261\% & 1.19\% \\ \hline
\multirow{2}{*}{\shortstack{\textbf{4} \\ $\bm{(M=6)}$}} & \textbf{$\bm{L_2}$ Error} & 0.056 & 0.033 & 0.030 & 0.027 & 0.028 & 0.026\\
{} & \textbf{Rel. $\bm{L_2}$ Error} & 7.02\% & 4.14\% & 3.81\% & 3.43\% & 3.47\%& 3.30\%\\ \hline
\multirow{2}{*}{\shortstack{\textbf{5} \\ $\bm{(M=7)}$}} & \textbf{$\bm{L_2}$ Error} & 0.235 & 0.122 & 0.074 & 0.065 & 0.061 & 0.059 \\
{} & \textbf{Rel. $\bm{L_2}$ Error} & 25.56\% & 13.33\% & 8.05\% & 7.11\% & 6.62\% & 6.39\%\\ \hline
\multirow{2}{*}{\shortstack{\textbf{6} \\ $\bm{(M=8)}$}} & \textbf{$\bm{L_2}$ Error} & 0.564 & 0.244 & 0.171 & 0.136 & 0.162 & 0.138 \\
{} & \textbf{Rel. $\bm{L_2}$ Error} & 56.57\% & 24.45\% & 17.20\% & 13.68\% & 16.30\% & 13.91\%\\
\bottomrule
\end{tabular}
\end{table}

\textbf{Example 3. Continuous Exact Solution with Discontinuity in Gradient.} Because of the piecewise nature of $\sigma(\bm{x})$, a common class of solutions to the interface problem consists of $C^0$ functions whose gradients are discontinuous. However, there are only a few setups for which an explicit analytical formula for such a solution is available. 

One such example is provided by Kyle Waller Bower~\cite{Bower_thesis}, which considers a 2D domain composed of two concentric circles. Let $\Omega_1$ be the outer circle of radius $r_1 = 1$, and $\Omega_2$ be the inner circle of radius $r_2 = \alpha$. Consider the problem such that the solution $u(r, \theta)$ satisfies
\begin{subequations}\label{C0_gradient_2D_problem}
\begin{empheq}[left=\empheqlbrace]{align}
 \Delta u &= 0 \quad \text{on $r \leq 1$, $r \neq \alpha$} \\
 % ========================
\partial_{\bm{r}}u &= \lambda \sin(m \theta) \quad \text{on $r = 1$} \\
 % ========================
 [u] &= 0 \quad \text{on $r = \alpha$}\\
 % ========================
 \sigma \partial_{\bm{r}^{-}} u &= \partial_{\bm{r}^{+}} u \quad \text{on $r = \alpha$},
\end{empheq}
\end{subequations}
in which $\sigma = \frac{\sigma_2}{\sigma_1}$ is the ratio of $\sigma(\bm{x})$ in disks of radius $\alpha$ and 1, and $m$ is a natural number. The analytical solution, obtained by separation of variables on the disk, is
\begin{equation}\label{C0_gradient_2D_analytical_soln}
    u(r, \theta) = 
    \begin{cases}
        \lambda A r^m \sin(m\theta), \quad & r \leq \alpha \\
        \lambda (Br^m + Cr^{-m}) \sin (m \theta), \quad & \alpha < r \leq 1
    \end{cases},
\end{equation}
in which $A = \frac{2}{m} \frac{1}{\alpha^{2m} (\sigma-1) + \sigma + 1}$, $B = \frac{A}{2} (\sigma + 1)$, and $C = B - \frac{1}{m}$.

\begin{figure}[htbp]
\centering
\begin{tabular}{@{}c@{}c c c c c}
% ============= WOI ==============
\multirow{3}{*}{
  \rotatebox[origin=c]{90}{%
    \parbox[c]{\dimexpr 0.5\columnwidth\relax}{\centering\vfill\textbf{WoI}\vfill}%
  }
} & 
\textbf{Ground Truth} & $\mathcal{W} = 10^6$ & $\mathcal{W} = 5 \times 10^6$ & $\mathcal{W} = 10^7$ & {} \\

% First Row: estimation
{} & \includegraphics[width=0.2\columnwidth]{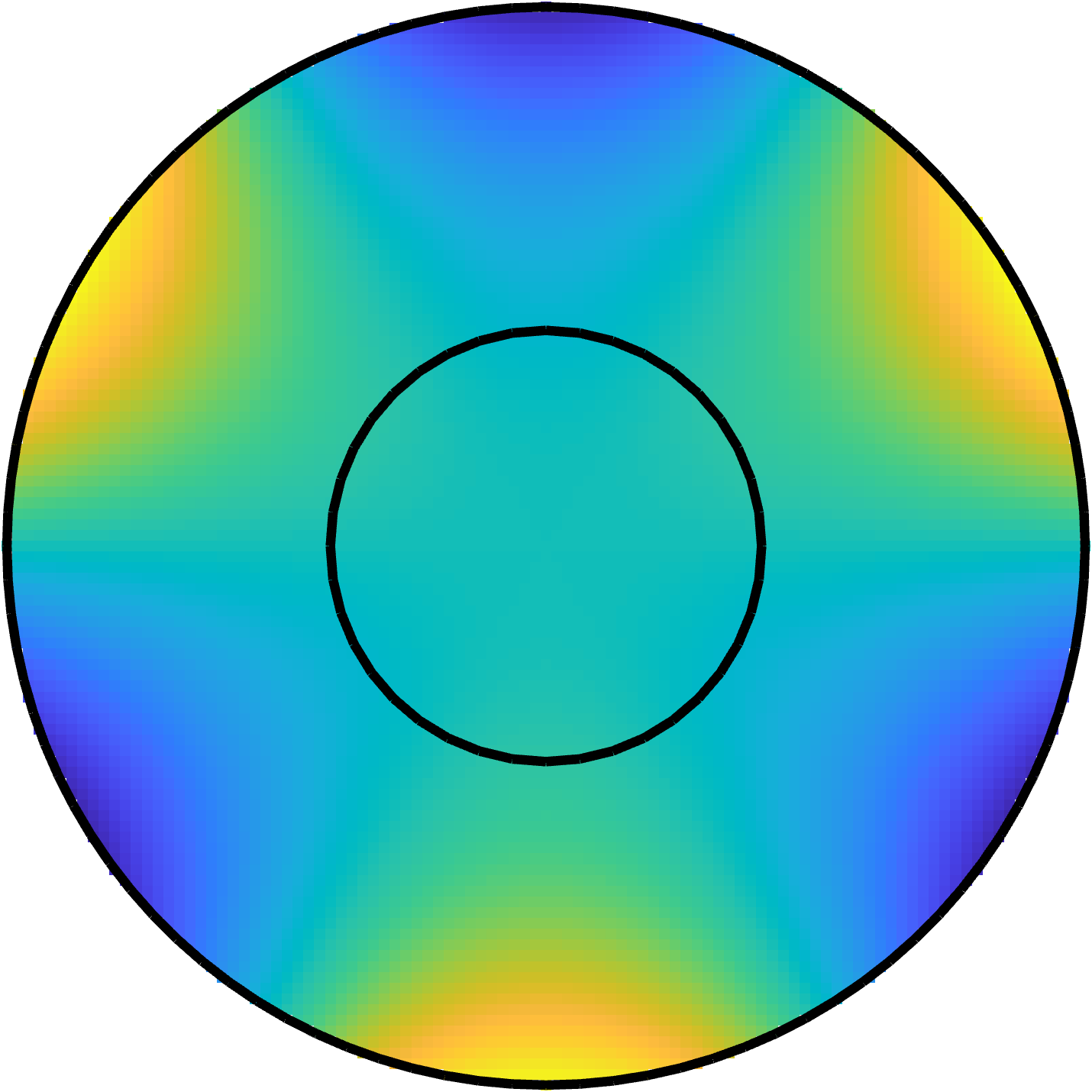} &
\includegraphics[width=0.2\columnwidth]{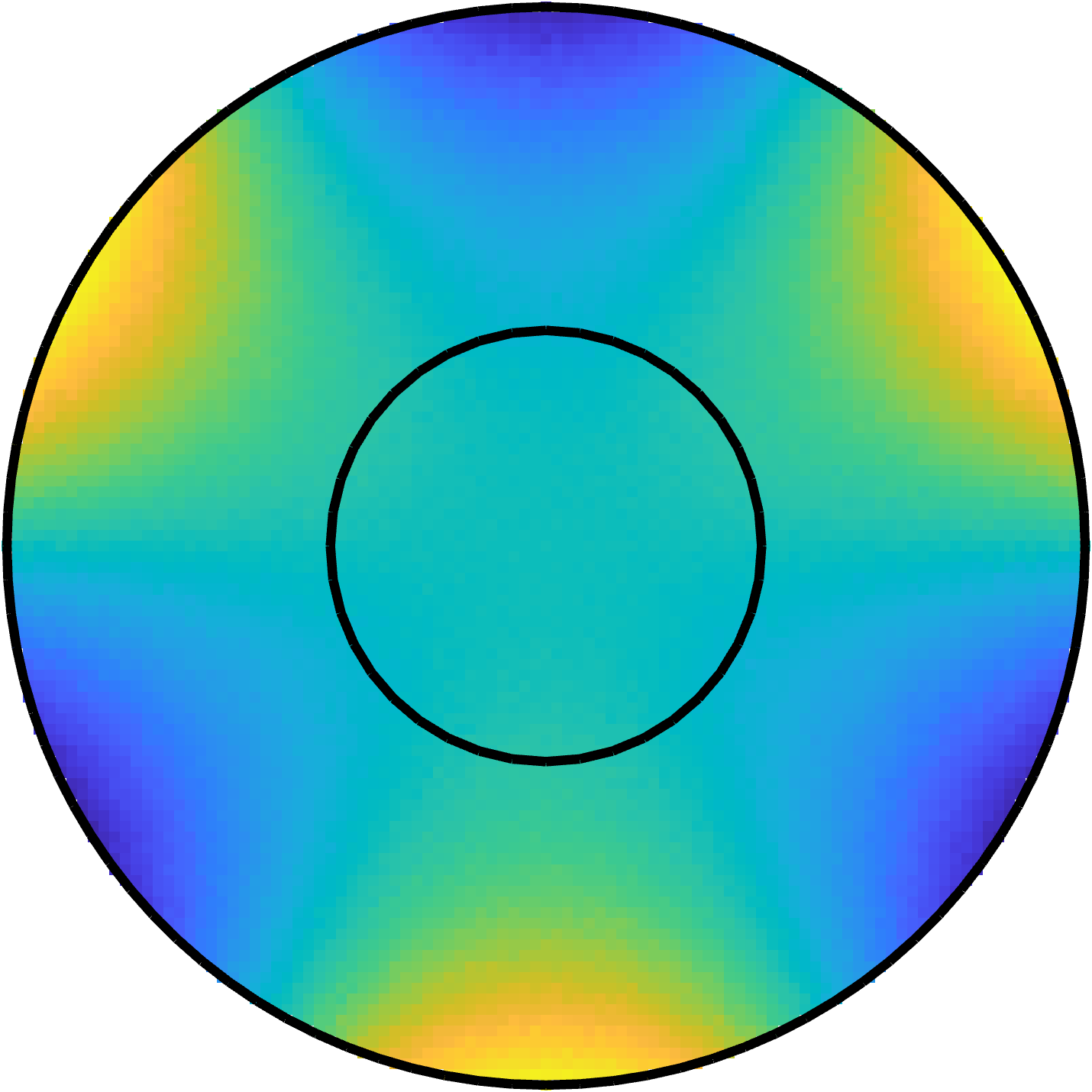}  &
\includegraphics[width=0.2\columnwidth]{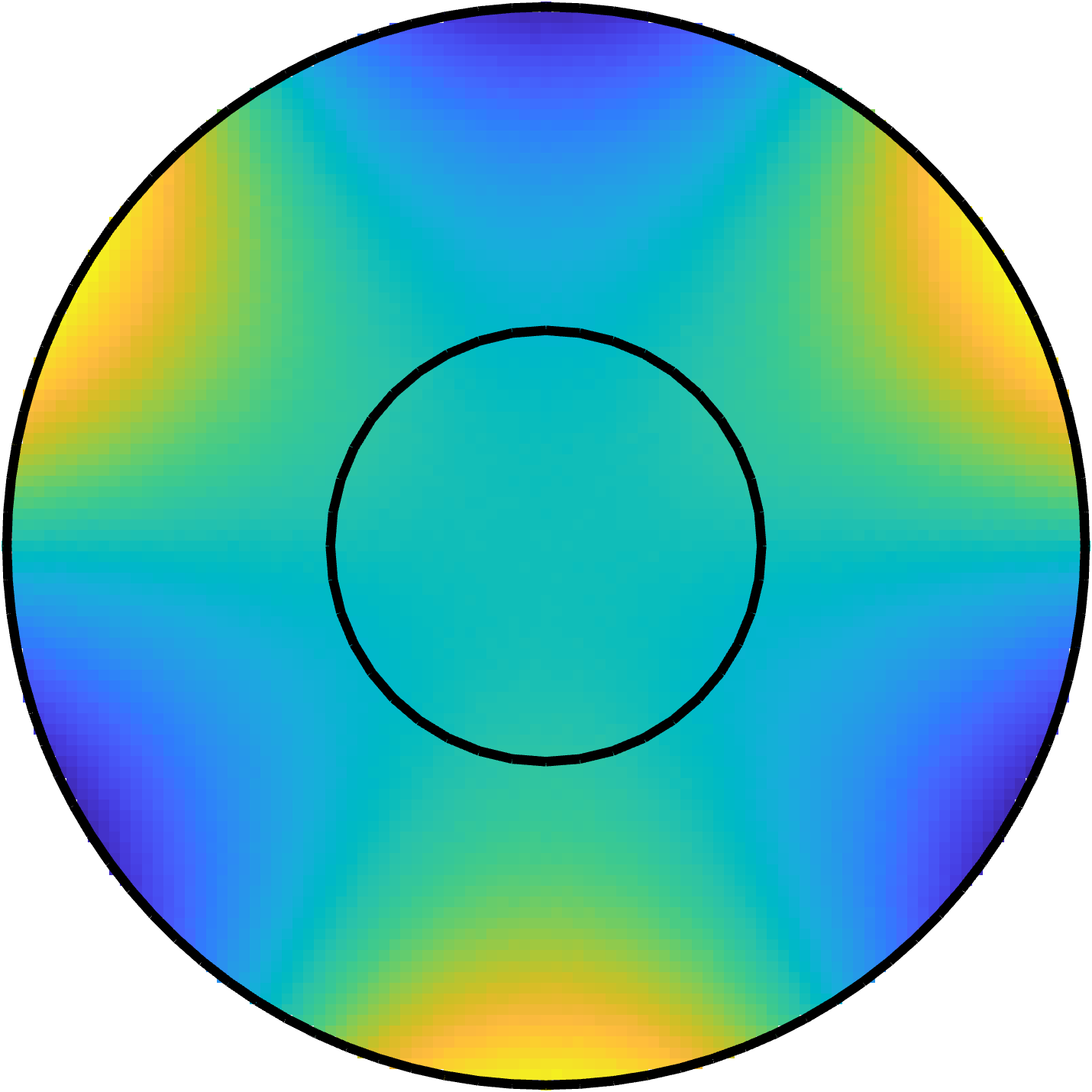} &
\includegraphics[width=0.2\columnwidth]{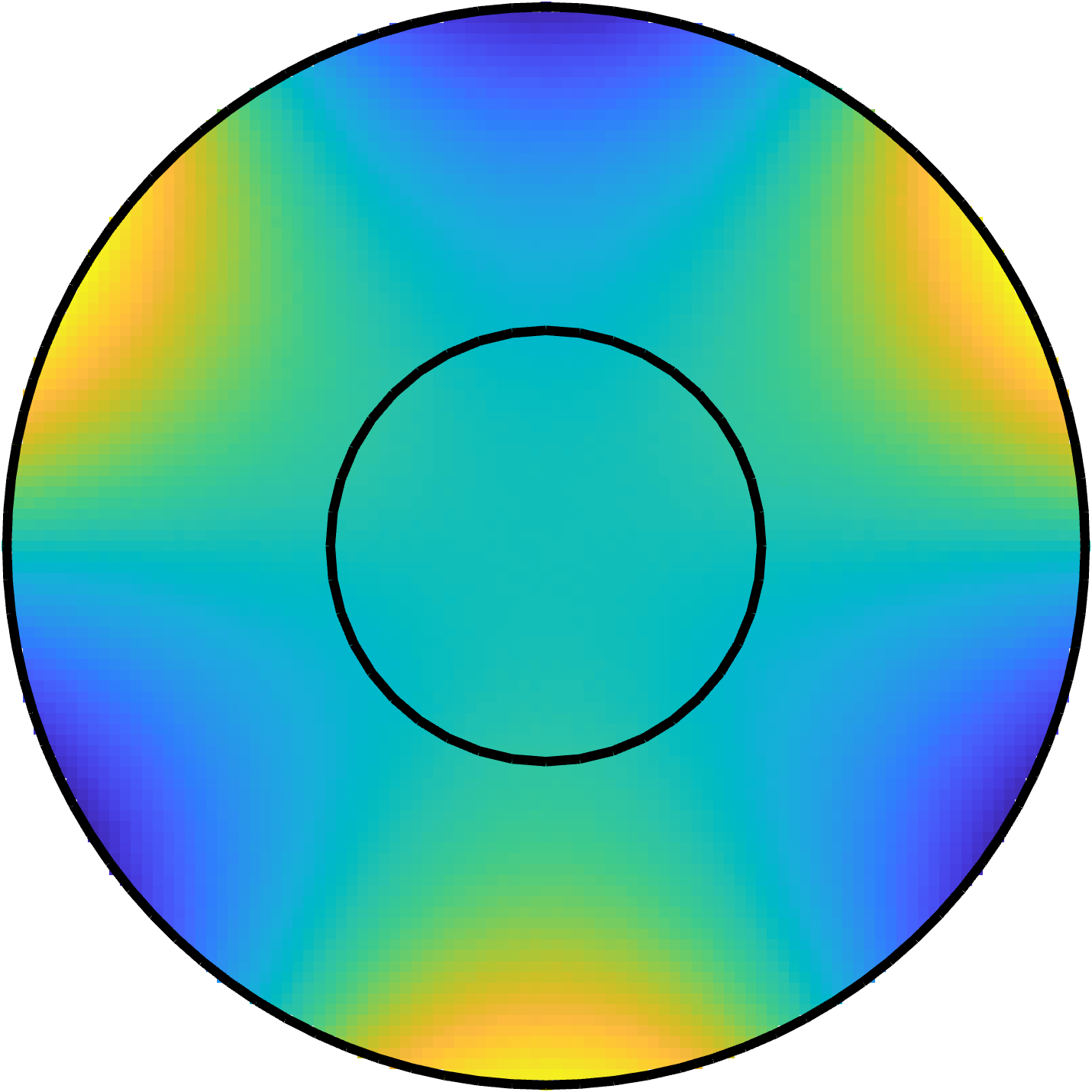} &
\includegraphics[width=0.03\columnwidth]{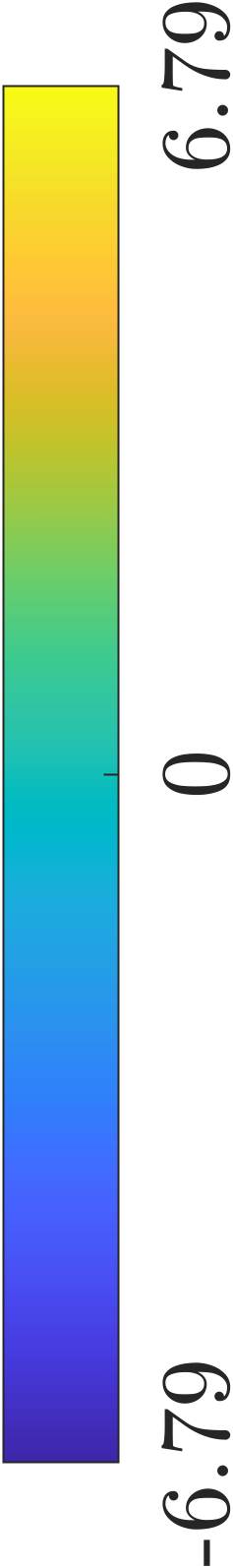} \\

% Second Row: pointwise estimation
{} & {} &
\includegraphics[width=0.2\columnwidth]{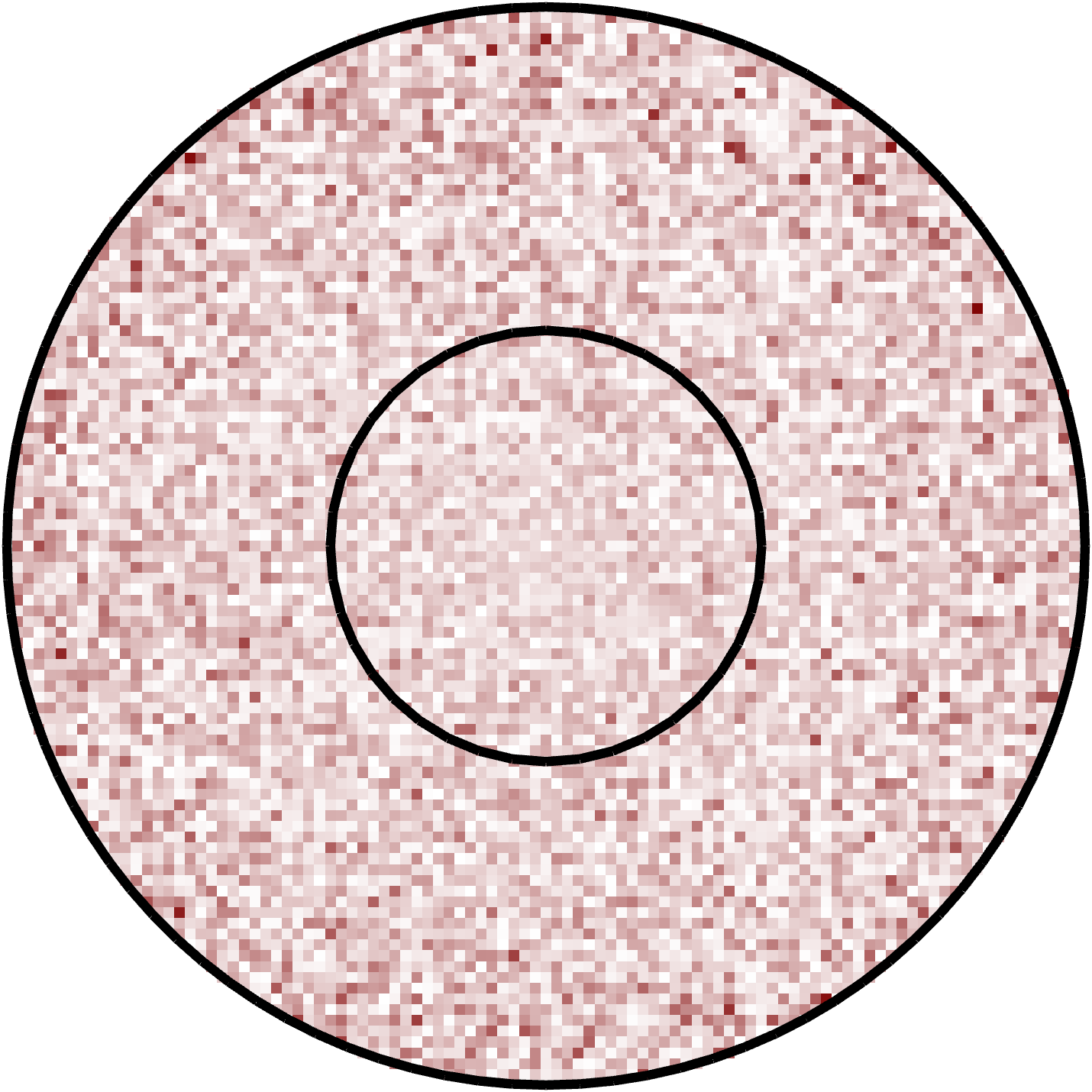} &
\includegraphics[width=0.2\columnwidth]{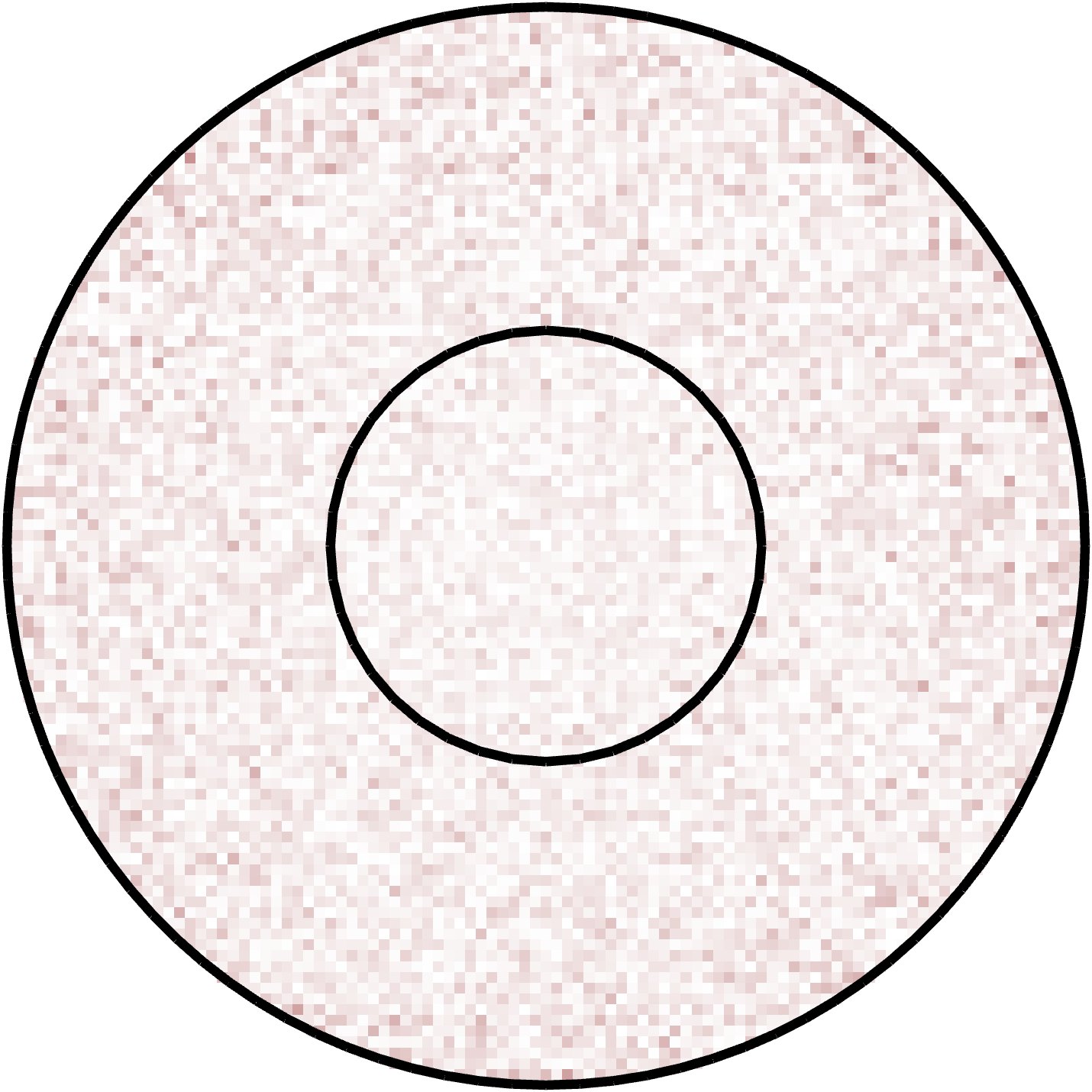} &
\includegraphics[width=0.2\columnwidth]{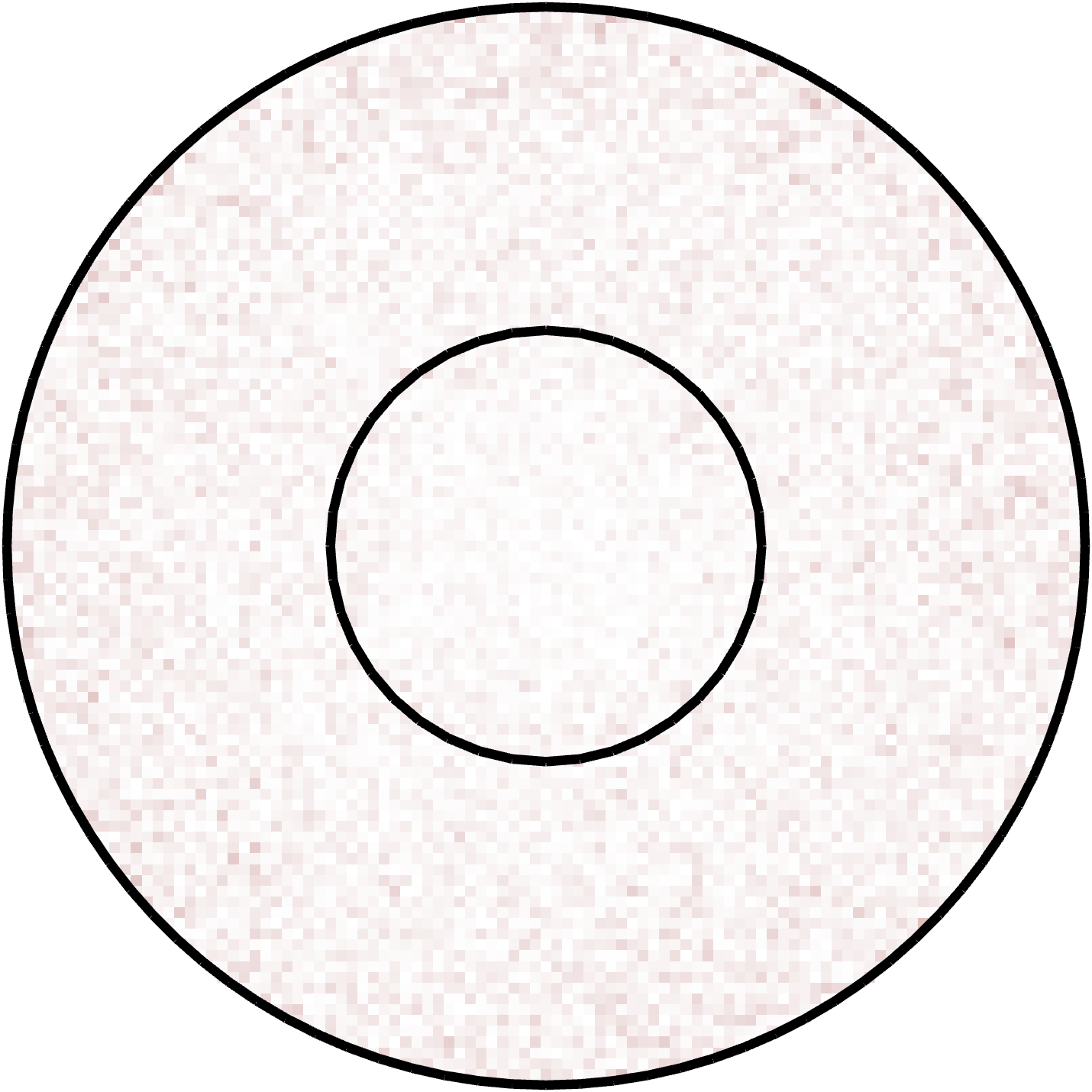} &
\includegraphics[width=0.03\columnwidth]{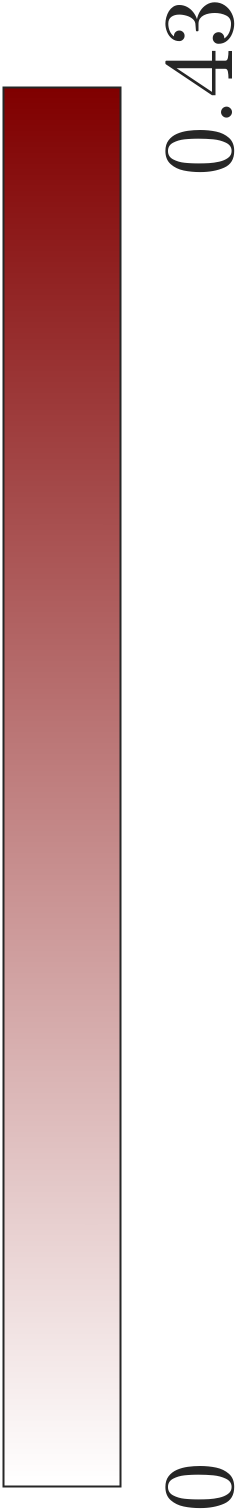}\\
{} & \textbf{$\bm{L_2}$ Error} & $0.107$ & $0.042$ & $0.025$ & {} \\
% relative L2 error
{} & \textbf{Relative $\bm{L_2}$ Error} & $4.51\%$ & $1.77\% $ & $1.07\% $ & {}\\
%===============================================
\hline \hline \\
%===============================================
% =========== variance-reduced WoI =============
\multirow{3}{*}{
  \rotatebox[origin=c]{90}{%
    \parbox[c]{\dimexpr 0.5\columnwidth\relax}{\centering\vfill\textbf{Var-reduced WoI}\vfill}%
  }
} & \textbf{Ground Truth} & $\mathcal{W} = 2 \times 10^6$ & $\mathcal{W} = 10^7$ & $\mathcal{W} = 2 \times 10^7$ & \\
% First Row: estimation
{} & \includegraphics[width=0.2\columnwidth]{figure/kink_in_gradient2D/N=2_kink_concentric_circle_ground_truth.png} &
\includegraphics[width=0.2\columnwidth]{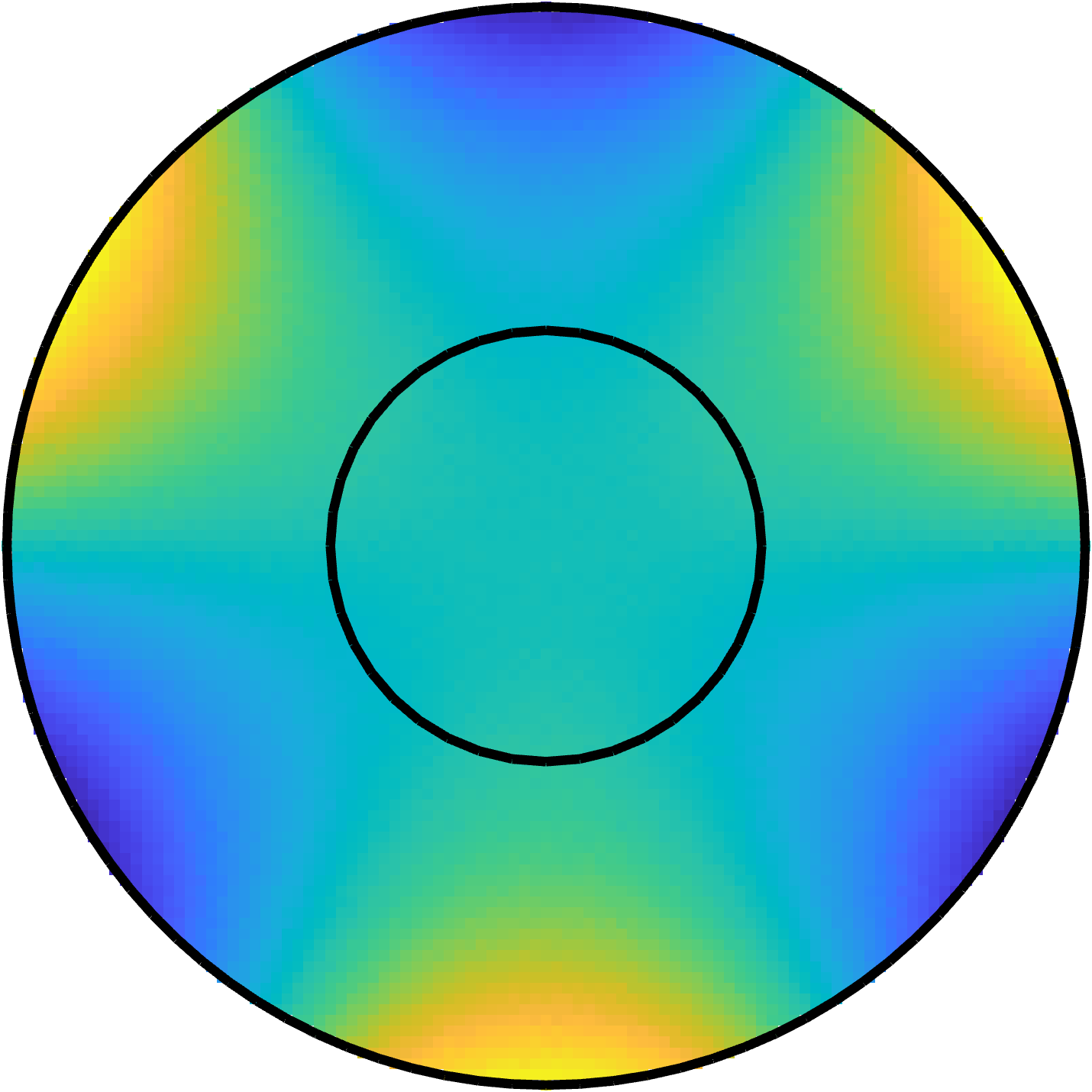}  &
\includegraphics[width=0.2\columnwidth]{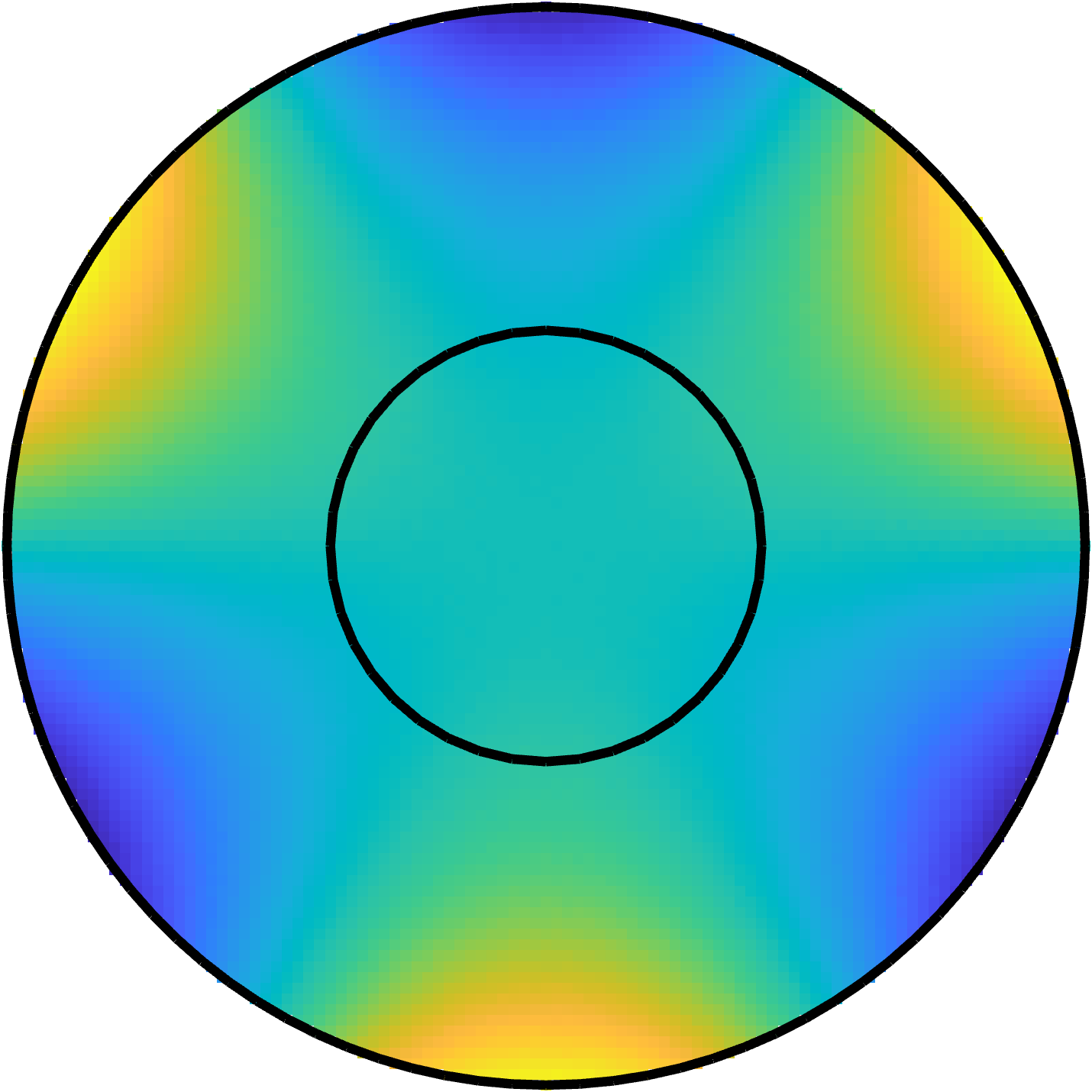} &
\includegraphics[width=0.2\columnwidth]{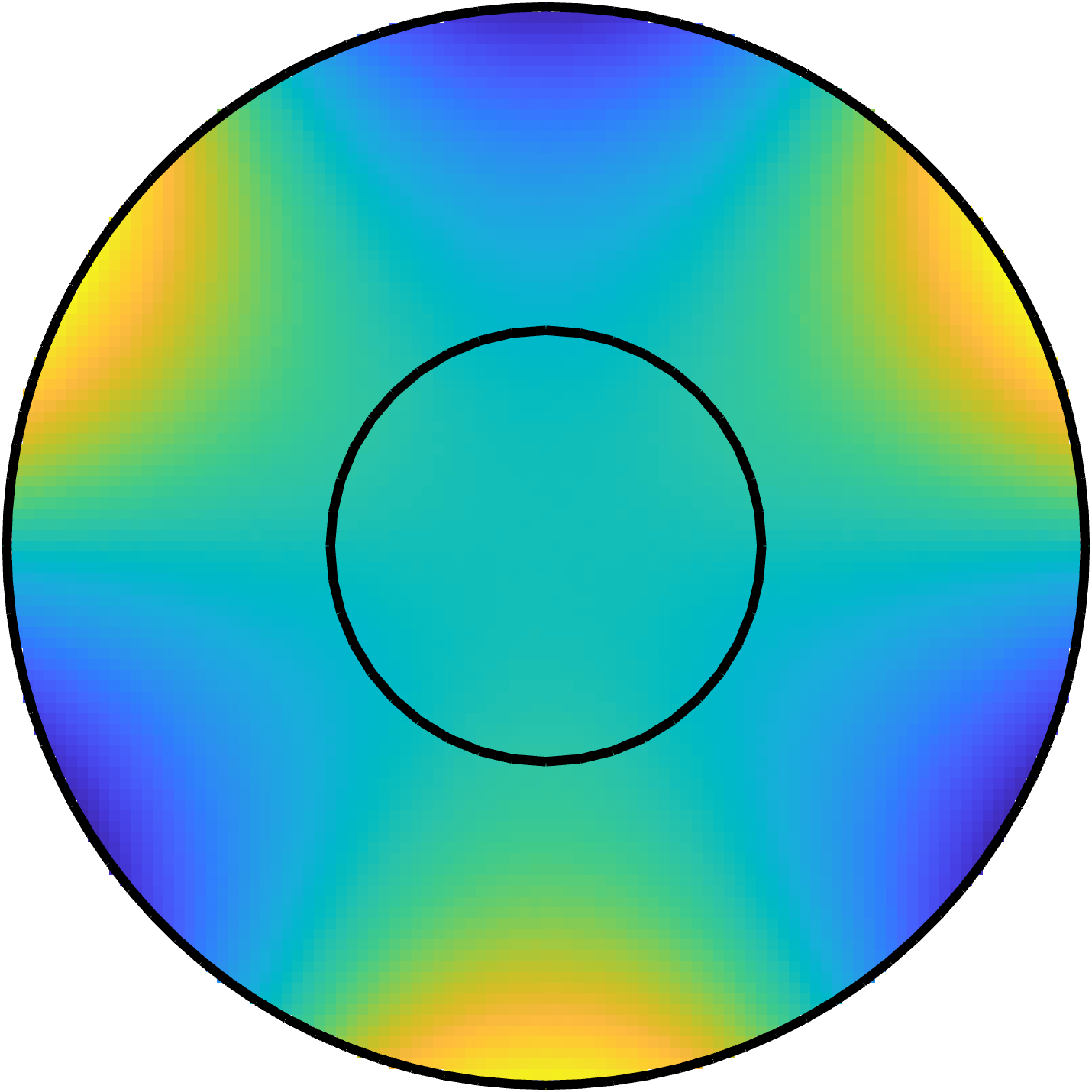} &
\includegraphics[width=0.03\columnwidth]{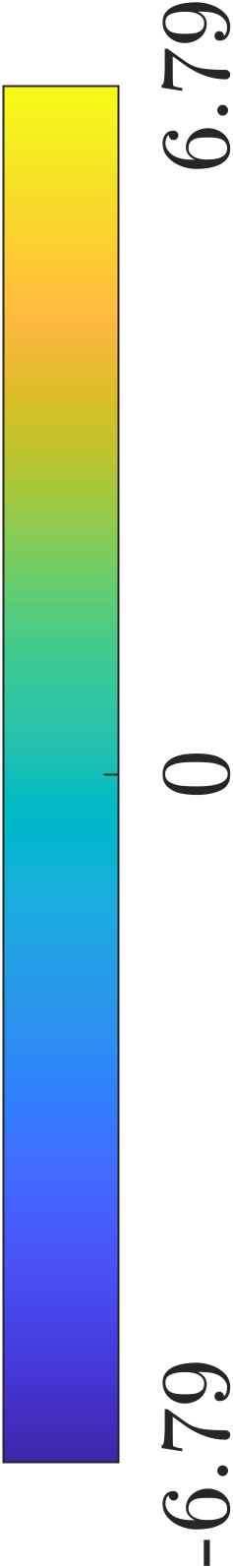}\\
% Second Row: pointwise estimation
{} & {} &
\includegraphics[width=0.2\columnwidth]{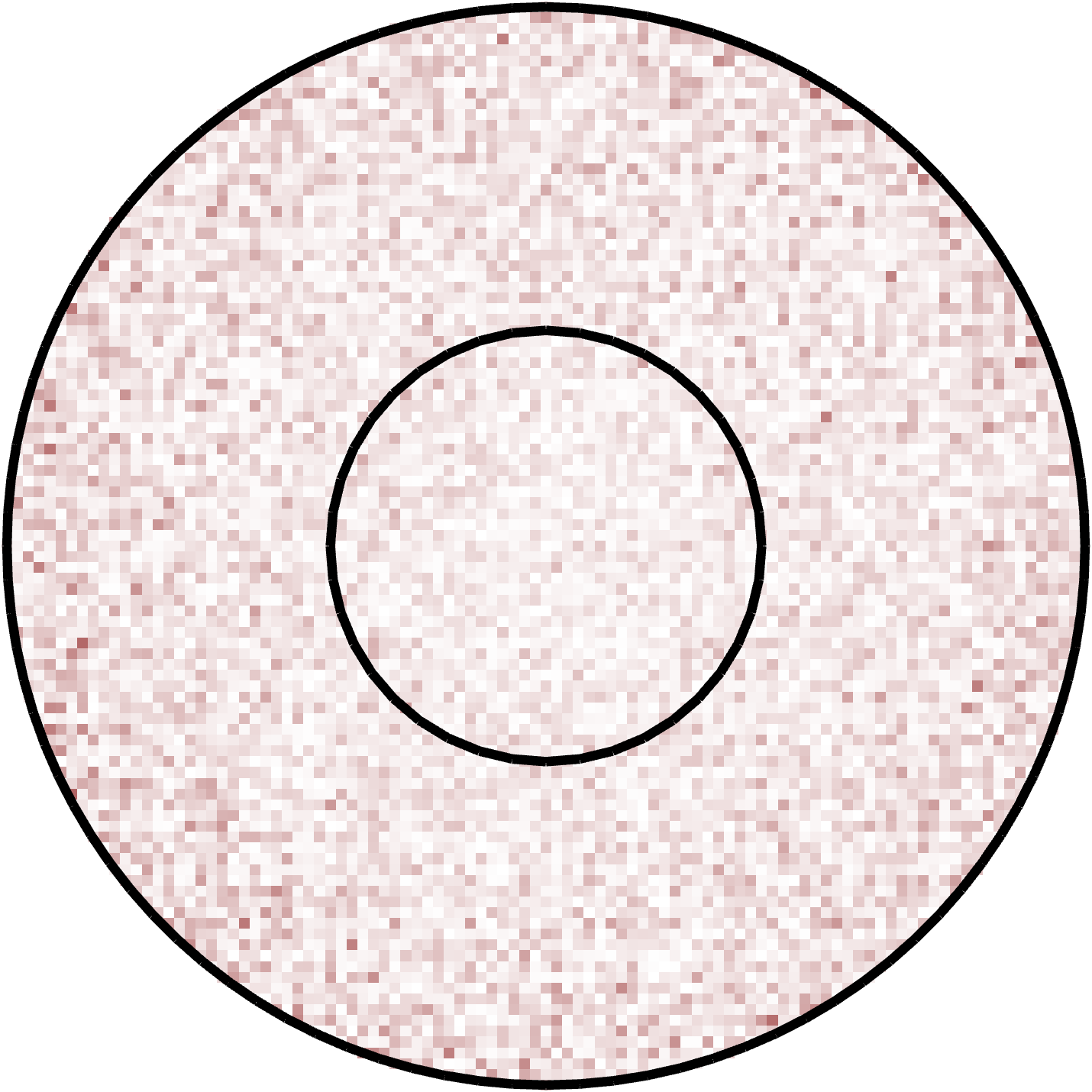} &
\includegraphics[width=0.2\columnwidth]{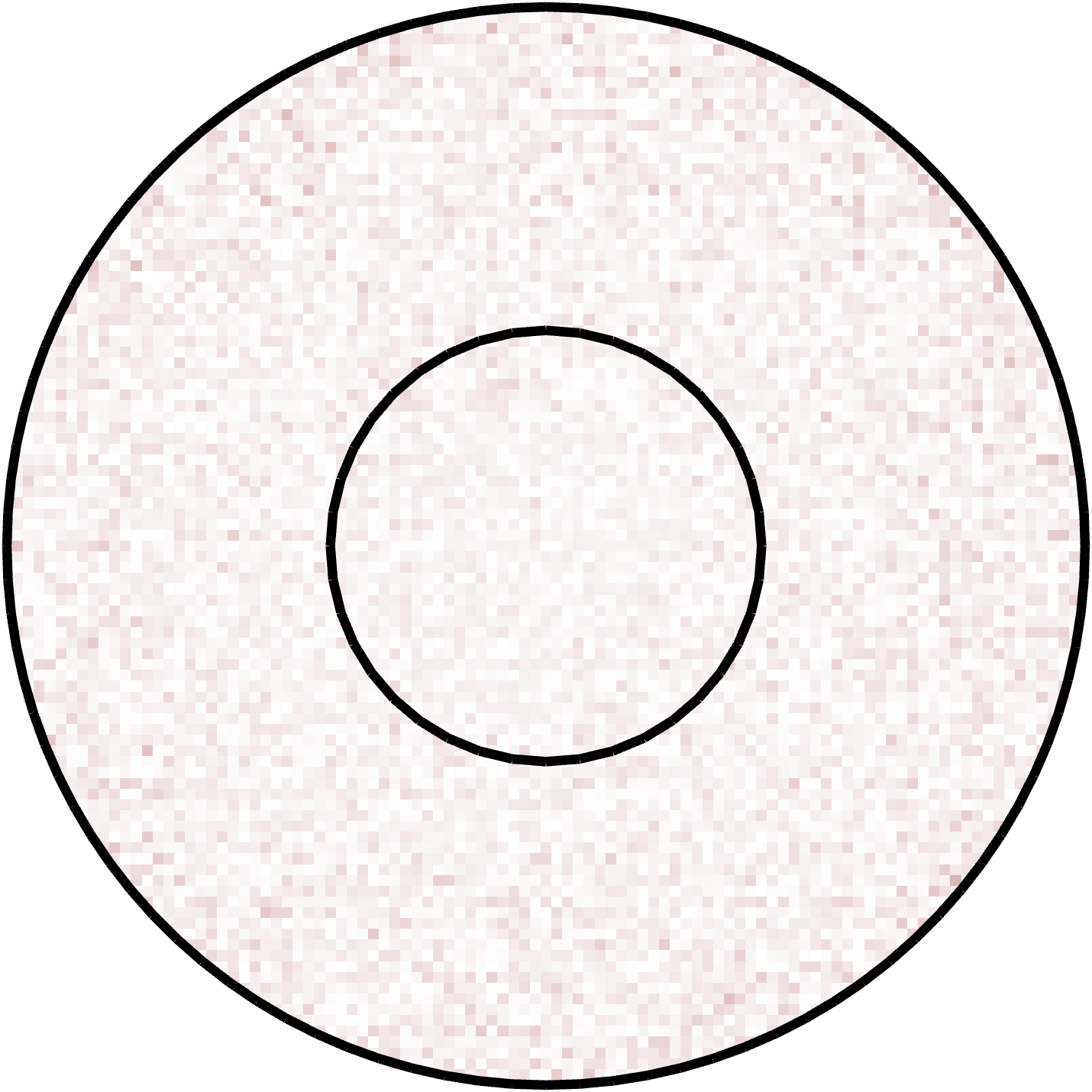} &
\includegraphics[width=0.2\columnwidth]{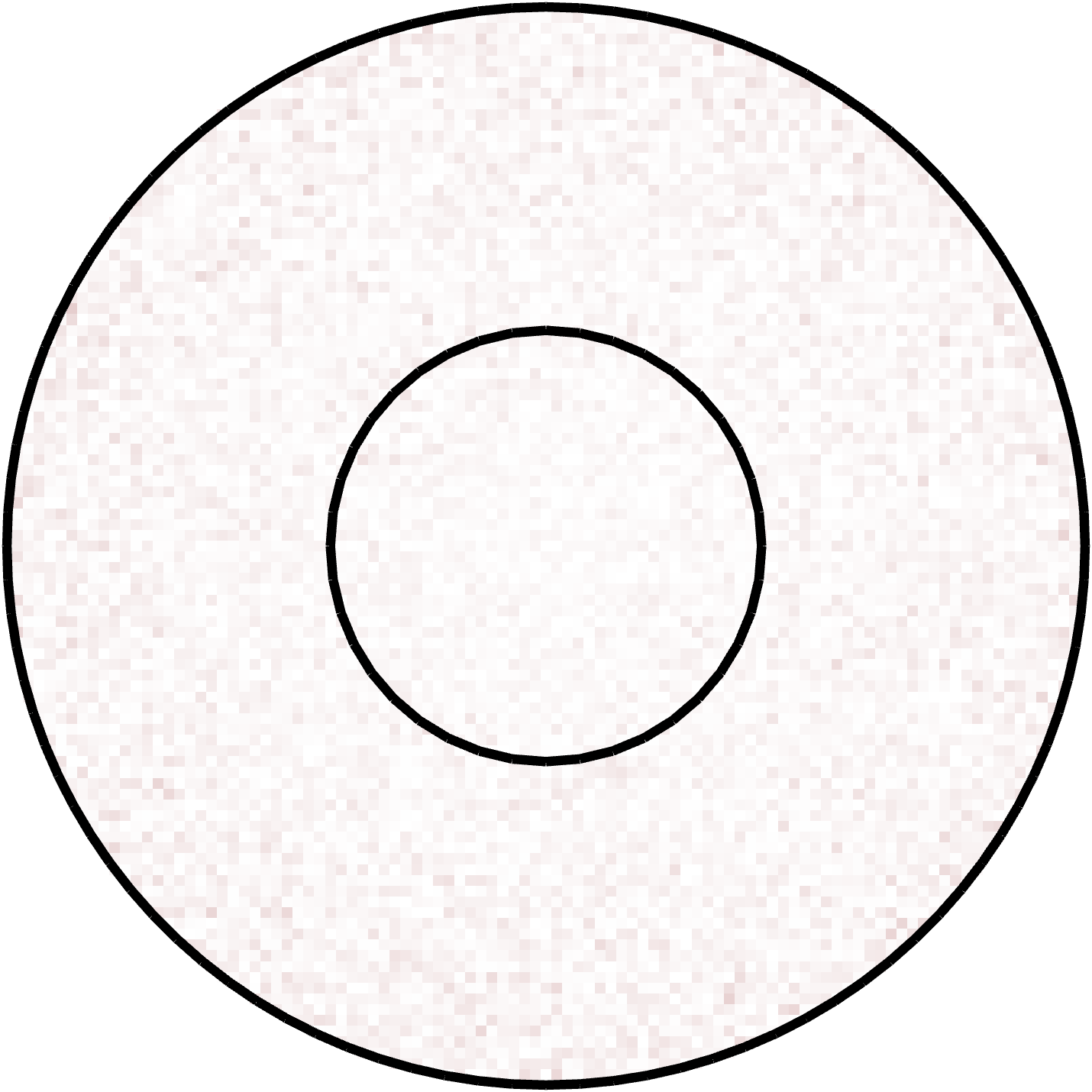} &
\includegraphics[width=0.03\columnwidth]{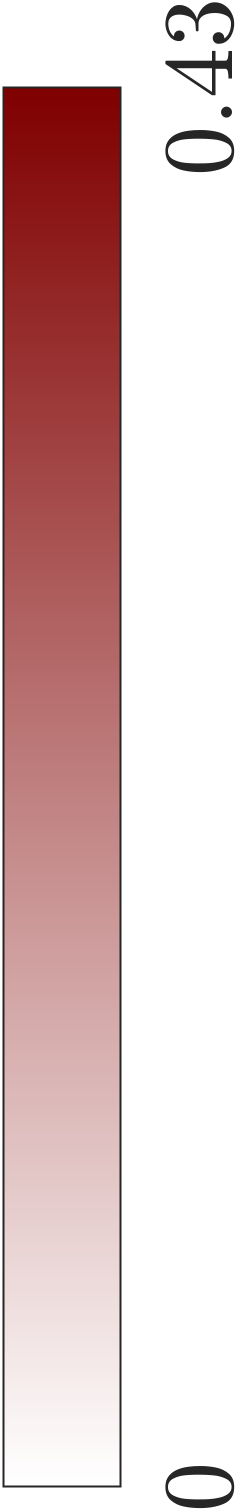}\\

% relative L2 error
{} & \textbf{$\bm{L_2}$ Error} & $0.060$ & $0.029$ & $0.019$ &\\
% relative L2 error
{} & \textbf{Relative $\bm{L_2}$ Error} & $2.55\%$ & $1.23\% $ & $0.80\%$ & \\
\end{tabular}
\caption{Solution to Eq.~\eqref{C0_gradient_2D_problem} with $\alpha = 0.4$, $\sigma = \frac{1}{3}$, $m = 3$, and $\lambda = 20$. In both halves of the figure, the first row shows the solution of WoI or variance-reduced WoI with different numbers of walkers. The second row plots the corresponding errors. Both estimators walk $M = 4$ steps.}\label{ex_woi_kink_gradient2D}
\end{figure}

We test our methods with $\alpha = 0.4$, $\sigma = \frac{1}{3}$, $m = 3$, and $\lambda = 20$ and present our results in Figure \ref{ex_woi_kink_gradient2D}.

\begin{figure}[htbp]
\centering
\begin{tabular}{@{}c@{}c c c c c}
% ============= WOI ==============
\multirow{3}{*}{
  \rotatebox[origin=c]{90}{%
    \parbox[c]{\dimexpr 0.5\columnwidth\relax}{\centering\vfill\textbf{On Boundary}\vfill}%
  }
} & 
\textbf{Neumann Data} & \textbf{Dirichlet Data} & \textbf{WoI} & \textbf{Var-reduced WoI} & {} \\

% First Row: estimation
{} & \includegraphics[width=0.2\columnwidth]{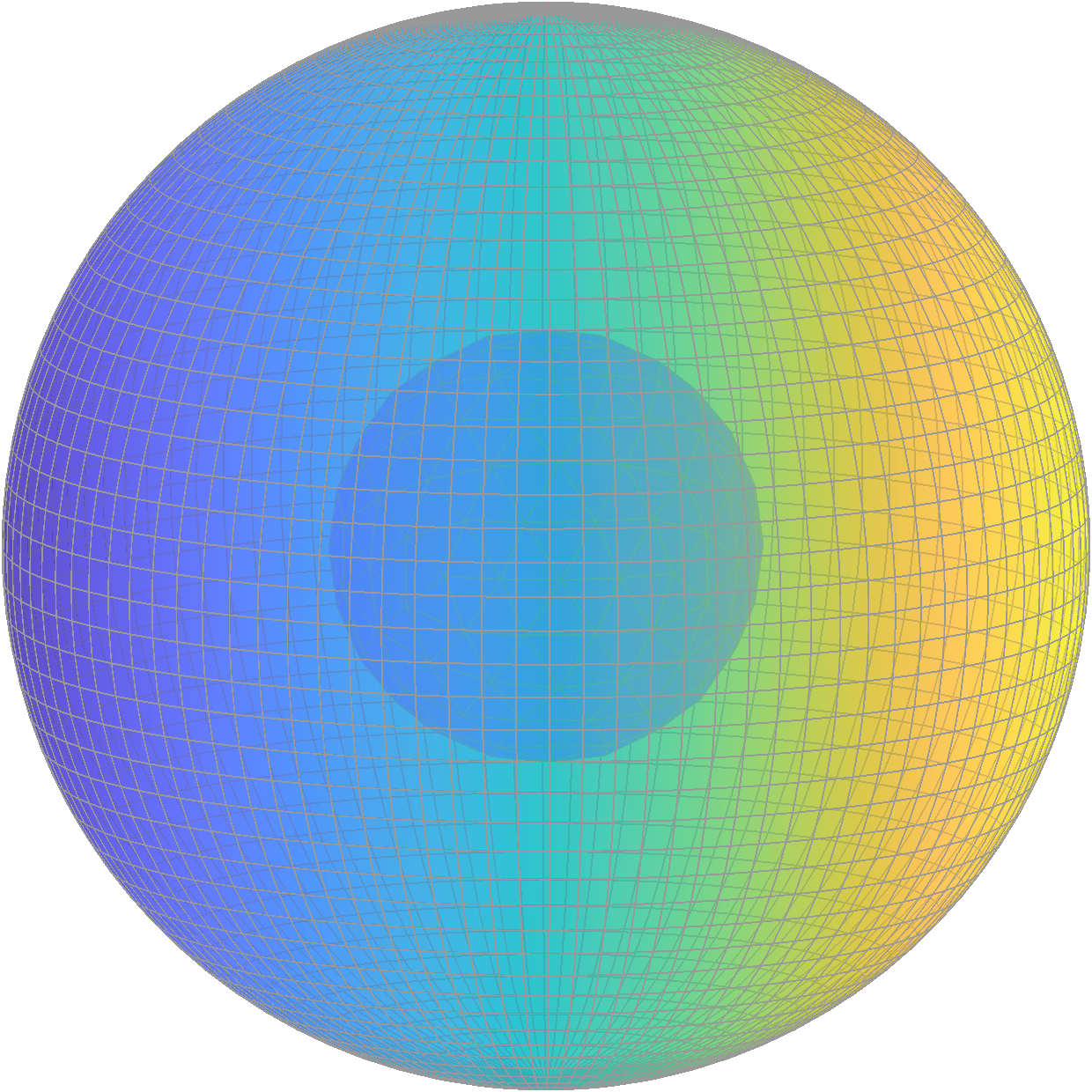} &
\includegraphics[width=0.2\columnwidth]{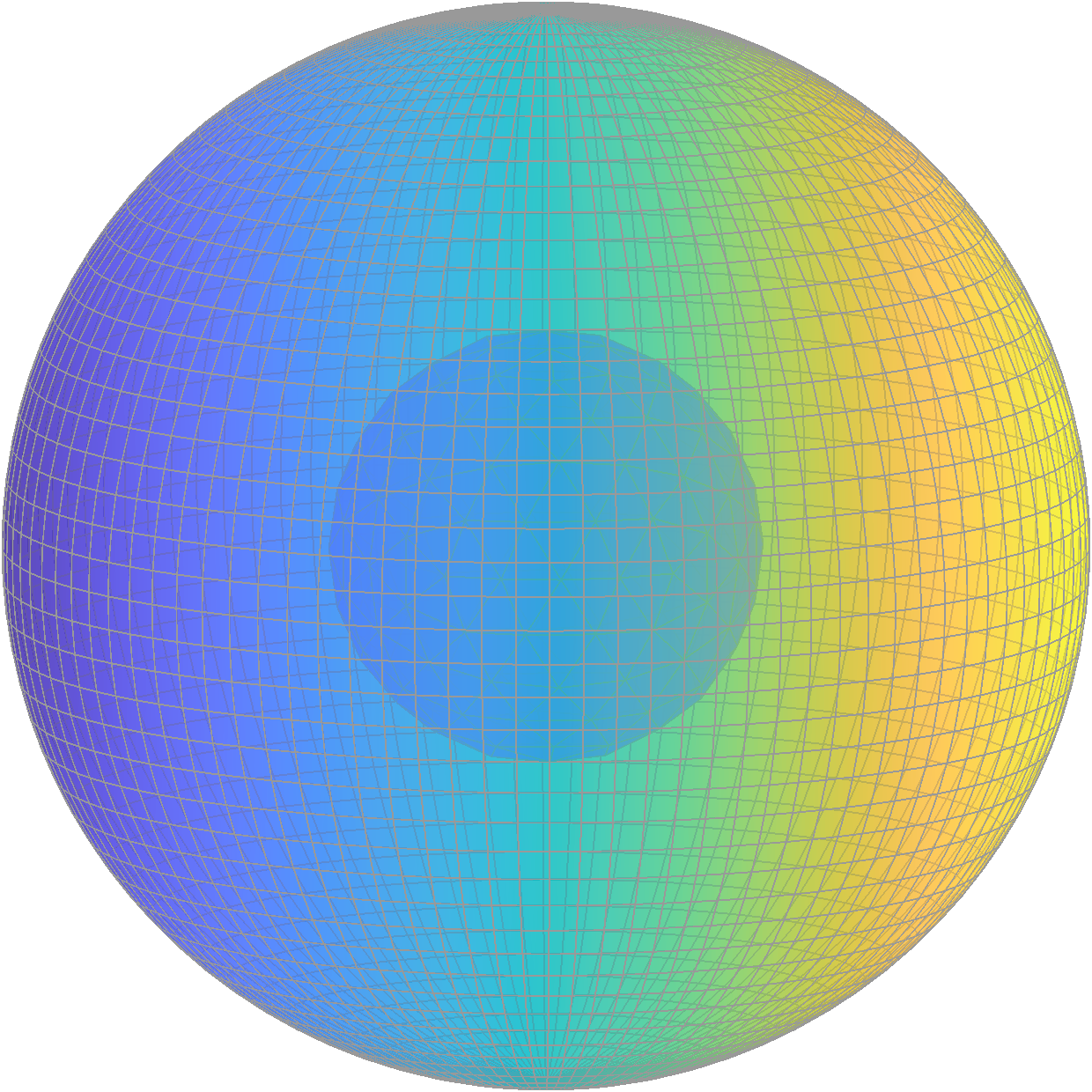}  &
\includegraphics[width=0.2\columnwidth]{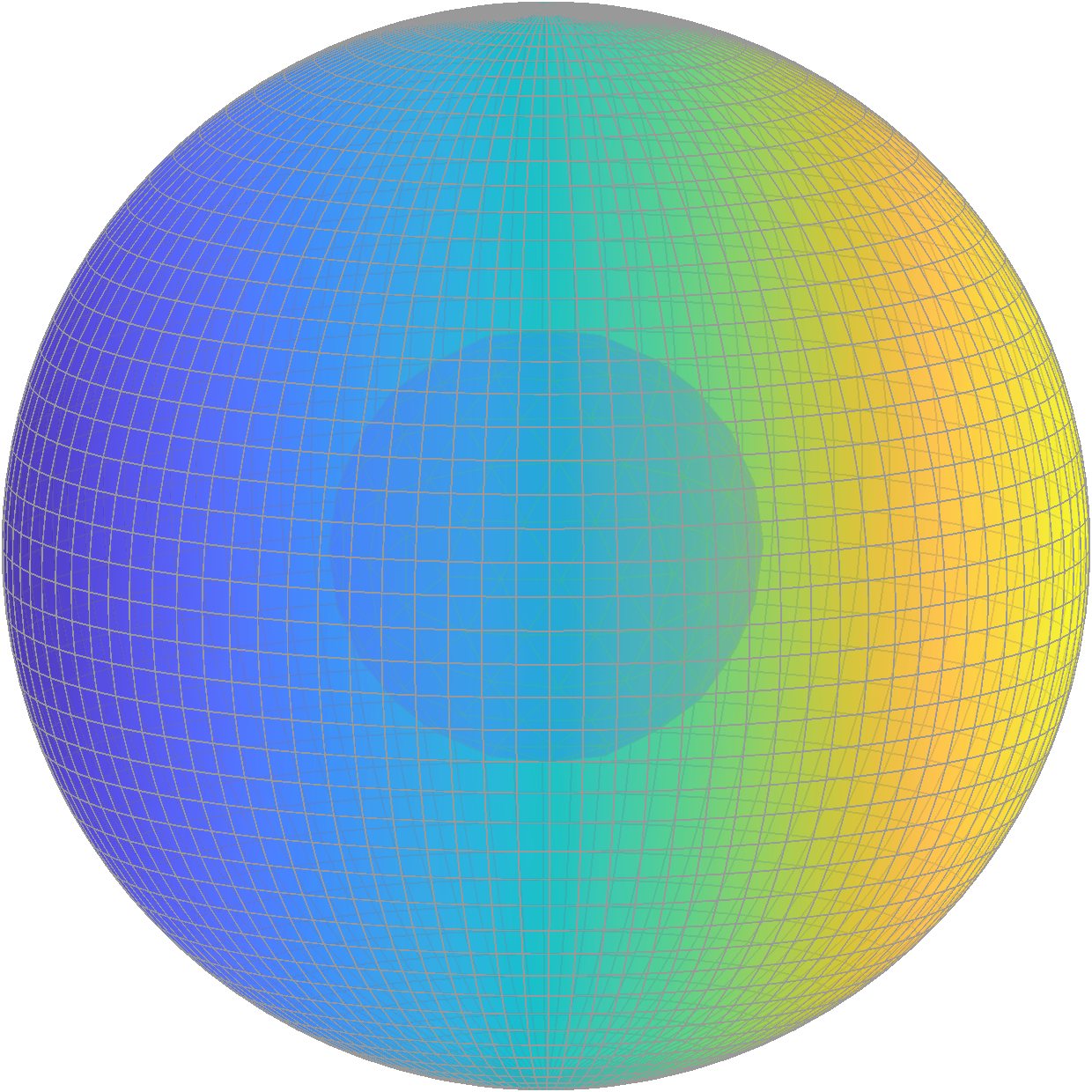} &
\includegraphics[width=0.2\columnwidth]{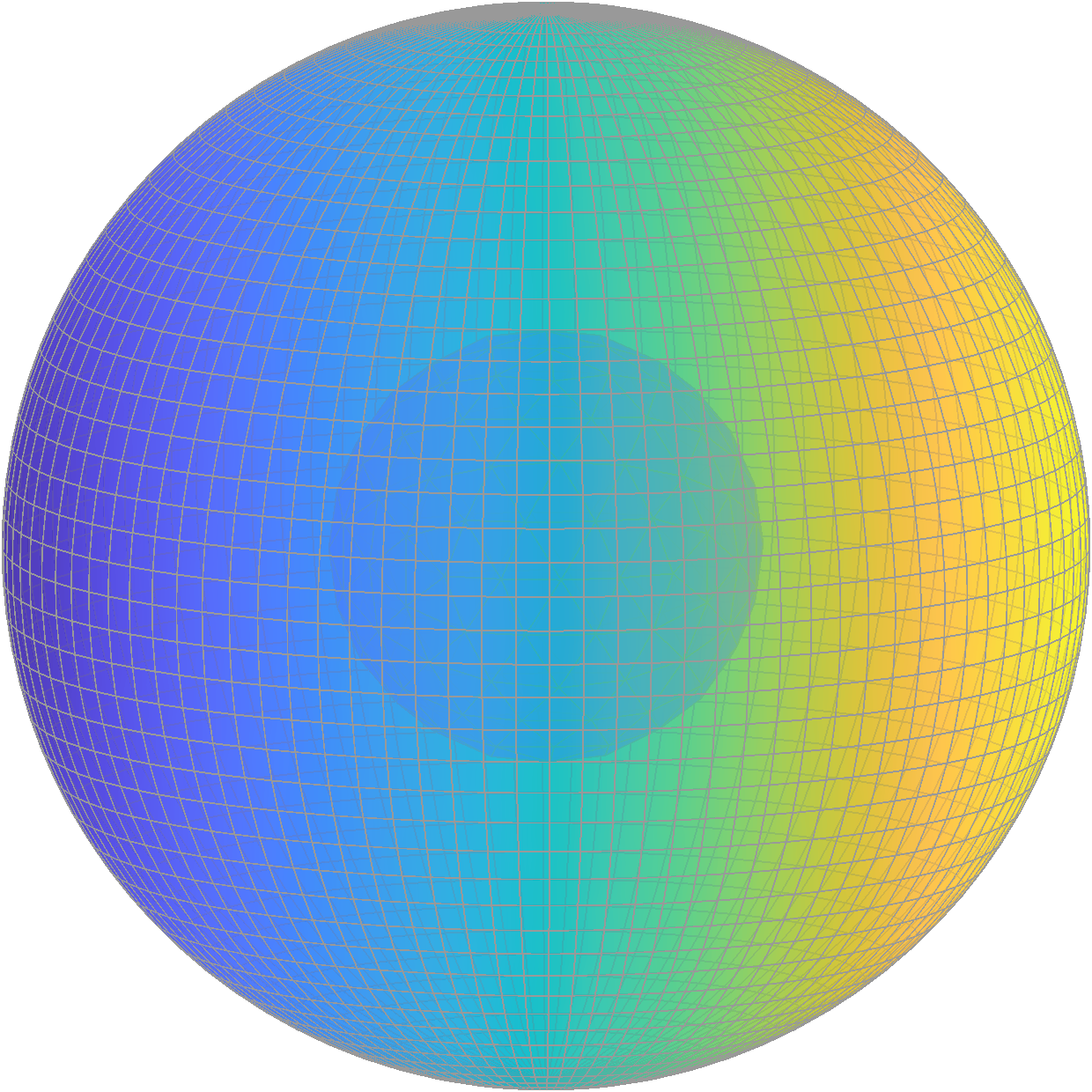} &
\includegraphics[width=0.03\columnwidth]{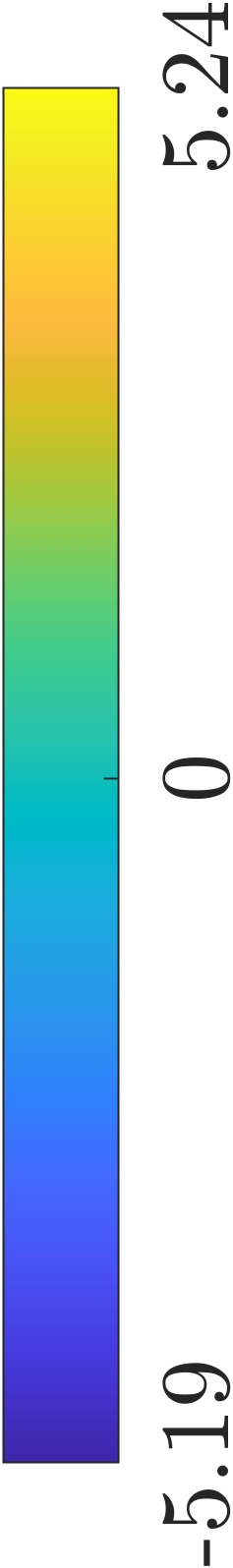} \\

% Second Row: pointwise estimation
{} & {} & {} &
\includegraphics[width=0.2\columnwidth]{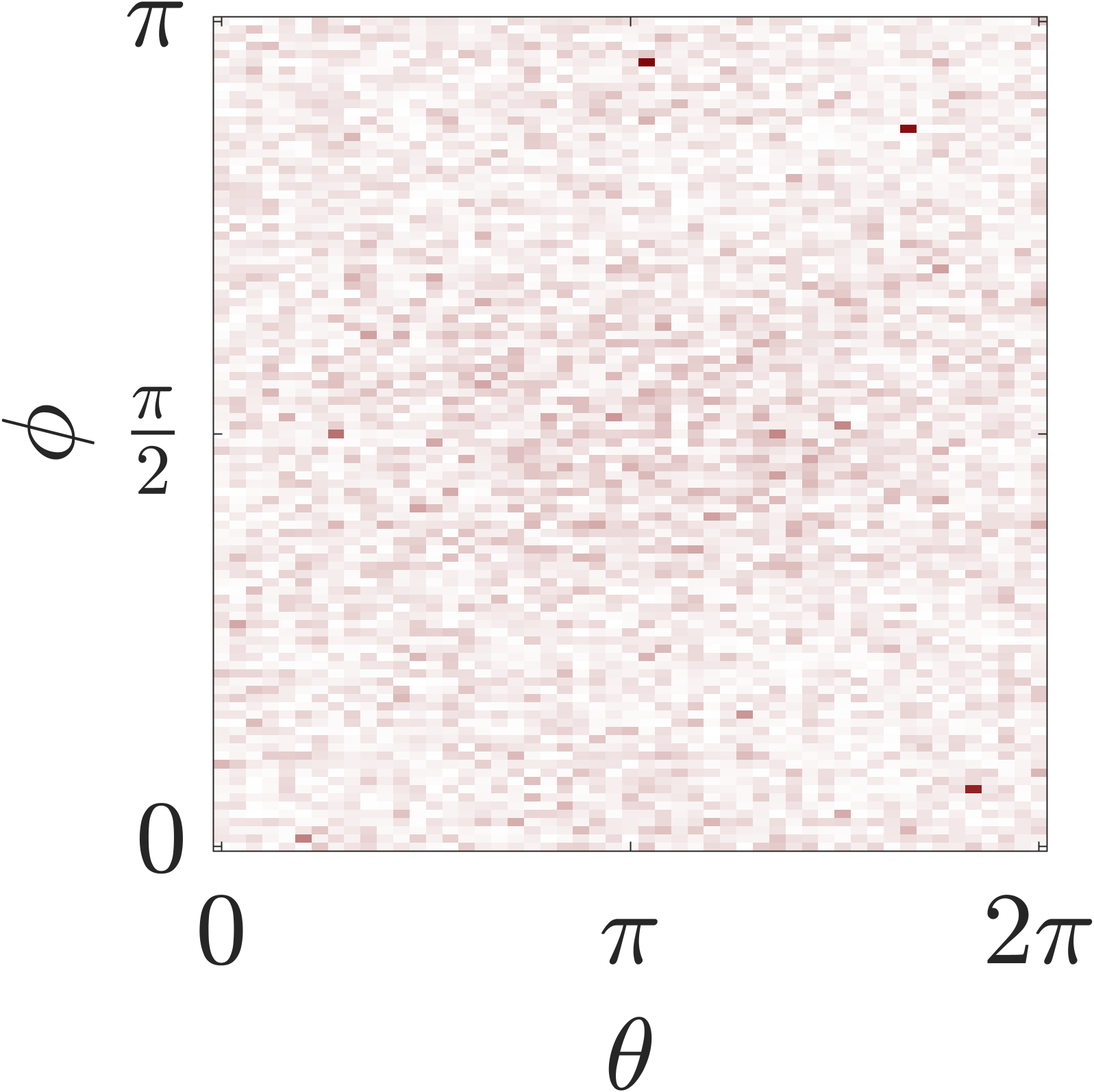} &
\includegraphics[width=0.2\columnwidth]{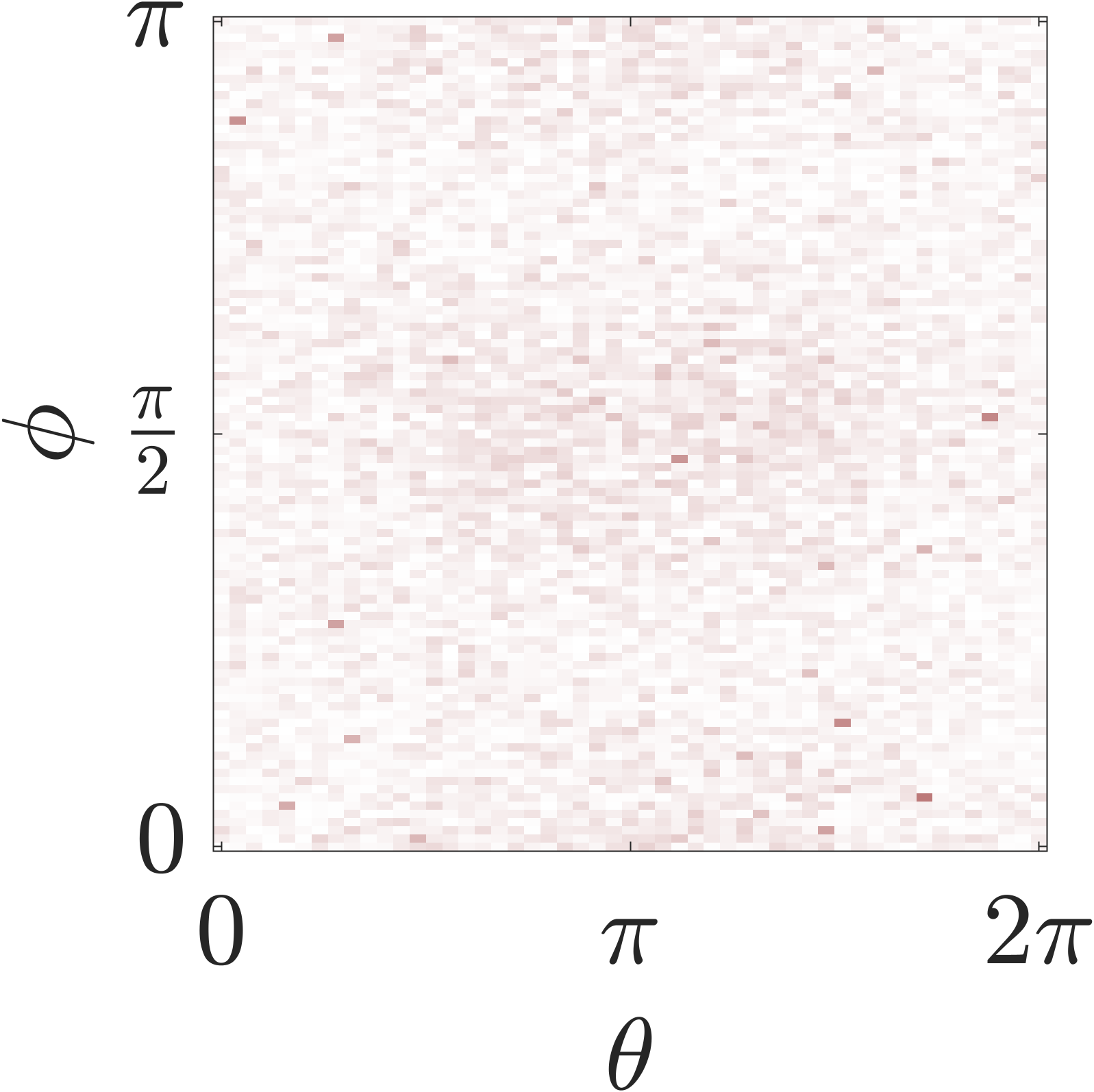} &
\includegraphics[width=0.034\columnwidth]{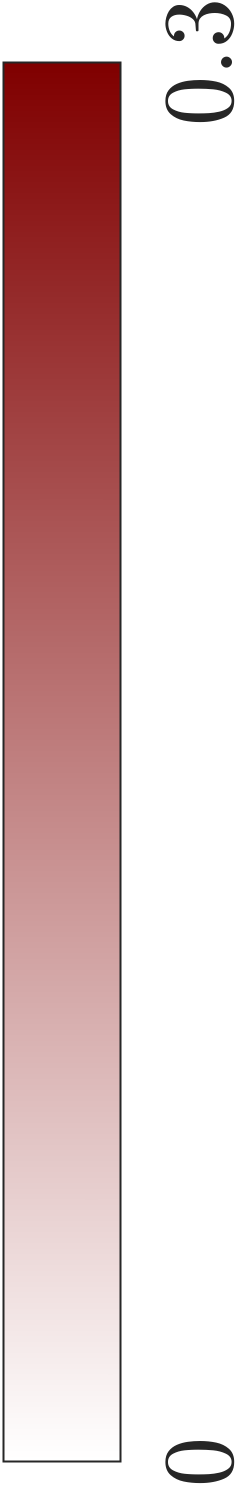}\\
{} & {} & \textbf{$\bm{L_2}$ Error} & $0.033$ & $0.021$ & {} \\
% relative L2 error
{} & {} & \textbf{Relative $\bm{L_2}$ Error} & $1.28\% $ & $0.83\% $ & {}\\
%===============================================
\hline \hline \\
%===============================================
% =========== variance-reduced WoI =============
\multirow{3}{*}{
  \rotatebox[origin=c]{90}{%
    \parbox[c]{\dimexpr 0.5\columnwidth\relax}{\centering\vfill\textbf{On Interface}\vfill}%
  }
} & \textbf{Neumann Data} & \textbf{Dirichlet Data} & \textbf{WoI} & \textbf{Var-reduced WoI} & {}  \\
% First Row: estimation
{} & \includegraphics[width=0.2\columnwidth]{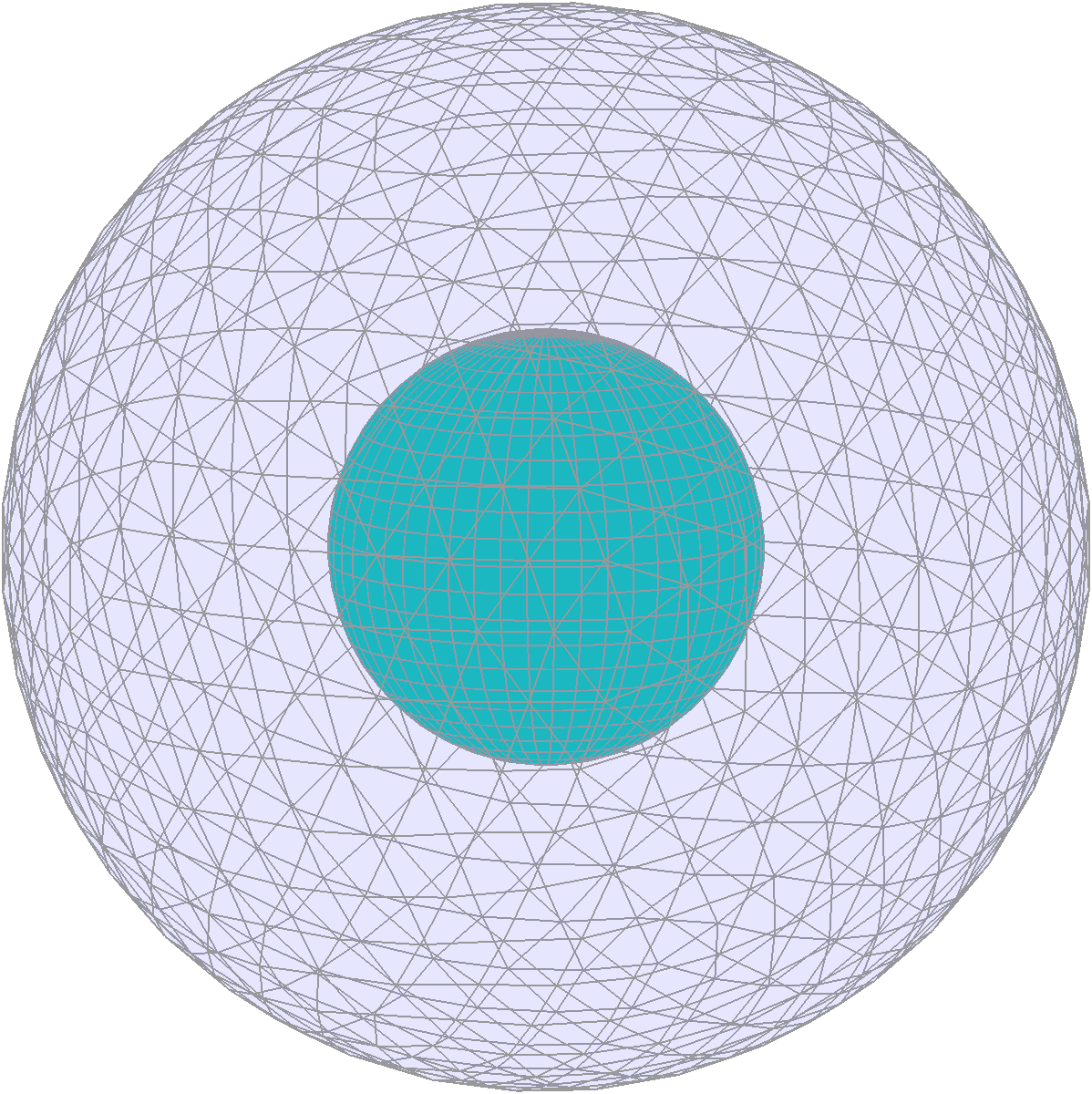} &
\includegraphics[width=0.2\columnwidth]{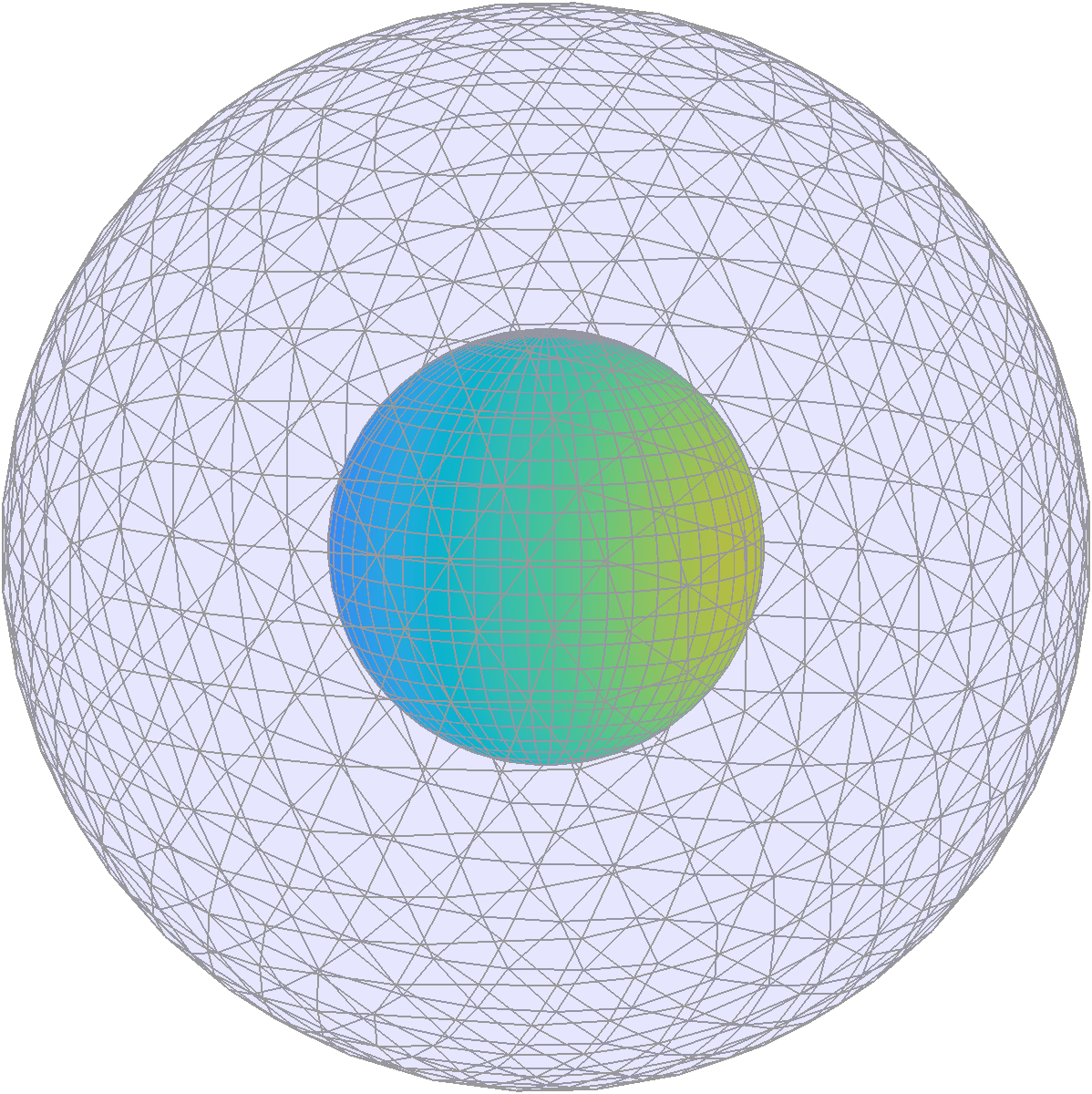}  &
\includegraphics[width=0.2\columnwidth]{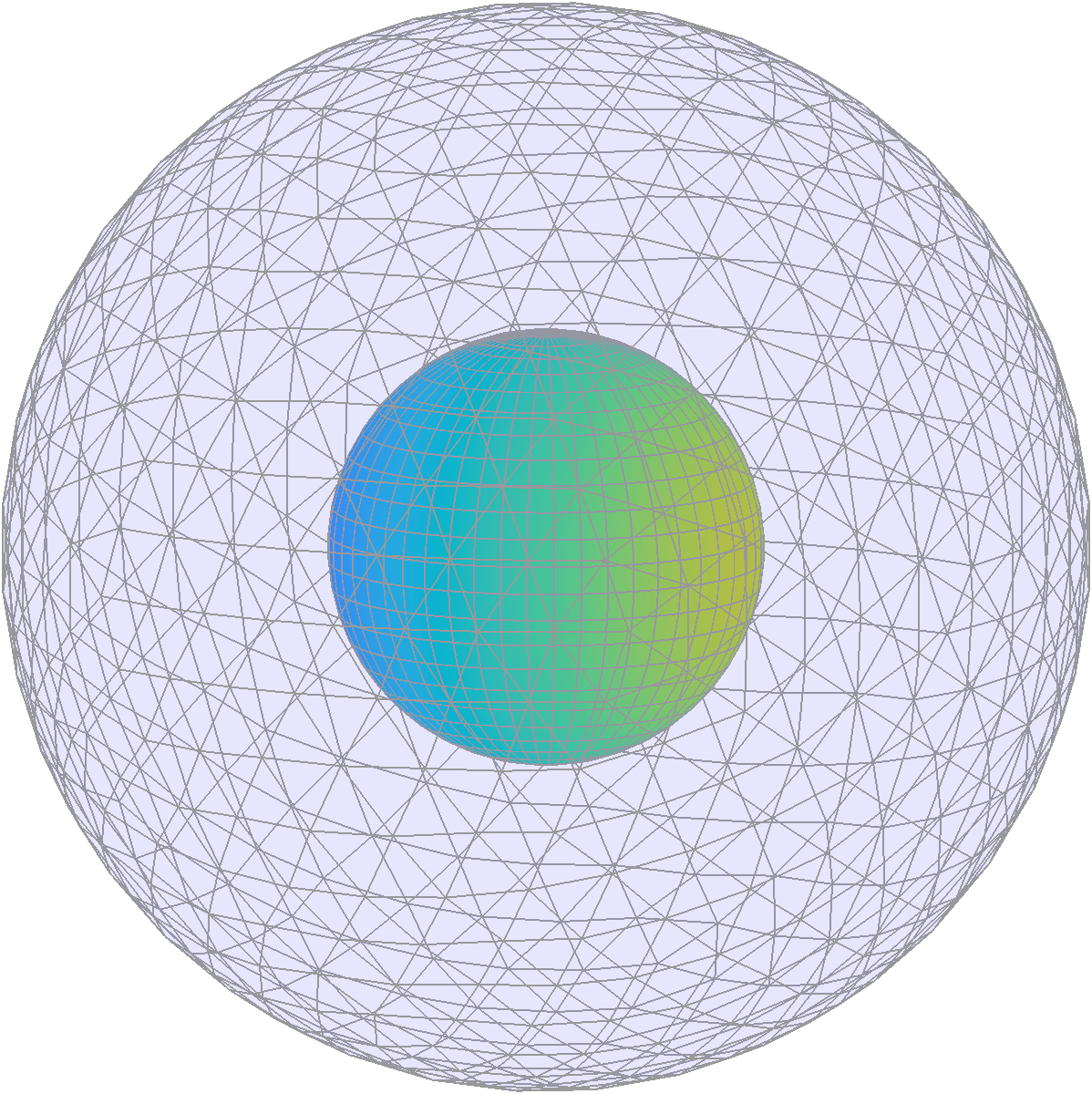} &
\includegraphics[width=0.2\columnwidth]{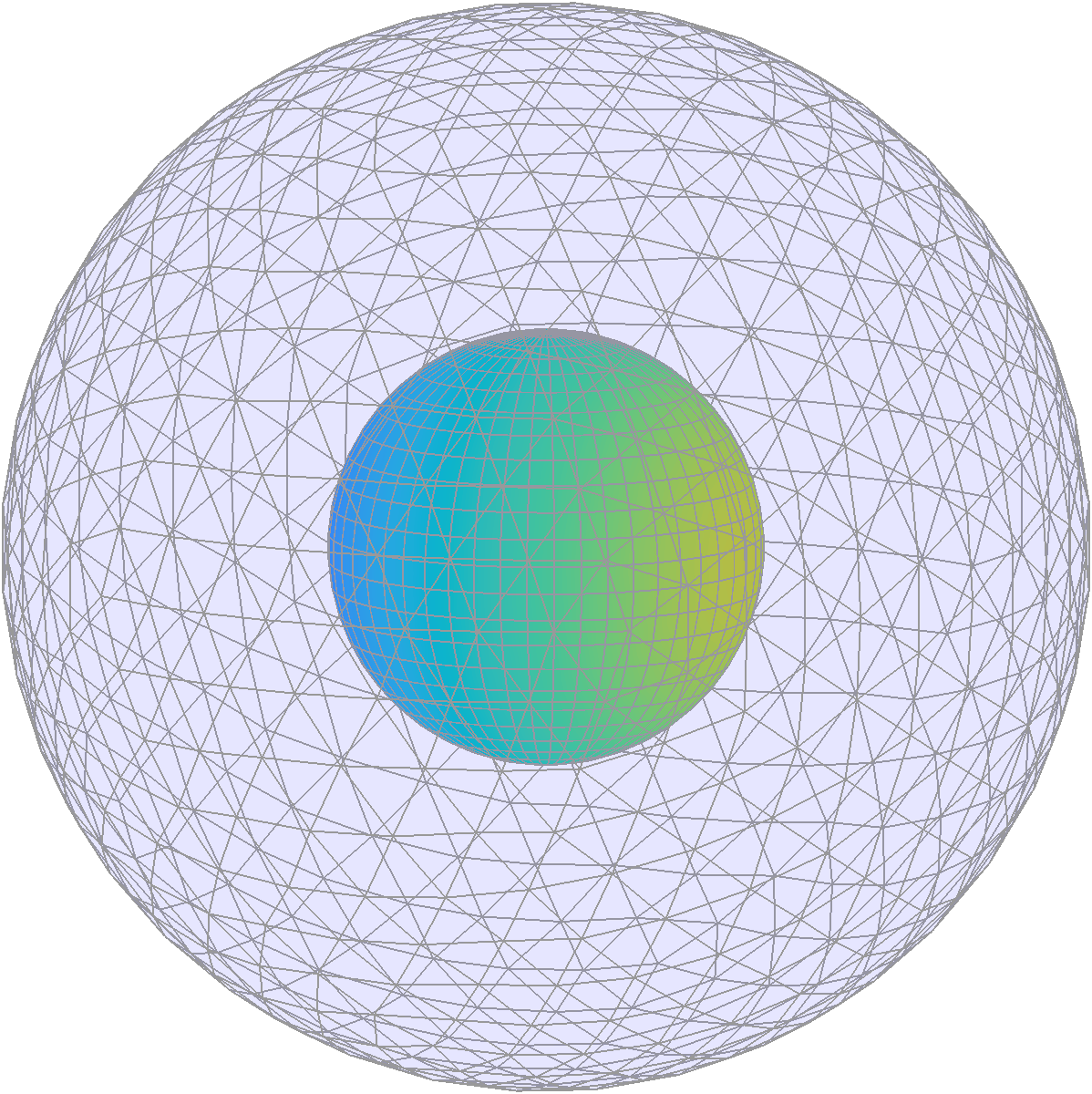} &
\includegraphics[width=0.03\columnwidth]{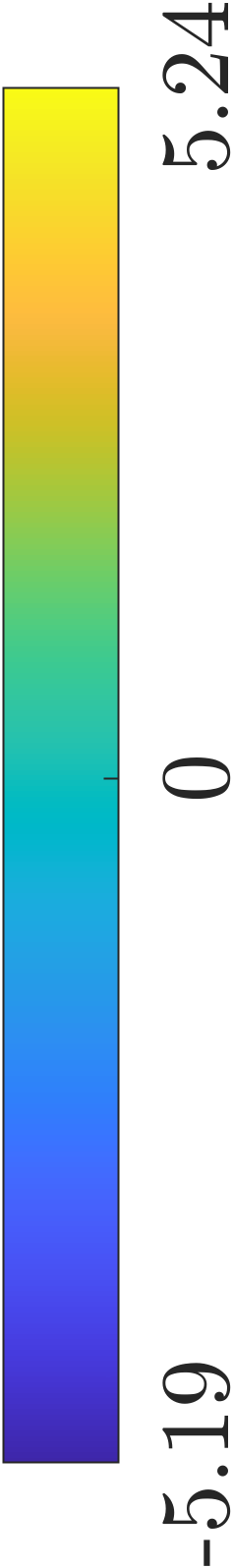}\\
% Second Row: pointwise estimation
{} & {} & {} &
\includegraphics[width=0.2\columnwidth]{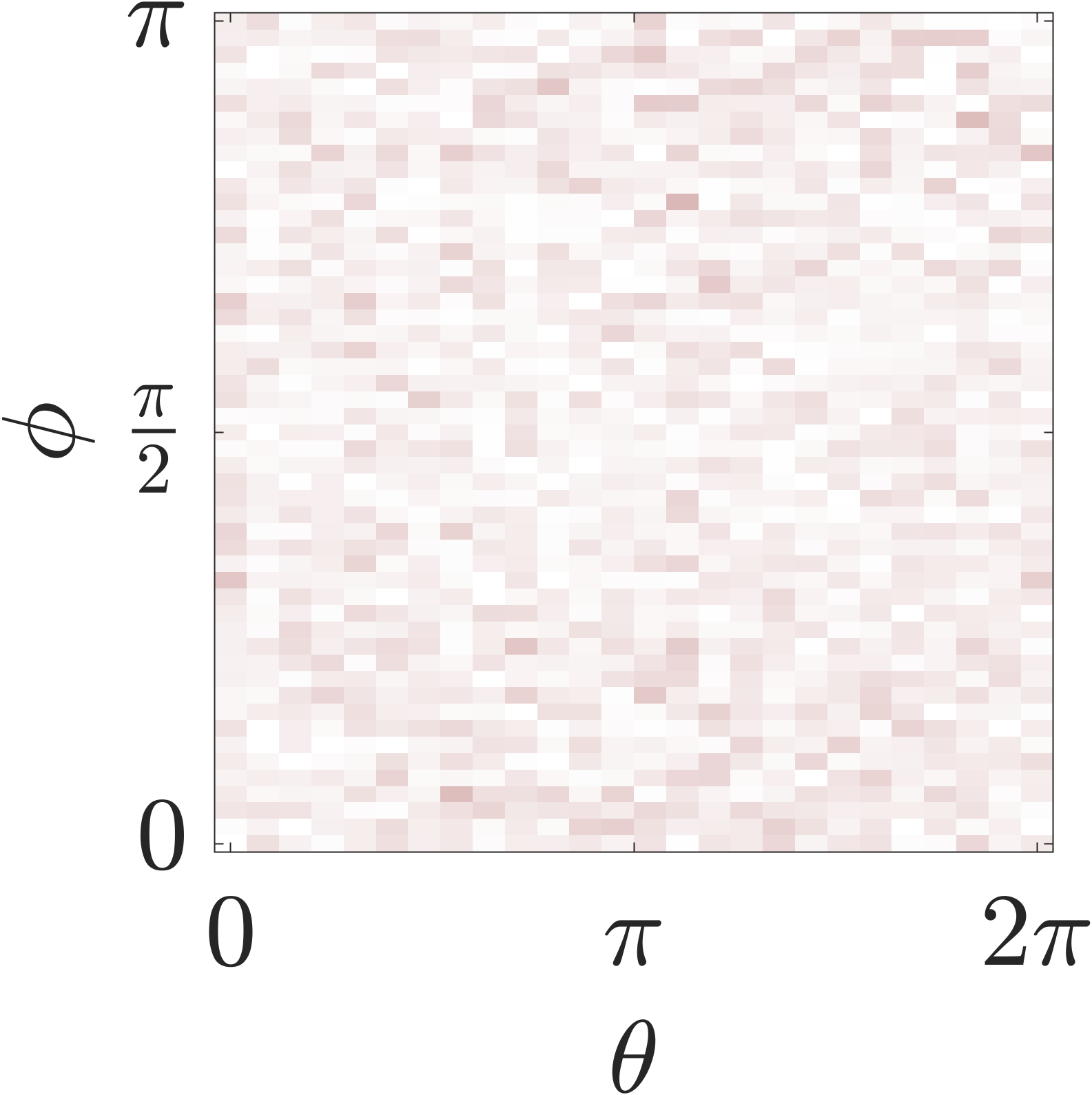} &
\includegraphics[width=0.2\columnwidth]{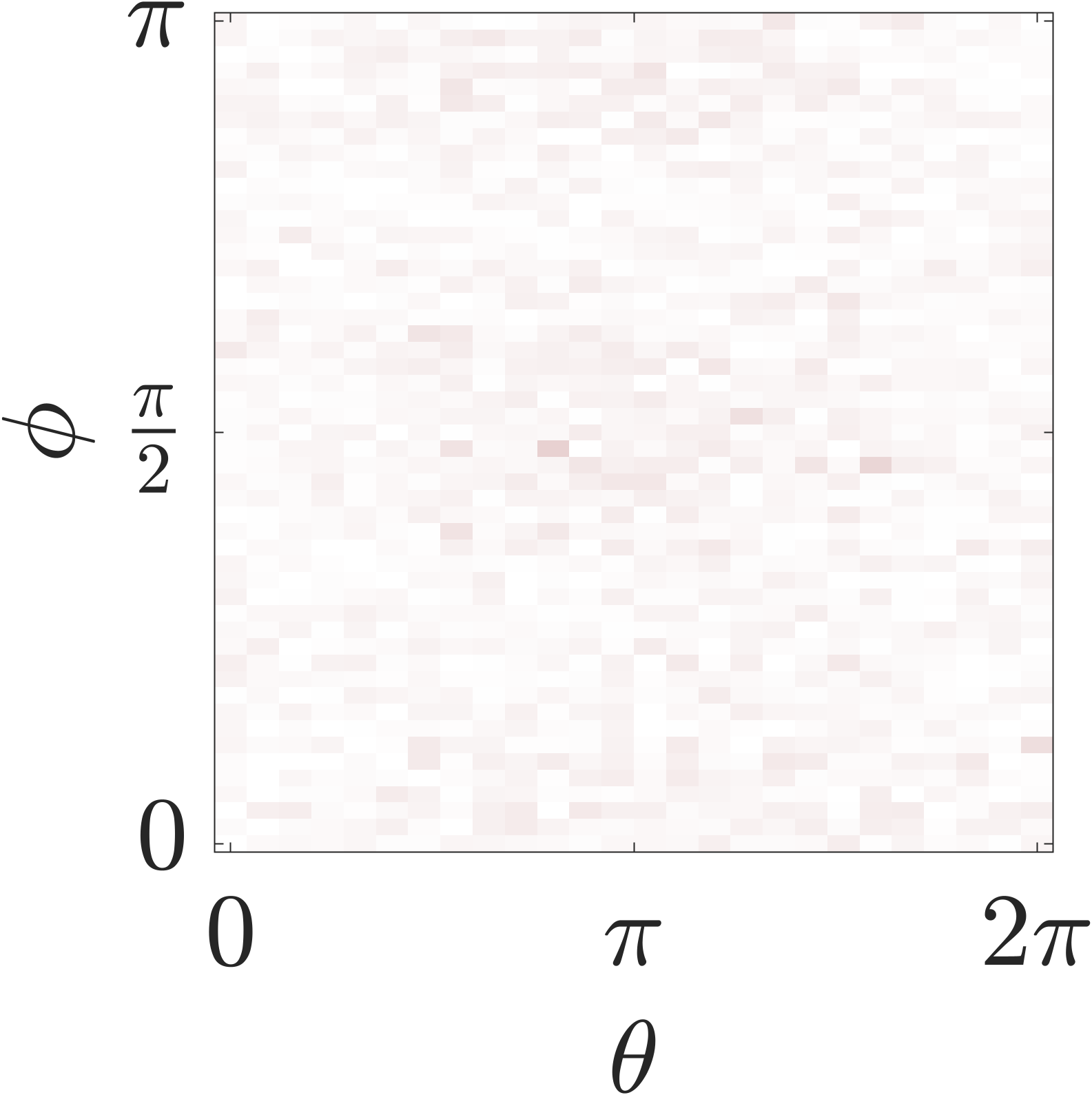} &
\includegraphics[width=0.034\columnwidth]{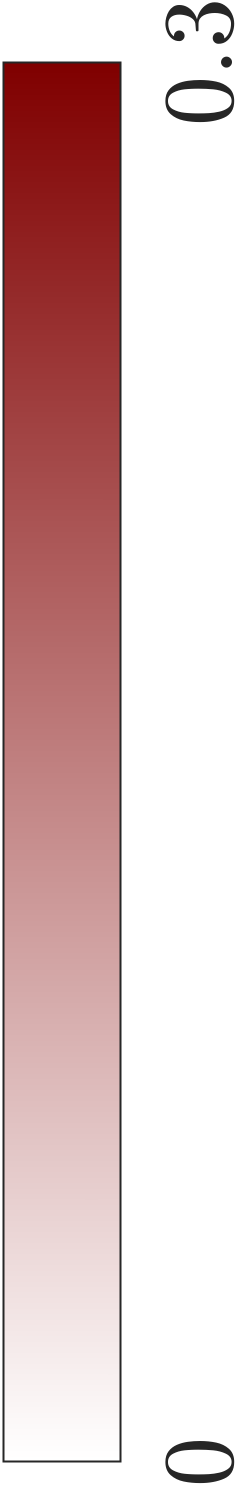}\\
{} & {} & \textbf{$\bm{L_2}$ Error} & $0.024$ & $0.011$ & {} \\
% relative L2 error
{} & {} & \textbf{Relative $\bm{L_2}$ Error} & $2.00\% $ & $0.93\% $ & {}\\
\end{tabular}
\caption{Solution to Eq.~\eqref{C0_gradient_3D_problem} with $\alpha = 0.4$, $\sigma = \frac{1}{2}$, and $\lambda = 5$. In both halves of the figure, the first row shows ground truth of the Neumann-to-Dirichlet (NtD) map, the NtD map obtained from WoI with $\mathcal{W} = 10^7$ walkers, and the NtD map obtained from variance-reduced WoI with $\mathcal{W} = 2 \times 10^7$ walkers. The second row plots the corresponding errors in the $\theta-\phi$ parametric space. Both estimators walk $M = 6$ steps. Since the spherical harmonic function $Y_1^1$ has property
$\sin(\phi) \cos(\theta) = \sqrt{\frac{4\pi}{3}}Y_1^1$,
the Neumann data on the boundary is proportional to the Dirichlet data.}\label{ex_woi_kink_gradient3D}
\end{figure}

Andrew Zheng~\cite{Zheng_thesis} provides yet another example in a domain defined by two concentric spheres. We let $\Omega_1$ be the outer sphere of radius 1, and $\Omega_2$ be the inner sphere of radius $\alpha$. The problem of interest requires $u(\rho, \theta, \phi)$ to satisfy
\begin{subequations}\label{C0_gradient_3D_problem}
\begin{empheq}[left=\empheqlbrace]{align}
 \Delta u &= 0 \quad \text{on $r \leq 1$, $r \neq \alpha$} \\
 % ========================
\partial_{\bm{r}}u &= \lambda \sin(\phi) \cos(\theta) \quad \text{on $r = 1$} \\
 % ========================
 [u] &= 0 \quad \text{on $r = \alpha$}\\
 % ========================
 \sigma \partial_{\bm{r}^{-}} u &= \partial_{\bm{r}^{+}} u \quad \text{on $r = \alpha$},
\end{empheq}
\end{subequations}
in which $\theta \in [0, 2\pi)$, $\phi \in [0, \pi)$, $\sigma = \frac{\sigma_2}{\sigma_1}$ is the ratio of $\sigma(\bm{x})$ in balls of radius $\alpha$ and 1. By separation of variables again, the analytical solution is
\begin{equation}\label{C0_gradient_3D_analytical_soln}
    u(r, \theta, \phi) = 
    \begin{cases}
        \lambda A_1^1 r Y_1^1 (\theta, \phi), \quad & r \leq \alpha \\
        \lambda (B_1^1 r + C_1^1 r^{-2}) 
        Y_1^1 (\theta, \phi), \quad & \alpha < r \leq 1
    \end{cases},
\end{equation}
in which $Y_1^1 (\theta, \phi)$ is the spherical harmonics function with $(m, \ell) = (1, 1)$, and $A_1^1$, $B_1^1$, $C_1^1$ satisfy the system
\begin{equation}
    \begin{bmatrix}
        0 & 1 & -2 \\
        \alpha & -\alpha & -\alpha^{-2} \\
        \sigma & -1 & 2\alpha^{-3}
    \end{bmatrix}
    \begin{bmatrix}
        A_1^1 \\ B_1^1 \\ C_1^1
    \end{bmatrix}
    =
    \begin{bmatrix}
        \sqrt{\frac{4\pi}{3}} \\ 0 \\ 0
    \end{bmatrix}.
    \notag
\end{equation}
In Figure \ref{ex_woi_kink_gradient3D}, we present the Neumann-to-Dirichlet (NtD) map defined by Eq.~\eqref{C0_gradient_3D_problem} by plotting the solution values on the boundary and the interface.
Unlike the family of WoS estimators, which typically leave a shell of inaccurate results near the boundary to stop the Brownian motion, and unlike the boundary integral methods that suffer from a near singular kernel when evaluating close-to-interface points, our results demonstrate that our WoI estimators maintain reasonable accuracy on the boundary and across the interfaces.

\begin{figure}[htbp]
\centering
\begin{tabular}{c c}
    \textbf{WoI Convergence} & \textbf{Var-reduced WoI Convergence} \\
    % u(x)
    \includegraphics[width=0.4\columnwidth, valign=c]{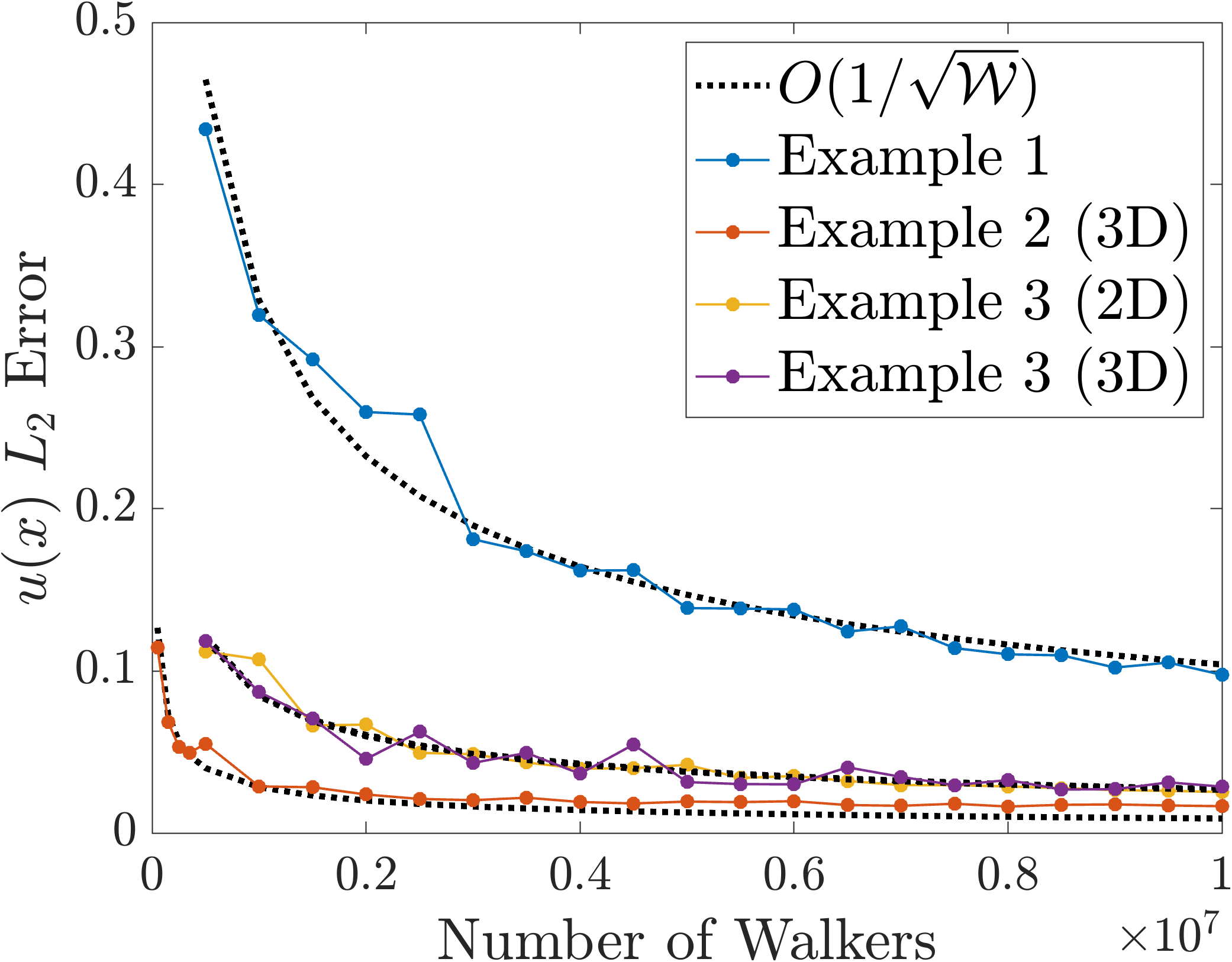} &
    \includegraphics[width=0.4\columnwidth, valign=c]{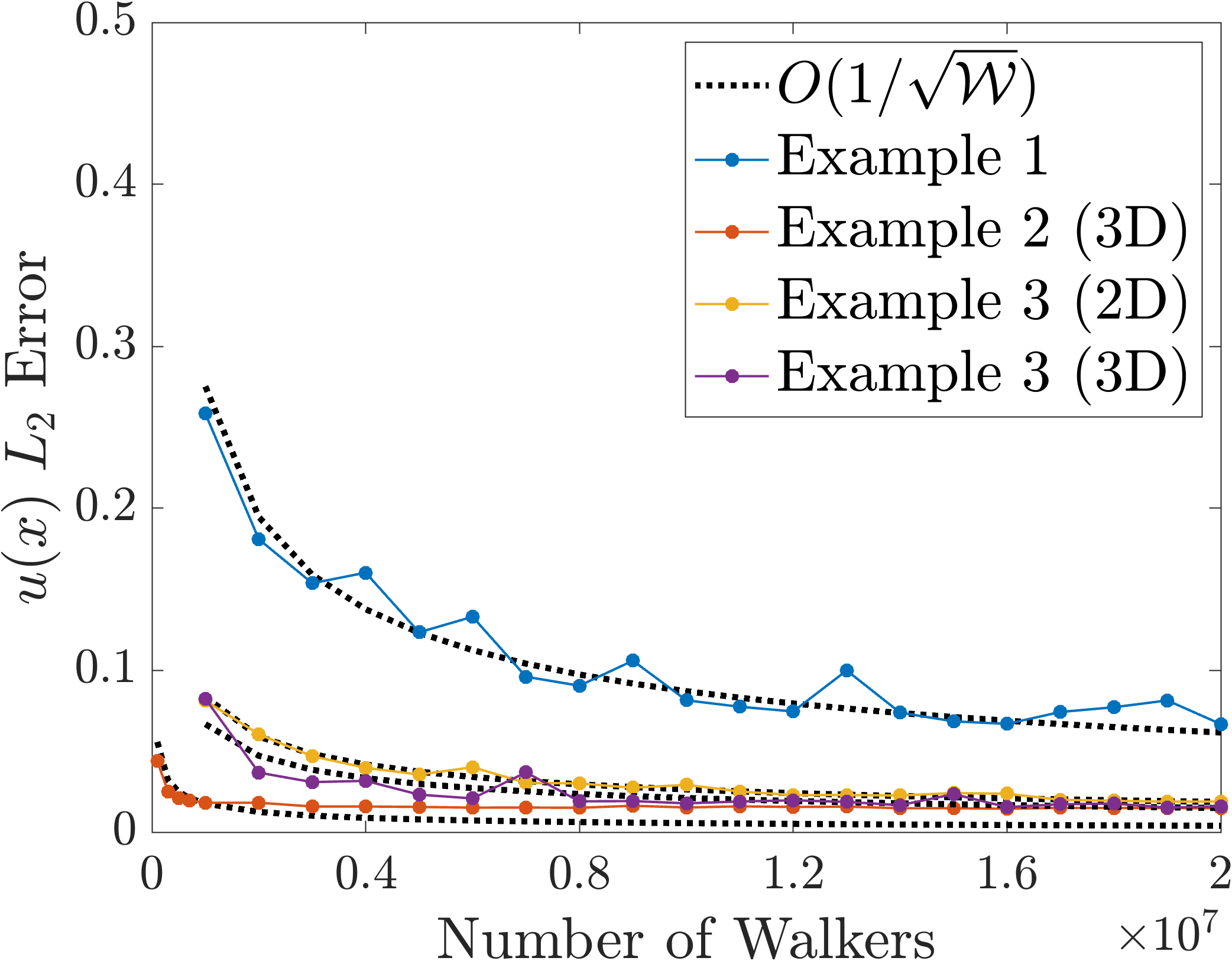} 
    \vspace{0.5cm} \\
    % gradient
    \includegraphics[width=0.4\columnwidth, valign=c]{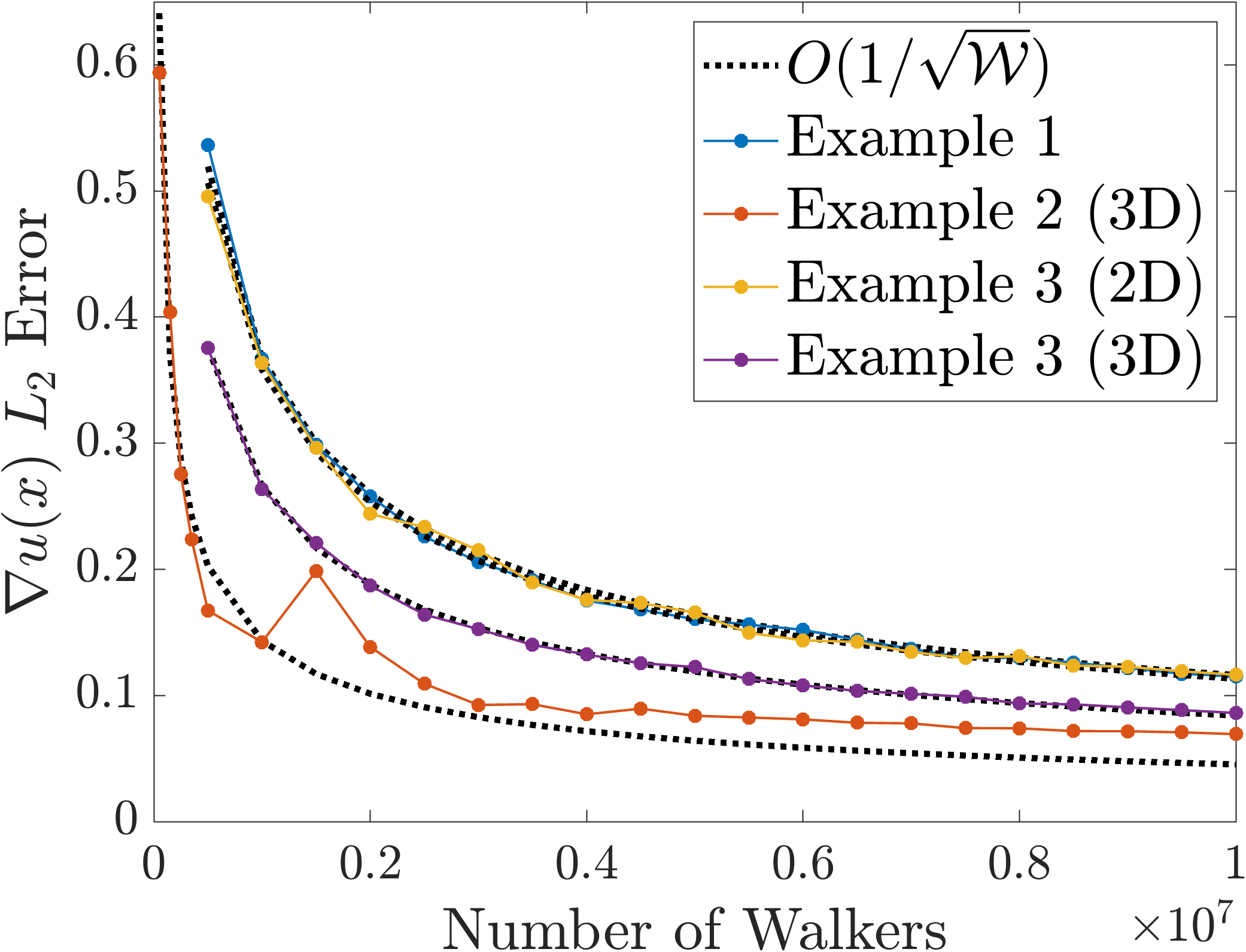} &
    \includegraphics[width=0.4\columnwidth, valign=c]{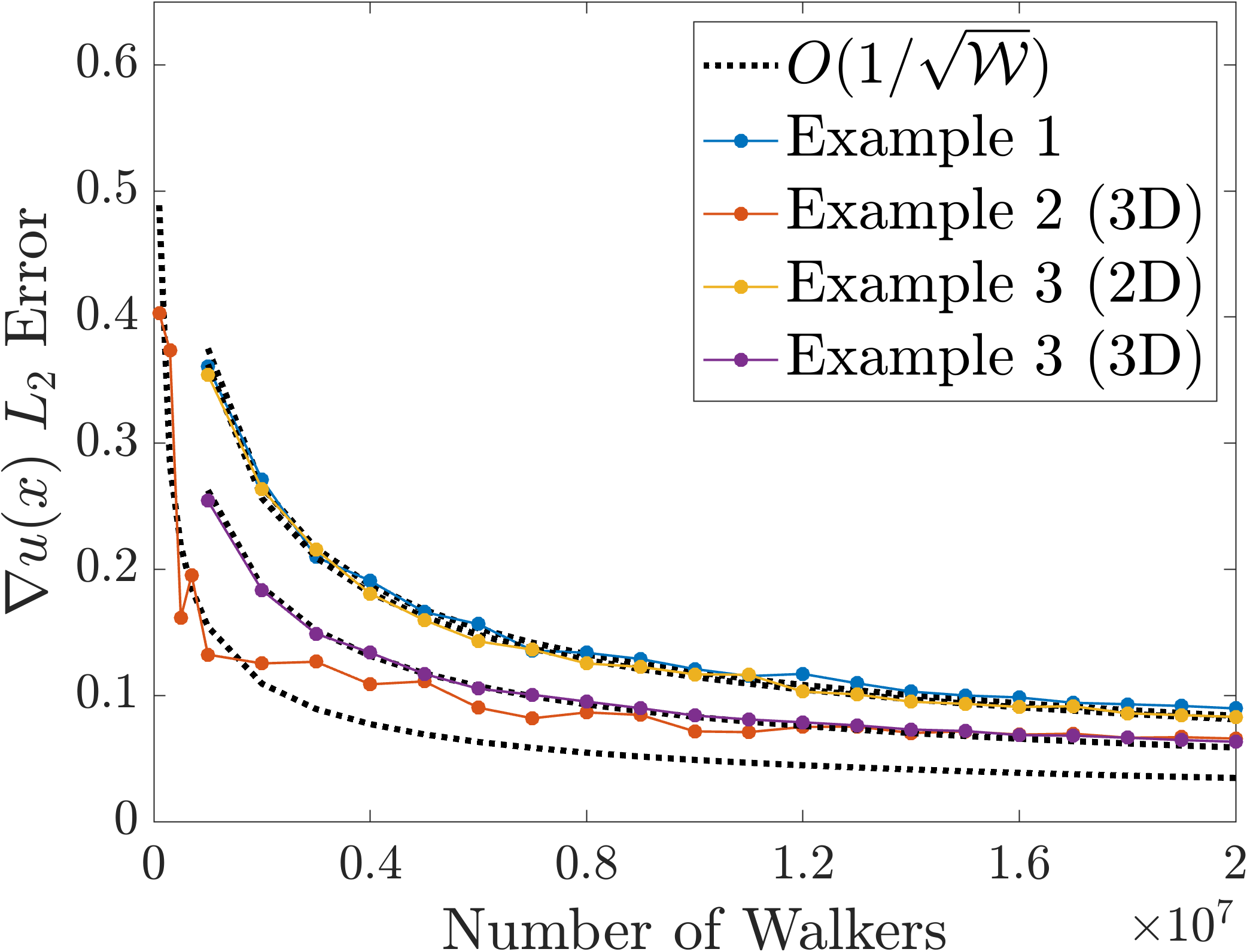}
\end{tabular}
\caption{$\mathcal{O}(1/\sqrt{\mathcal{W}})$ convergence rate, where $\mathcal{W}$ is the number of walkers, of both WoI and variance-reduced WoI estimators and their corresponding gradient estimators. The unexpected tail behavior in the convergence plots for Example 2 (3D) arises because the solution and its gradient are nearly zero over much of the domain, which causes Monte Carlo sampling errors to dominate.}\label{ex_conv}
\end{figure}

\textit{Convergence.} Just like any Monte Carlo estimator, our WoI and variance-reduced WoI estimators have a $\mathcal{O}(1 / \sqrt{\W})$ convergence rate, where $\W$ denotes the number of walkers. The same convergence rate holds for the gradient estimators once boundary and interface query points are excluded. We illustrate the convergence behaviour for the three previous examples in Figure \ref{ex_conv}. Since the gradient estimators are inaccurate at query points located on the boundary and the interface, we verify the convergence behaviour of the problem defined by Eq.~\eqref{C0_gradient_3D_analytical_soln} with $7000$ randomly sampled query points. 

\textit{Neural Network Representation.} The WoI and variance-reduced WoI estimators have two remaining drawbacks. First, these two Monte Carlo estimators, by design, produce pointwise solutions, which can be restrictive when a continuous form is needed. The error plots presented so far indicate that both the WoI estimator and the variance-reduced WoI estimator exhibit high-frequency errors. However, high-frequency errors are typically considered undesirable artifacts, as the solutions to many physical systems and real-world problems often exhibit smooth behaviour.

Incorporating a neural network addresses both issues. The universal approximation property of neural networks~\cite{Kratsios2021UniversalApproximation} guarantees a continuous representation of the solution. The neural network learns a coordinate-to-value map, serving as a continuous surrogate for our Monte Carlo estimators, which enables the solution to be evaluated at arbitrary query points without re-running the estimators. In addition, empirical and theoretical studies~\cite{Rahaman_2019, Bietti_Mairal_2019, Cao_2021, Fridovich‑Keil_2022} from the machine learning community point out that deep neural networks exhibit spectral bias, fitting low-frequency components first before other complex components. 

Here, we consider a simple multilayer perceptron (MLP) with the $tanh$ activation function as the low-pass filter for the errors. Our MLP consists of 6 layers: an input layer whose dimension matches that of the interface problem, four hidden layers of width 64, and a 1-dimensional output layer. We denote the neural network representation of $u(\bm{x})$ by $u_\theta(\bm{x})$, where $\theta$ denotes the trainable parameters. The loss function is defined as the weighted sum of the mean squared errors over the interior and boundary training points,
\begin{equation}\label{MLP_loss}
    \mathcal{L}(\theta) = \frac{\mu}{\mathcal{N}} \sum_{\bm{x} \in \mathcal{X}} |u_\theta(\bm{x}) - \hat{u}(\bm{x})|^2 + \frac{1}{\mathcal{N}_\text{bc}} \sum_{\bm{x} \in \mathcal{X}_\text{bc}} |\partial_{\bm{n}}u_\theta(\bm{x}) - b_1(\bm{x})|^2,  
\end{equation}
in which $\mathcal{X}$ and $\mathcal{X}_\text{bc}$ are the sets of training points in $\Omega$ and on $\bdr{\Omega}$, respectively, and $\mathcal{N}$ and $\mathcal{N}_\text{bc}$ are their corresponding cardinalities. The proposed two-term data-driven loss is more straightforward to balance than a PINN-style loss \cite{HE2022114358, WU2024113217} since the number of terms in the latter design scales linearly as the number of interfaces grows. For all our examples, we train the MLP for 15000 epochs using the ADAM optimizer with $\mu = 2$ and a learning rate of $10^{-4}$.

\begin{figure}[htbp]
\centering
\begin{tabular}{@{}c@{} c c c c c}
% ============= WOI ==============
\multirow{4}{*}{
  \rotatebox[origin=c]{90}{%
    \parbox[c]{\dimexpr 0.65\columnwidth\relax}{\centering\vfill\textbf{Example 3 (2D)}\vfill}%
  }
} & \textbf{WoI Solution} & \raisebox{-.25\height}{\shortstack{\textbf{Learned Solution} \\ \textbf{(over training points)}}} & {} & \raisebox{-.25\height}{\shortstack{\textbf{Learned Error} \\ \textbf{(training error)}}} & {} \\
% First Row: estimation
{} & \includegraphics[width=0.25\columnwidth]{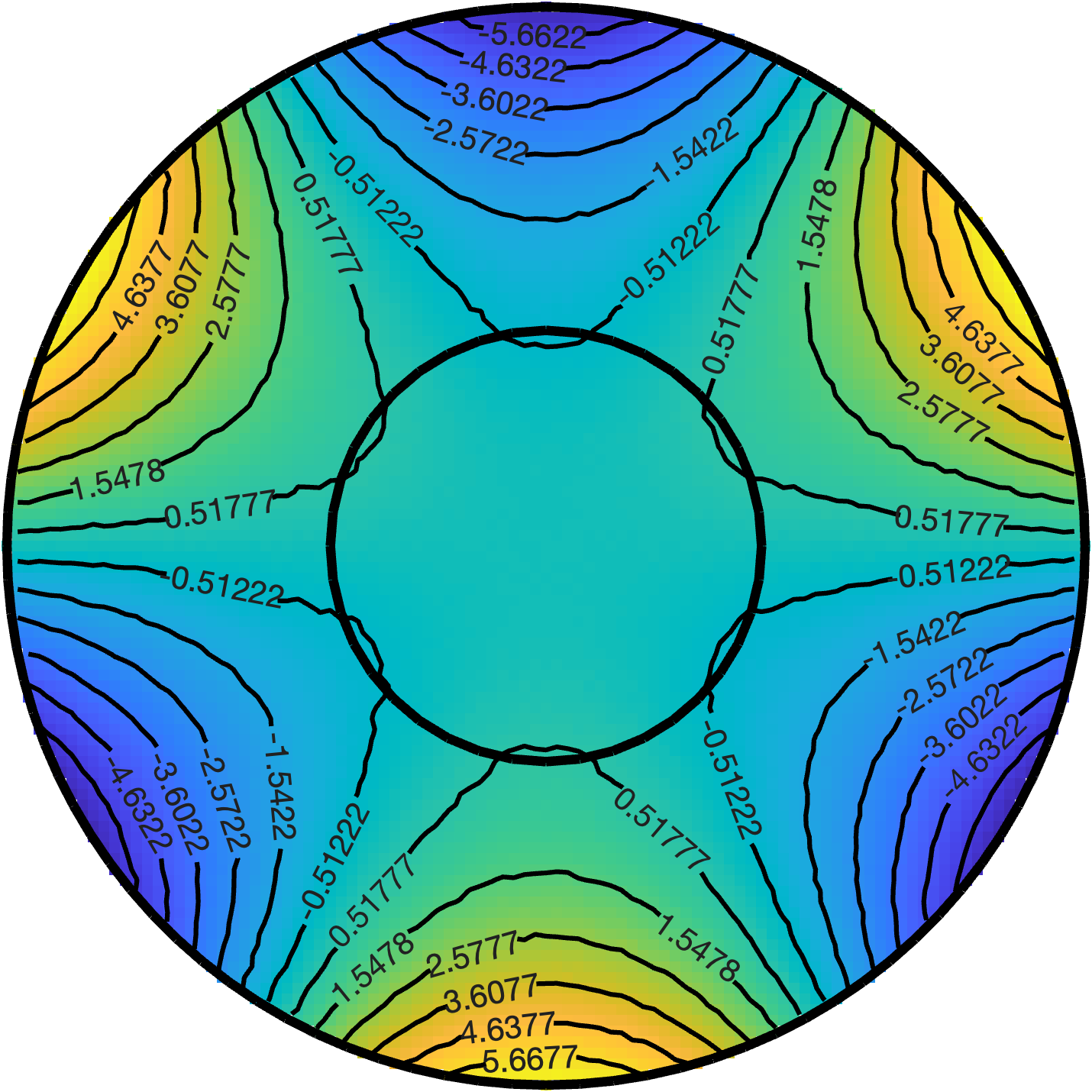} &
\includegraphics[width=0.25\columnwidth]{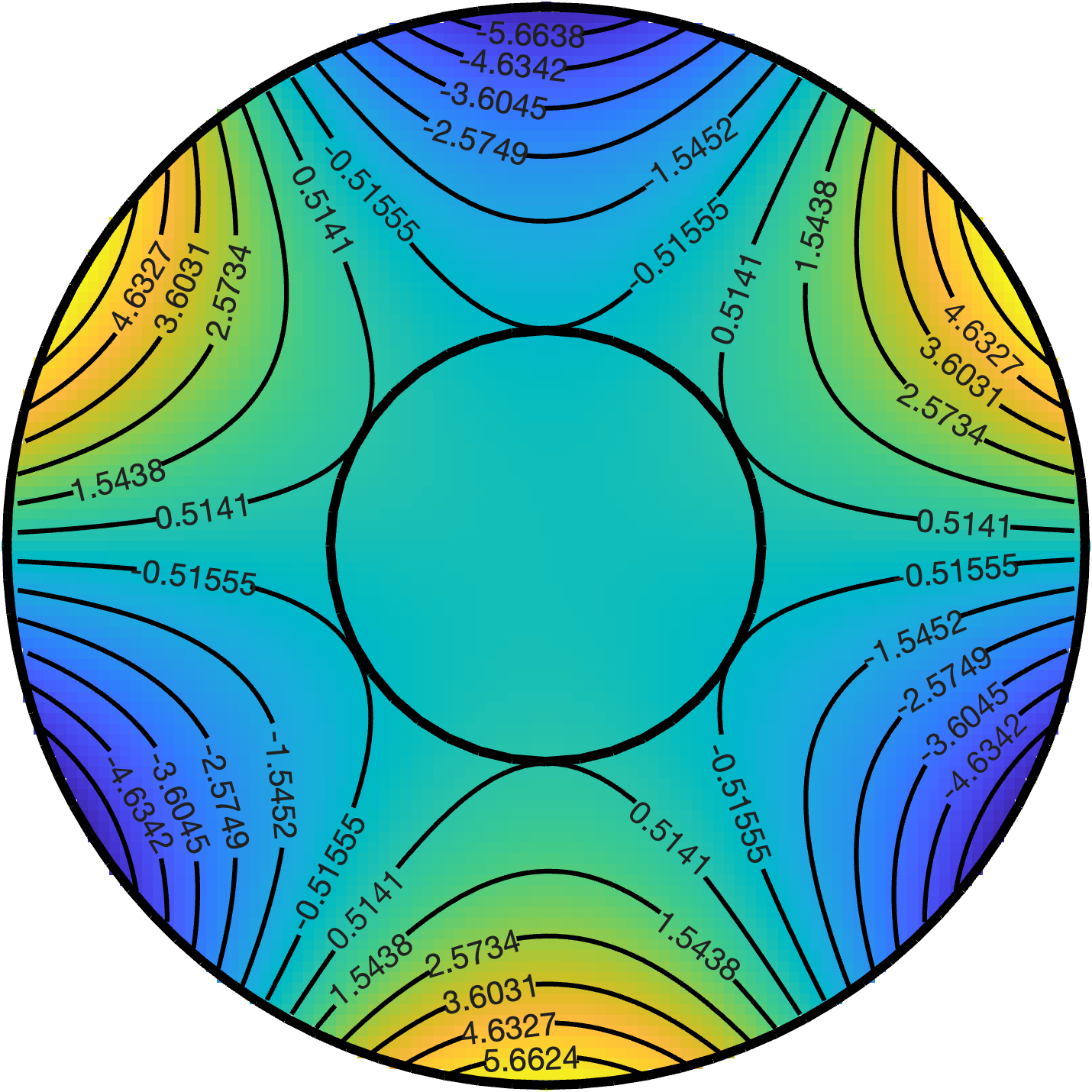}  &
\includegraphics[width=0.038\columnwidth]{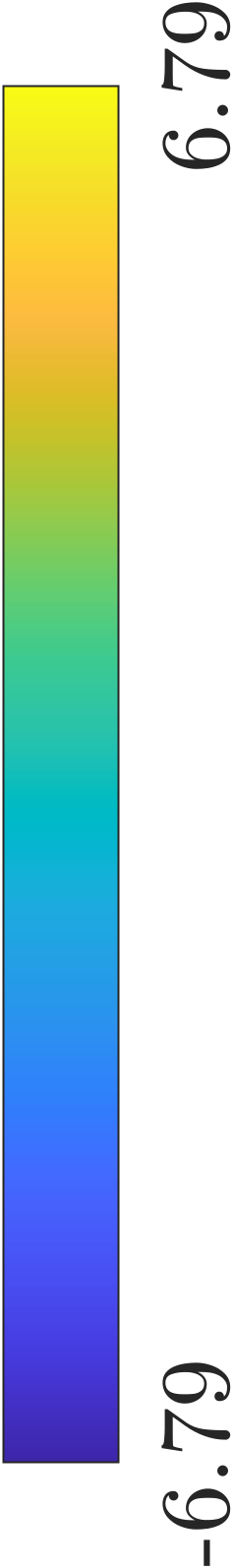} &
\includegraphics[width=0.25\columnwidth]{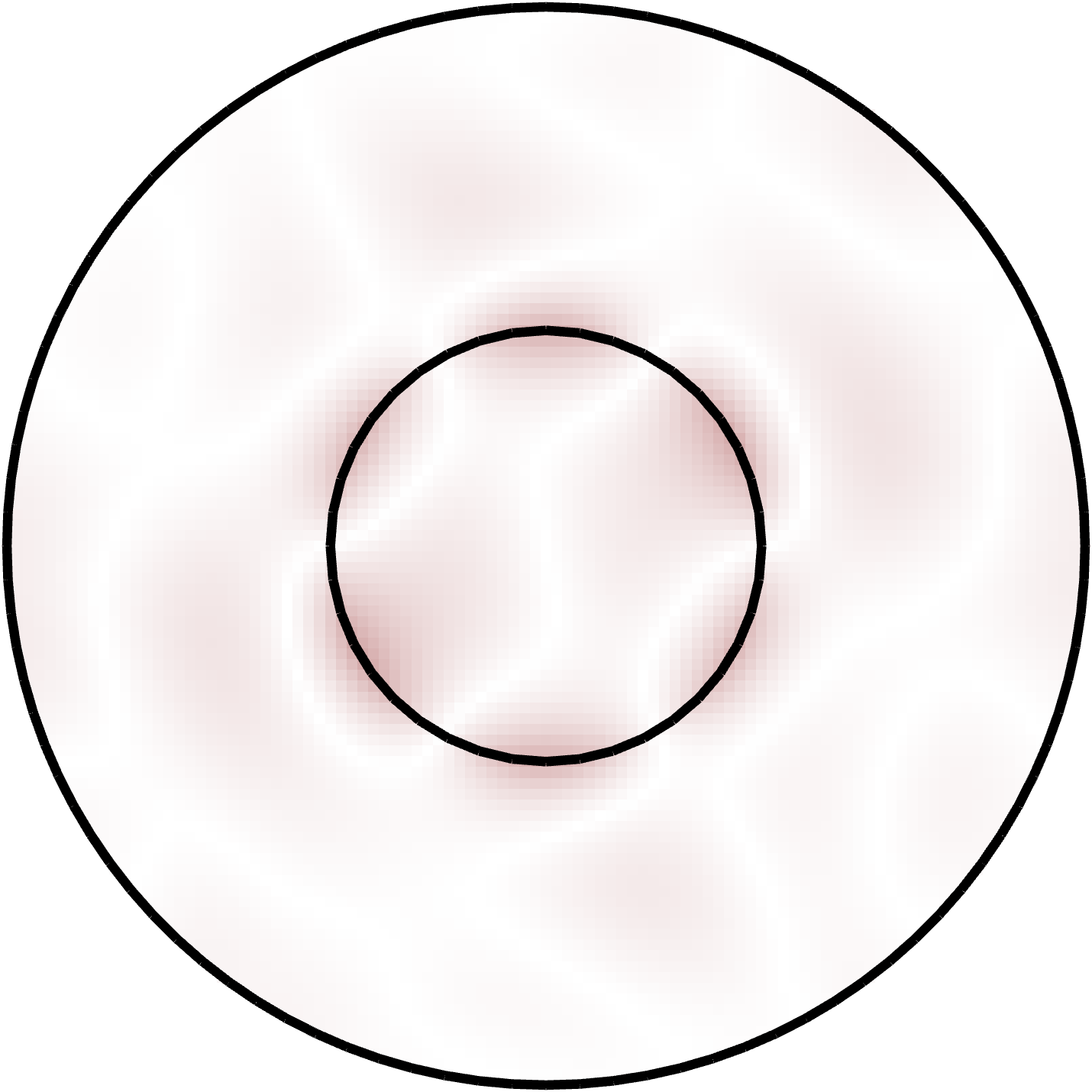} &
\includegraphics[width=0.04\columnwidth]{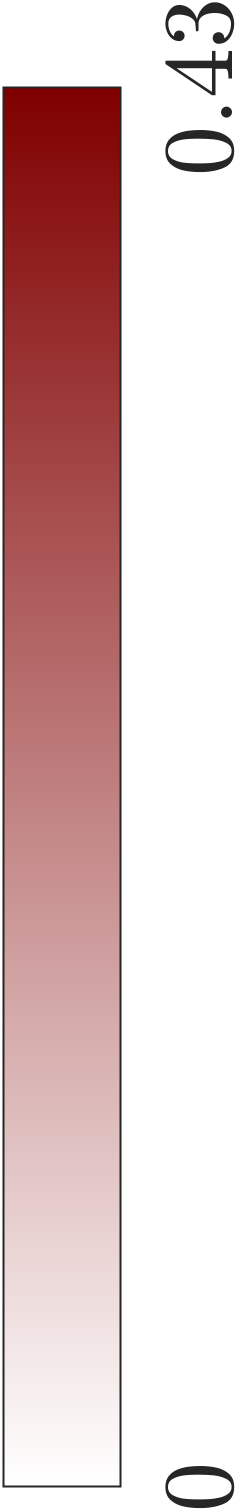} \vspace{0.3cm} \\

{} & \raisebox{-.25\height}{\shortstack{\textbf{Var-reduced} \\ \textbf{WoI Solution}}} & \raisebox{-.25\height}{\shortstack{\textbf{Learned Solution} \\ \textbf{(over training points)}}} & {} & \raisebox{-.25\height}{\shortstack{\textbf{Learned Error} \\ \textbf{(training error)}}} & {} \\
% Second Row: pointwise estimation
{} & \includegraphics[width=0.25\columnwidth]{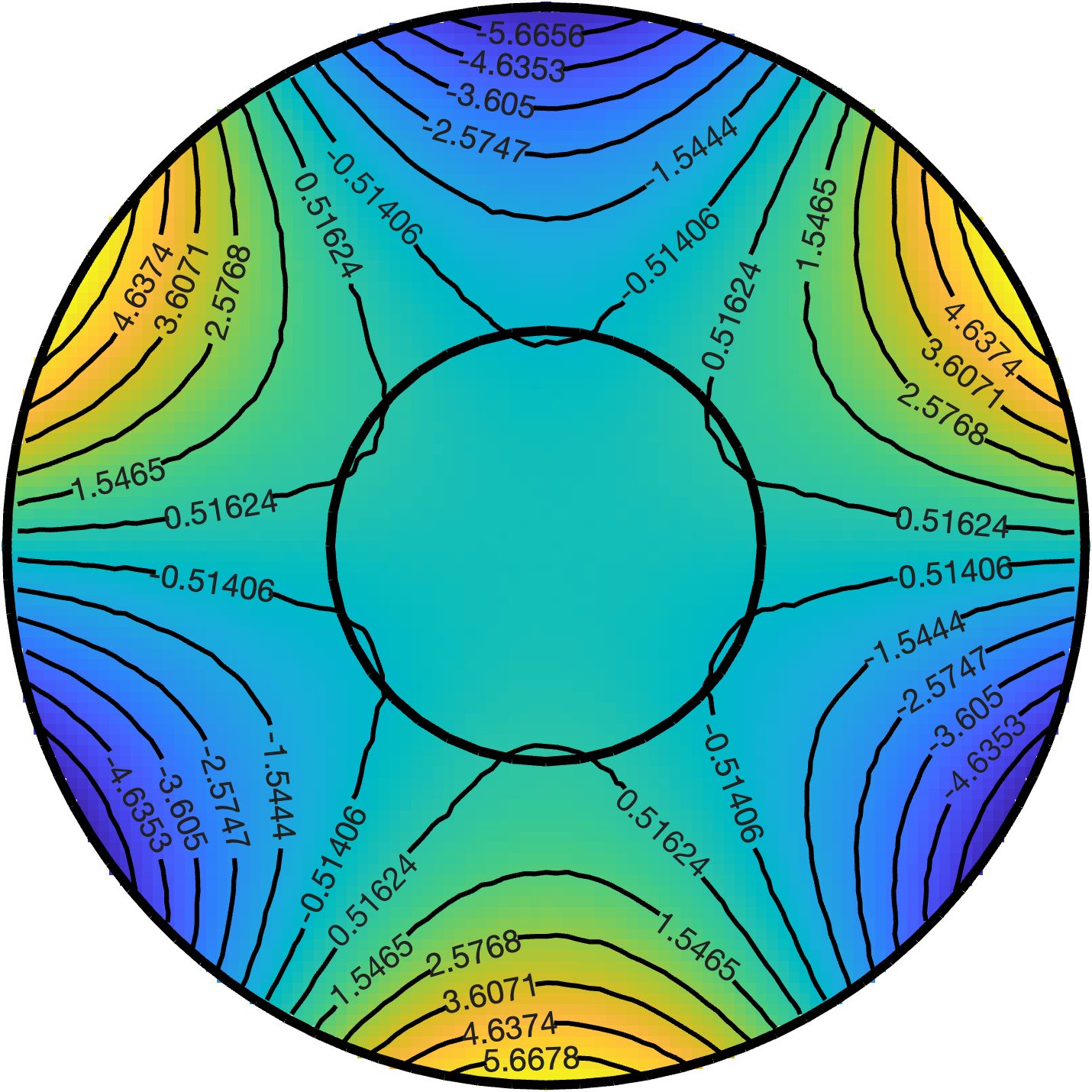} &
\includegraphics[width=0.25\columnwidth]{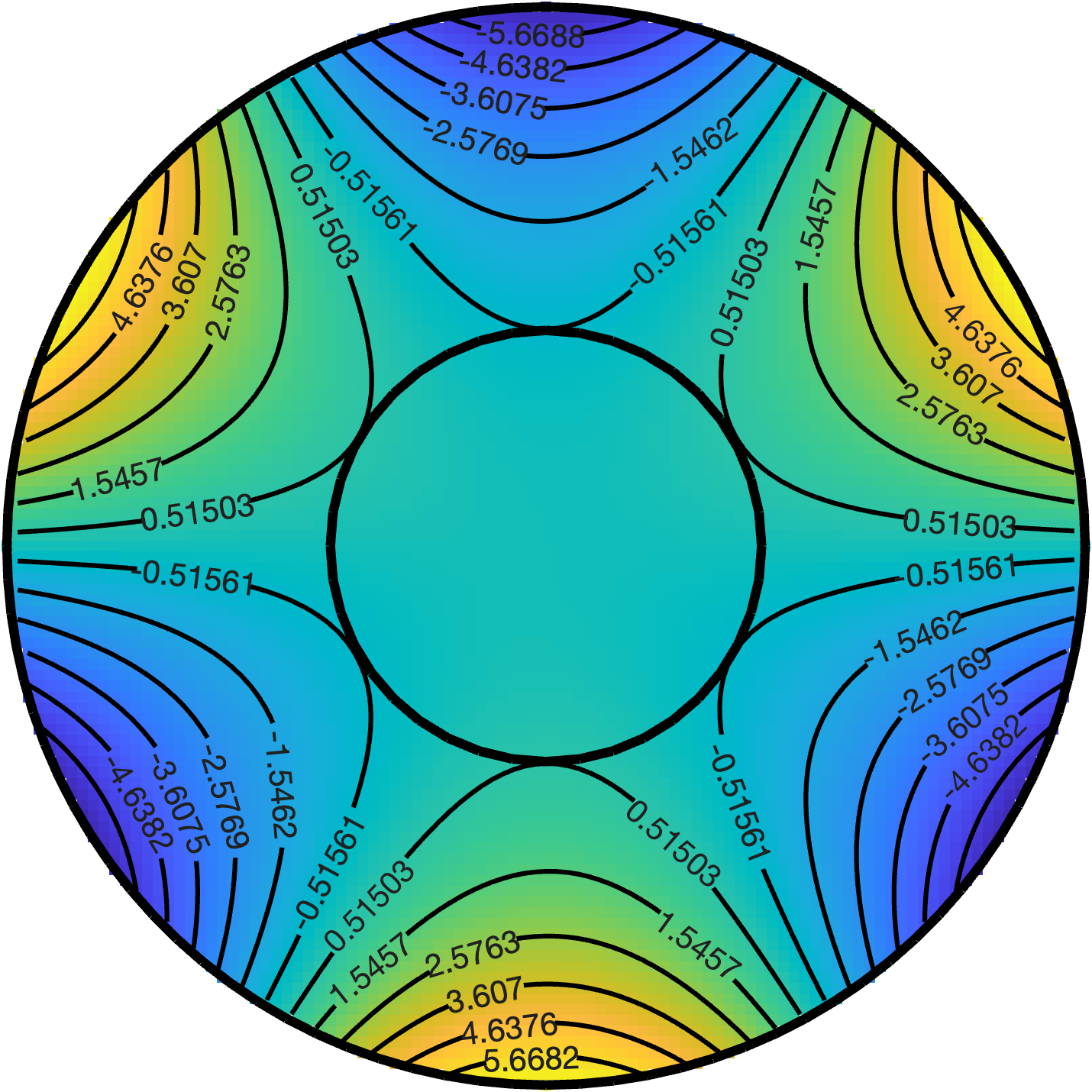}  &
\includegraphics[width=0.038\columnwidth]{figure/kink_in_gradient2D_nn/estimation_cb.png} &
\includegraphics[width=0.25\columnwidth]{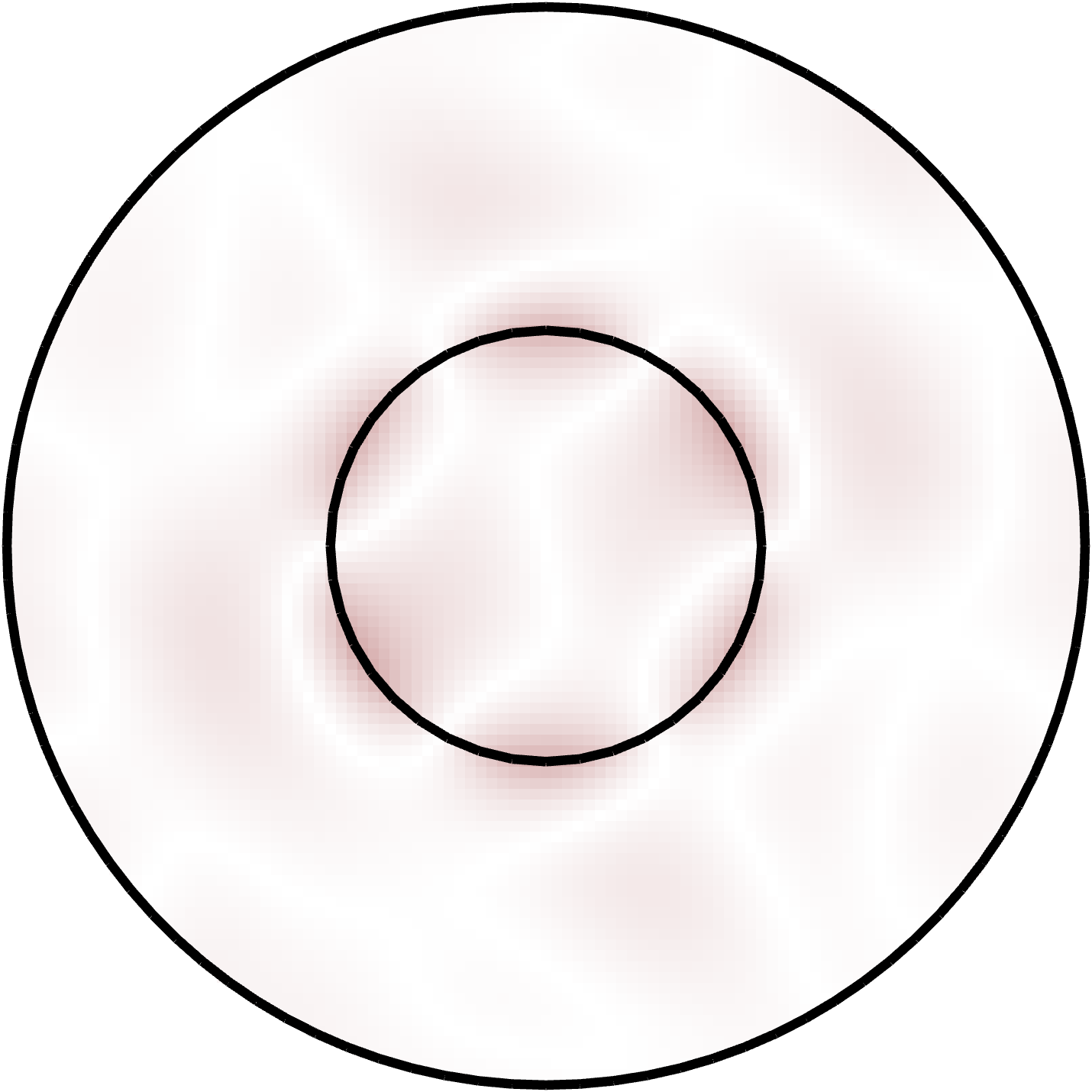} &
\includegraphics[width=0.04\columnwidth]{figure/kink_in_gradient2D_nn/error_cb.png}
\end{tabular}
\caption{Solution and error of Eq.~\eqref{C0_gradient_2D_analytical_soln} learned by MLP assumed by the coordinates of training points. The contour lines in both the WoI solution and the variance-reduced WoI solution exhibit noticeable waviness, whereas those in the learned solution are smooth. In contrast to the pixelated error pattern in Figure \ref{ex_woi_kink_gradient2D}, The error in the learned solution changes smoothly and is concentrated along the interface, where the gradient is $C^0$. This error distribution reflects the frequency content present in the ground-truth solution. For consistency, we use the same color scale for errors as in Figure \ref{ex_woi_kink_gradient2D} for the purpose of comparing the learned error and the high-frequency Monte Carlo error.}\label{kink_gradient2D_smooth_error}
\end{figure}

\begin{figure}[htbp]
\centering
\setlength{\tabcolsep}{3pt}
\begin{tabular}{@{}c@{} c c c c c}
% ============= WOI ==============
\multirow{4}{*}{
  \rotatebox[origin=c]{90}{%
    \parbox[c]{\dimexpr 0.65\columnwidth\relax}{\centering\vfill\textbf{Example 3 (3D) -- Boundary}\vfill}%
  }
} & \textbf{WoI Solution} & \raisebox{-.25\height}{\shortstack{\textbf{Learned Solution} \\ \textbf{(over training points)}}} & {} & \raisebox{-.25\height}{\shortstack{\textbf{Learned Error} \\ \textbf{(training error)}}} & {} \\

% First Row: estimation
{} & \includegraphics[width=0.26\columnwidth]{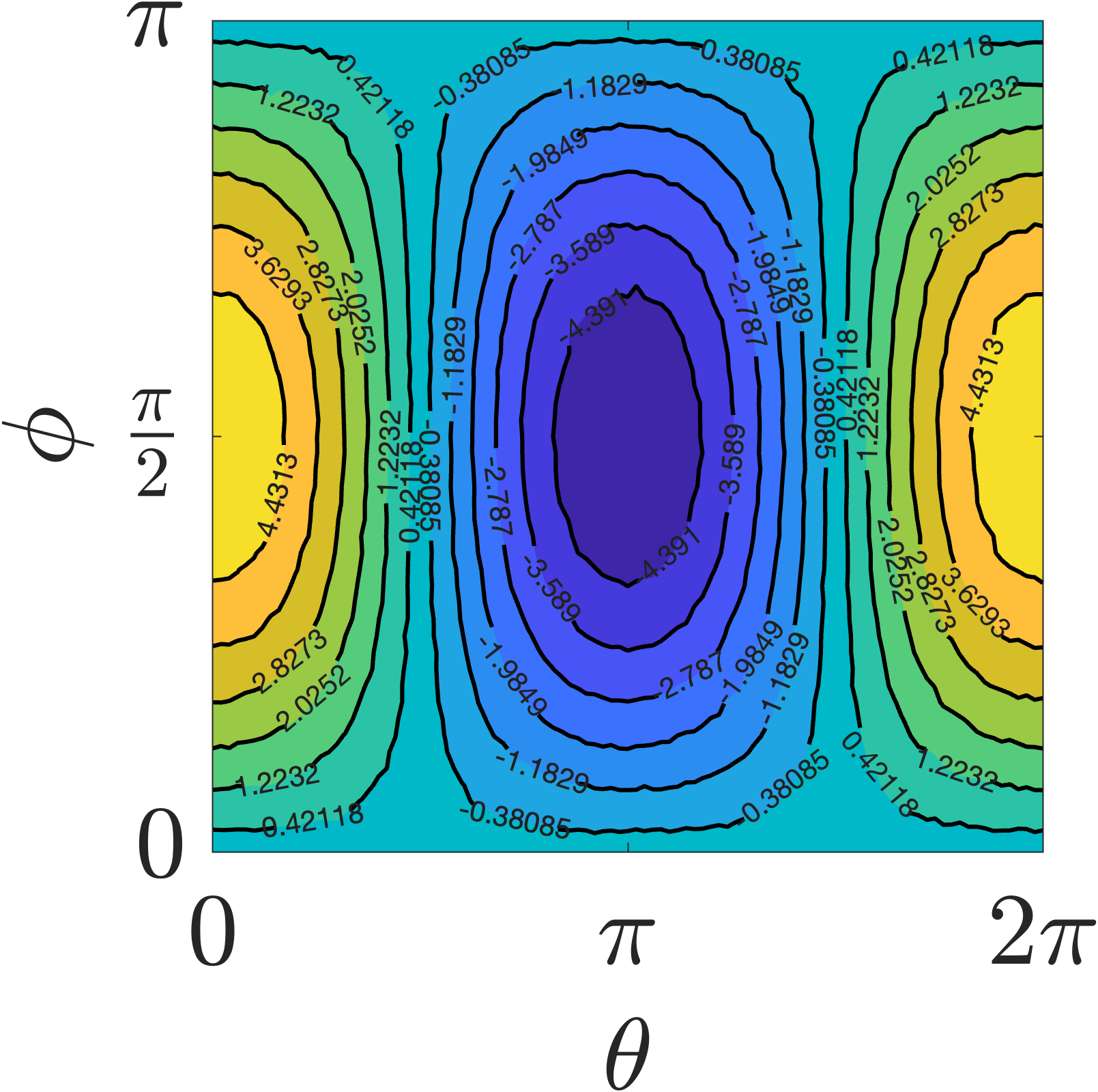} &
\includegraphics[width=0.26\columnwidth]{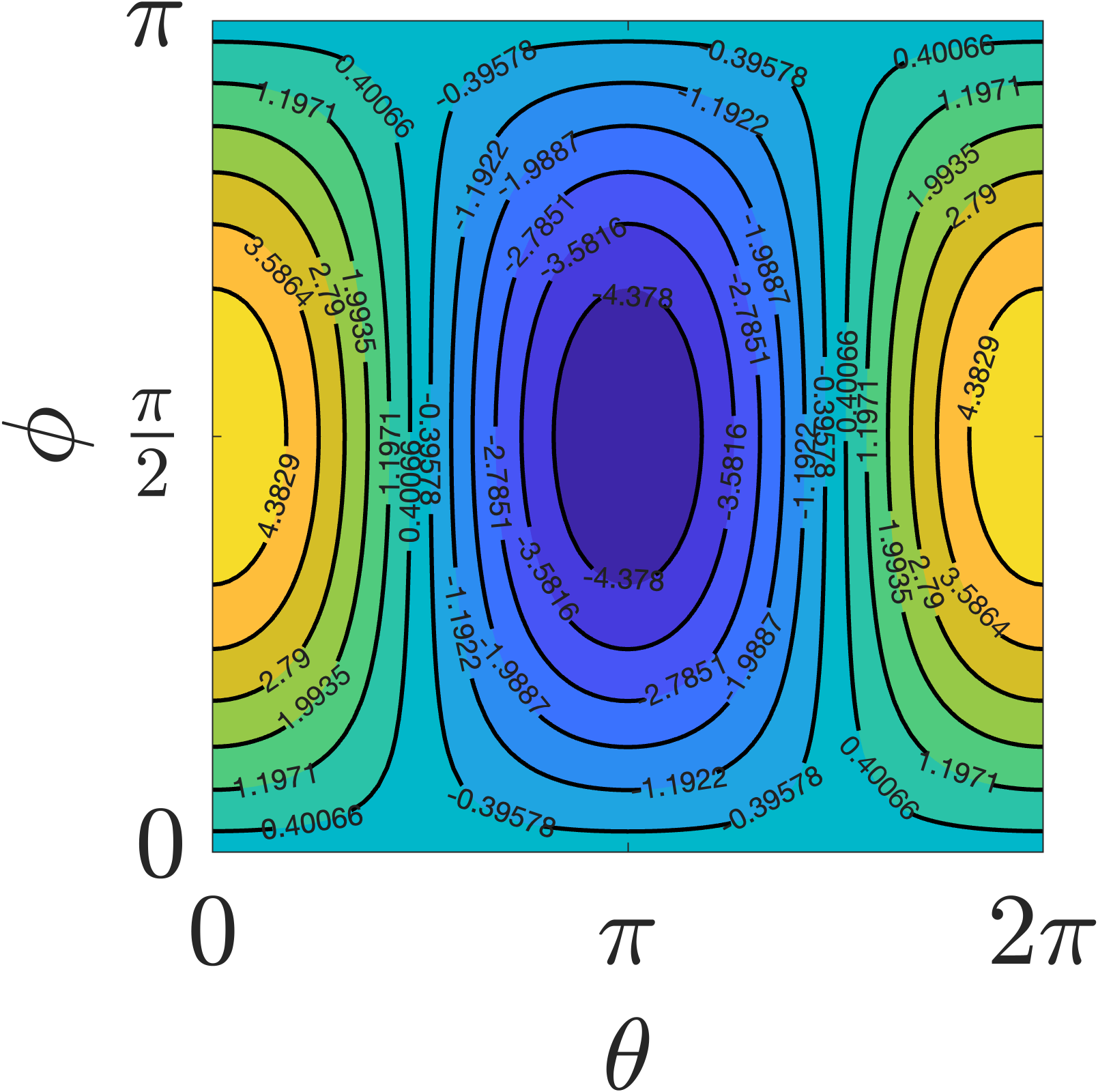}  &
\includegraphics[width=0.04\columnwidth]{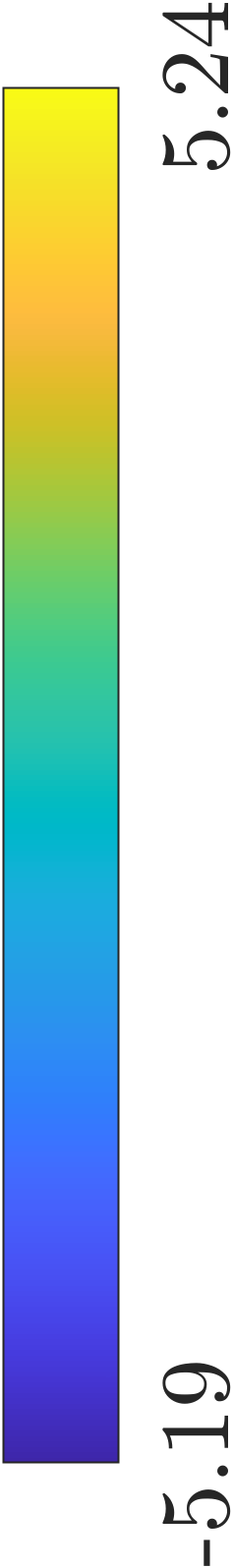} &
\includegraphics[width=0.26\columnwidth]{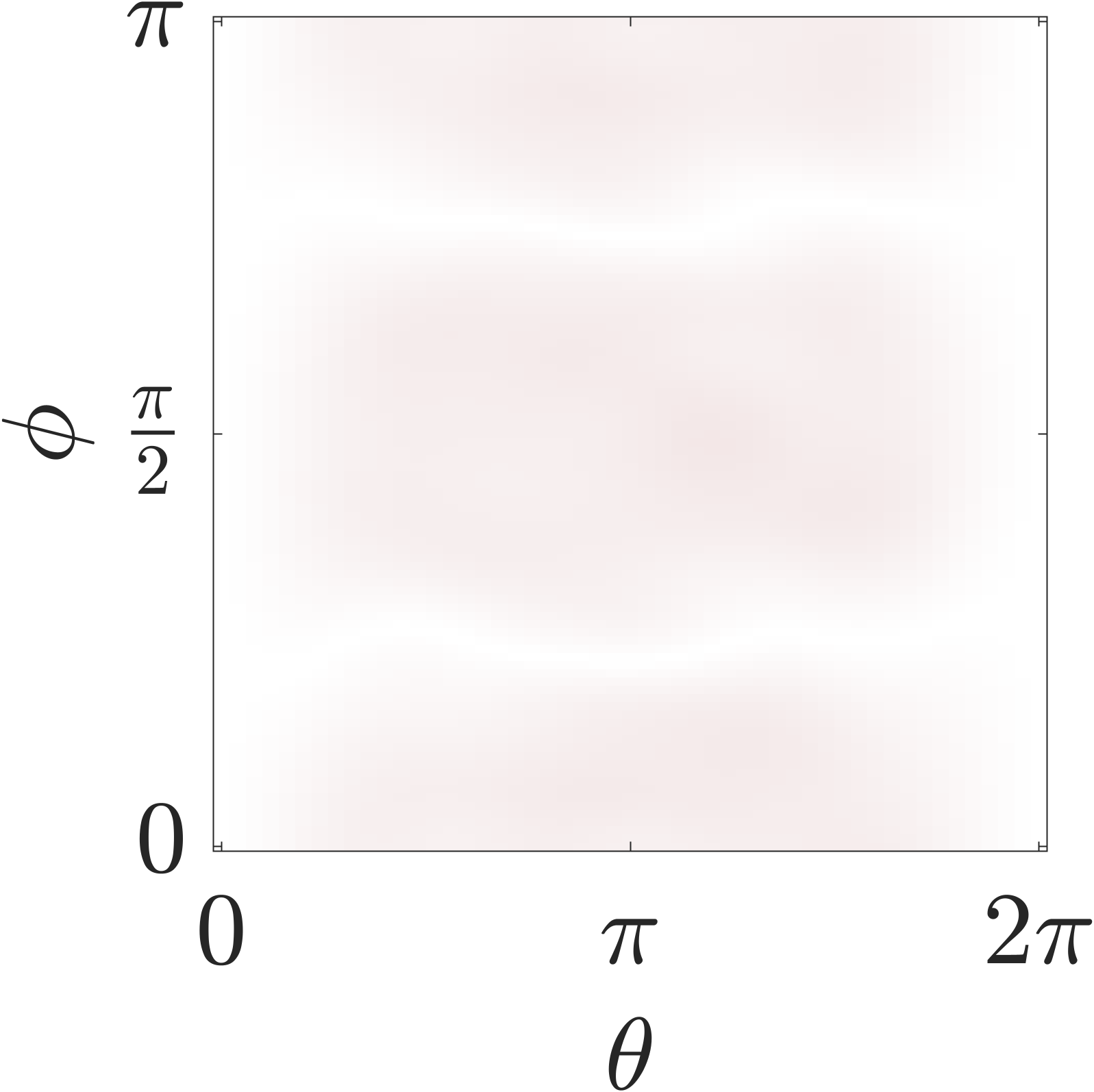} &
\includegraphics[width=0.042\columnwidth]{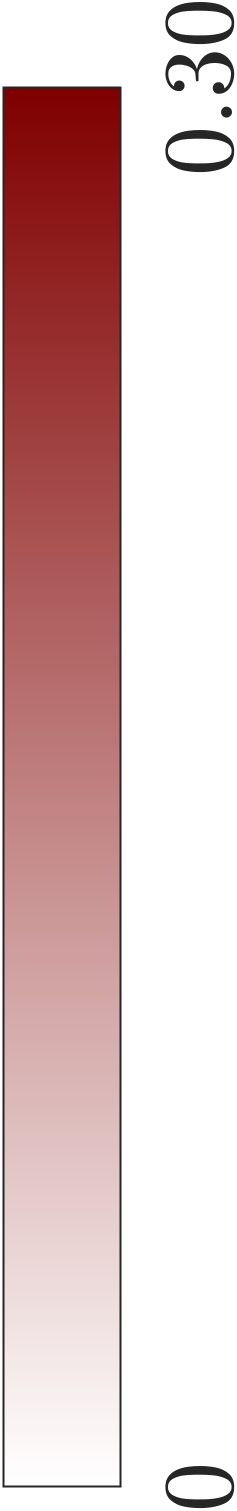} \\

{} & \raisebox{-.25\height}{\shortstack{\textbf{Var-reduced} \\ \textbf{WoI Solution}}} & \raisebox{-.25\height}{\shortstack{\textbf{Learned Solution} \\ \textbf{(over training points)}}} & {} & \raisebox{-.25\height}{\shortstack{\textbf{Learned Error} \\ \textbf{(training error)}}} & {} \\
% Second Row: pointwise estimation
{} & \includegraphics[width=0.26\columnwidth]{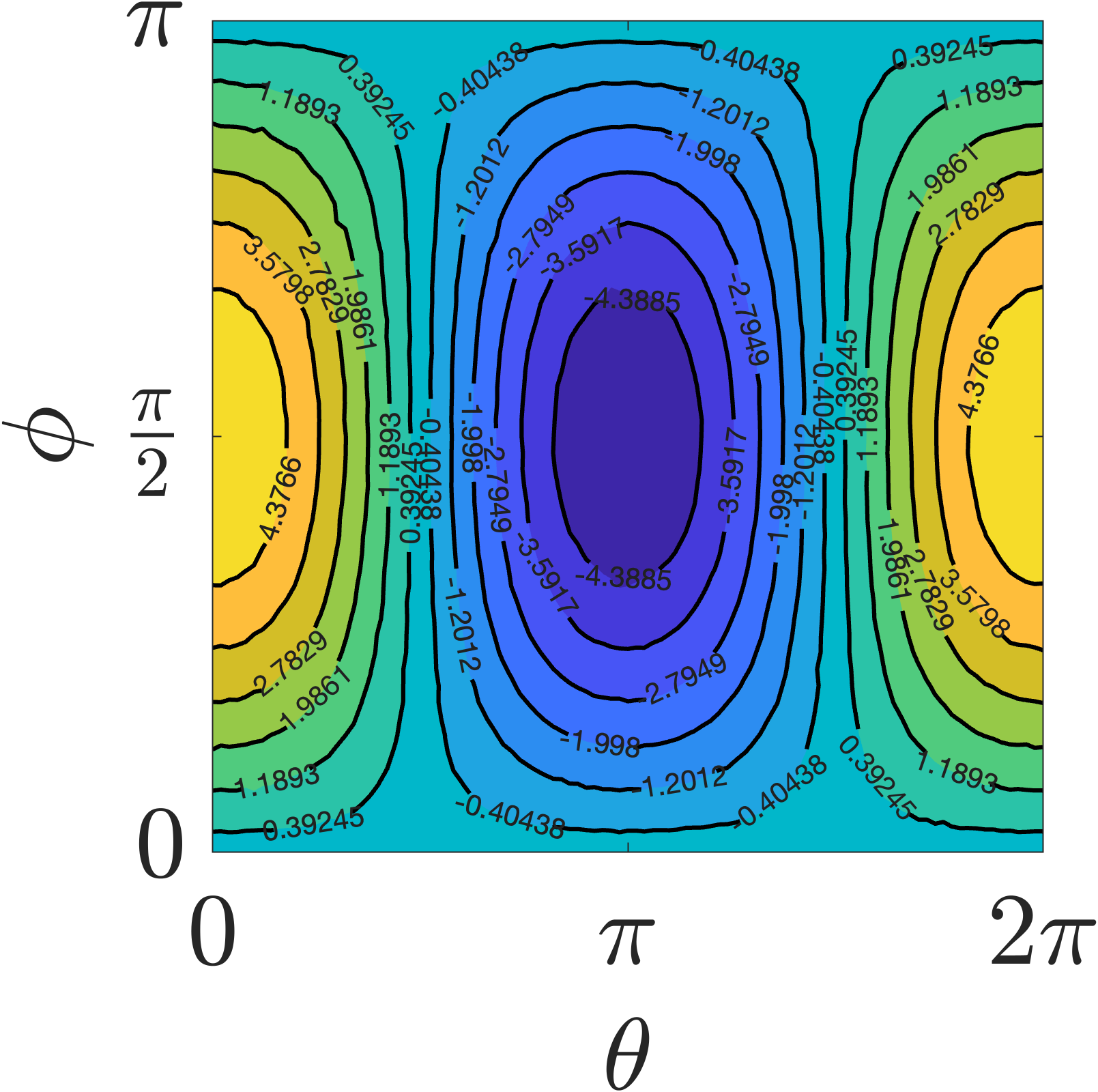} &
\includegraphics[width=0.26\columnwidth]{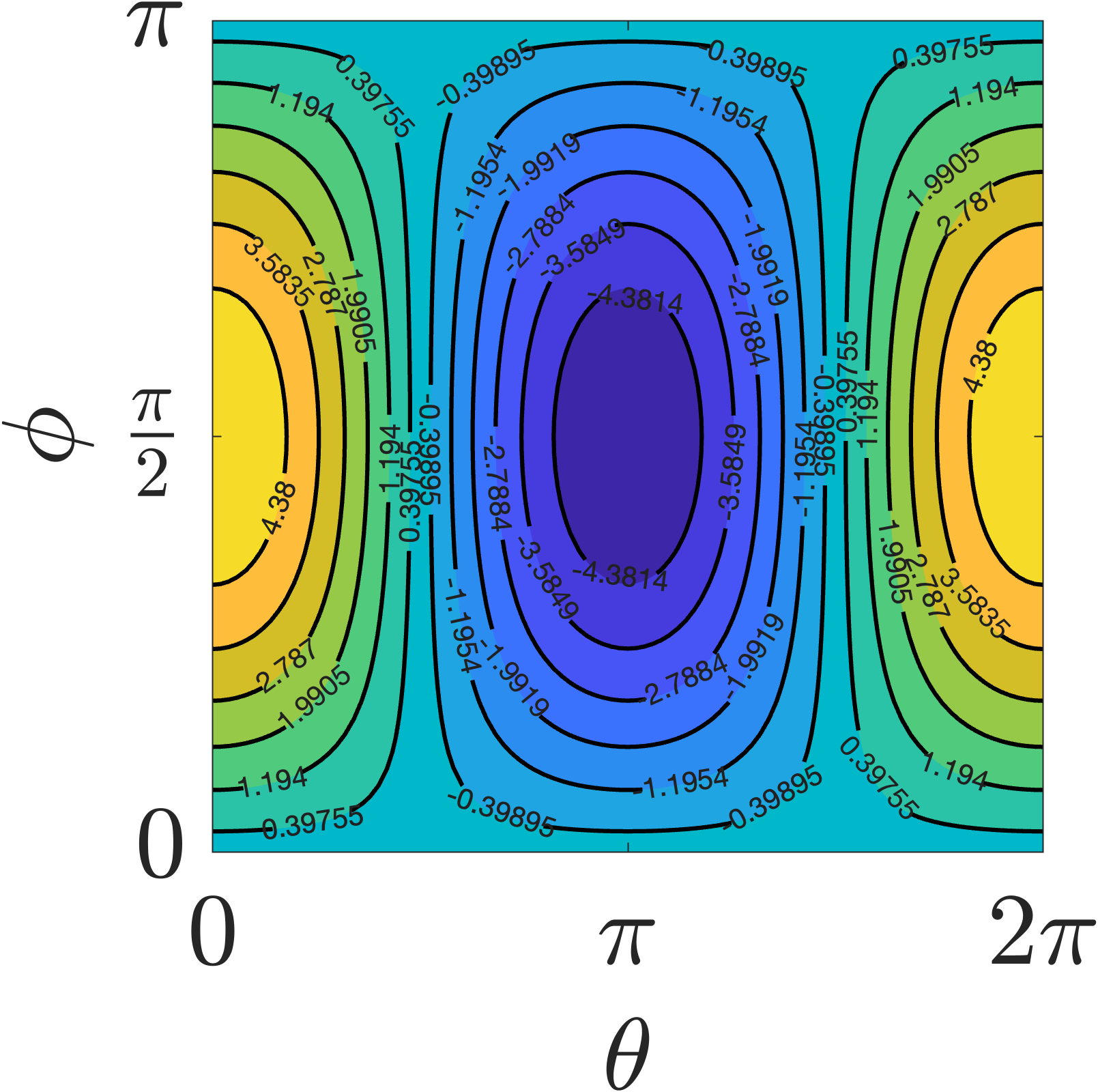}  &
\includegraphics[width=0.04\columnwidth]{figure/kink_in_gradient3D_nn/estimation_cb.png} &
\includegraphics[width=0.26\columnwidth]{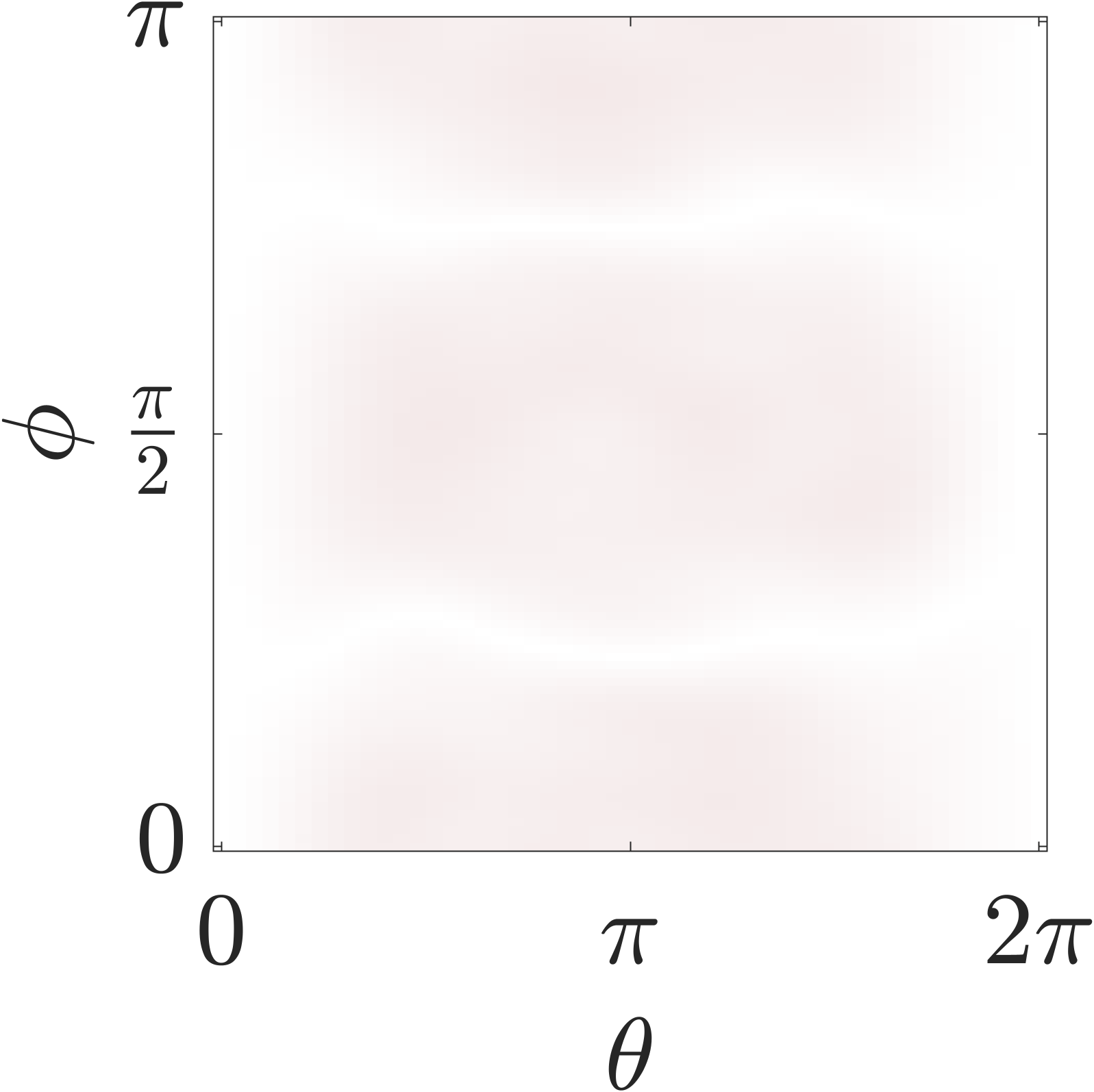} &
\includegraphics[width=0.042\columnwidth]{figure/kink_in_gradient3D_nn/error_cb.png}
\end{tabular}
\caption{Solution and error of Eq.~\eqref{C0_gradient_3D_analytical_soln} on the boundary of the domain learned by MLP assumed by the coordinates of training points. The boundary of the spherical domain is presented in $\theta-\phi$ parametric space. The contour lines in both the WoI solution and the variance-reduced WoI solution are wiggly, whereas those in the learned solutions change smoothly. The error in the learned solution change gently comparing with the ones in Figure \ref{ex_woi_kink_gradient3D}. For consistency, we use the same color scale for error as in Figure~\ref{ex_woi_kink_gradient3D} for the purpose of comparing the learned error and the high-frequency Monte Carlo error.}\label{kink_gradient3D_smooth_error}
\end{figure}

Figure \ref{kink_gradient2D_smooth_error} and Figure \ref{kink_gradient3D_smooth_error} present the learned solutions and learned errors to Eq.~\eqref{C0_gradient_2D_problem} and Eq.~\eqref{C0_gradient_3D_problem}, respectively. 
For both interface problems, interior training points are sampled on a uniform grid within the circular (or spherical) domain, while boundary training points are uniformly distributed in angular coordinates along the boundary. For Eq.~\eqref{C0_gradient_2D_problem}, we use $\mathcal{N} = 7845$ points in the domain and $\mathcal{N}_\text{bc} = 500$ points on the boundary. For Eq.~\eqref{C0_gradient_2D_problem}, we use $\mathcal{N} = 8326$ points in the domain and $\mathcal{N}_\text{bc} = 5151$ boundary points. The set of test points contains 8000 randomly sampled points within the domain. In both figures, the learned solutions and errors are plotted at the coordinates of the training points. Table \ref{nn_training_test_error_table} reports the $L_2$ and relative $L_2$ errors of the learned solutions, using training data generated by WoI and variance-reduced WoI, evaluated on both training and test points.

\begin{table}
\centering
\captionsetup{justification=centering}
\caption{MLP training and test errors with respect to the true solution. For cases without an analytical ground truth, the numerical solution is used as the reference.}
\label{nn_training_test_error_table}
\renewcommand{\arraystretch}{1.5}
\begin{tabular}{c c c c c c} \toprule
\textbf{Problem} & \makecell{\textbf{Data} \\ \textbf{Generator}} & \makecell{$\bm{L_2}$ \textbf{Training} \\ \textbf{Error}} & \makecell{\textbf{Relative} $\bm{L_2}$ \\ \textbf{Training Error}} & \makecell{$\bm{L_2}$ \\ \textbf{Test Error}} & \makecell{\textbf{Relative} $\bm{L_2}$ \\ \textbf{Test Error}}\\ \hline 
Ex.3 (2D) & WoI & 0.0288 & 1.21\% & 0.0289 & 1.23\% \\
Ex.3 (2D) & \makecell{Var-reduced \\ WoI} & 0.029 & 1.22\% & 0.029 & 1.21\% \\
Ex.3 (3D) & WoI & 0.015 & 0.592\% & 0.022 & 1.01\% \\
Ex.3 (3D) & \makecell{Var-reduced \\ WoI} & 0.015 & 0.580\% & 0.021 & 0.98\% \\
Ex.4 & WoI & 0.073 & 3.39\% & 0.096 & 3.51\% \\
Ex.4 & \makecell{Var-reduced \\ WoI} & 0.056 & 2.64\% & 0.093 & 3.41\% \\
\bottomrule
\end{tabular}
\end{table}

Before closing this part of the discussion, we note that numerous other neural network architectures with various activation functions can serve as high-frequency error filters, such as the variant of ResNet~\cite{HAN2020109672}. Yet, identifying the optimal neural network design is beyond the scope of this work.

\textit{More Applications.} Having established the effectiveness of our WoI framework on examples with known ground truth, we now turn to problems rooted in physics and motivated by real-world applications. These problems typically lack analytical solutions, and in some cases, numerical boundary integral solvers fail due to memory constraints. Whenever a numerical solution is available, we compare both the Monte Carlo solution and the learned solution from our framework against it. However, in cases where the boundary integral solver cannot produce a solution, we demonstrate that our framework can still output reasonable and physically meaningful results.

\textbf{Example 4. The Electrical Conductivity Problem.} The electrical conductivity problem models the electrical potential induced by a current in a domain of a conductive medium, given a spatially varying electrical conductivity. The problem serves as the foundational forward problem for a wide range of inverse problems, including Direct Current Resistivity (DCR), Electrical Resistivity Tomography (ERT), and Electrical Impedance Tomography (EIT). We interpret $\sigma(\bm{x})$ as the electrical conductivity, $u(\bm{x})$ as the electrical potential, and consequently, $-\nabla u(\bm{x})$ as the charge density. Since conductivity is an intrinsic material property, it is natural to model $\sigma(\bm{x})$ as a piecewise constant function.

We consider a domain with a cow-shaped insulator enclosed by a spherical conductor. Let $\Omega_1$ be the spherical conductor of radius 2 and $\Omega_2$ be the cow interface, where $\sigma_1 = 1.2$ and $\sigma_2 = 0$. The interface problem of interest, described by spherical coordinates, is
\begin{subequations}\label{EIT}
\begin{empheq}[left=\empheqlbrace]{align}
 \Delta u &= 0 \quad \text{in $\Omega \ \backslash \ \bdr{\Omega}_2$}  \\
 % ========================
 \sigma_1 \partial_{\bm{r}}u &= A \sin^9 \phi \cos^9 \theta \quad \text{on $\bdr{\Omega_1}$} \\
 % ========================
 [u](r, \theta, \phi) &= 0 \quad \text{on $\bdr{\Omega}_2$}\\
 % ========================
 [\sigma \partial_{\bm{r}} u](r, \theta, \phi) &= 0 \quad \text{on $\bdr{\Omega}_2$},
\end{empheq}
\end{subequations}
in which $\phi \in [0, \pi]$, $\theta \in [0, 2 \pi]$, and $A = 10$. 

\begin{figure}[htbp]
\centering
\begin{minipage}[l]{\textwidth}
\begin{tabular}{c c c c}
    % Row 1
    \hspace{0.8cm} \textbf{Domain} & \textbf{A Slice of Domain} & \textbf{Numerical Solution} & {} \\
    \hspace{0.8cm} \includegraphics[width=0.25\columnwidth, valign=c]{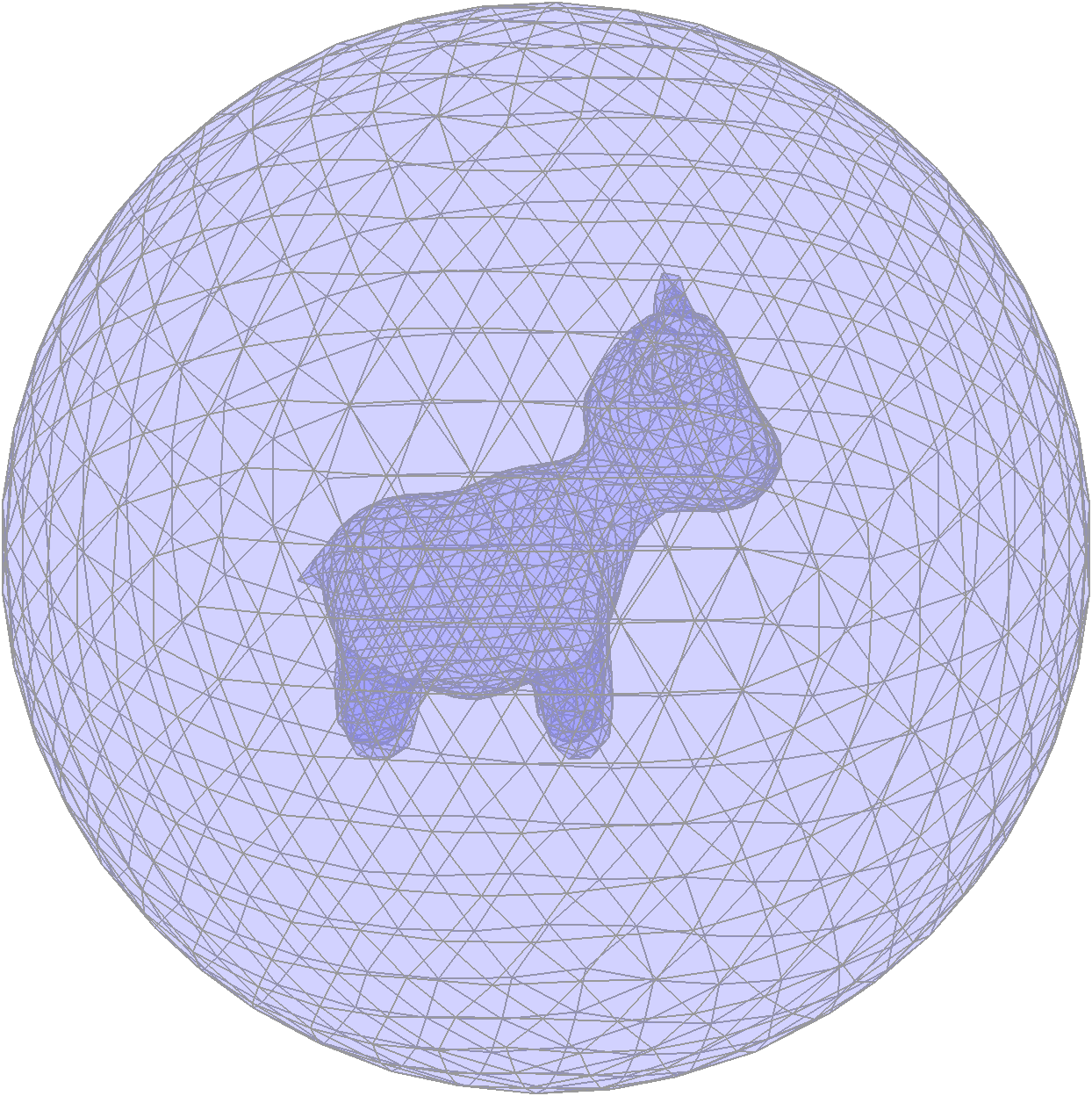} & 
    \includegraphics[width=0.25\columnwidth, valign=c]{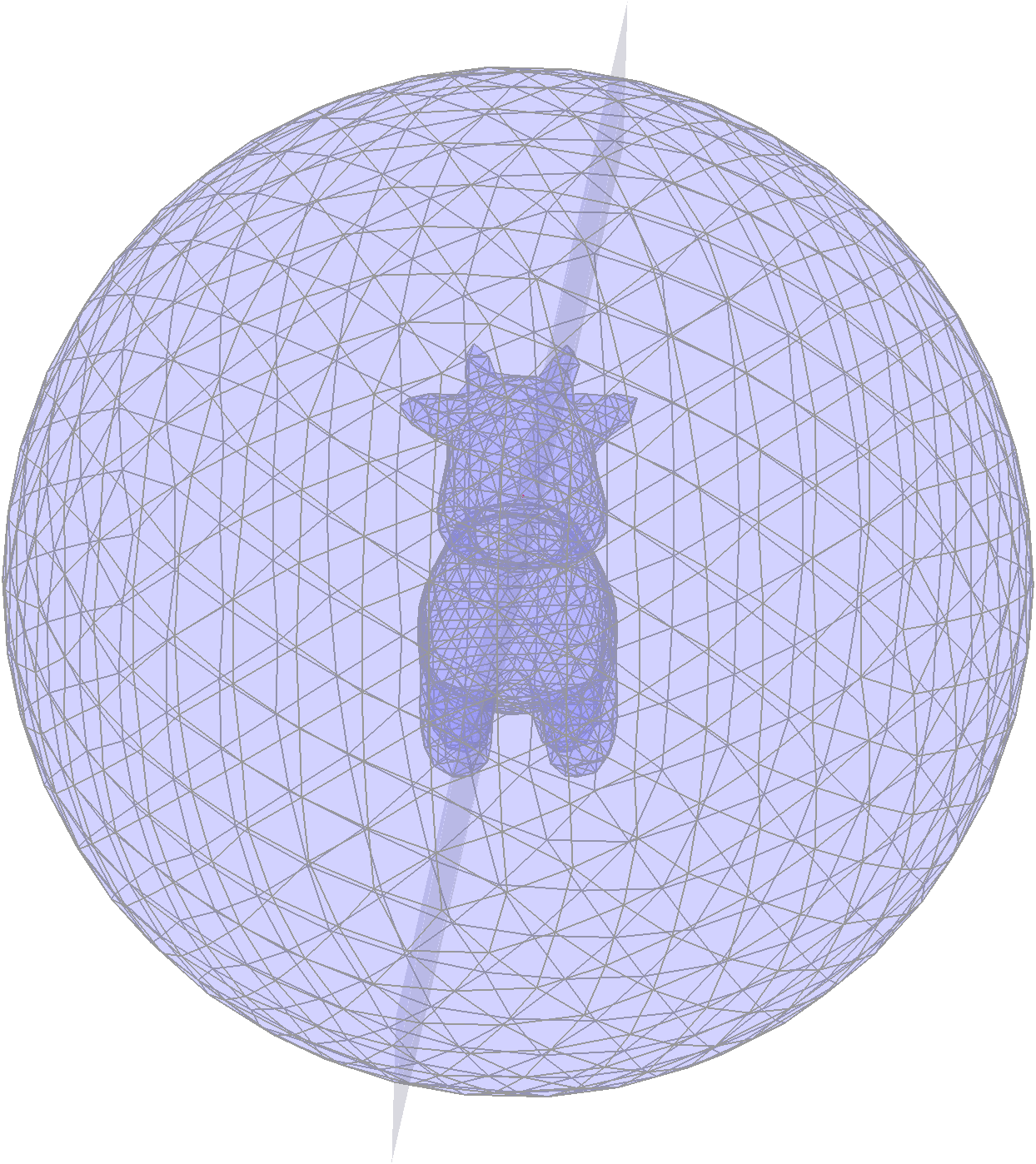} &
    \includegraphics[width=0.25\columnwidth, valign=c]{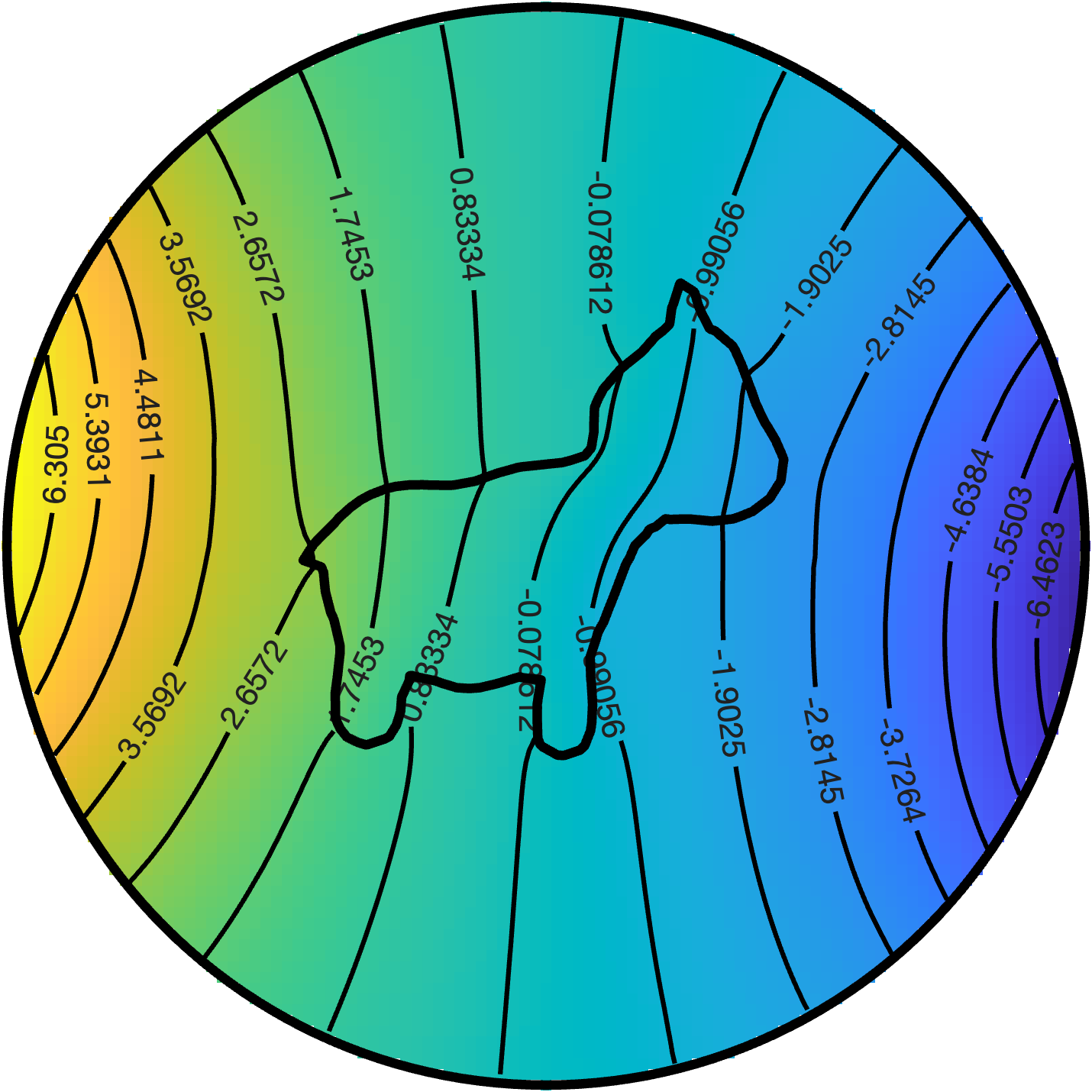} &
    \includegraphics[width=0.04\columnwidth, valign=c]{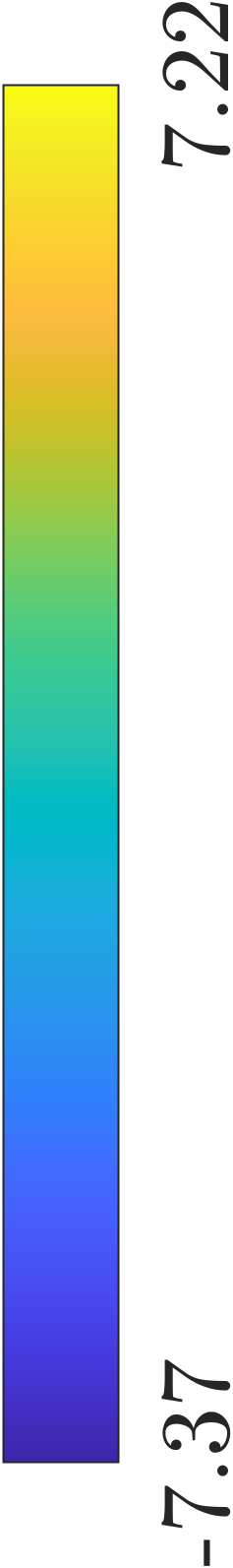}
\end{tabular}
\end{minipage}
\\ [0.5cm]
\begin{minipage}[l]{\textwidth}
\setlength{\tabcolsep}{2pt}
\begin{tabular}{c c c c c c}
    % WoI
    {} & \textbf{WoI Solution} & \raisebox{-.25\height}{\shortstack{\textbf{Learned Solution} \\ \textbf{(on test points)}}} & \textbf{WoI Error} & \raisebox{-.25\height}{\shortstack{\textbf{Learned Error} \\ \textbf{(test error)}}} & {} \vspace{0.3cm} \\
    % First Row: estimation
    \rotatebox[origin=t]{90}{\centering \textbf{WoI Framework}} & \includegraphics[width=0.21\columnwidth, valign=c]{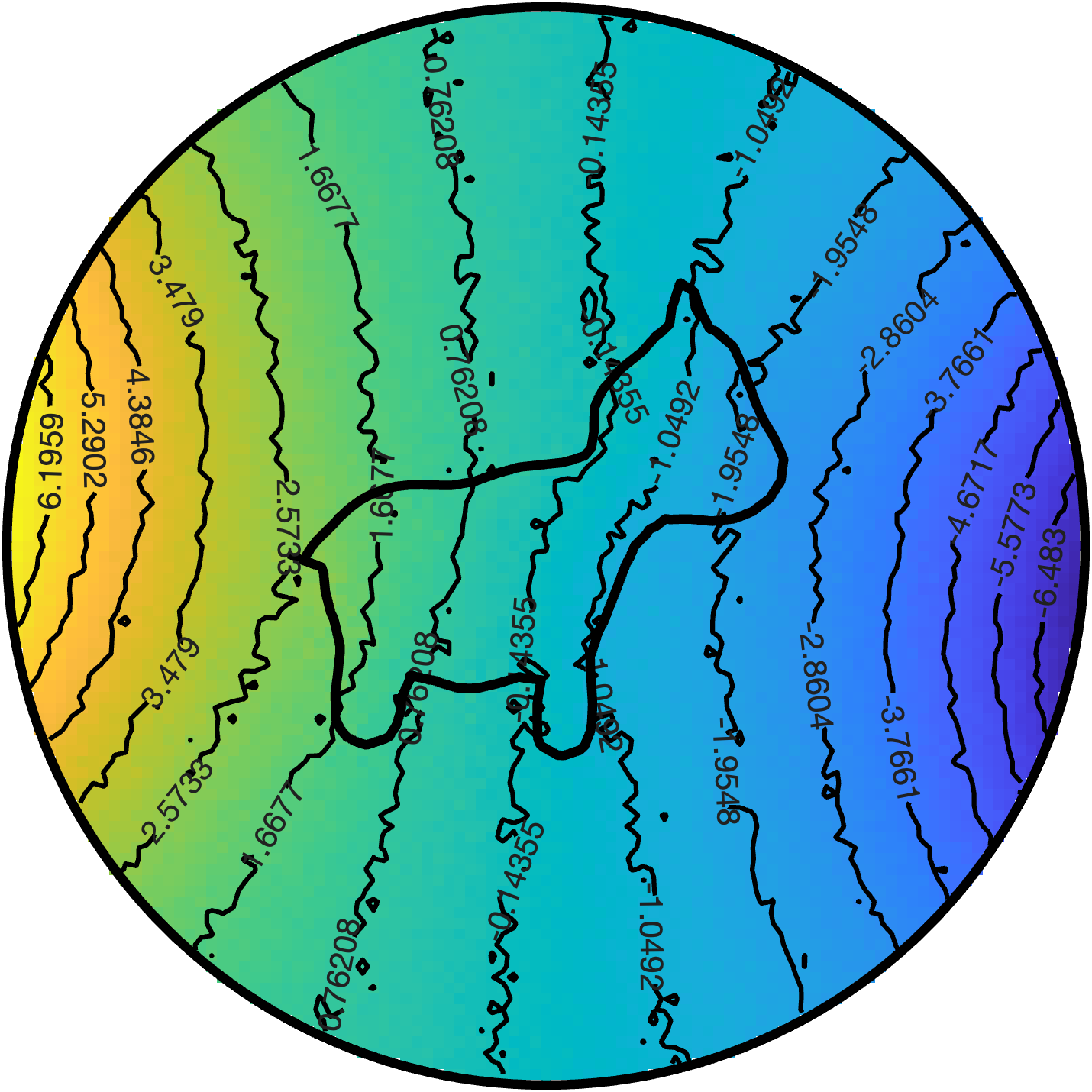}  &
    \includegraphics[width=0.21\columnwidth, valign=c]{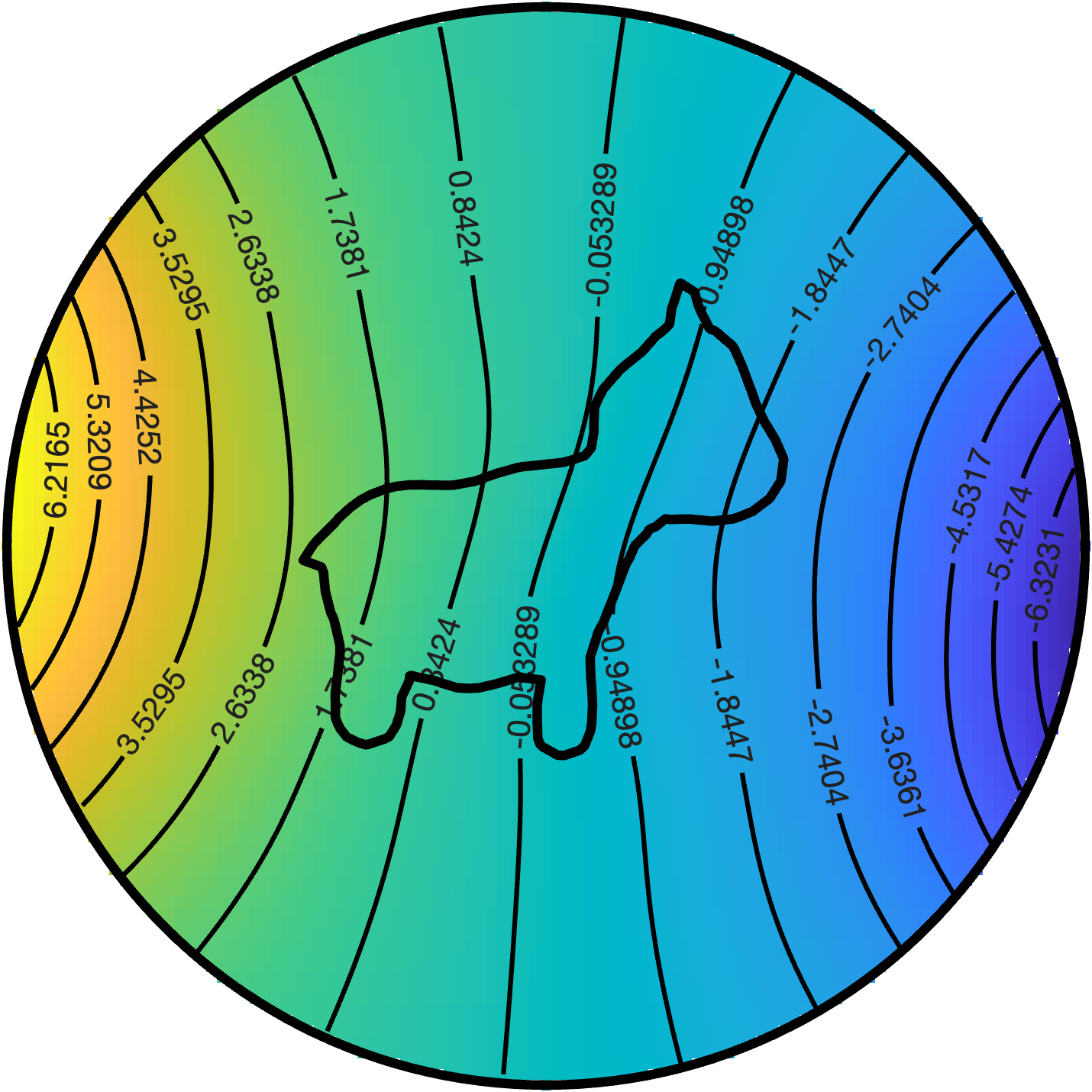} &
    \includegraphics[width=0.21\columnwidth, valign=c]{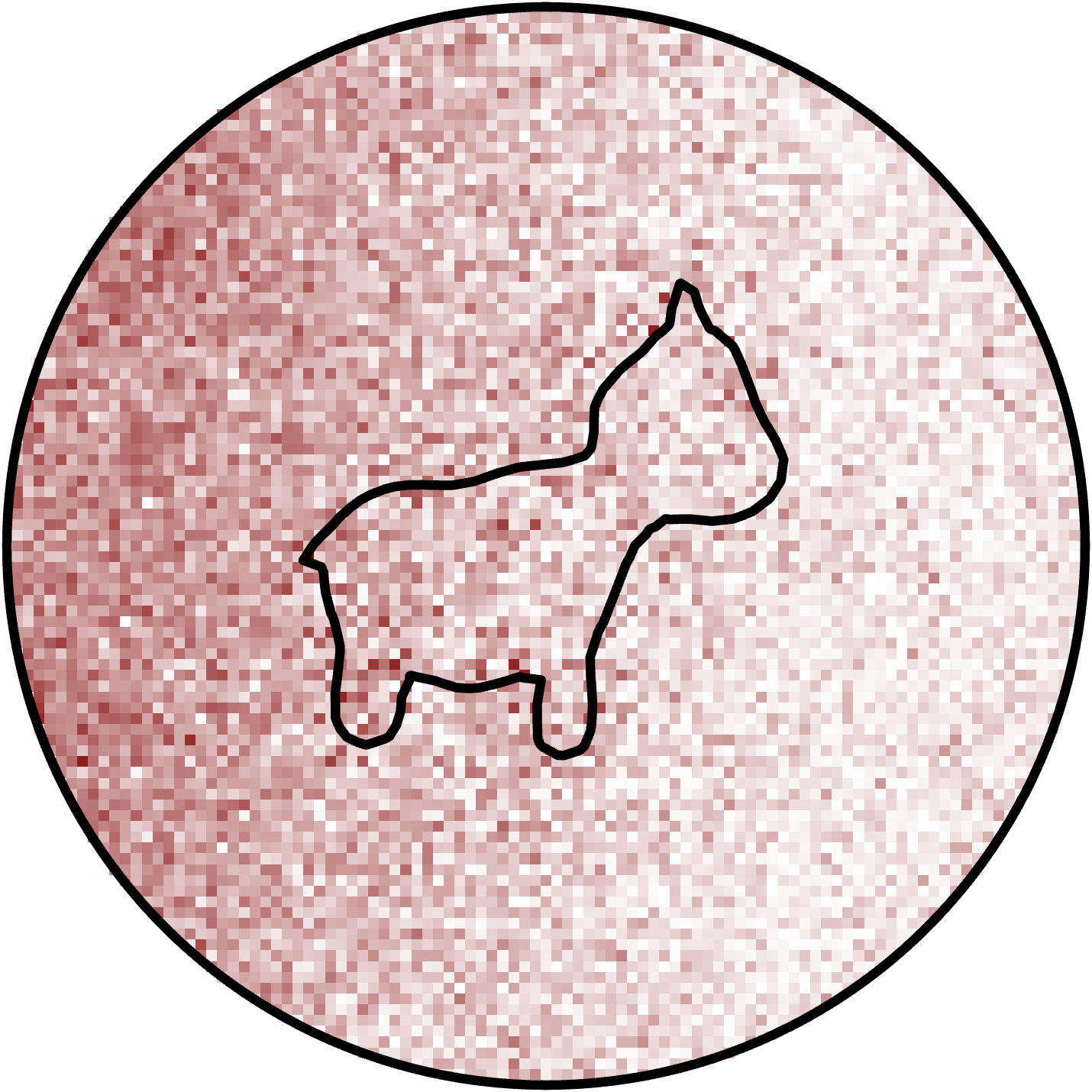} &
    \includegraphics[width=0.21\columnwidth, valign=c]{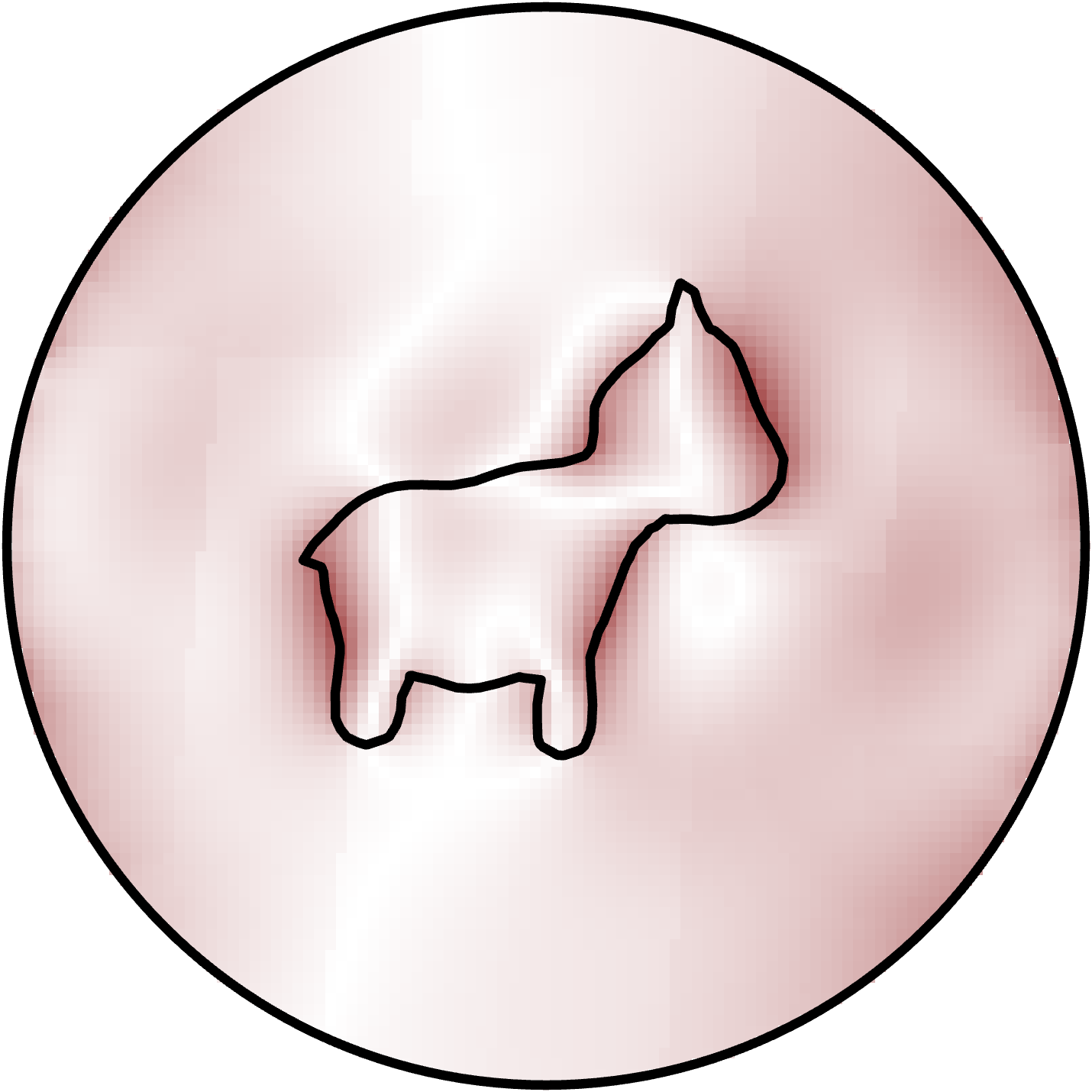} &
    \includegraphics[width=0.04\columnwidth, valign=c]{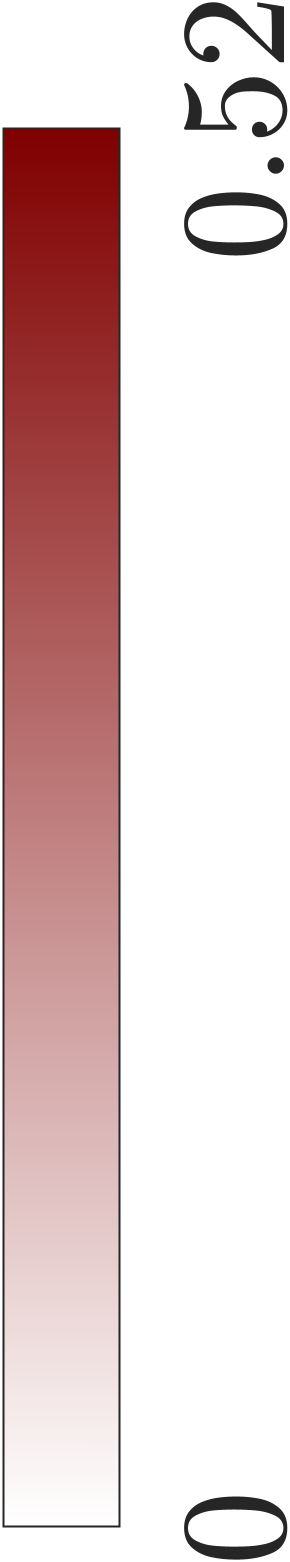} \vspace{0.1cm} \\
    % L2 difference
    {} & \multicolumn{2}{c}{\textbf{$\bm{L_2}$ Difference}} & $0.150$ & $0.096$ & {} \\
    % relative L2 difference
    {} & \multicolumn{2}{c}{\textbf{Relative $\bm{L_2}$ Difference}} & $5.51\%$ & $3.51\%$ & {} 
    \vspace{0.5cm}\\
    % var-reduced WoI
    {} & \shortstack{\textbf{Var-reduced WoI} \\ \textbf{Solution}} & \shortstack{\textbf{Learned Solution} \\ \textbf{(on test points)}} & \shortstack{\textbf{Var-reduced WoI} \\ \textbf{Error}} & \shortstack{\textbf{Learned Error} \\ \textbf{(test error)}} & {} \\
    % estimation and error
    \rotatebox[origin=t]{90}{\centering \shortstack{\textbf{Var-reduced WoI} \\ \textbf{Framework}}} & \includegraphics[width=0.21\columnwidth, valign=c]{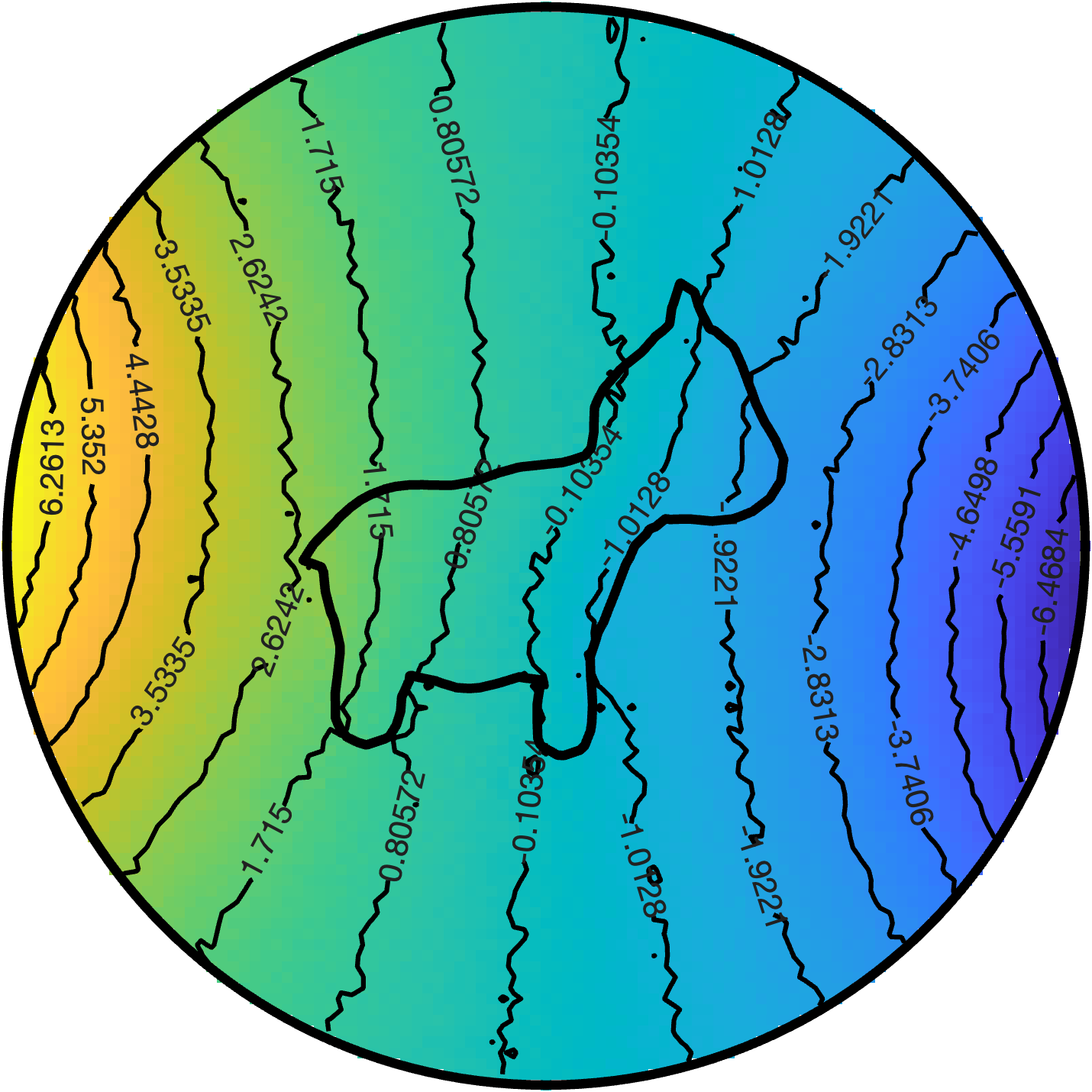}  &
    \includegraphics[width=0.21\columnwidth, valign=c]{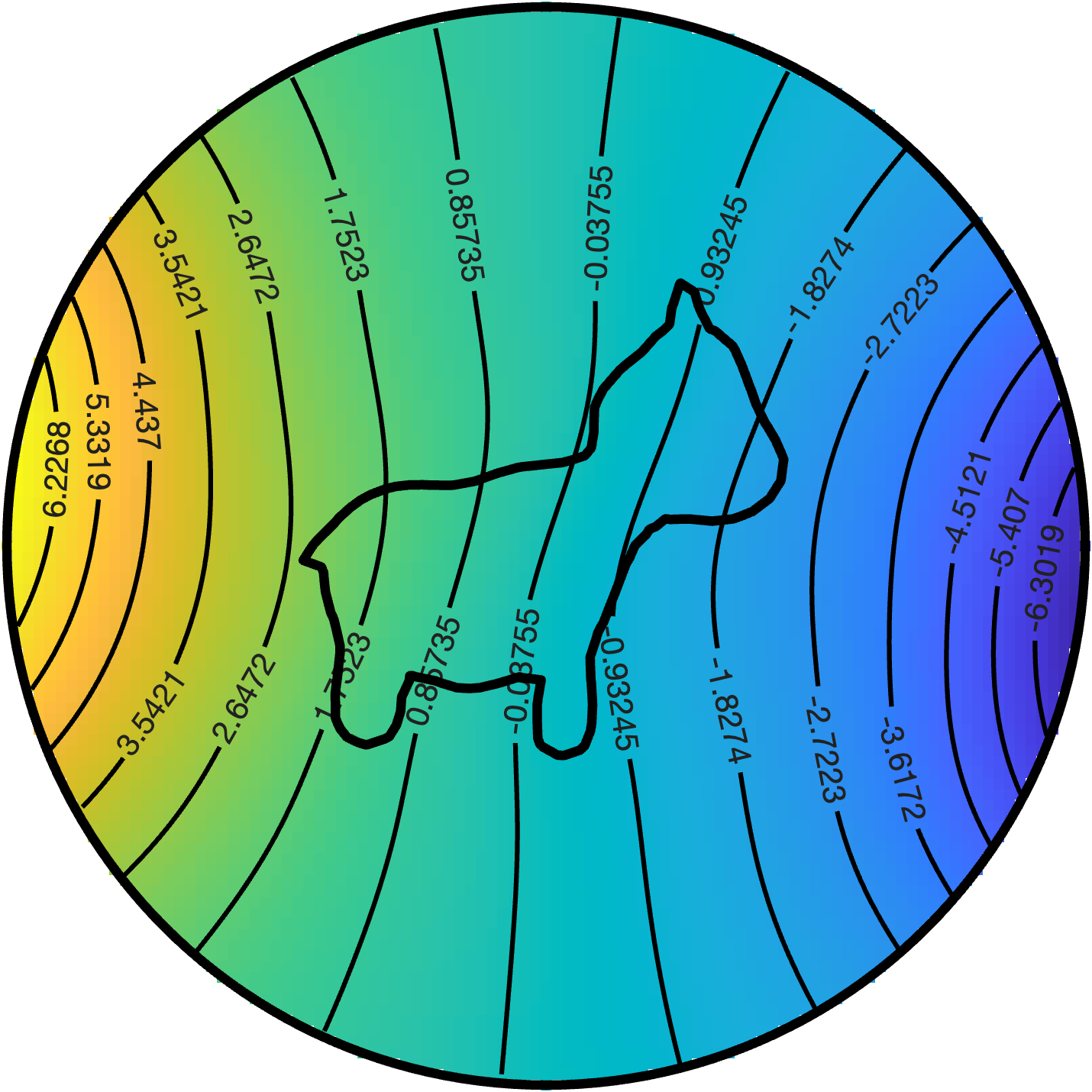} &
    \includegraphics[width=0.21\columnwidth, valign=c]{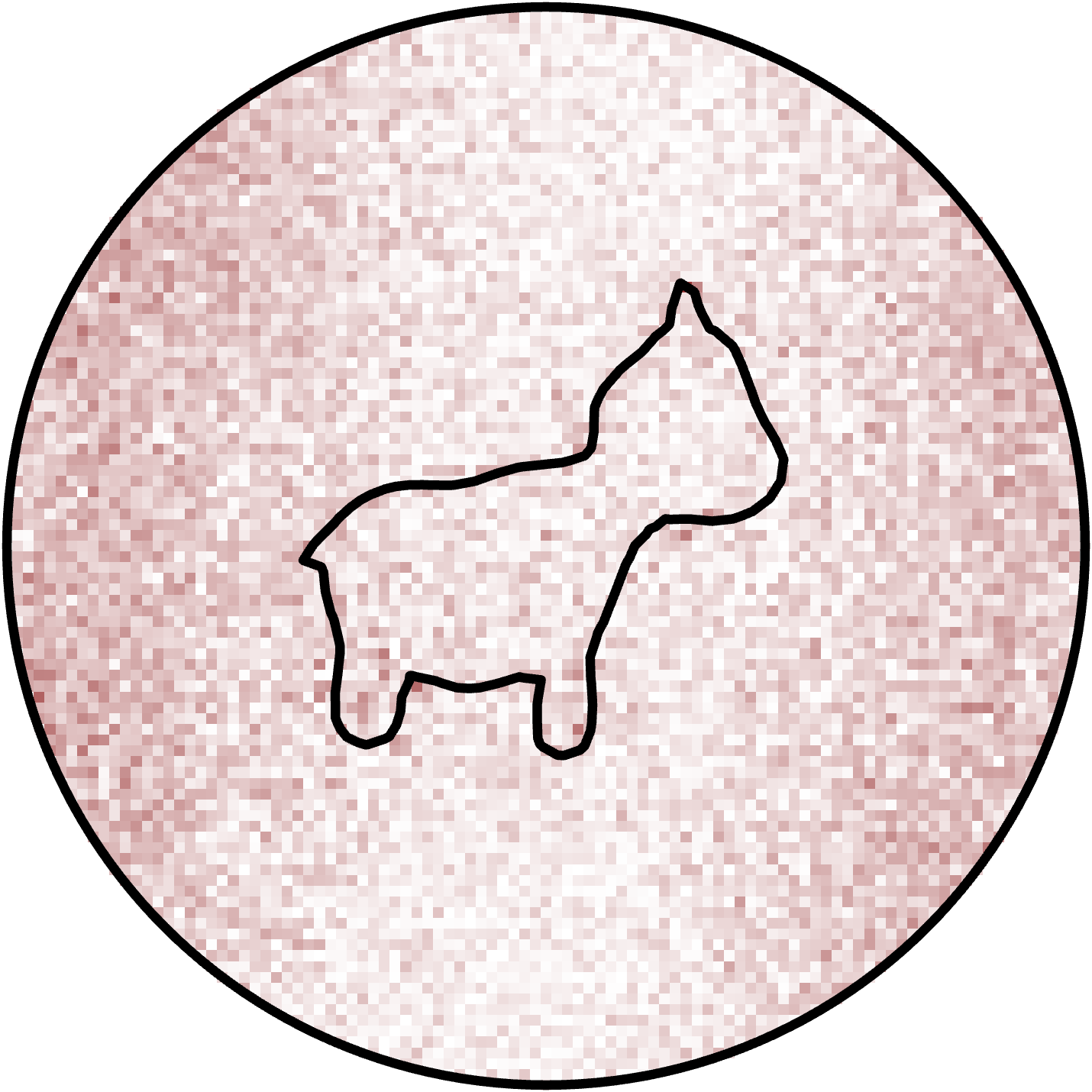} &
    \includegraphics[width=0.21\columnwidth, valign=c]{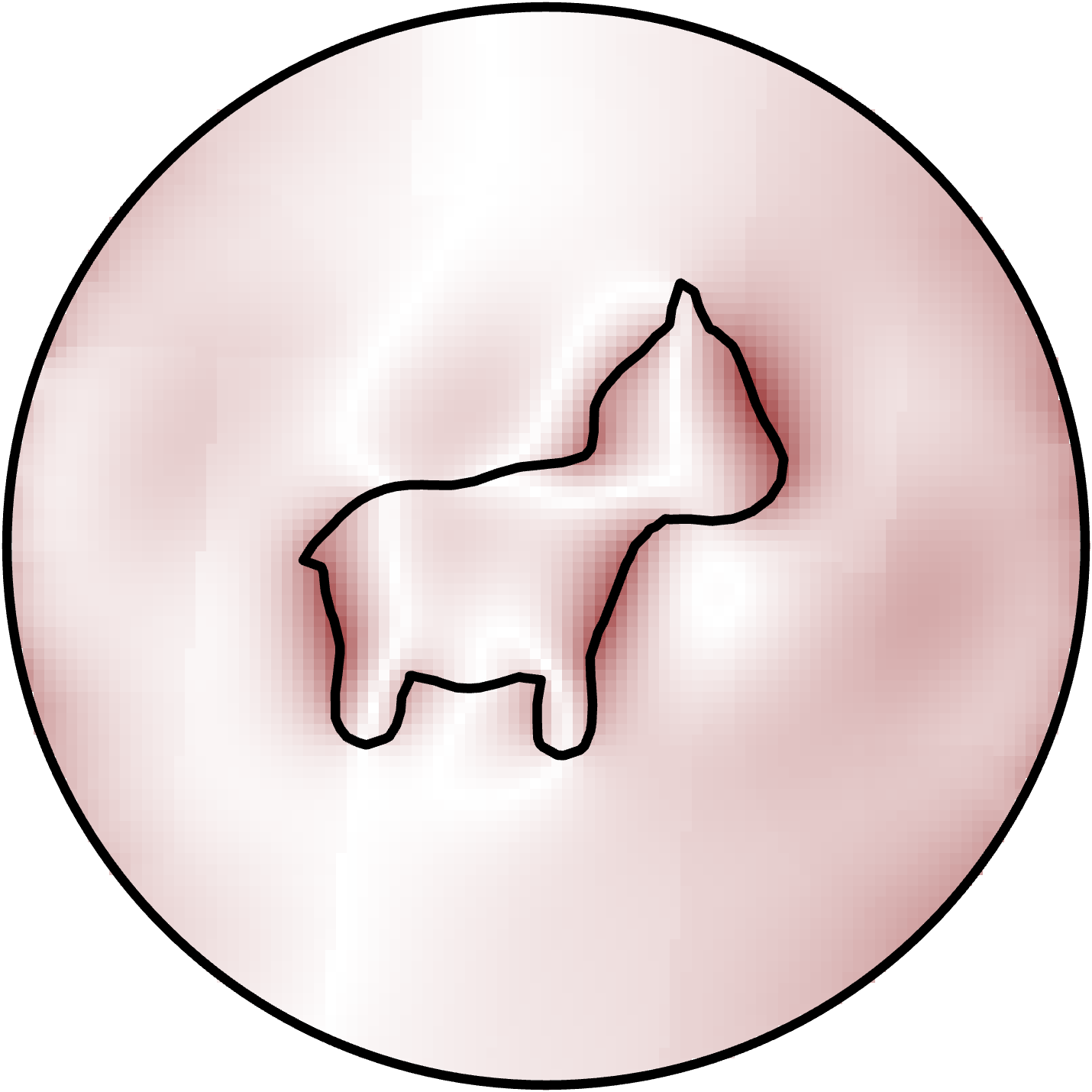} &
    \includegraphics[width=0.04\columnwidth, valign=c]{figure/EIT/nn/error_cb_large_text.png} \vspace{0.1cm} \\
    % L2 difference
    {} & \multicolumn{2}{c}{\textbf{$\bm{L_2}$ Difference}} & $0.094$ & $0.093$ & {} \\
    % relative L2 difference
    {} & \multicolumn{2}{c}{\textbf{Relative $\bm{L_2}$ Difference}} & $3.46\%$ & $3.41\%$ & {}
\end{tabular}
\end{minipage}
\caption{The first row shows the domain, a slice of the domain we show our estimation on, and the numerical solution obtained by a boundary integral equation method solver~\cite{Zheng_thesis}. The second row shows the solution and error estimated by WoI, and the learned solution and error with training data generated by the WoI estimator. The third row shows the solution and error estimated by variance-reduced WoI, and the learned solution and error with training data generated by the variance-reduced WoI estimator.}\label{ex_EIT}
\end{figure}

Since Eq.~\eqref{EIT} does not admit an analytical solution, we compute a numerical solution using a boundary integral equation method solver developed by Andrew Zheng for comparison~\cite{Zheng_thesis}. The neural network is trained using $\mathcal{N} = 7200$ random points within the domain and $\mathcal{N}_\text{bc} = 1000$ points on the boundary. Figure~\ref{ex_EIT} visualizes the domain, a representative slice used to display the solutions, the numerical solution, the WoI and variance-reduced WoI solutions, the learned solutions, and their differences from the numerical solution. Under this setup, the learned solution evaluates the generalization behaviour of the MLP.

\textbf{Example 5. Groundwater Flow.} Groundwater flow through subsurface environments is governed by Darcy’s law. In natural settings, the subsurface typically consists of multiple geological materials, such as sand, soil, and rock, each exhibiting distinct hydraulic conductivities. In this context, $u(\bm{x})$ represents the hydraulic head, and $\sigma(\bm{x})$ denotes the hydraulic conductivity. The model relates the volume flow rate of the fluid, $\sigma(\bm{x}) \nabla u(\bm{x})$, proportionally to the gradient of the hydraulic head and the properties of the porous medium.

We consider a steady-state Darcy flow problem in a hill-shaped domain, where the elevation is higher on the right and lower on the left. Beneath the hill surface are six distinct rocky regions, each characterized by its own hydraulic conductivity.
The domain is defined as $\Omega = [-A_1, A_1] \times [-A_2, A_2] \times [-A_3, h(x)]$, in which $h(x)$ describes the hill-shaped upper surface and is defined as
\[
h(x) = 
\begin{cases}
    kx & \quad \text{if $x \in (-a, a)$} \\
    ka & \quad \text{if $x \in [a, A_1]$} \\
    -ka & \quad \text{if $x \in [-A_1, -a]$} \\
\end{cases}.
\]
The hydraulic head satisfies the interface problem,
\begin{subequations}\label{groundwater}
\begin{empheq}[left=\empheqlbrace]{align}
 \Delta u &= 0 \quad \text{in $\Omega \ \backslash \ \bigcup_{i=2}^7\bdr{\Omega}_i$}  \\
 % ========================
 \sigma_1 \partial_{\bm{n}}u &= 
 \begin{cases}
     \sin(\pi x  / A_1) \cos(\pi y  / 2A_1) & \quad \text{if $z = h(x)$} \\
     0 & \quad \text{otherwise}
 \end{cases}
 \quad \text{on $\bdr{\Omega_1}$}, \\
 % ========================
 [u](\bm{x}) &= 0 \quad \text{on $\bdr{\Omega}_i$, $\forall i \in \{2, 3, 4, 5, 6, 7\}$}\\
 % ========================
 [\sigma \partial_{\bm{n}} u](\bm{x}) &= 0 \quad \text{on $\bdr{\Omega}_i$, $\forall i \in \{2, 3, 4, 5, 6, 7\}$}
\end{empheq}
\end{subequations}
that enforces continuity of head across material boundaries. The boundary conditions are interpreted as representing rainfall on the hilltop and water extraction in a residential area located downhill. 

\begin{figure}[htbp]
\centering
\begin{minipage}[l]{\textwidth}
\begin{tabular}{c c c}
    % Row 1: 2 plots
    \rotatebox[origin=t]{90}{\centering \textbf{Domain}} & \quad \quad
    \includegraphics[height=3.5cm, valign=c]{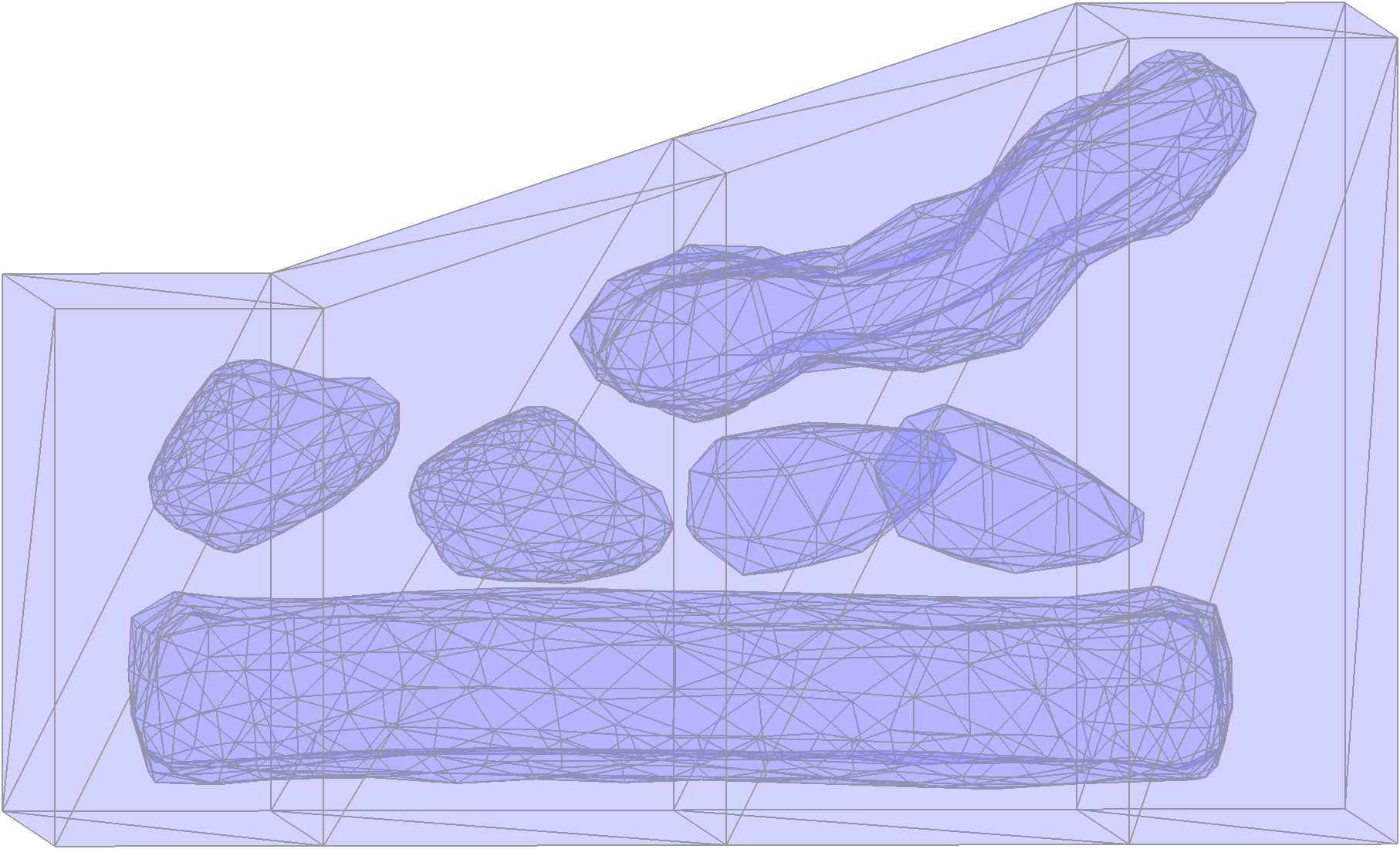} & \quad \quad
    \includegraphics[height=3.7cm, valign=c]{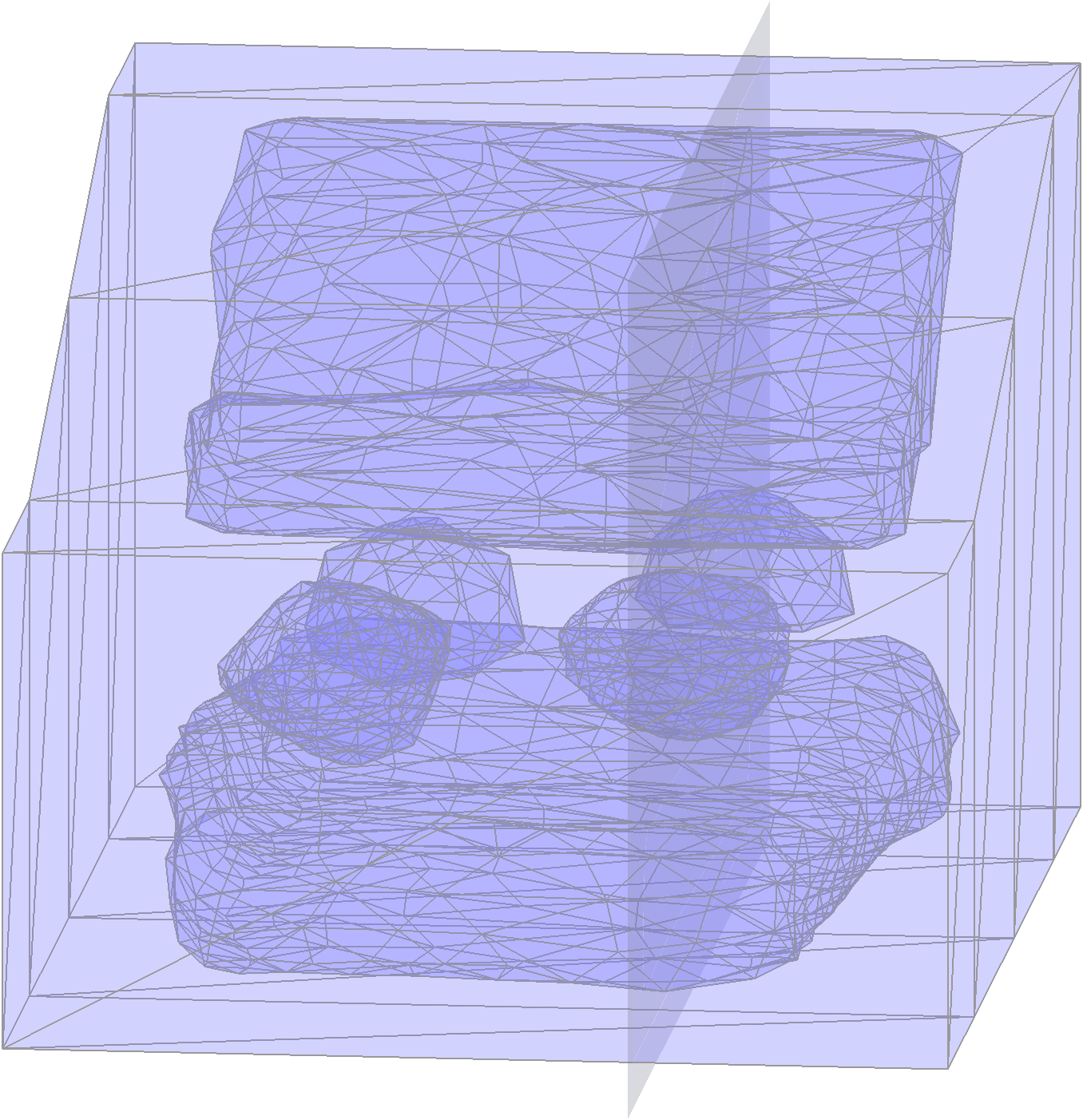}
\end{tabular}
\end{minipage}
\\ [0.5cm]
\begin{minipage}[l]{\textwidth}
\begin{tabular}{c c c c c}
    % number of walkers
    {} & $\mathcal{W} = 10^7$ & $\mathcal{W} = 2 \times 10^7$ & $\mathcal{W} = 4 \times 10^7$ & {} \\
    % WoI: estimation
    \rotatebox[origin=c]{90}{\centering \textbf{WoI}} &
    \includegraphics[width=0.25\columnwidth, valign=c]{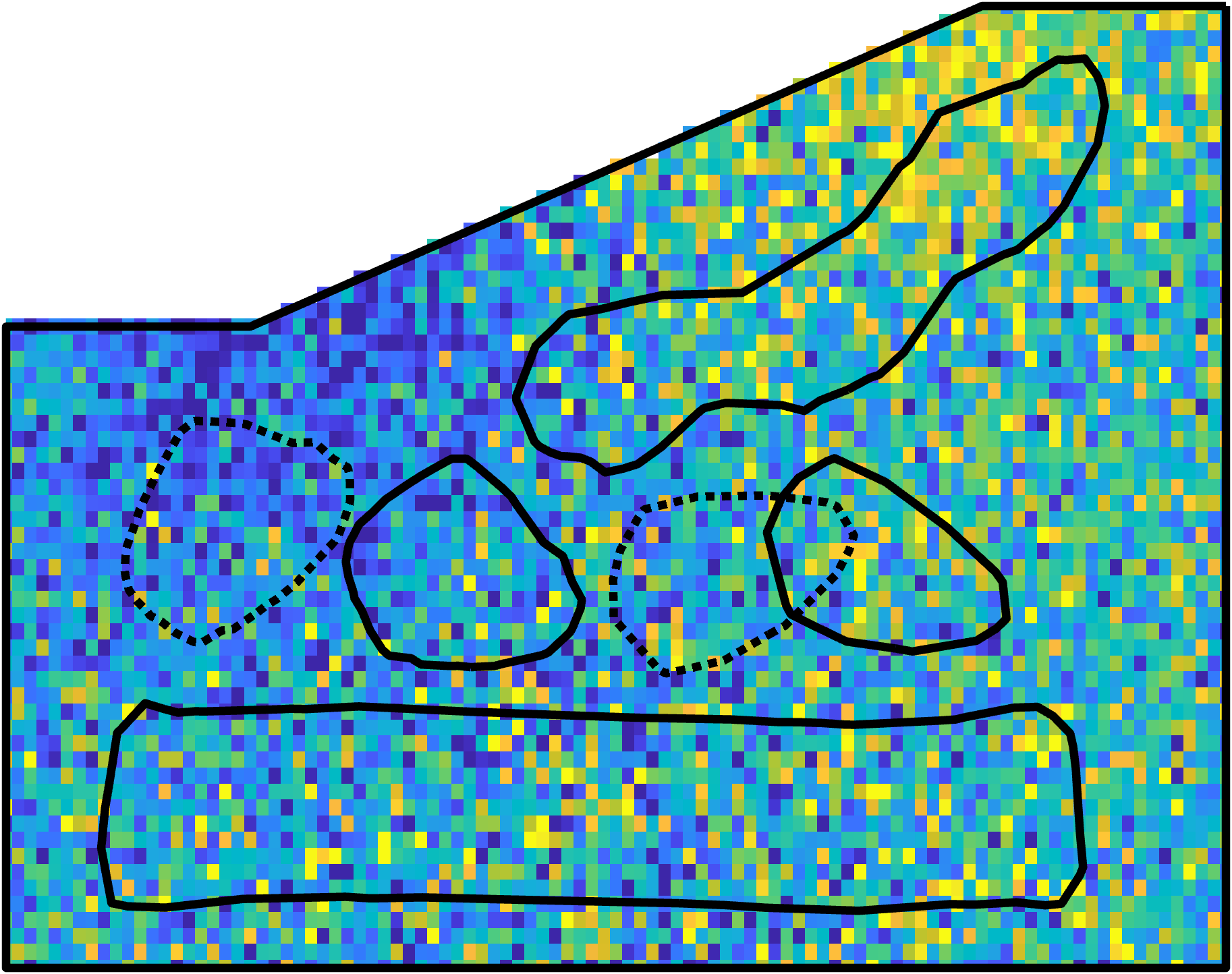} &
    \includegraphics[width=0.25\columnwidth, valign=c]{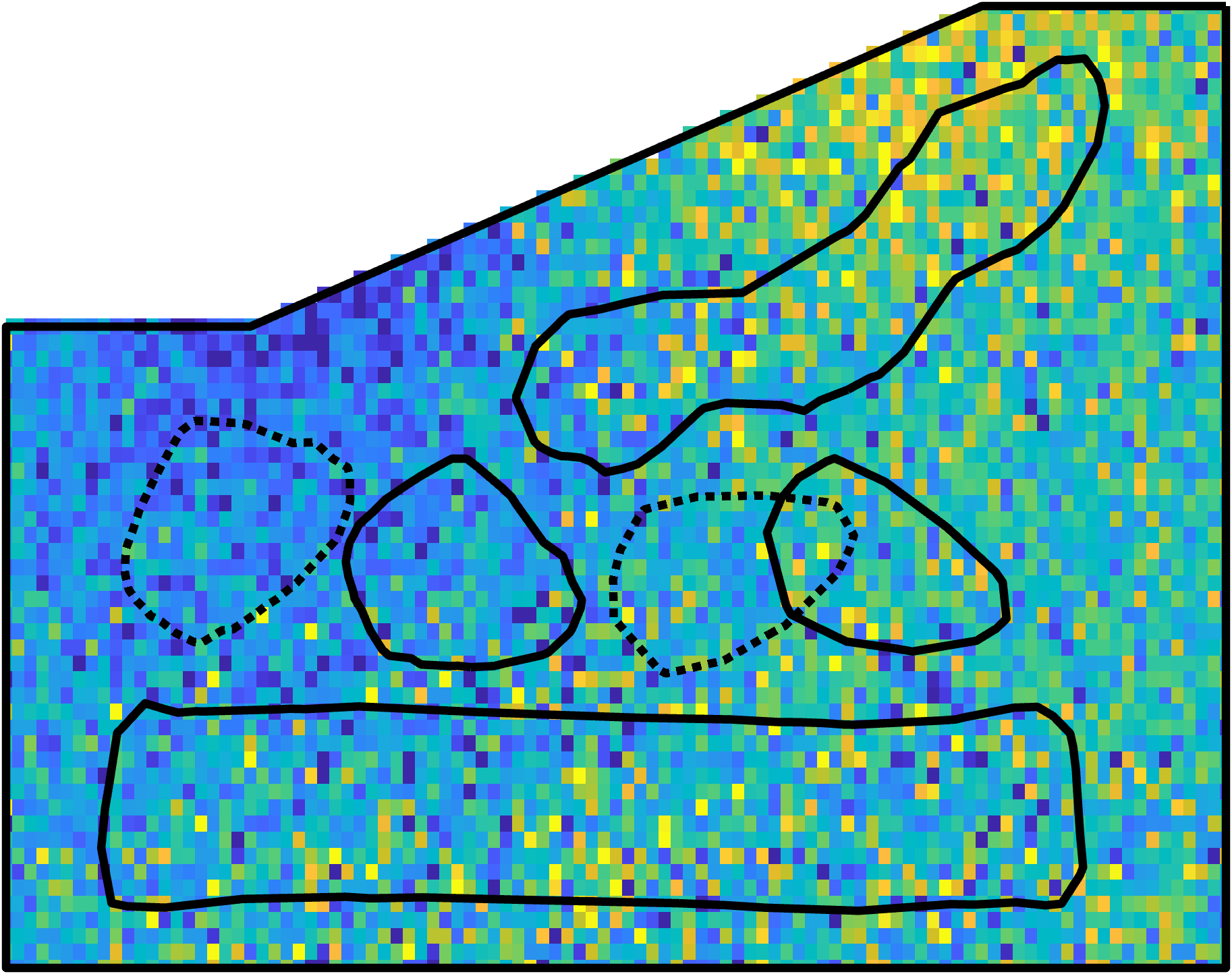}  &
    \includegraphics[width=0.25\columnwidth, valign=c]{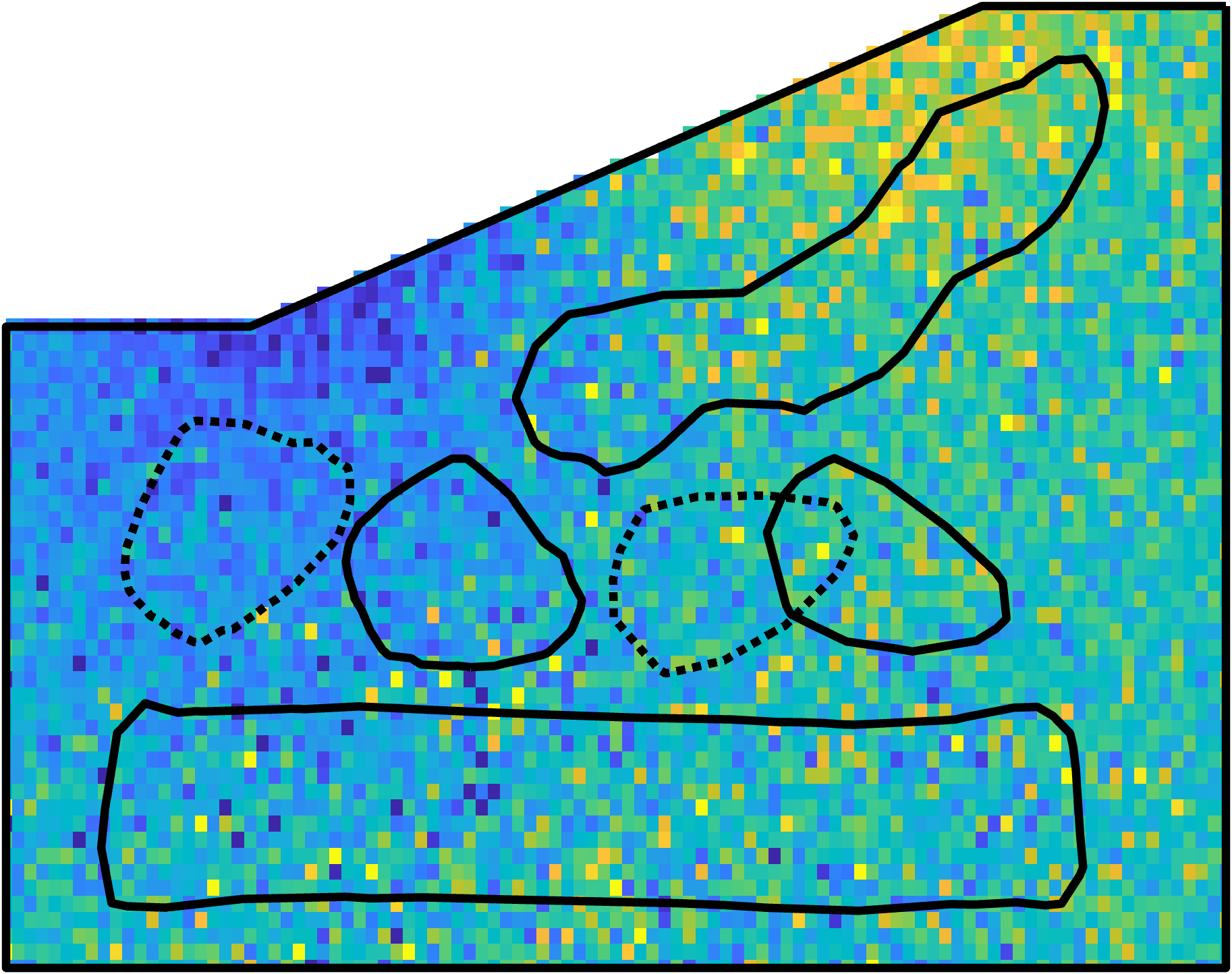} &
   \includegraphics[width=0.034\columnwidth, valign=c]{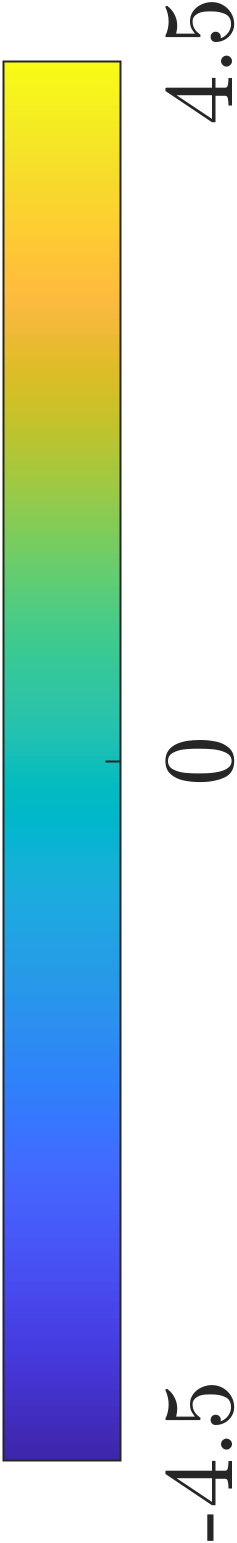}
   \vspace{0.5cm} \\ 
    % number of walkers
    {} & $\mathcal{W} = 2 \times 10^7$ & $\mathcal{W} = 4 \times 10^7$ & $\mathcal{W} = 8 \times 10^7$ & {} \\
    % second Row: estimation
    \rotatebox[origin=t]{90}{\centering \textbf{Var-reduced WoI}} &
    \includegraphics[width=0.25\columnwidth, valign=c]{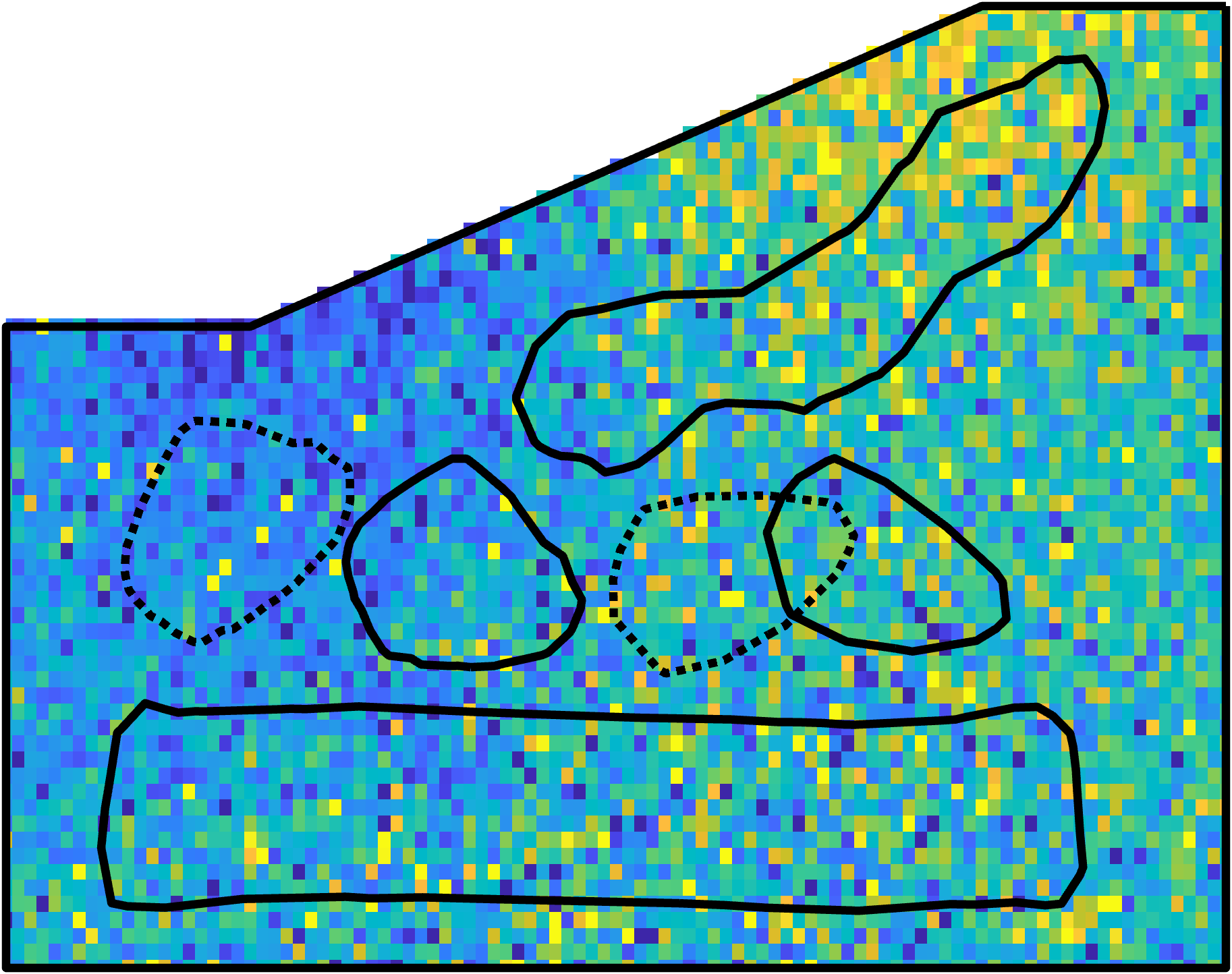} &
    \includegraphics[width=0.25\columnwidth, valign=c]{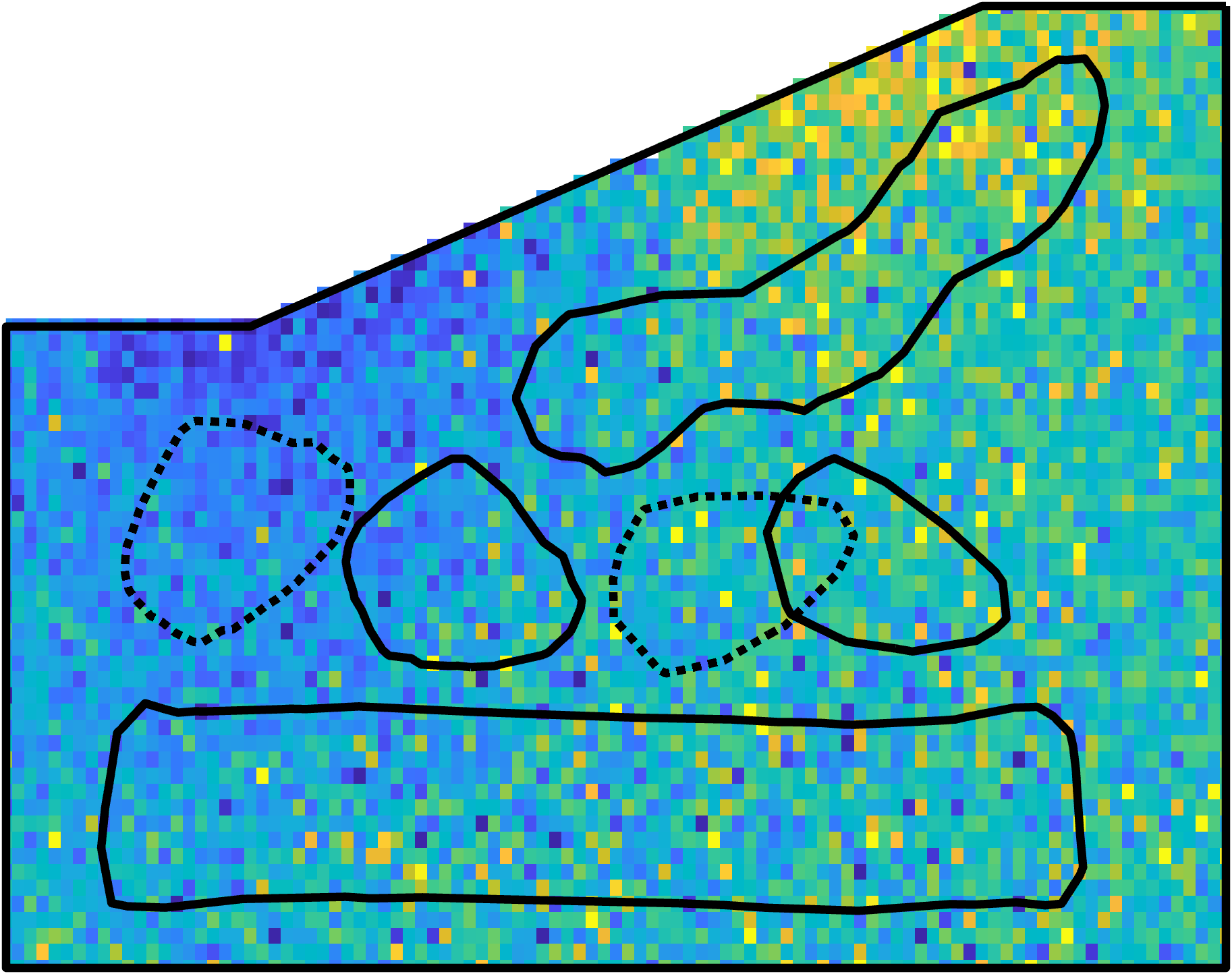}  &
    \includegraphics[width=0.25\columnwidth, valign=c]{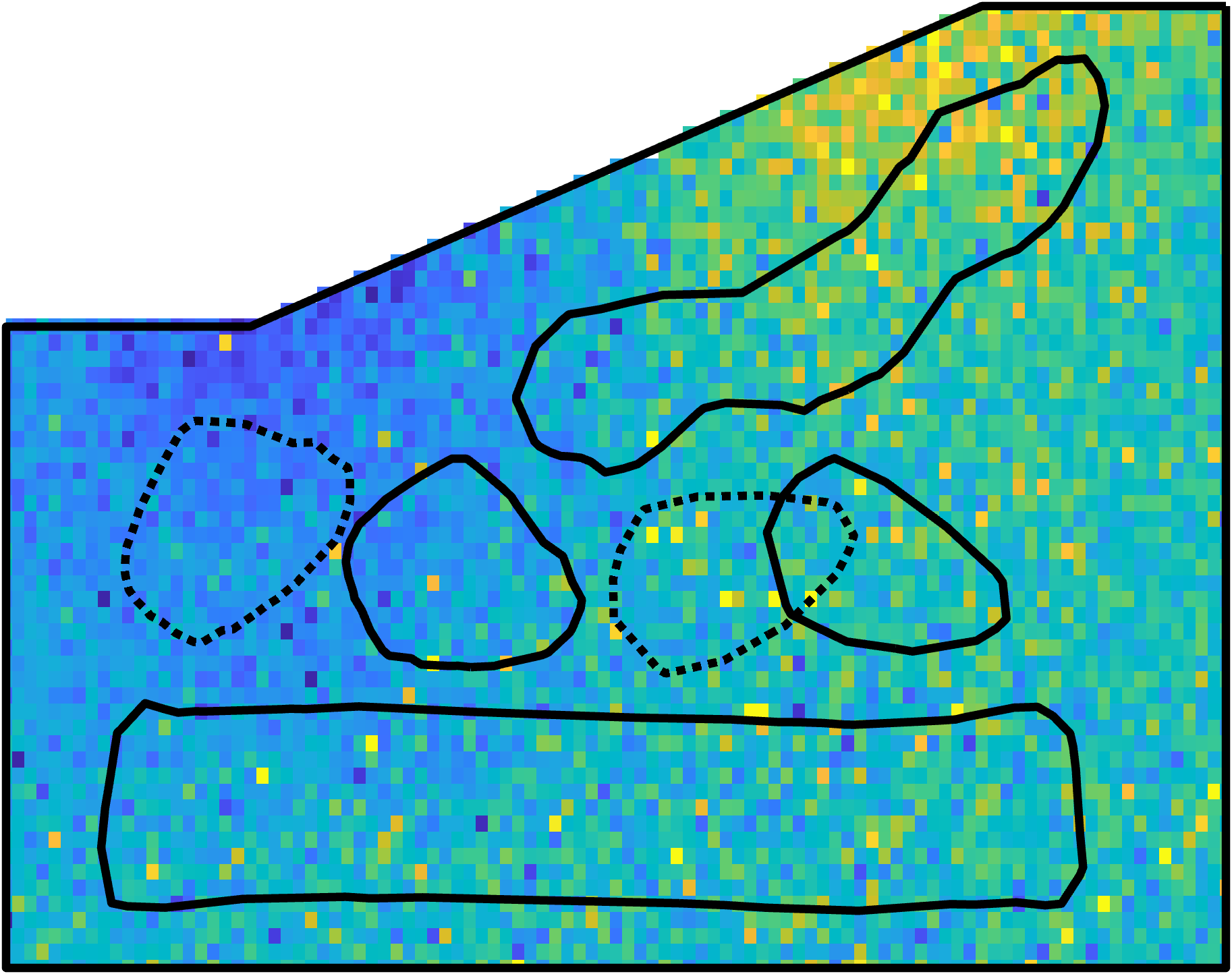} &
    \includegraphics[width=0.034\columnwidth, valign=c]{figure/ground_water/estimation_cb.png} \\
\end{tabular}
\end{minipage}
 \\[0.5cm]  % adds 1 cm vertical space
\begin{minipage}[l]{\textwidth}
\begin{tabular}{c c c c}
    {} & \textbf{WoI Data} & \textbf{Var-reduced WoI Data} \\
    % Row 1: 2 plots
    \rotatebox[origin=t]{90}{\centering \shortstack{\textbf{Learned} \\ \textbf{Hydraulic Head}}} & \quad \quad \quad
    \includegraphics[width=0.28\columnwidth, valign=c]{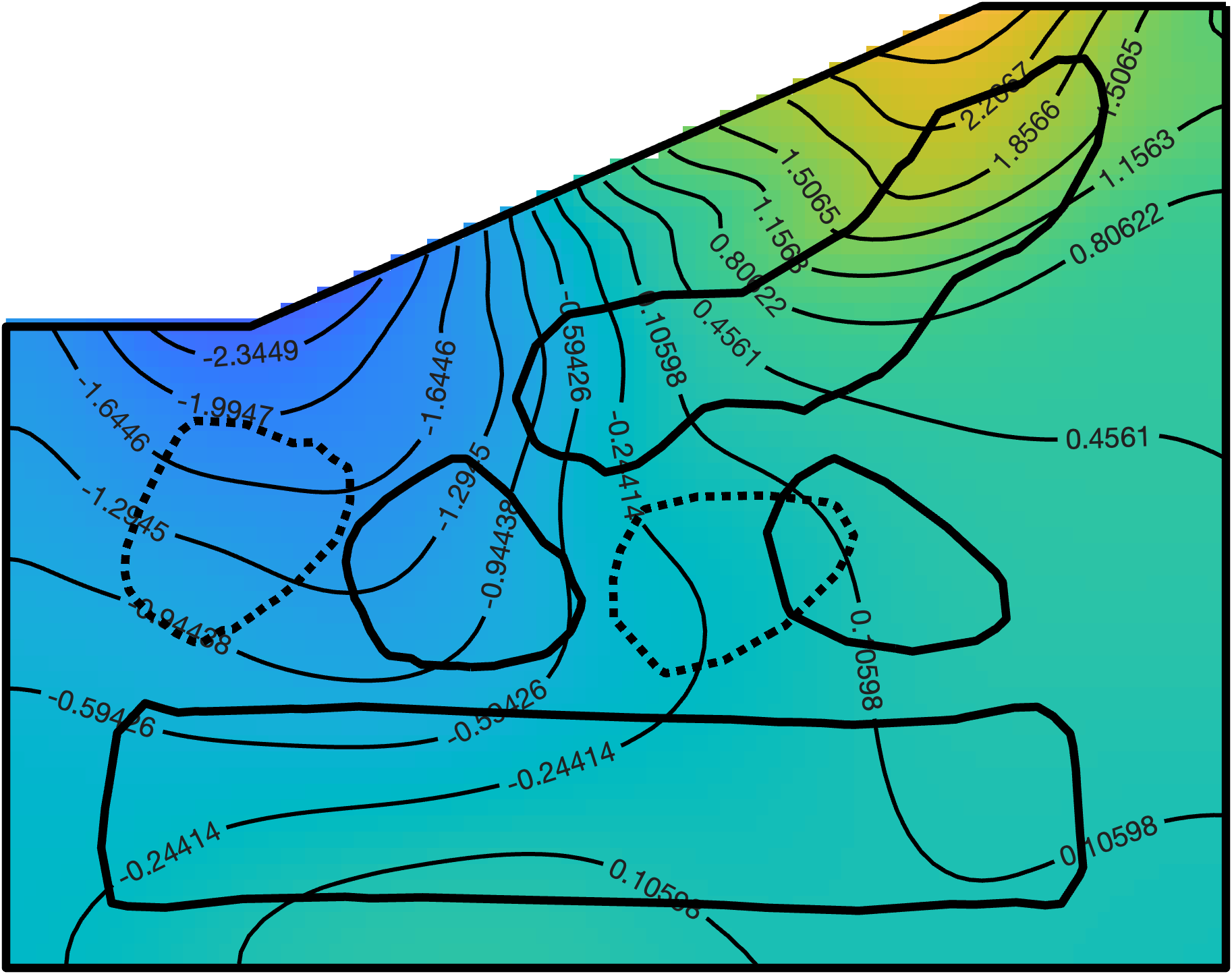} & \quad \quad \quad
    \includegraphics[width=0.28\columnwidth, valign=c]{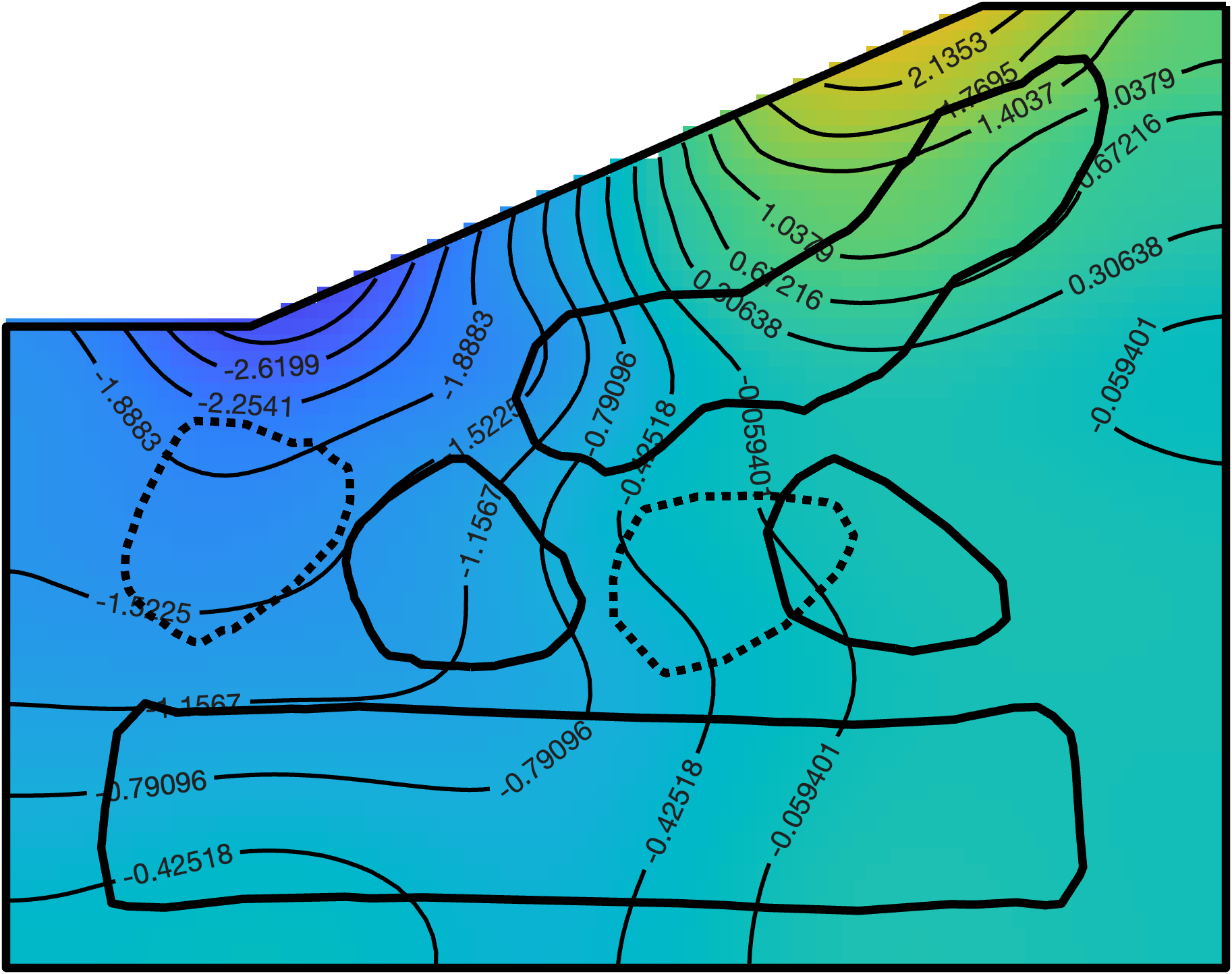}
    & \hspace{0.75cm} 
    \includegraphics[width=0.034\columnwidth, valign=c]{figure/ground_water/estimation_cb.png}
\end{tabular}
\end{minipage}
\caption{WoI and variance-reduced WoI solution to Eq.~\eqref{groundwater} and learned solution using data generated by both Monte Carlo methods. Both estimators walk $N = 5$ steps. The rocks with dotted outline do not intersect with the cutting plane. The hill is $\Omega_1$, the tilted rock at top is $\Omega_2$, the four rocks in the middle are $\Omega_3, \Omega_4, \Omega_5, \Omega_6$, the bedrock at the bottom is $\Omega_7$. We set $\sigma_1 = 0.2$, $\sigma_2 = 0.05$, $\sigma_3 = \sigma_4 = \sigma_5 = \sigma_6 = 0.03$, $\sigma_7 = 0.005$ with unit cm/s.}\label{ex_ground_water}
\end{figure}

We solve Eq.~\eqref{groundwater} with $A_1 = A_3 = 2$, $A_2 = 1.5$, $a = 1.2$ and $k = \frac{1}{3}$. While the numerical solver~\cite{Zheng_thesis} fails to produce a solution due to the number of interfaces, our framework successfully computes a descent solution. We present our WoI solution, variance-reduced WoI solution, and learned solutions over a slice of the domain in Figure \ref{ex_ground_water}. The neural network is trained using $\mathcal{N} = 7200$ randomly sampled in the domain and $\mathcal{N}_\text{bc} = 1200$ points on the boundary. The learned solution is, therefore, assumed by a separate set of test points.

\section{Conclusion}\label{conclusion}
Herein presented was a grid-free, highly parallelizable Monte Carlo estimator,  WoI, that has $\mathcal{O}(1/\sqrt{\mathcal{W}})$ convergence rate for solving an elliptic interface problem. We demonstrated that our estimator has uniform accuracy up to the boundary, close to the interfaces, and in the interior of the domain in contrast to other methods for this problem. Moreover, the estimator generalizes naturally to higher dimensional problems that are intractable to grid-based methods. Also, we can readily compute the gradient of the solution with almost no extra computational cost. 
Finally, we introduced the WoI framework in which our Monte Carlo estimator acts as a data generator for training a neural network.

Several promising directions remain open for future research. While we have proposed one effective PDF in Eq.~\eqref{naive_woi_p0} and Eq.~\eqref{naive_woi_p}, alternative PDFs that consider the geometric properties of the boundary and the interfaces may lead to improved efficiency or accuracy. Also, it would be interesting to explore whether and how our Monte Carlo estimator can be extended to other boundary conditions, such as Dirichlet, Robin, or mixed. Similarly, seeking the possibility of adapting our estimator to more general elliptic PDEs that include advection and reaction terms is another compelling direction for future research. Study of the network architectures that can represent solutions to these elliptic interface problems and be effectively trained by WoI data are of particular importance. Recent work~\cite{liu2025} has introduced a neural field construction capable of representing continuous functions with discontinuous gradients. Integrating such representation into our framework may improve the accuracy of the learned solutions.

% ===================================
\section*{Acknowledgments}
We express our gratitude to Andrew Zheng for his contribution in generating the numerical solution in Example 4. We also thank the \href{https://www.mathcha.io/}{Mathcha} team for providing an excellent platform that facilitated the creation of our diagrams. Additionally, we acknowledge the authors of \href{https://github.com/alecjacobson/gptoolbox}{gptoolbox} and the \href{https://www.mathworks.com/matlabcentral/fileexchange/35623-tree-data-structure-as-a-matlab-class}{MATLAB tree structure} for their toolboxes, which helped us build some of our code.
% ===================================
% Reference
%\newpage
\addcontentsline{toc}{section}{References}
\bibliographystyle{plain}
\noindent
\bibliography{ref}

% ===================================
% Appendices
\appendix
\section{A Proof of Theorem~\ref{naive_woi_theorem}} \label{Append-proof_of_naive_woi_theorem}
We prove Theorem \ref{naive_woi_theorem} in this appendix.
\begin{theorem}
Given a query point $\bm{x}$ and a fixed set of $h_0, h_1, \dots, h_i \in \{1, \dots, N\}$. Let $\{Y_0, Y_1, \dots, Y_i, \dots\}$ be a Markov chain starting from $\bm{x}$ such that $Y_i \in \bdr{\Omega_{h_i}}$. Let $p_0(\bm{x}, \bm{y})$ and $p(\bm{x}, \bm{y})$ be any probability density functions (PDFs) denoting the transitional distribution of $\bm{y}$ given $\bm{x}$. Define a family of functions
\begin{equation}
    \adj{Q}_{(h_i, \dots, h_1, h_0)} = 
    \begin{cases}
        \frac{\beta_{h_0}(Y_0)}{p_0(\bm{x}, Y_0)} \quad & \text{if $i = 0$}\\
         \adj{Q}_{(h_{i-1}, \dots, h_1, h_0)} \frac{\adj{K}_{h_i, h_{i-1}}(Y_i, Y_{i-1})}{p(Y_{i-1}, Y_i)}\quad & \text{if $i \geq 1$}.\\
    \end{cases}
\end{equation}
Then, $\adj{Q}_{(h_i, \dots, h_1, h_0)} \Phi(\bm{x}, Y_i)$ is an unbiased estimator for 
\[
\int_{\bdr{\Omega}_{h_i}} \Phi(\bm{x}, \bm{y}) 
(
\adj{\K}_{h_i, h_{i-1}} \cdots \adj{\K}_{h_2, h_1}\adj{\K}_{h_1, h_0} \beta_{h_0}
) \ d A_{\bm{y}}.
\]
i.e. 
\[
\int_{\bdr{\Omega}_{h_i}} \Phi(\bm{x}, \bm{y}) 
(\adj{\K}_{h_i, h_{i-1}} \cdots \adj{\K}_{h_2, h_1}\adj{\K}_{h_1, h_0} \beta_{h_0}
) \ d A_{\bm{y}} 
= \mathbb{E} \bigg[\adj{Q}_{(h_i, \dots, h_1, h_0)} \Phi(\bm{x}, Y_i) \bigg].
\]
\end{theorem}
\textit{Proof.} Unwrapping the recursive definition of $\adj{Q}_{(h_i, \dots h_1, h_0)}$, we have
\begin{equation}
    \adj{Q}_{(h_i, \dots h_1, h_0)} = \frac{\prod_{m=1}^i \adj{K}_{h_m, h_{m-1}}(Y_m, Y_{m-1})}{\prod_{m=1}^{i} p(Y_{m-1}, Y_m)} \frac{\beta_{h_0}(Y_0)}{p_0(\bm{x}, Y_0)} \notag
\end{equation}
Let $\rho_{Y_0, Y_1, \dots Y_i, \dots}(\bm{y}_0, \bm{y}_1, \dots, \bm{y}_i)$ be the joint distribution of $Y_0, Y_1, \dots, Y_i$. By chain rule,
\begin{align}
    \rho_{Y_0, Y_1, \dots Y_i}(\bm{y}_0, \bm{y}_1, \dots, \bm{y}_i)
    &=
    p_0(\bm{x}, \bm{y}_0) \prob(\bm{y}_1 | \bm{y}_0) \prob(\bm{y}_2 | \bm{y}_0, \bm{y}_1) \dots \prob(\bm{y}_i | \bm{y}_0, \bm{y}_1 \dots, \bm{y}_{i-1}) \notag \\
    &=
    p_0(\bm{x}, \bm{y}_0) \prob(\bm{y}_1 | \bm{y}_0) \prob(\bm{y}_2 | \bm{y}_1) \dots \prob(\bm{y}_i | \bm{y}_{i-1}) \notag \\
    &=
    p_0(\bm{x}, \bm{y}_0) p(\bm{y}_0, \bm{y}_1) p(\bm{y}_1 , \bm{y}_2) \dots p(\bm{y}_{i-1}, \bm{y}_i) \notag \\
    &=
    p_0(\bm{x}, \bm{y}_0) \prod_{m=1}^{i} p(\bm{y}_{m-1}, \bm{y}_m), \notag
\end{align}
where the second equality holds since ${Y_0, Y_1, \dots, Y_i, \dots}$ is a Markov chain, and the third equality follows immediately from the definition of transitional probability, $p(\bm{x}, \bm{y}) = \prob(\bm{y} | \bm{x})$. By definition of expected value,
\begin{align}
    & \mathbb{E} \bigg[\adj{Q}_{(h_i, \dots h_1, h_0)} \Phi(\bm{x}, Y_i) \bigg] \notag \\
    =
    & \int_{\bdr{\Omega_{h_i}}} \dots \int_{\bdr{\Omega_{h_1}}} \int_{\bdr{\Omega_{h_0}}} \adj{Q}_{(h_i, \dots h_1, h_0)} \Phi(\bm{x}, \bm{y}_i) \rho_{Y_0, Y_1, \dots Y_i}(\bm{y}_0, \bm{y}_1, \dots, \bm{y}_i) \ d\ell_{\bm{y}_{h_0}} \ d\ell_{\bm{y}_{h_1}} \dots \ d\ell_{\bm{y}_{h_i}}.
\end{align}
Substitute $\adj{Q}_{(h_i, \dots h_1, h_0)}$ and $\rho_{Y_0, Y_1, \dots Y_i}(\bm{y}_0, \bm{y}_1, \dots, \bm{y}_i)$, then we have
\begin{align}
    & \mathbb{E} \bigg[\adj{Q}_{(h_i, \dots h_1, h_0)} \Phi(\bm{x}, Y_i) \bigg] \notag \\
    =
    &\int_{\bdr{\Omega_{h_i}}} \dots \int_{\bdr{\Omega_{h_1}}} \int_{\bdr{\Omega_{h_0}}}
    \begin{aligned}[t]
        & \frac{\prod_{m=1}^i \adj{K}_{h_m, h_{m-1}}(\bm{y}_m, \bm{y}_{m-1})}{\prod_{m=1}^{i} p(\bm{y}_{m-1}, \bm{y}_m)} \frac{\beta_{h_0}(\bm{y}_0)}{p_0(\bm{x}, \bm{y}_0)}
        \Phi(\bm{x},  \bm{y}_i) \notag \\
        &p_0(\bm{x}, \bm{y}_0) \prod_{m=1}^{i} p(\bm{y}_{m-1}, \bm{y}_m) 
    \ d\ell_{\bm{y}_{h_0}} \ d\ell_{\bm{y}_{h_1}} \dots \ d\ell_{\bm{y}_{h_i}}
    \end{aligned} \notag \\
    =
    &\int_{\bdr{\Omega_{h_i}}} \dots \int_{\bdr{\Omega_{h_1}}} \int_{\bdr{\Omega_{h_0}}} 
    \prod_{m=1}^i \adj{K}_{h_m, h_{m-1}}(\bm{y}_m, \bm{y}_{m-1}) \beta_{h_0}(\bm{y}_0) 
    \Phi(\bm{x},  \bm{y}_i) \ d\ell_{\bm{y}_{h_0}}
    \ d\ell_{\bm{y}_{h_1}} \dots \ d\ell_{\bm{y}_{h_i}} \notag \\
    =
    &\int_{\bdr{\Omega}_{h_i}} \Phi(\bm{x}, \bm{y}) 
    (\adj{\K}_{h_i, h_{i-1}} \cdots \adj{\K}_{h_2, h_1}\adj{\K}_{h_1, h_0} \beta_{h_0}
    ) \ \lint{y}, \notag 
\end{align}
as required. \hfill \qed
\section{A Proof of Theorem~\ref{woi_theorem}} \label{Append-proof_of_woi_theorem}
We prove Theorem \ref{woi_theorem} in this appendix.

\begin{theorem}
    Let $H_0, H_1, \dots, H_i$ be independent discrete random variables with p.d.f
    \begin{equation}
        \prob(H_i = n) = 
            \begin{cases}
                \frac{1}{N} & \quad \text{$\forall n \in \{1, \ldots, N\}$, when $i = 0$} \\
                \frac{|\alpha_n|}{||\bm{\alpha}||_1} & \quad \text{$\forall n \in \{1, \ldots, N\}$, when $i \geq 1$}
            \end{cases}, \notag
    \end{equation}
    where $||\cdot||_1$ is the $\ell_1$ norm. It can be shown that 
    \begin{align}
        &\sum_{h_{i}, \dots, h_1, h_0=1}^{N}
        \int_{\bdr{\Omega}_{h_i}} \Phi(\bm{x}, \bm{y}) 
        \bigg(
        (\alpha_{h_i} \adj{\K}_{h_i, h_{i-1}}) \cdots (\alpha_{h_2}\adj{\K}_{h_2, h_1} ) (\alpha_{h_1}\adj{\K}_{h_1, h_0}) \beta_{h_0}
        \bigg) \ \lint{y} \notag \\
        =
        & N ||\bm{\alpha}||_1^i
        \mathbb{E}_{H_i, \dots, H_0} \bigg[
        \prod_{m=1}^i \sign(\alpha_{H_m})
         \int_{\bdr{\Omega}_{H_i}}
         \Phi(\bm{x}, \bm{y})
         \bigg(\adj{\K}_{H_i, H_{i-1}} \dots \adj{\K}_{H_1, H_0} \beta_{H_0} \lint{y}
        \bigg) \bigg]. \notag
    \end{align}
\end{theorem}

\textit{Proof.} Since $H_0, H_1, \dots H_M$ are independent, their joint probability density is
\begin{align}
    \prob(H_0 = h_0, H_1 = h_1, \dots, H_i = h_i)
    &= 
    \prob(H_0 = h_0)  \bigg(\prod_{m=1}^i \prob (H_i = h_i)\bigg) \notag \\
    &= 
    \frac{1}{N}\bigg(\prod_{m=1}^i \frac{|\alpha_{h_m}|}{||\bm{\alpha}||_1}\bigg). \notag
\end{align}
Observe that $\alpha_{h_i} = \sign(\alpha_{h_i}) |\alpha_{h_i}|$ holds for all $h_i \in \{1, \ldots, N\}$. Then,
\begin{align}
    &\sum_{h_{i}, \dots, h_1, h_0=1}^{N}
    \int_{\bdr{\Omega}_{h_i}} \Phi(\bm{x}, \bm{y}) 
    \bigg(
    (\alpha_{h_i} \adj{\K}_{h_i, h_{i-1}}) \cdots (\alpha_{h_2}\adj{\K}_{h_2, h_1} ) (\alpha_{h_1}\adj{\K}_{h_1, h_0}) \beta_{h_0}
    \bigg) \ \lint{y} \notag \\
    =
    &\sum_{h_{i}, \dots, h_1, h_0=1}^{N}
    \bigg(\prod_{m=1}^i \sign(\alpha_{h_m}) |\alpha_{h_m}| \bigg)
    \int_{\bdr{\Omega}_{h_i}} \Phi(\bm{x}, \bm{y}) 
    (\adj{\K}_{h_i, h_{i-1}} \cdots \adj{\K}_{h_2, h_1} \adj{\K}_{h_1, h_0} \beta_{h_0}) \ \lint{y} \notag \\
    =
    &N ||\bm{\alpha}||^i \sum_{h_{i}, \dots, h_1, h_0=1}^{N}
    \begin{aligned}[t]
        &\bigg(\frac{1}{N}\prod_{m=1}^i \frac{|\alpha_{h_m}|}{||\bm{\alpha}||} \bigg)
        \bigg( \prod_{m=1}^i \sign(\alpha_{h_m}) \bigg) \\
        &\int_{\bdr{\Omega}_{h_i}} \Phi(\bm{x}, \bm{y}) 
    (\adj{\K}_{h_i, h_{i-1}} \cdots \adj{\K}_{h_2, h_1} \adj{\K}_{h_1, h_0} \beta_{h_0}) \ \lint{y}
    \end{aligned} \notag \\
    =
    &N ||\bm{\alpha}||^i \sum_{h_{i}, \dots, h_1, h_0=1}^{N}
    \begin{aligned}[t]
        &\prob(H_0 = h_0, H_1 = h_1, \dots, H_i = h_i)
        \bigg( \prod_{m=1}^i \sign(\alpha_{h_m}) \bigg) \\
        &\int_{\bdr{\Omega}_{h_i}} \Phi(\bm{x}, \bm{y}) 
    (\adj{\K}_{h_i, h_{i-1}} \cdots \adj{\K}_{h_2, h_1} \adj{\K}_{h_1, h_0} \beta_{h_0}) \ \lint{y}
    \end{aligned} \notag \\
    =
    & N ||\bm{\alpha}||_1^i
    \mathbb{E}_{H_i, \dots, H_0} \bigg[
    \prod_{m=1}^i \sign(\alpha_{H_m})
     \int_{\bdr{\Omega}_{H_i}}
     \Phi(\bm{x}, \bm{y})
     \bigg(\adj{\K}_{H_i, H_{i-1}} \dots \adj{\K}_{H_1, H_0} \beta_{H_0} \lint{y}
    \bigg) \bigg], \notag
\end{align}
where the last step follows immediately from the definition of the expected value.
\end{document}